%% file: sigma26-037.tex
\pdfoutput=1
\RequirePackage{ifpdf}
\ifpdf 
\documentclass[pdftex]{sigma}
\else
\documentclass{sigma}
\fi

\numberwithin{equation}{section}

\newtheorem{Theorem}{Theorem}[section]
\newtheorem*{Theorem*}{Theorem}
\newtheorem*{Claim*}{Claim}
\newtheorem{Corollary}[Theorem]{Corollary}
\newtheorem{Lemma}[Theorem]{Lemma}
\newtheorem{Proposition}[Theorem]{Proposition}
\newtheorem{Conjecture}[Theorem]{Conjecture}

\theoremstyle{definition}
\newtheorem{Definition}[Theorem]{Definition}
\newtheorem{Remark}[Theorem]{Remark}

\usepackage{mathrsfs}
\usepackage{tikz}
\usepackage{enumitem}
\usepackage{subcaption}
\usepackage[all,cmtip]{xy}
\usepackage{bigdelim}

\newcommand{\circled}[1]{\scalebox{0.8}{\textcircled{\fontsize{8}{9}\selectfont \raisebox{-0.1pt}{#1}}}}

\newcommand{\circledd}[1]{\scalebox{0.8}{\textcircled{\fontsize{6}{7}\selectfont \raisebox{0.5pt}{#1}}}}

\newcommand*\circleddd[1]{\scalebox{0.4}{\tikz[baseline=(char.base)]{
 \node[shape=circle,draw,inner sep=1pt,very thick] (char) {#1};}}}

\newcommand*\circlee[1]{\scalebox{0.6}{\tikz[baseline=(char.base)]{
 \node[shape=circle,draw,inner sep=1pt,very thick] (char) {#1};}}}

\newcommand{\tablearrow}[2]{$\hspace{0,5mm} \bullet \hspace{-1,3mm} \overset{\circled{#1}}{\overset{{\rm #2}}{\longrightarrow}} \hspace{-1,3mm}$}

\newcommand{\tablearrowd}[2]{$\hspace{0,5mm} \bullet \hspace{-1,3mm} \overset{\circledd{#1}}{\overset{{\rm #2}}{\longrightarrow}} \hspace{-1,3mm}$}

\newcommand{\tablearrowdd}[2]{$\hspace{0,5mm} \bullet \hspace{-1,3mm} \overset{\circleddd{#1}}{\overset{{\rm #2}}{\longrightarrow}} \hspace{-1,3mm}$}

\newcommand{\vbar}[0]{\textcolor{gray}{$\, |$}}

\newcommand{\til}[1]{\widetilde{#1}}

\newcommand{\ol}[1]{\overline{#1}}

\definecolor{pinky}{rgb}{1.0, 0, 1.0}

\definecolor{realpurple}{rgb}{1.0, 0, 1.0}
\definecolor{grass}{rgb}{0, 0.56, 0}
\definecolor{yellowi}{HTML}{fba358}

\newcommand{\red}[1]{{\color{red}#1}}
\newcommand{\blue}[1]{{\color{blue}#1}}

\usepackage[labelfont=bf,labelsep=period]{caption}

\begin{document}

\allowdisplaybreaks

\newcommand{\arXivNumber}{2405.14727}

\renewcommand{\PaperNumber}{037}

\FirstPageHeading

\ShortArticleName{Quantized Geodesic Lengths for Teichm\"uller Spaces: Algebraic Aspects}

\ArticleName{Quantized Geodesic Lengths for Teichm\"uller Spaces: \\ Algebraic Aspects}

\Author{Hyun Kyu KIM}

\AuthorNameForHeading{H.K.~Kim}

\Address{School of Mathematics, Korea Institute for Advanced Study (KIAS),\\
85 Hoegi-ro, Dongdaemun-gu, Seoul 02455, South Korea}
\Email{\mail{hkim@kias.re.kr}, \mail{hyunkyu87@gmail.com}}
\URLaddress{\url{https://sites.google.com/site/hyunkyukimmath/}}

\ArticleDates{Received June 09, 2025, in final form March 13, 2026; Published online April 15, 2026}

\Abstract{In 1980's H.~Verlinde suggested to construct and use a quantization of Teich\-m\"uller spaces to construct spaces of conformal blocks for the Liouville conformal field theory. This suggestion led to a mathematical formulation by Fock in 1990's and later by Fock, Goncharov and Shen, called the modular functor conjecture, based on the Chekhov--Fock quantum Teichm\"uller theory. In 2000's, Teschner combined the Chekhov--Fock version and the Kashaev version of quantum Teichm\"uller theory to construct a solution to a modified form of the conjecture. We embark on a direct approach to the conjecture based on the Chekhov--Fock(--Goncharov) theory. We construct quantized trace-of-monodromy along simple loops via Bonahon and Wong's quantum trace maps developed in 2010's, and investigate algebraic structures of them, which will eventually lead to construction and properties of quantized geodesic length operators. We show that a special recursion relation used by Teschner is satisfied by the quantized trace-of-monodromy, and that the quantized trace-of-monodromy for disjoint loops commute in a certain strong sense.}


\Classification{18M20; 57K31; 57K20; 13F60; 81R60; 46L65}

\section{Introduction}
\label{sec:introduction}

\subsection{The modular functor conjecture}

Let $\mathfrak{S}$ be a punctured oriented surface of finite type. One can think of $\mathfrak{S}$ as $\Sigma\setminus\mathcal{P}$, where $\Sigma$ is a~closed oriented surface and $\mathcal{P}$ is a finite subset of $\Sigma$, whose members are referred to as {\it punctures} of $\mathfrak{S} = \Sigma \setminus \mathcal{P}$. Let $\mathfrak{L}$ be a measure space, playing the role of the space of labels. A~version of a (2-dimensional) {\it modular functor} is an assignment
\[
\bigl(\mathfrak{S}, \vec{l}\;\bigr) \mapsto \mathscr{H}_{\mathfrak{S},\vec{l}}
\]
that assigns a Hilbert space \smash{$\mathscr{H}_{\mathfrak{S},\vec{l}}$} to each pair $\bigl(\mathfrak{S}, \vec{l}\;\bigr)$ consisting of a punctured oriented surface $\mathfrak{S} = \Sigma\setminus\mathcal{P}$ and a tuple of labels at punctures \smash{$\vec{l} = (l_p)_{p\in \mathcal{P}} \in \mathfrak{L}^\mathcal{P}$}, such that this assignment satisfies a certain set of axioms. Among such axioms is functoriality, which in particular requires that $\smash{\mathscr{H}_{\mathfrak{S},\vec{l}}}$ should carry the structure of a projective unitary representation of the {\it mapping class group} ${\rm MCG}(\mathfrak{S}) = {\rm Diff}^+(\mathfrak{S})/{\rm Diff}(\mathfrak{S})_0$ of $\mathfrak{S}$, the group of orientation-preserving self-diffeomorphisms of~$\mathfrak{S}$ considered up to isotopy, or more precisely, of the subgroup ${\rm MCG}(\mathfrak{S})_{\vec{l}}$ of ${\rm MCG}(\mathfrak{S})$ preserving~\smash{$\vec{l}$}.

Another pivotal axiom is the compatibility under cutting and gluing of the surfaces. Let $\gamma$~be an unoriented simple loop in $\mathfrak{S}$; let us further assume that it is an {\em essential} loop, i.e., is not null-homotopic. Cutting $\mathfrak{S}$ along $\gamma$ yields a surface $\mathfrak{S}\setminus\gamma$ with two (new) boundary holes. By shrinking these two holes to punctures, we see that $\mathfrak{S}\setminus\gamma$ is diffeomorphic to a new punctured surface $\mathfrak{S}'$ whose set of punctures can be identified as $\mathcal{P} \sqcup\{p_1,p_2\}$, with $p_1$ and $p_2$ being the newly introduced punctures. One denotes an element of \smash{$\mathfrak{L}^{\mathcal{P} \sqcup \{p_1,p_2\}}$} as \smash{$\vec{l} \sqcup\{\chi_1,\chi_2\}$}, which in particular assigns $\chi_i$ to $p_i$. The cutting/gluing axiom states that there must be a natural isomorphism
\begin{align}
\label{eq:intro_decomposition}
\mathscr{H}_{\mathfrak{S},\vec{l}} \ \to \ \int_{\mathfrak{L}}^\oplus \mathscr{H}_{\mathfrak{S}',\vec{l} \sqcup\{\chi, \ol{\chi}\}} \, \mathrm{d}\chi,
\end{align}
where $\ol{\chi}$ denotes the image of $\chi \in \mathfrak{L}$ under a suitable involutive operation on $\mathfrak{L}$. We note that in a traditional formulation, one may require that the Hilbert spaces $\mathscr{H}_{\cdot,\cdot}$ should be finite-dimensional and that $\mathfrak{L}$ should be finite, in which case the direct integral in the right-hand side becomes a direct sum. The adjective `natural' for the decomposition isomorphism means an appropriately defined equivariance under actions of certain algebraic structures, such as the mapping class groups.

A modular functor is a mathematical model for the `spaces of conformal blocks' for surfaces, which are crucial ingredients in the study of 2-dimensional conformal field theories in physics; see \cite{BK, MS, Segal} and also~\cite{T}. In \cite{Verlinde}, H.~Verlinde suggested that for the Liouville conformal field theory, the Hilbert space $\mathscr{H}_{\mathfrak{S},\vec{l}}$ in the modular functor formulation should play the role of the Hilbert space of quantum states for the {\it Teichm\"uller space} of $\mathfrak{S}$. The Teichm\"uller space of $\mathfrak{S}$, which we denote by $\mathscr{T}(\mathfrak{S})$, is the space of all complex structures on $\mathfrak{S}$, or equivalently, the space of all complete hyperbolic metrics on $\mathfrak{S}$, considered up to pullback by self-diffeomorphisms of $\mathfrak{S}$ isotopic to the identity (Definition~\ref{def:Teichmuller_space}). It is known that $\mathscr{T}(\mathfrak{S})$ is a contractible smooth manifold (with corners; see Corollary~\ref{cor:with_corners}), equipped with the Weil--Petersson Poisson structure. The suggestion of Verlinde led to the initial development of the quantum theory of Teichm\"uller spaces in 1990's by Kashaev \cite{Kash98} and independently by Chekhov and Fock~\cite{F97, CF}. The approaches by these two groups are similar but different, and the latter approach later generalized to the quantum theory of cluster varieties by Fock and Goncharov \cite{FG09} (see also \cite{JLSS} for a further development).

A quantization of $\mathscr{T}(\mathfrak{S})$ should associate to each punctured surface $\mathfrak{S}$ a Hilbert space $\mathscr{H}_\mathfrak{S}$ that is a projective unitary representation of the mapping class group ${\rm MCG}(\mathfrak{S})$, together with a~collection of operators on $\mathscr{H}_\mathfrak{S}$ which quantize certain functions on the Teichm\"uller space~$\mathscr{T}(\mathfrak{S})$. As advocated in a recent work of the author \cite{KS}, the task of determining a suitable collection of functions on $\mathscr{T}(\mathfrak{S})$ to be quantized should already be viewed as a crucial part of the quantization problem. In various aspects, one of the most important and natural functions on $\mathscr{T}(\mathfrak{S})$ is the {\em trace-of-monodromy} function $f_\gamma$ (or trace-of-holonomy, or Wilson loop) along a loop $\gamma$ in $\mathfrak{S}$ considered up to isotopy, especially when $\gamma$ is an essential simple loop. It is well known via the Poincar\'e--Koebe uniformization theorem that $\mathscr{T}(\mathfrak{S})$ can be identified with one of the two components of ${\rm Hom}^{{\rm df},\mathfrak{S}}(\pi_1(\mathfrak{S}), {\rm PSL}_2(\mathbb{R}))/{\rm PSL}_2(\mathbb{R})$, the space of all group homomorphisms $\rho\colon \pi_1(\mathfrak{S}) \to {\rm PSL}_2(\mathbb{R})$ satisfying some conditions (Section~\ref{subsec:Teichmuller_space}), considered up to overall conjugation by elements in ${\rm PSL}_2(\mathbb{R})$. The trace-of-monodromy function $f_\gamma$ is defined as
\[
f_\gamma([\rho]) = |{\rm trace}(\rho([\gamma]))|,
\]
when $[\gamma]$ denotes the element of $\pi_1(\mathfrak{S})$ represented by $\gamma$; here we are being somewhat less careful in keeping track of a basepoint in $\mathfrak{S}$. This function is related to another natural function $l_\gamma$ on~$\mathscr{T}(\mathfrak{S})$, the {\em geodesic length function}, by
\[
f_\gamma = 2\cosh(l_\gamma/2) = {\rm e}^{\frac{1}{2} l_\gamma} + {\rm e}^{-\frac{1}{2} l_\gamma}.
\]
For a {\em peripheral} loop $\gamma$ around a puncture $p$, i.e., a small loop around~$p$, we have $l_\gamma=0$, and for a non-peripheral loop $\gamma$, the value of $l_\gamma$ is the hyperbolic length of the unique geodesic loop $\gamma'$ homotopic to $\gamma$.

Verlinde's suggestion on the construction of spaces of conformal blocks using a sought-for quantum Teichm\"uller theory can be interpreted as follows. The label space $\mathfrak{L}$ should be \smash{$(\mathbb{R}_{\ge 0})^\mathcal{P}$}, the set of all tuples \smash{$\vec{l} = (l_p)_{p\in \mathcal{P}}$} of non-negative real numbers $l_p \in \mathbb{R}_{\ge 0}$, where the involution $\chi \mapsto \ol{\chi}$ is identity. We should be able to quantize the trace-of-monodromy functions $f_\gamma$ along simple loops $\gamma$, to a self-adjoint operator \smash{${\bf f}_\gamma^\hbar = {\bf f}_{\gamma;\mathfrak{S}}^\hbar$} on \smash{$\mathscr{H}_{\mathfrak{S},\vec{l}}$}, the {\em quantized trace-of-monodromy operator}, depending on the real quantum parameter $\hbar$. In particular, \smash{${\bf f}^\hbar_\gamma$} must have simple spectrum $[2,\infty)$, so that there is a unique self-adjoint operator \smash{${\bf l}^\hbar_\gamma$}, called the {\em quantized geodesic length operator}, having simple spectrum $[0,\infty)$ and satisfying
\[
{\bf f}^\hbar_\gamma = 2\cosh\bigl({\bf l}^\hbar_\gamma/2\bigr).
\]
We require that \smash{${\bf f}^\hbar_\gamma$} for a peripheral loop $\gamma$ around a puncture $p$ should act as the scalar $2\cosh\bigl(l_p/2\bigr)$. Moreover, the decomposition isomorphism in equation~\eqref{eq:intro_decomposition} should be obtained from the spectral decomposition of \smash{${\bf f}^\hbar_\gamma$} or of \smash{${\bf l}^\hbar_\gamma$}, where the $\chi$-th slice \smash{$\mathscr{H}_{\mathfrak{S}',\vec{l} \sqcup\{\chi, \chi\}}$} corresponds~to the value $\chi$ in the spectrum $[0,\infty)$ of \smash{${\bf l}^\hbar_\gamma$}. For any simple loop $\eta$ in $\mathfrak{S}$ disjoint from $\gamma$, the action of the operator \smash{${\bf f}^\hbar_{\gamma;\mathfrak{S}}$} on \smash{$\mathscr{H}_{\mathfrak{S},\vec{l}}$} should correspond to that of the operator \smash{${\bf f}^\hbar_{\gamma; \mathfrak{S}'}$} on \smash{$\mathscr{H}_{ \mathfrak{S}',\vec{l} \sqcup\{\chi, \chi\}}$} via the decomposition isomorphism in equation~\eqref{eq:intro_decomposition}. Here is a rough summary of necessary ingredients:
\begin{enumerate}[label={\rm (V\arabic*)}]\itemsep=0pt
\item\label{V1} For each punctured surface $\mathfrak{S} = \Sigma\setminus\mathcal{P}$ (with $|\mathcal{P}|\ge 1$, and $\mathfrak{S}$ having negative Euler characteristic) and for each $\vec{l} \in (\mathbb{R}_{\ge 0})^\mathcal{P}$, construct a projective unitary representation
\[
\rho^\hbar_{\mathfrak{S},\vec{l};{\rm MCG}}\colon \ {\rm MCG}(\mathfrak{S})_{\vec{l}} ~\to~ {\rm U}\bigl(\mathscr{H}_{\mathfrak{S},\vec{l}}\bigr) = \bigl\{\mbox{unitary operators on $\mathscr{H}_{\mathfrak{S},\vec{l}}$}\bigr\}.
\]

\item\label{V2} For each simple loop $\gamma$ in $\mathfrak{S}$, construct a self-adjoint operator \smash{${\bf f}^\hbar_\gamma$} on \smash{$\mathscr{H}_{\mathfrak{S},\vec{l}}$} having simple spectrum $[2,\infty)$, depending only on the isotopy class of $\gamma$, quantizing the function $f_\gamma$ in a suitable sense and satisfying the equivariance under the mapping class group action, in the sense that, for each $g \in {\rm MCG}(\mathfrak{S})_{\vec{l}}$,
\[
{\bf f}^\hbar_{g.\gamma} = \rho^\hbar_{\mathfrak{S},\vec{l};{\rm MCG}}(g) \, {\bf f}^\hbar_{\gamma} \, \bigl(\rho^\hbar_{\mathfrak{S},\vec{l};{\rm MCG}}(g)\bigr)^{-1}.
\]

\item\label{V3} For two simple loops $\gamma$ and $\eta$ in $\mathfrak{S}$ that are disjoint, the operators ${\bf f}^\hbar_\gamma$ and ${\bf f}^\hbar_\eta$ {\em strongly} commute with each other, i.e., all spectral projections commute.
\end{enumerate}
The {\em modular functor conjecture} refers to the problem of proving equation~\eqref{eq:intro_decomposition} through the above suggestion \ref{V1}, \ref{V2} and \ref{V3} (see Conjecture \ref{conjecture:sequel}), using some version of quantum Teichm\"uller theory \cite{F97, FG09, Ian, T, Verlinde}; the name of the conjecture seems to have first appeared in \cite{F97}. Recently, Goncharov and Shen formulated a more complete version of the modular functor conjecture \cite{GS19}.

\subsection[Direct approach using Chekhov--Fock--Goncharov quantum Teichm\"uller space and Bonahon--Wong quantum trace]{Direct approach using Chekhov--Fock--Goncharov quantum \\ Teichm\"uller space and Bonahon--Wong quantum trace}

The early results on quantum Teichm\"uller theory by Kashaev \cite{Kash98} and by Chekhov and Fock~\cite{CF}, the latter having developed further later by Fock and Goncharov \cite{FG09}, were solely on \ref{V1}, but did not solve \ref{V2}--\ref{V3}. By combining the combinatorial ingredients of these works, Teschner~\cite{T} built a new quantum Teichm\"uller theory to solve all of \ref{V1}, \ref{V2} and \ref{V3}, claiming to have settled a modified version of the modular functor conjecture. Although this work \cite{T} contains several original ideas, it has some unsatisfactory aspects. First, it is not clear which classical space is quantized by the operators in this construction. Second, the canonicity of the constructed quantized trace-of-monodromy operators ${\bf f}^\hbar_\gamma$ is not sufficiently investigated. Third, the proofs involving operators should be made more careful and rigorous, as the issue of domains can become a subtle matter, e.g., as in \cite{Ruijsenaars}, which is addressed briefly in Section~\ref{subsec:toward_spectral_properties}. Therefore, a more geometric and direct solution to the conjecture, as suggested in~\cite{F97, FG09}, based on the Chekhov--Fock--Goncharov quantum Teichm\"uller theory \cite{F97, CF, FG09}, which has the most straightforward connection to the usual Teicm\"uller space among known quantization schemes, has still been elusive.

At the time of \cite{F97, CF, FG09}, perhaps the biggest hurdle for such a direct approach was that there was no known canonical construction of quantized trace-of-monodromy operator ${\bf f}^\hbar_\gamma$ in full generality, and only partial results were established. The heart of this problem is of algebraic nature, which we briefly describe now. The Chekhov--Fock--Goncharov quantum Teichm\"uller theory is based on special coordinate charts on the Teichm\"uller space $\mathscr{T}(\mathfrak{S})$ developed in 1980's \cite{F93,F97,P12, Thurston, Thurston86} (see Section~\ref{subsec:shear}), each associated to the choice of an {\em ideal triangulation}~$\Delta$ of~$\mathfrak{S}$, which is a collection of {\em ideal} arcs in $\mathfrak{S}$ running between punctures $\mathcal{P}$ that divide $\mathfrak{S}$ into ideal triangles. For any chosen $\Delta$, per each ideal arc $i$ of $\Delta$ is associated a positive real valued coordinate function $X_i$ on $\mathscr{T}(\mathfrak{S})$ (to be more precise, on a certain branched covering space~$\mathscr{T}^+(\mathfrak{S})$), called the exponentiated {\em shear coordinate} function. When one chooses another ideal triangulation $\Delta'$, the corresponding shear coordinates~$X_i'$ are related to the shear coordinates~$X_i$ for~$\Delta$ by certain rational formulas, being examples of the cluster $\mathscr{X}$-mutations in Fock and Goncharov's theory of cluster varieties \cite{F97, FG06, FG09}. It is known that the trace-of-monodromy function $f_\gamma$ is an integer coefficient Laurent polynomial in the square-roots of these functions $\sqrt{X_i}$, $i\in \Delta$, per each $\Delta$. As a quantum deformed algebra for this coordinate system, a special quantum torus algebra $\mathcal{Z}^\omega_\Delta$ called the (square-root) {\it Chekhov--Fock algebra} is suggested \cite{BL07, BW, FG09, Hiatt, Liu}, which is an associative algebra generated by $Z_i^{\pm 1}$, $i \in \Delta$, with the relations
\[
Z_i Z_j = \omega^{2\varepsilon_{ij}} Z_j Z_i, \qquad \forall i, j \in \Delta,
\]
where $\omega = {\rm e}^{\pi {\rm i} \hbar/4}$ is a quantum parameter, and $\varepsilon_{ij} \in \{-2,-1,0,1,2\}$ counts the number of appearances of the arcs $i$ and $j$ in a same ideal triangle of $\Delta$, with sign (Definition~\ref{def:varepsilon_Delta}). As~$\omega\to1$, each variable $Z_i$ should recover the function $\sqrt{X_i}$. Per each pair of ideal triangulations~$\Delta$ and~$\Delta'$, a quantum coordinate change isomorphism is constructed (Proposition~\ref{prop:Theta})
\[
\Theta^\omega_{\Delta,\Delta'} \colon \ {\rm Frac}( (\mathcal{Z}^\omega_{\Delta'})_{\rm bl} ) \to {\rm Frac}\bigl( \bigl(\mathcal{Z}^\omega_\Delta\bigr)_{\rm bl}\bigr)
\]
between the skew fields of fractions of certain `balanced' subalgebras (Definition~\ref{def:balanced_algebra}) of the Chekhov--Fock algebras \cite{BL07, F97, CF, FG09, Hiatt, KLS}. This recovers the classical coordinate change as $\omega\to 1$ and satisfies the consistency equations $\Theta^\omega_{\Delta,\Delta'} \circ \Theta^\omega_{\Delta',\Delta''} = \Theta^\omega_{\Delta,\Delta''}$.

The problem of algebraically quantizing the trace-of-monodromy function $f_\gamma$ along a simple loop $\gamma$ then aims to construct an element $f^\omega_{\gamma,\Delta} \in \bigl(\mathcal{Z}^\omega_\Delta\bigr)_{\rm bl} \subset \mathcal{Z}^\omega_\Delta$, a suitable (noncommutative) quantum Laurent polynomial expression recovering $f_{\gamma}$ as $\omega\to 1$, per each ideal triangulation~$\Delta$, such that these elements for different $\Delta$ are compatible with each other in the sense that
\[
f^\omega_{\gamma,\Delta} = \Theta^\omega_{\Delta,\Delta'} \bigl(f^\omega_{\gamma,\Delta'}\bigr).
\]
This problem, sometimes referred to as the quantum ordering problem \cite{CF2} (but is not completely solved there), is solved (say, in \cite{AK}) by using Bonahon and Wong's {\em quantum trace maps} \cite{BW, Le19} (see Proposition~\ref{def:skein_algebra})
\[
{\rm Tr}^\omega_\Delta = {\rm Tr}^\omega_{\Delta;\mathfrak{S}} \colon \ \mathcal{S}^\omega(\mathfrak{S}) \to \bigl(\mathcal{Z}^\omega_\Delta\bigr)_{\rm bl} \subset \mathcal{Z}^\omega_\Delta
\]
associated to each ideal triangulation $\Delta$ of a punctured surface $\mathfrak{S}$, where $\mathcal{S}^\omega(\mathfrak{S})$ is the {\em $($Kauffman bracket$)$ skein algebra} of $\mathfrak{S}$ \cite{Pr,Turaev} (see also \cite{Muller}), which is generated by isotopy classes of framed links living in the 3-manifold $\mathfrak{S} \times (-1,1)$, modulo the {\em skein} relations (Definition~\ref{def:skein_algebra}). The~most crucial properties are that, first, the naturality holds in the sense of ${\rm Tr}^\omega_\Delta = \Theta^\omega_{\Delta,\Delta'} \circ {\rm Tr}^\omega_{\Delta'}$, and second, for a simple loop $\gamma$, if we denote by $K_\gamma$ the lift of $\gamma$ at a `constant elevation' in $\mathfrak{S} \times \{0\}$ equipped with a constant upward vertical framing, then ${\rm Tr}^\omega_\Delta([K_\gamma])$ recovers $f_\gamma$ as $\omega\to 1$. Thus, we are led to the following natural solution to the quantum ordering problem
\[
f^\omega_{\gamma,\Delta} := {\rm Tr}^\omega_\Delta([K_\gamma]) \in \mathcal{Z}^\omega_\Delta.
\]
This element is what we call the {\em $($algebraic$)$ quantized trace-of-monodromy} for $\gamma$ with respect to the ideal triangulation $\Delta$ (Definition~\ref{def:quantized_trace-of-monodromy}). See the recent work of Mandel and Qin \cite{MQ} which adds to the canonicity of this $f^\omega_{\gamma,\Delta}$.

However, this only deals with the algebraic aspect of the quantization, and one needs to represent these elements as operators on a Hilbert space. In \cite{ F97,CF, FG09, Kim_irreducible}, a representation $\rho^\hbar_\Delta$ of the Chekhov--Fock algebra $\mathcal{Z}^\omega_\Delta$ on a Hilbert space $\mathscr{H}_\Delta = L^2(\mathbb{R}^n)$ (for some $n$) is constructed, and moreover, per each pair of ideal triangulations $\Delta,\Delta'$ a unitary operator ${\bf K}^\hbar_{\Delta,\Delta'} \colon \mathscr{H}_{\Delta'} \to \mathscr{H}_\Delta$ that `intertwines' these representations in a suitable consistent sense is constructed. Applying this representation, one would obtain an operator $\rho^\hbar_\Delta\bigl(f^\omega_{\gamma,\Delta}\bigr)$ on the Hilbert space $\mathscr{H}_\Delta$, which is supposed to be the operator that quantizes the function $f_\gamma$. However, the operators assigned by the representation $\rho^\hbar_\Delta$ are not defined on the whole Hilbert space $\mathscr{H}_\Delta$ but only on a dense subspace, notably the {\em Schwartz space} $\mathscr{S}_\Delta$ of Fock and Goncharov \cite[Definition~5.2]{FG09}, which is roughly the common maximal domain of all operators $\rho^\hbar_\Delta\bigl(f^\omega_{\gamma,\Delta}\bigr)$ of our concern (see equation~\eqref{eq:Schwartz}).

In order to use the known results above to attack \ref{V2} and \ref{V3}, which are steps toward the modular functor conjecture, one should first show that
\begin{align}
\label{eq:intro_essential_uniqueness}
\mbox{$\rho^\hbar_\Delta\bigl(f^\omega_{\gamma,\Delta}\bigr)$ has a unique (or at least, canonical) self-adjoint extension in $\mathscr{H}_\Delta$,}
\end{align}
call it the quantized trace-of-monodromy operator ${\bf f}^\hbar_{\gamma,\Delta}$, then show that ${\bf f}^\hbar_{\gamma,\Delta}$ has simple spectrum $[2,\infty)$, and finally show the strong commutativity for two disjoint loops as in \ref{V3}.

The major goal of the present paper is to embark on the study of the above problem of studying the analytic properties of the operators $\rho^\hbar_\Delta\bigl(f^\omega_{\gamma,\Delta}\bigr)$, which seem to have never been undertaken in a rigorous manner in the literature. In the present paper, we will collect algebraic preliminaries, i.e., study algebraic properties of the quantized trace-of-monodromy $f^\omega_{\gamma,\Delta} \in \mathcal{Z}^\omega_\Delta$, which we will use in the sequel paper \cite{sequel} to prove the desired analytic properties of the corresponding operators $\rho^\hbar_\Delta\bigl(f^\omega_{\gamma,\Delta}\bigr)$; we give a rough outline of the latter task in Section~\ref{sec:consequences}.

To emphasize the value of the program set out by the current paper and the sequel \cite{sequel} once more, on the one hand, the algebraic quantized trace-of-monodromy $f^\omega_{\gamma,\Delta}$ that fits the Chekhov--Fock--Goncharov `cluster' theoretic formulation satisfying the relevant naturality under change of $\Delta$ and whose classical geometric meaning is clear was not known at the early stage of quantum Teichm\"uller theory (\cite{F97, CF}, not even in \cite{FG09}). It follows only from the work of Bonahon and Wong \cite{BW}, but the analytic properties of the corresponding operator \smash{$\rho^\hbar_\Delta\bigl(f^\omega_{\gamma,\Delta}\bigr)$} is yet to be studied, in order to obtain the sought-for operator \smash{${\bf f}^\hbar_{\gamma,\Delta}$} or \smash{${\bf f}^\hbar_\gamma$}. On the other hand, Teschner \cite{T} has suggested a totally different system of operators which would play the role of \smash{${\bf f}^\hbar_\gamma$}. However, the classical geometric meaning of these operators or their relevance to cluster theoretic formulation is not clear, and the analytic proofs in \cite{T} seem to require a more careful justification.

\subsection{Summary of results}

To efficiently describe our results, let us introduce a notation first. For an ideal triangulation $\Delta$, define $N_\Delta$ to be the free $\mathbb{Z}$-module generated by $\{ v_i \mid i \in \Delta\}$, equipped with a skew symmetric $\mathbb{Z}$-bilinear form $\langle \cdot, \cdot \rangle$ such that $\langle v_i, v_j\rangle = \varepsilon_{ij}$ (Definition~\ref{def:N_Delta}). For each $v \in N_\Delta$, define the corresponding {\em Weyl-ordered Laurent monomial} $Z_v \in \mathcal{Z}^\omega_\Delta$ as follows: write $v = \sum_{j=1}^r a_j v_{i_j}$ for some sequence of ideal arcs $i_1,i_2,\dots,i_r$ and integers $a_1,\dots,a_r$, then
\[
Z_v := \omega^{- \sum_{j<k} a_j a_k \varepsilon_{i_j \, i_k} } Z_{i_1}^{a_1} Z_{i_2}^{a_2} \cdots Z_{i_r}^{a_r},
\]
which can be thought of as a standard multiplicative normalization of $Z_{i_1}^{a_1} Z_{i_2}^{a_2} \cdots Z_{i_r}^{a_r}$. It is straightforward to verify that
\begin{gather*}
Z_0 = 1, \qquad
Z_{v_i} = Z_i, \quad \forall i \in \Delta, \qquad
Z_{-v} = (Z_v)^{-1}, \\
Z_v Z_w = \omega^{\langle v,w\rangle} Z_{v+w}, \quad \forall v,w\in N_\Delta.
\end{gather*}

Our results can be largely divided into two. One is on the structure of the quantized trace-of-monodromy $f^\omega_{\gamma,\Delta} \in \mathcal{Z}^\omega_\Delta$, for each single essential simple loop $\gamma$ in $\mathfrak{S}$. Per each $\gamma$, we find an~ideal triangulation $\Delta$ such that $f^\omega_{\gamma,\Delta}$ has a structure that is simple enough, or more precisely, suitable for our purpose mentioned in the previous subsection. The following is our first main theorem.
\begin{Theorem}[algebraic structure of quantized trace-of-monodromy along a single loop; Theorem~\ref{thm:AS}]
\label{thm:main1}
Let $\mathfrak{S}$ be a punctured surface admitting an ideal triangulation. Let $\gamma$ be an essential simple loop in $\mathfrak{S}$. Then the following holds:
\begin{enumerate}[label={\rm (\arabic*)}]\itemsep=0pt
\item\label{thm1_1_item1} If $\gamma$ is a peripheral loop around a puncture, then for any ideal triangulation $\Delta$ one has
\[
f^\omega_{\gamma,\Delta} = Z_v + Z_{-v}
\]
for some $v\in N_\Delta$.

\item\label{thm1_1_item2} If $\gamma$ is not peripheral and there exists an ideal triangulation $\Delta$ such that $\gamma$ meets $\Delta$ exactly at two points, one on $a\in \Delta$ and another on $b\in \Delta$, then
\[
f^\omega_{\gamma,\Delta} = Z_{v_a+v_b} +Z_{-v_a-v_b} + Z_{\epsilon(v_a-v_b)}
\]
for some $\epsilon \in \{+,-\}$.

\item\label{thm1_1_item3} If none of the items {\rm \ref{thm1_1_item1}} and {\rm \ref{thm1_1_item2}} holds, then $\gamma$ is a separating loop that cuts out a subsurface of genus $g\ge 1$ with one boundary component and no puncture, and there exists an ideal triangulation $\Delta$ of $\mathfrak{S}$, essential non-peripheral simple loops $\gamma_1$ and $\gamma_2$, and elements $v_1$ and $v_2$ of $N_\Delta$ satisfying $\langle v_1,v_2\rangle =\pm 4$, such that the three curves $\gamma,\gamma_1$ and $\gamma_2$ cut out a~three-holed sphere, and the {\em Teschner recursion relation} $($Definition~{\rm \ref{def:Teschner_triple_level1}--\ref{def:Teschner_triple_level2})} holds for the ordered triple $(\gamma,\gamma_1,\gamma_2)$:
\[
f^\omega_{\gamma,\Delta} = Z_{ v_1 - v_2} + Z_{v_2 - v_1} + Z_{v_1 + v_2} + Z_{v_1} f^\omega_{\gamma_1,\Delta} + Z_{v_2} f^\omega_{\gamma_2,\Delta},
\]
where each $Z_{v_j}$ commutes with each $f^\omega_{\gamma_k,\Delta}$ in the following strong sense: for each $j=1,2$ and each Laurent monomial term $Z_v$ appearing in $f^\omega_{\gamma_k,\Delta}$, $k=1,2$, one has $\langle v_j, v\rangle=0$.

\end{enumerate}

\end{Theorem}
In the item \ref{thm1_1_item3}, we also show that $f^\omega_{\gamma_1,\Delta}$ and $f^\omega_{\gamma_2,\Delta}$ commute with each other in a certain sense that is stronger than the mere algebraic commutativity $f^\omega_{\gamma_1,\Delta}f^\omega_{\gamma_2,\Delta}=f^\omega_{\gamma_2,\Delta}f^\omega_{\gamma_1,\Delta}$. We remark that the above notion of the Teschner recursion relation is inspired by a recursion formula \cite[equation~(15.5) of Definition~6]{T} used by Teschner in his recursive construction of operators for loops. Teschner used this recursion formula to argue the sought-for analytic properties \ref{V2} of the operators, and we expect that our version of the Teschner recursion should also be used to prove the analytic properties of the operators $\rho^\hbar_\Delta\bigl(f^\omega_{\gamma,\Delta}\bigr)$; namely, equation~\eqref{eq:intro_essential_uniqueness} and \ref{V2} for~$\gamma_1$ and~$\gamma_2$ should be used to prove those properties for $\gamma$.
In the sequel \cite{sequel}, which is briefly previewed in Section~\ref{sec:consequences}, we deal with these analytic statements in a rigorous manner.

In the meantime, in \cite{sequel} and Section~\ref{sec:consequences}, we also justify the usefulness of the statements of~the items \ref{thm1_1_item1}--\ref{thm1_1_item2} for proofs of analytic properties. For \ref{thm1_1_item1}, the operator \smash{$\rho^\hbar_\Delta\bigl(f^\omega_{\gamma,\Delta}\bigr)$} is unitarily equivalent to the multiplication operator ${\rm e}^x + {\rm e}^{-x}$ on a dense subspace of $L^2(\mathbb{R},{\rm d}x)$, tensored with identity on $L^2\bigl(\mathbb{R}^{n-1}\bigr)$. This operator is well known and is straightforward to study. For~\ref{thm1_1_item2}, \smash{$\rho^\hbar_\Delta\bigl(f^\omega_{\gamma,\Delta}\bigr)$} is unitarily equivalent to the operator \smash{${\rm e}^{8 \pi {\rm i} \hbar \frac{{\rm d}}{{\rm d}x}} + {\rm e}^{-8 \pi {\rm i} \hbar \frac{{\rm d}}{{\rm d}x}} + {\rm e}^x$} on a dense subspace of $L^2(\mathbb{R},{\rm d}x)$ (like equation~\eqref{eq:our_Hermite_subspace}), tensored with identity; on $L^2(\mathbb{R},{\rm d}x)$ it acts as $\varphi(x) \mapsto \varphi(x+8\pi {\rm i} \hbar) + \varphi(x-8\pi {\rm i} \hbar) + {\rm e}^x \varphi(x)$ for nice $\varphi$. This operator is taken to be the standard model of the operator ${\bf f}^\omega_{\gamma,\Delta}$ by Kashaev already in~\cite{Kash00}, and its analytic properties have been rigorously studied in the work of Takhtajan and Faddeev~\cite{FT}, which implies the properties we want.

Combining, a rough plan for analytic proofs in \cite{sequel} is as follows. The loops falling to~\ref{thm1_1_item1} or~\ref{thm1_1_item2} are manageable as explained above. For other loops $\gamma$, in the Teschner recursion \ref{thm1_1_item3}, the `complexities' of $\gamma_1$ and $\gamma_2$ can be viewed as being less than that of $\gamma$, so that the Teschner recursion lets us reduce the complexity recursively, until we reach the class of loops falling into the item \ref{thm1_1_item2} whose complexity is manageable with regard to the analytic proofs.

In order to show Theorem~\ref{thm:main1}, we first establish a certain topological classification (Proposition~\ref{prop:loop_classification}) of all essential simple loops $\gamma$ in $\mathfrak{S}$, in relation to their relative positions with respect to ideal triangulations. Using this classification, for each $\gamma$, we find a suitable ideal triangulation $\Delta$ and directly compute the quantized trace-of-monodromy $f^\omega_{\gamma,\Delta}$, and verify Theorem~\ref{thm:main1}. Along the way, we develop some tools to use in such computations (Lemma~\ref{lem:AP}, Corollary~\ref{cor:AP}), which has the following useful corollary.
\begin{Proposition}[Proposition~\ref{prop:term-by-term_Weyl-ordered_f_omega_gamma_Delta}]
\label{prop:intro_term-by-term}
Let $\mathfrak{S}$ be a punctured surface admitting an ideal triangulation. Let $\gamma$ be an arbitrary essential simple loop in $\mathfrak{S}$. Then there exists an ideal triangulation~$\Delta$ of $\mathfrak{S}$ such that the quantized trace-of-monodromy can be written as
\[
f^\omega_{\gamma,\Delta} = \sum_{v\in N_\Delta'} Z_v,
\]
for some finite subset $N_\Delta'$ of $N_\Delta$. In particular, $f^\omega_{\gamma,\Delta}$ is a sum of some mutually distinct Weyl-ordered Laurent monomials $Z_v$, and therefore is termwise Weyl-ordered.
\end{Proposition}
We note that the termwise Weyl-ordered expression has been the first attempt of constructing a quantized trace-of-monodromy \cite{F97, CF2} before the work of Bonahon and Wong \cite{BW}. Only by this Proposition~\ref{prop:intro_term-by-term} (together with \cite{BW}) is this first attempt fully justified (see Section~\ref{subsec:quantized_trace-of-monodromy})!

Our second result is on a certain algebraic version of strong commutativity of quantized trace-of-monodromy $f^\omega_{\xi_1,\Delta}$ and $f^\omega_{\xi_2,\Delta}$ for any pair of disjoint essential simple loops $\xi_1,\xi_2$. By~construction and from the basic property of Bonahon and Wong's quantum trace \cite{BW}, we already know that these two commute with each other, as elements of the algebra $\mathcal{Z}^\omega_\Delta$. However, in view of \ref{V3}, we seek for a stronger notion of commutativity that would imply the strong commutativity of corresponding operators assigned by representations. We develop one such notion in Definitions~\ref{def:algebraic_SC1}--\ref{def:algebraic_SC2}; we say that elements $X$ and $Y$ of $\mathcal{Z}^\omega_\Delta$ {\em algebraically strongly commute} if there exist $\mathbb{Z}$-submodules $N_1$ and $N_2$ of the lattice $N_\Delta$ such that
\begin{enumerate}[label={\rm (\arabic*)}]\itemsep=0pt
\item $X = \sum_{v\in N_1} a_v Z_v$ and $Y = \sum_{v\in N_2} b_v Z_v$ for some $a_v,b_v \in \mathbb{Z}\bigl[\omega^{\pm 1}\bigr]$ that are zero for all but finitely many $v$,

\item $N_1 \perp N_2$, in the sense that $\langle v,w\rangle=0$ for all $v\in N_1$ and $w\in N_2$, and

\item the $\mathbb{R}$-spans $(N_j)_\mathbb{R} = N_j \otimes_\mathbb{Z} \mathbb{R} \subset N_\Delta \otimes_\mathbb{Z} \mathbb{R} = (N_\Delta)_\mathbb{R}$, equipped with the $\mathbb{R}$-bilinearly extended form $\langle \cdot,\cdot\rangle$, satisfy the condition $(N_1)_\mathbb{R} \cap (N_2)_\mathbb{R} \subset (N_\Delta)_\mathbb{R}^\perp$, where $(N_\Delta)_\mathbb{R}^\perp = \bigl\{ v\in (N_\Delta)_\mathbb{R} \mid \langle v,w\rangle=0, ~ \forall w\in (N_\Delta)_\mathbb{R}\bigr\}$.
\end{enumerate}
We will justify in the sequel \cite{sequel} that this condition indeed implies the strong commutativity of corresponding operators. However, it turns out that for some pairs of loops $\xi_1,\xi_2$, one cannot find an ideal triangulation $\Delta$ such that $f^\omega_{\xi_1,\Delta}$ and $f^\omega_{\xi_2,\Delta}$ algebraically strongly commute as above. In fact, for such pairs of loops, we employ the notion of Teschner recursion relation again. This is a~sensible idea because it is expected that the Teschner recursion relation for a triple of loops $(\gamma,\gamma_1,\gamma_2)$ with respect to the ideal triangulation $\Delta$ can be used to prove that the self-adjoint extension of \smash{$\rho^\hbar_\Delta\bigl(f^\omega_{\gamma,\Delta}\bigr)$} \big(i.e., \smash{${\bf f}^\hbar_{\gamma,\Delta}$}\big) strongly commutes with the self-adjoint extension of \smash{$\rho^\hbar_\Delta\bigl(f^\omega_{\gamma_j,\Delta}\bigr)$} \big(i.e., \smash{${\bf f}^\hbar_{\gamma_j,\Delta}$}\big), for $j=1,2$. Proof of this strong commutativity of operators is argued by Teschner~\cite{T} in his setting, and will be rigorously established in~\cite{sequel} for our setting, which is previewed in Section~\ref{sec:consequences}. For the present paper, we focus just on the algebraic statements. Our second main theorem is as follows.
\begin{Theorem}[algebraic strong commutativity of quantized trace-of-monodromy for disjoint loops; Theorem~\ref{thm:algebraic_commutativity}]
\label{thm:intro_second_main}
Let $\mathfrak{S}$ be a punctured surface admitting an ideal triangulation. Let~$\xi_1$ and~$\xi_2$ be any disjoint and mutually non-isotopic essential simple loops in~$\mathfrak{S}$. Then,
\begin{enumerate}[label={\rm (\arabic*)}]\itemsep=0pt
\item if one of $\xi_1$ and $\xi_2$ is a peripheral loop around a puncture, then $f^\omega_{\xi_1,\Delta}$ and $f^\omega_{\xi_2,\Delta}$ algebraically strongly commute, for any ideal triangulation~$\Delta$ of~$\mathfrak{S}$;

\item if there is an ideal triangulation $\Delta$ of $\mathfrak{S}$ such that the set of ideal triangles intersecting~$\xi_1$ and that for $\xi_2$ do not have a common ideal triangle, then $f^\omega_{\xi_1,\Delta}$ and $f^\omega_{\xi_2,\Delta}$ algebraically strongly commute;

\item if none of the above holds, then there exists an ideal triangulation $\Delta$ of $\mathfrak{S}$ and an essential simple loop $\xi_3$ in $\mathfrak{S}$ such that one of the triples $(\xi_1,\xi_2,\xi_3)$ and $(\xi_2,\xi_1,\xi_3)$ satisfies the Teschner recursion relation with respect to $\Delta$.
\end{enumerate}
\end{Theorem}
See Theorem~\ref{thm:algebraic_commutativity} for a more detailed statement. To prove this theorem, we establish a certain classification statement for pairs of disjoint essential simple loops (Proposition~\ref{prop:classification_of_pair_of_loops}), in relation to relative positions with respect to ideal triangulations. We then prove by direct computation that some triples of loops do satisfy the Teschner recursion relation (Section~\ref{sec:on_commutativity}).

In Section~\ref{sec:consequences}, we present a preview of the sequel \cite{sequel} about representations and some analytic arguments toward a proof of \ref{V2}--\ref{V3}, our partial suggestion for how to use the results of the present paper and the sequel \cite{sequel} to actually attack the modular functor conjecture, as well as some more future research directions. We note in particular that \ref{V2}--\ref{V3} can be proved for each genus zero punctured surface essentially by using the results of the current paper.

We hope that, even without the implications on the analytic properties of the operators, the results of the present paper will be useful for the understanding of the (algebraic) quantized trace-of-monodromy $f^\omega_{\gamma,\Delta}$ and of the Bonahon--Wong quantum trace ${\rm Tr}^\omega_\Delta$ in general.

We are informed that after the first version of the present article was put online, Schrader and Shapiro \cite{SS25} proved a certain algebraic version of the modular functor conjecture for quantum higher Teichm\"uller theory formulated by Goncharov and Shen \cite{GS19}. It would be interesting to investigate the relationship to our results in the case of quantum Teichm\"uller theory.

\section{Quantization of Teichm\"uller spaces}

In the present section, we establish the basic setting, as a preparation for later sections. We~review known constructions, based on which we suggest a quantized trace-of-monodromy along a simple loop (Definition~\ref{def:quantized_trace-of-monodromy}), which is an algebraic quantization of the trace-of-monodromy function. Using it, in Section~\ref{subsec:on_quantized_length_operators}, we organize an approach to the modular functor conjecture, and formulate a goal to seek for (Conjecture \ref{conjecture:sequel}), involving the suggested quantized geodesic length operator. We recall a known tool for computation of quantized trace-of-monodromy, and state one useful new result (Proposition~\ref{prop:term-by-term_Weyl-ordered_f_omega_gamma_Delta}).

\subsection{Teichm\"uller space, trace-of-monodromy and length functions}
\label{subsec:Teichmuller_space}

We first fix notations about surfaces. Our main motivation is about punctured surfaces obtained by removing finitely many points from a closed oriented surface. But we also deal with more general type of surfaces as follows (see \cite{FG06}, and, e.g., also \cite{Le19}). Most of the contents of the current subsection and the next one are review materials and can be found for example in \cite{P12}.
\begin{Definition}
\label{def:generalized_marked_surface}
A {\em marked surface} is a pair $(\Sigma,\mathcal{P})$, where $\Sigma$ is a compact oriented smooth real 2-dimensional manifold $\Sigma$ with possibly empty boundary $\partial \Sigma$, and $\mathcal{P}$ is a {\it nonempty} finite subset of $\Sigma$ such that each boundary component of~$\Sigma$ meets~$\mathcal{P}$.
\begin{itemize}\itemsep=0pt
\item An element of $\mathcal{P}$ is called a {\em marked point}. An element of $\mathcal{P}$ that lies in the interior $\mathring{\Sigma} = \Sigma\setminus\partial \Sigma$ of $\Sigma$ is called a {\em puncture}. We say $(\Sigma,\mathcal{P})$ is a {\em punctured surface} if $\partial \Sigma = \varnothing $.

\item We identify $(\Sigma,\mathcal{P})$ with the surface
\[
\mathfrak{S} = \Sigma\setminus \mathcal{P}
\]
whenever convenient; the marked points are then thought of as `ideal' points of $\mathfrak{S}$.

\item An isomorphism between marked surfaces $(\Sigma,\mathcal{P})$ and $(\Sigma',\mathcal{P}')$ is an orientation-preserving diffeomorphism $\Sigma \to \Sigma'$ inducing a bijection $\mathcal{P} \to \mathcal{P}'$. We say two marked surfaces are isomorphic to each other if there exists an isomorphism between them.
\end{itemize}
\end{Definition}

A basic classical object of our interest is the following.
\begin{Definition}
\label{def:Teichmuller_space}
Let $\mathfrak{S} = \Sigma \setminus \mathcal{P}$ be a punctured surface. When $\Sigma$ is of genus $0$, we assume $|\mathcal{P}|\ge 3$.

The {\em Teichm\"uller space} $\mathscr{T}(\mathfrak{S})$ of $\mathfrak{S}$ is the set of all complete hyperbolic metrics $g$ on $\mathfrak{S}$, considered up to isotopy, i.e., up to pullback by self-diffeomorphisms of $\mathfrak{S}$ that are isotopic to the identity diffeomorphism. Denote by $[g]$ the element of $\mathscr{T}(\mathfrak{S})$ represented by $g$.
\end{Definition}
It is well known that $\mathscr{T}(\mathfrak{S})$ is in bijection with one of the two components of
\begin{align}
\label{eq:Hom_space}
{\rm Hom}^{{\rm df},\mathfrak{S}}(\pi_1(\mathfrak{S}), {\rm PSL}_2(\mathbb{R}))/{\rm PSL}_2(\mathbb{R}),
\end{align}
which is the moduli space of all group homomorphisms $\rho : \pi_1(\mathfrak{S}) \to {\rm PSL}_2(\mathbb{R})$ that are faithful (i.e., injective) and have discrete images such that the quotient $\mathbb{H}^2 /\pi_1(\mathfrak{S})$ of the upper half plane is diffeomorphic to $\mathfrak{S}$. Two such homomorphisms $\rho$ and $\rho'$ are identified if $\exists \sigma \in {\rm PSL}_2(\mathbb{R})$ such that $\forall [\gamma] \in \pi_1(\mathfrak{S})$, $\rho'([\gamma]) = \sigma \, \rho([\gamma]) \, \sigma^{-1}$.

Given a point $[g] \in \mathscr{T}(\mathfrak{S})$, let $[\rho]$ be the corresponding point of the space in equation~\eqref{eq:Hom_space}. We say that $\rho$ is a {\em monodromy representation} of $[g]$. For an oriented loop $\gamma$ in $\mathfrak{S}$ representing an element $[\gamma] \in \pi_1(\mathfrak{S})$, the matrix $\rho([\gamma]) \in {\rm PSL}_2(\mathbb{R})$ which is determined by $[\rho]$ or $[g]$ up to conjugation, is referred to as the {\em monodromy} along $\gamma$ with respect to $[\rho]$ or $[g]$. Note that
\begin{align*}
f_{[\gamma]} \colon \ \mathscr{T}(\mathfrak{S}) \to \mathbb{R}, \qquad [\rho] \mapsto |{\rm trace}(\rho([\gamma]))|
\end{align*}
is a well-defined function on $\mathscr{T}(\mathfrak{S})$ which depends on $[\gamma]$, but not on the orientation of $\gamma$. We call~$f_{[\gamma]}$ the {\em trace-of-monodromy} function along an unoriented loop $\gamma$ or its homotopy class $[\gamma]$. This function plays a crucial role in the study of algebraic structures of $\mathscr{T}(\mathfrak{S})$, and particularly in the present paper. We remark that, throughout the paper we mostly work with unoriented loops; we may equip them with arbitrary orientations when necessary, without explicitly saying~so.

For the free homotopy class $[\gamma]$ of any (oriented) loop $\gamma$ in $\mathfrak{S}$, the infimum of the hyperbolic length of a curve in the class $[\gamma]$ yields a well-defined function on $\mathscr{T}(\mathfrak{S})$, which we denote by
\[
l_{[\gamma]} \colon \ \mathscr{T}(\mathfrak{S}) \to \mathbb{R}_{\ge 0}
\]
and call the {\em length function} along $[\gamma]$ or along $\gamma$. It is known that ${\rm e}^{l_{[\gamma]}/2} + {\rm e}^{-l_{[\gamma]}/2}$ coincides with $|{\rm trace}(\rho([\gamma]))|$, where the latter is a function written in terms of the monodromy representations~$[\rho]$. Namely, the trace-of-monodromy function $f_{[\gamma]}$ and the length function $l_{[\gamma]}$ are related to each other by the equation
\begin{align}
\label{eq:trace_of_monodromy_as_length}
f_{[\gamma]} = 2\cosh\bigl(l_{[\gamma]}/2\bigr).
\end{align}
In particular, the values of $f_{[\gamma]}$ lie in the interval $[2,\infty)$, and we have $l_{[\gamma]} = 2 \cosh^{-1} \bigl( f_{[\gamma]}/2\bigr)$ when we consider $\cosh$ as a function from $[0,\infty)$ to $[1,\infty)$. Hence, at the classical level, dealing with $f_{[\gamma]}$ is equivalent to dealing with $l_{[\gamma]}$.

\begin{Definition}\label{def:essential_and_peripheral}
Let $\gamma$ be a simple loop in a punctured surface $\mathfrak{S} = \Sigma\setminus\mathcal{P}$ that is {\em essential}, i.e., is not null-homotopic. We say $\gamma$ is a {\em peripheral} loop around a puncture $p \in \mathcal{P}$, if $\gamma$ is contractible to $p$ in $\Sigma$, or equivalently, if $\gamma$ is isotopic to a small loop surrounding $p$.
\end{Definition}
For each puncture $p \in \mathcal{P}$, let $\gamma_p$ be a peripheral loop around $p$. We say that $p$ is a {\em cusp} with respect to $[g] \in \mathscr{T}(\mathfrak{S})$ (or to $[\rho]$) if $l_{[\gamma_p]}([g])=0$ (or $f_{[\gamma_p]}([\rho])=2$), and a {\em funnel} if $l_{[\gamma_p]}([g])>0$ (or $f_{[\gamma_p]}([\rho])>2$).

We believe that the space that one would like to quantize to realize the suggestion of Verlinde \cite{Verlinde} and to settle the modular functor conjecture \cite{F97,FG09} is the following subspace of the Teichm\"uller space, with the prescribed lengths of peripheral loops.
\begin{Definition}
\label{def:relative_Teichmuller_space}
Let $\mathfrak{S} = \Sigma\setminus \mathcal{P}$ be a punctured surface as in Definition~\ref{def:Teichmuller_space}. For any $\vec{l} = (l_p)_{p \in \mathcal{P}} \in (\mathbb{R}_{\ge 0})^\mathcal{P}$, we define the {\em relative $($or slice$)$ Teichm\"uller space} $\mathscr{T}(\mathfrak{S})_{\vec{l}}$ as
\begin{align*}
\mathscr{T}(\mathfrak{S})_{\vec{l}} = \bigl\{[g] \in \mathscr{T}(\mathfrak{S}) \mid l_{[\gamma_p]} = l_p, ~ \forall p \in \mathcal{P} \bigr\}.
\end{align*}
\end{Definition}
We postpone the discussion on geometric structures of $\mathscr{T}(\mathfrak{S})$ and $\mathscr{T}(\mathfrak{S})_{\vec{l}}$ to later subsections.

\subsection{Thurston's shear coordinate systems}
\label{subsec:shear}

We recall another version of the Teichm\"uller space that appears in the known results on quantization of Teichm\"uller spaces.
\begin{Definition}[\cite{F97}]
\label{def:enhanced_Teichmuller_space}
Let $\mathfrak{S} = \Sigma\setminus \mathcal{P}$ be a punctured surface; when $\Sigma$ is of genus $0$, we assume $|\mathcal{P}|\ge 3$. Let $\mathscr{T}(\mathfrak{S})$ be the Teichm\"uller space of $\mathfrak{S}$ defined in Definition~\ref{def:Teichmuller_space}. For $[g] \in \mathscr{T}(\mathfrak{S})$, a~{\em framing} $\beta$ on $[g]$ is the choice of a sign $\beta(p) \in \{+,-\}$ per each puncture $p\in \mathcal{P}$ that is a funnel with respect to $[g]$.

Define the {\em enhanced Teichm\"uller space} $\mathscr{T}^+(\mathfrak{S})$ as the space parametrizing all pairs $([g], \beta)$ consisting of a point $[g]$ of $\mathscr{T}(\mathfrak{S})$ and a framing $\beta$ on $[g]$.
\end{Definition}
\begin{Remark}
See \cite{A16} for the definition of $\mathscr{T}^+(\mathfrak{S})$ for a marked surface.
\end{Remark}
So, naturally we have the framing-forgetting map
\begin{align}
\label{eq:framing-forgetting_map}
\mathscr{T}^+(\mathfrak{S}) \to \mathscr{T}(\mathfrak{S})
\end{align}
which is a branched $2^{|\mathcal{P}|}$-fold covering. For a loop $\gamma$ in $\mathfrak{S}$ that is not null-homotopic, the trace-of-monodromy function $f_{[\gamma]} \in C^\infty( \mathscr{T}(\mathfrak{S}))$ and the length function $l_{[\gamma]} \in C^\infty( \mathscr{T}(\mathfrak{S}))$ yield corresponding functions on $\mathscr{T}^+(\mathfrak{S})$ via pullback under the map in equation~\eqref{eq:framing-forgetting_map}, which we still denote by the same symbols and call by the same names.

We will see that the enhanced Teichm\"uller space $\mathscr{T}^+(\mathfrak{S})$ is a smooth manifold equipped with the {\it Weil--Petersson} Poisson structure which is compatible with the (Atiyah--Bott--)Goldman symplectic structure \cite{Goldman} on the ${\rm SL}_2$-character variety of $\mathfrak{S}$, and possesses a system of special coordinate charts that are suitable for the quantization problem. These coordinate charts are enumerated by the ideal triangulations of $\mathfrak{S}$.
\begin{Definition}
\label{def:triangulable}
A marked surface $(\Sigma,\mathcal{P})$ is called {\em triangulable} if it is not isomorphic to one of the following:
\begin{enumerate}[label={\rm (\arabic*)}]\itemsep=0pt
\item A {\em monogon}: $\Sigma$ is diffeomorphic to a closed disc, and $\mathcal{P}$ consists of a single point lying in~$\partial \Sigma$.

\item A {\em bigon}: $\Sigma$ is diffeomorphic to a closed disc, and $\mathcal{P}$ consists of two points lying in $\partial \Sigma$.

\item A sphere with less than three punctures: $\Sigma$ is diffeomorphic to the 2-sphere $S^2$, and $|\mathcal{P}|\le 2$.
\end{enumerate}
\end{Definition}
In particular, in Definitions \ref{def:Teichmuller_space}, \ref{def:relative_Teichmuller_space} and \ref{def:enhanced_Teichmuller_space}, the hypothesis for $\mathfrak{S}$ can be phrased as $\mathfrak{S}$ being a triangulable punctured surface.
\begin{Definition}
\label{def:ideal_triangulation}
Let $(\Sigma,\mathcal{P})$ be a marked surface, with $\mathfrak{S} = \Sigma \setminus \mathcal{P}$. Let $[0,1] \subset \mathbb{R}$ be the closed interval.
\begin{enumerate}[label={\rm (\arabic*)}]\itemsep=0pt
\item An {\em ideal arc} in $(\Sigma,\mathcal{P})$ (or in $\mathfrak{S}$) is an immersion $\alpha \colon [0,1] \to \Sigma$ such that $\alpha(0),\alpha(1) \in \mathcal{P}$, whose restriction to the interior $\alpha|_{(0,1)} \colon (0,1) \to \Sigma$ is an embedding. An isotopy of ideal arcs is an isotopy within the class of ideal arcs.

\item A {\em boundary arc} in $(\Sigma,\mathcal{P})$ is an ideal arc in $(\Sigma,\mathcal{P})$ whose image lies in $\partial \Sigma$.

\item When $(\Sigma,\mathcal{P})$ is triangulable, an {\em ideal triangulation} $\Delta$ of $(\Sigma,\mathcal{P})$ is a maximal collection of ideal arcs in $(\Sigma,\mathcal{P})$ such that the members of $\Delta$ do not intersect except possibly at their endpoints, and no two members of $\Delta$ are isotopic.
\end{enumerate}
\end{Definition}
We often identify an ideal arc $\alpha \colon [0,1]\to\Sigma$ with its image $\alpha([0,1]) \subset \Sigma$, and call $\alpha((0,1))$ the {\em interior} of this ideal arc. An ideal triangulation $\Delta$ is often considered up to simultaneous isotopy of its members. By applying an isotopy, we assume that each ideal arc of $\Delta$ isotopic to a boundary arc is actually a boundary arc. When $\Delta$ is viewed as the corresponding subset $\bigcup_{i \in \Delta} i$ of $\Sigma$, the closure $\ol{T}$ in $\Sigma$ of a connected component $T$ of the complement $\Sigma \setminus \Delta$ is called an {\em ideal triangle} of $\Delta$. Note that $\ol{T} \setminus T$ is a union of three or two ideal arcs of $\Delta$; we call these ideal arcs the {\em sides} of $\ol{T}$. In case when $\ol{T}$ has only two sides, we say that the ideal triangle $\ol{T}$ is {\em self-folded}, and the ideal arc of $\Delta$ whose interior is contained in the interior of $\ol{T}$ is called the {\em self-folded arc} of $\ol{T}$ and of $\Delta$. Outside this paragraph, ideal triangles of an ideal triangulation $\Delta$ will be denoted by symbols such as $T$; that is, this time $T$ is a closed subset of $\Sigma$ and contains its sides. The combinatorial information of an ideal triangulation $\Delta$ is encoded in the following matrix.
\begin{Definition}[the exchange matrix for an ideal triangulation \cite{FST}]
\label{def:varepsilon_Delta}
Let $\Delta$ be an ideal triangulation of a triangulable marked surface $\mathfrak{S} = \Sigma \setminus \mathcal{P}$. Define the integer valued $\Delta\times \Delta$ matrix
\[
\varepsilon = \varepsilon^\Delta = (\varepsilon_{ij})_{i,j\in \Delta} = \bigl(\varepsilon^\Delta_{ij}\bigr)_{i,j\in \Delta},
\]
called the {\em exchange matrix for $\Delta$}, as follows: first, define
\begin{align*}
a_{ij} & = \sum_{T\colon \mbox{\tiny ideal triangles of $\Delta$}} a^T_{ij}, \\
a^T_{ij} & = \begin{cases}
 \hphantom{-}1 & \mbox{if $T$ is not self-folded, and its sides are $i,j,k$ appearing in clockwise order,} \\
 -1 & \mbox{if $T$ is not self-folded, and its sides are $j,i,k$ appearing in clockwise order,} \\
 \hphantom{-}0 & \mbox{otherwise}.
 \end{cases}
\end{align*}
Now, let $\pi \colon \Delta \to \Delta$ be given by
\begin{align}
\label{eq:pi}
\pi(i) = \begin{cases}
 i & \mbox{if $i$ is not a self-folded arc}, \\
 j & \begin{aligned} & \text{if $i$ is a self-folded arc of an ideal triangle $T$} \\
 & \text{and $j$ is the side of $T$ not equal to $i$}.\end{aligned}
 \end{cases}
\end{align}
Define
\begin{align*}
\varepsilon_{ij} := a_{\pi(i) \pi(j)} - a_{\pi(j) \pi(i)}.
\end{align*}
\end{Definition}
In particular, when $i$, $j$ are not self-folded arcs, we have $\varepsilon_{ij} = a_{ij} - a_{ji}$.

The following proposition is about the promised special coordinate systems on $\mathscr{T}^+(\mathfrak{S})$.
\begin{Proposition}[Thurston \cite{Thurston86}, Fock \cite{F93}]\label{prop:X_i}
Let $\mathfrak{S} = \Sigma\setminus \mathcal{P}$ be a triangulable punctured surface. Let $\Delta$ be an ideal triangulation of $\mathfrak{S}$. For each ideal arc $i \in \Delta$, there exists a function
\[
X_i = X_{i;\Delta} \colon \ \mathscr{T}^+(\mathfrak{S}) \to \mathbb{R}_{>0}
\]
called the {\em $($modified$)$ exponentiated shear coordinate} function, such that
\[
\mathscr{T}^+(\mathfrak{S}) \to (\mathbb{R}_{>0})^\Delta \colon \ ([g],\beta) \mapsto (X_i([g],\beta))_{i\in \Delta}
\]
is a real analytic bijection, such that the Weil--Petersson Poisson structure on $\mathscr{T}^+(\mathfrak{S})$ is given~by
\begin{align}
\label{eq:Poisson_brackets_of_shear_coordinates}
\{X_i, X_j\} = \varepsilon_{ij} X_i X_j \qquad \forall i,j \in \Delta.
\end{align}
\end{Proposition}

\begin{figure}
\centering
\begin{subfigure}[b]{0.4\textwidth}\centering
\hspace{0mm} \scalebox{0.7}{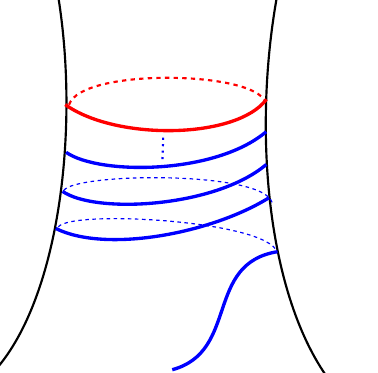}
\caption*{(1) positive framing}
\end{subfigure}
\qquad
\begin{subfigure}[b]{0.4\textwidth}\centering
\scalebox{0.75}{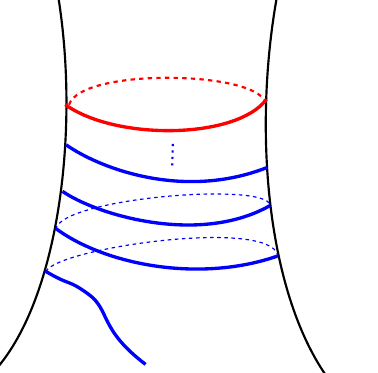}
\caption*{(2) negative framing}
\end{subfigure}
\caption{Framing and spiraling direction at a funnel end for puncture $p$.}
\label{fig:spiraling_direction}
\end{figure}

First, we briefly recall the construction of the unmodified exponentiated shear coordinate function $\til{X}_i$. Let $([g],\beta) \in \mathscr{T}^+(\mathfrak{S})$. Then the metric $g$ makes $\mathfrak{S} = \Sigma \setminus \mathcal{P}$ a complete hyperbolic surface $(\mathfrak{S},g)$. Stretch each ideal arc $i$ of $\Delta$ to the unique simple geodesic in its isotopy class; when an endpoint of $i$ is a puncture $p$ that is a funnel, we stipulate that $i$ spirals toward the unique simple geodesic loop surrounding $p$, with the spiraling direction specified by the framing data $\beta$ as depicted in Figure~\ref{fig:spiraling_direction}. The universal cover of $(\mathfrak{S},g)$ is the upper half plane $\mathbb{H}^2$ equipped with the standard hyperbolic metric. One can lift $\Delta$ to $\til{\Delta}$ in the universal cover; let $\til{i}$ be any of the lifts of $i$ in $\mathbb{H}^2$. Then $\til{i}$ is the diagonal of a unique ideal quadrilateral in $\mathbb{H}^2$ formed by the two ideal triangles of $\til{\Delta}$ in $\mathbb{H}^2$ having $\til{i}$ as a side. Denoting the four ideal vertices of this quadrilateral by $p_1,p_2,p_3,p_4 \in \mathbb{RP}^1 = \partial \mathbb{H}^2$ as below in the figure in the right, the coordinate is defined by the formula in the left
\[
\til{X}_i([g],\beta) = - \frac{p_1-p_2}{p_1-p_4} \frac{p_3-p_4}{p_3-p_2} \qquad \raisebox{-0,4\height}{\scalebox{0.9}{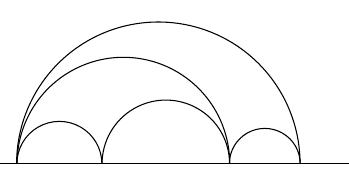}}
\]
Finally, we define the modified exponentiated shear coordinate functions as follows (\cite[Definition~9.2]{AB}, \cite[Section~10.7]{FG06}):
\[
X_i := \begin{cases}
 \til{X}_i & \text{if $i$ is not a self-folded arc}, \\
 \til{X}_i \til{X}_{\pi(i)} & \text{if $i$ is a self-folded arc, with $\pi$ as in equation~\eqref{eq:pi}}.
 \end{cases}
\]
We note that $\til{X}_i$ is the original definition of shear coordinates by Thurston, which is an incarnation of Fock and Goncharov's cluster $\mathscr{X}$-coordinates on their moduli space $\mathscr{X}_{{\rm PGL}_2,\mathfrak{S}}$. We use the modified version which suits better for the cluster variety structure, which in part is indicated by the Poisson structure in equation~\eqref{eq:Poisson_brackets_of_shear_coordinates}, and in part will be dealt with soon (Proposition~\ref{prop:shear_coordinate_change}).

We now give a description of $\mathscr{T}(\mathfrak{S})$ by embedding it into $\mathscr{T}^+(\mathfrak{S})$ as a subset defined by a system of inequalities written in terms of the exponentiated shear coordinates. We find it convenient to work with the unexponentiated shear coordinates:
\begin{align}
\label{eq:log_shear}
x_i = x_{i;\Delta} := \log(X_{i;\Delta}) \colon \ \mathscr{T}^+(\mathfrak{S}) \to \mathbb{R}, \qquad
\til{x}_i := \log\bigl(\til{X}_i\bigr).
\end{align}
In particular, $\mathscr{T}^+(\mathfrak{S}) \to \mathbb{R}^\Delta$, $([g],\beta) \mapsto (x_i)_{i \in \Delta}$, is a real analytic bijection, and the Poisson bracket can be presented as $\{x_i, x_j\} = \varepsilon_{ij}$. The following is well known.
\begin{Lemma}[\cite{CF, FG07}]
\label{lem:puncture_function_as_sum}
Let $\mathfrak{S} = \Sigma \setminus \mathcal{P}$ be a triangulable punctured surface. Let $\Delta$ be an ideal triangulation of $\mathfrak{S}$. For each puncture $p \in \mathcal{P}$, the length function $l_{[\gamma_p]}$ on $\mathscr{T}^+(\mathfrak{S})$ along a~peripheral loop $\gamma_p$ around $p$ satisfies
\[
l_{[\gamma_p]} = \left| x_p \right|,
\]
where
\[
x_p = x_{p;\Delta} := \sum_{i \in \Delta} \theta^i_p \, \til{x}_i,
\]
where $\theta_p^i \in \{0,1,2\}$ is the valence of $i$ at $p$.
\end{Lemma}
That is, if none of the endpoints of $i$ is $p$, then $\theta_p^i=0$. If only one endpoint of $i$ is $p$, then $\theta_p^i=1$. If both endpoints of $i$ are $p$, then $\theta_p^i=2$. We can also define a modified valence \smash{$\til{\theta}^i_p$} to express $x_p$ as \smash{$\sum_{i\in \Delta} \til{\theta}^i_p \, x_i$}, but let us not bother. So $x_p([g],\beta) =0$ if and only if $p$ is a cusp with respect to $[g]$. Moreover, when $x_p([g],\beta) \neq 0$, its sign determines the framing data of $\beta$ at $p$.
\begin{Lemma}[\cite{CF,FG07}]
Let $\mathfrak{S} = \Sigma \setminus \mathcal{P}$ be a triangulable punctured surface, and $\Delta$ an ideal triangulation of $\mathfrak{S}$. Suppose $x_p([g],\beta) \neq 0$ for some $p\in \mathcal{P}$ and $([g],\beta) \in \mathscr{T}^+(\mathfrak{S})$. Then the sign~$\beta(p)$ assigned by the framing $\beta$ at $p$ coincides with the sign of $x_p([g],\beta)$.
\end{Lemma}

As a section of the framing-forgetting map in equation~\eqref{eq:framing-forgetting_map}, consider the embedding
\[
\iota \colon \mathscr{T}(\mathfrak{S}) \to \mathscr{T}^+(\mathfrak{S}) \colon \ [g] \mapsto ([g], \beta_0),
\]
where $\beta_0$ denotes the all-plus framing, i.e., $\beta_0(p)=+$ for all $p \in \mathcal{P}$ when $p$ is a funnel with respect to $[g]$. It is then easy to observe the following characterization of the image of $\iota$.
\begin{Proposition}
Let $\mathfrak{S} = \Sigma \setminus \mathcal{P}$ be a triangulable punctured surface, and $\Delta$ an ideal triangulation of $\mathfrak{S}$. Then
\[
\iota(\mathscr{T}(\mathfrak{S})) = \{ ([g],\beta) \in \mathscr{T}^+(\mathfrak{S}) \mid x_{p;\Delta}([g],\beta) \ge 0, \, \forall p \in \mathcal{P} \}.
\]
\end{Proposition}
One can be more general, and for each sign sequence $\vec{\epsilon} = (\epsilon_p)_{p \in \mathcal{P}} \in \{+,-\}^\mathcal{P}$, one can embed~$\mathscr{T}(\mathfrak{S})$ into the subset defined by $\epsilon_p x_p \ge 0$, $\forall p \in \mathcal{P}$. These subsets for different $\vec{\epsilon}$ can be thought of as giving a decomposition of $\mathscr{T}^+(\mathfrak{S})$. One may also want to recall the fact that for each $\Delta$, the set of puncture functions $\{ x_p \mid p \in \mathcal{P} \}$, each of which is a linear function on $\mathbb{R}^\Delta \approx \mathscr{T}^+(\mathfrak{S})$, is linearly independent (e.g., \cite[Lemma~8]{BL07}, \cite[Lemma~3.2]{Kim_irreducible}). Consequently, the following hold.
\begin{Corollary}\label{cor:with_corners}
The Teichm\"uller space $\mathscr{T}(\mathfrak{S})$ is a smooth manifold with corners, of dimension $6g-6+3|\mathcal{P}|$.
The Teichm\"uller space $\mathscr{T}(\mathfrak{S})_{\vec{l}}$ is a smooth manifold diffeomorphic to $\mathbb{R}^{6g-6+2|\mathcal{P}|}$.
\end{Corollary}
This is one reason why it is better to deal with $\mathscr{T}^+(\mathfrak{S})$ instead of $\mathscr{T}(\mathfrak{S})$, while keeping track of the values of the constraint function $x_p$ at punctures.

A crucial property that the exponentiated shear coordinates enjoy is the coordinate change formula associated to the change of ideal triangulations. To describe this phenomenon efficiently, we choose to work with {\em labeled} ideal triangulations, each of which is an ideal triangulation $\Delta$ together with a bijection $\Delta \to I$ to a fixed index set $I$, which labels the members of $\Delta$ by elements of $I$. We still use just the same notation $\Delta$, when the underlying labelings are clear. There are two elementary changes of ideal triangulations, which generate all changes.
\begin{Definition}
Let $\mathfrak{S}= \Sigma\setminus \mathcal{P}$ be a triangulable marked surface. Let $\Delta$ and $\Delta'$ be labeled ideal triangulations of $\mathfrak{S}$.

We say $\Delta$ and $\Delta'$ are related by a {\em flip at the arc $k \in I$}, denoted by $\mu_k$, if the underlying ideal triangulations (considered up to isotopy) differ by exactly one arc, namely one labeled by $k$, and if the labelings are related by the natural bijection between the underlying ideal triangulations. We write $\Delta' = \mu_k(\Delta)$.

We say $\Delta$ and $\Delta'$ are related by a {\em label permutation} $P_\sigma$ associated to a bijection $\sigma \colon I \to I$, if the underlying ideal triangulations are the same and the labelings differ by $\sigma$, in the sense that if an arc of $\Delta$ is labeled by $i\in I$, then the same arc in $\Delta'$ is labeled by $\sigma(i)\in I$. We write $\Delta' = P_\sigma(\Delta)$.
\end{Definition}

\begin{Proposition}[{\cite[Fact 1.24]{P12}}, \cite{Thurston79}]
Let $\mathfrak{S} = \Sigma\setminus \mathcal{P}$ be a triangulable marked surface. Any two labeled ideal triangulations are related by a finite sequence of flips and label permutations.
\end{Proposition}
The algebraic relations among elementary changes are well known; see, e.g., \cite{FST, Kim_ratio}.

To investigate the coordinate change formulas for the exponentiated shear coordinates, it suffices to consider the flip and label permutation.
\begin{Proposition}[\cite{AB, BL07, F97, FST, Liu}]
\label{prop:shear_coordinate_change}
Let $\mathfrak{S} = \Sigma\setminus \mathcal{P}$ be a triangulable marked surface. Let~$\Delta$ and $\Delta'$ be labeled ideal triangulations of $\mathfrak{S}$.
\begin{enumerate}[label={\rm (\arabic*)}]\itemsep=0pt
\item If $\Delta' = \mu_k(\Delta)$ for $k\in I$, then the exponentiated shear coordinates $\{X_i \mid i \in I\}$ for $\Delta$ and those $\{X_i'\mid i \in I\}$ for $\Delta'$ are related by
\[
X'_i =\begin{cases}
 X_k^{-1} & \mbox{if } i=k, \\
 X_i \bigl(1+X_k^{-{\rm sign}(\varepsilon_{ik})}\bigr)^{-\varepsilon_{ik}} & \mbox{if } i\neq k,
 \end{cases}
\]
where $\varepsilon = \varepsilon^\Delta=(\varepsilon_{ij})_{i,j\in I}$ is the exchange matrix for $\Delta$ (here, ${\rm sign}(a)=a/|a|$ if $a\neq 0$). Moreover, the exchange matrix $\varepsilon' = \varepsilon^{\Delta'}=(\varepsilon'_{ij})_{i,j\in I}$ for $\Delta'$ is related to $\varepsilon^\Delta$ via the following formula:
\begin{align*}
\varepsilon'_{ij}
= \begin{cases}
 -\varepsilon_{ij} & \mbox{if } k\in \{i,j\}, \\
 \varepsilon_{ij} + \frac{1}{2}(\varepsilon_{ik}|\varepsilon_{kj}| + |\varepsilon_{ik}|\varepsilon_{kj}) & \mbox{if } k\notin \{i,j\}.
 \end{cases}
\end{align*}

\item If $\Delta' = P_\sigma(\Delta)$ for a permutation $\sigma$ of $I$, then $X'_{\sigma(i)} = X_i$, $\forall i\in I$. Moreover, the exchange matrices for $\Delta$ and $\Delta'$ are related by
\begin{align*}
\varepsilon'_{\sigma(i) \sigma(j)} = \varepsilon_{ij}, \qquad
X'_{\sigma(i)} = X_i \quad \forall i \in I.
\end{align*}
\end{enumerate}
\end{Proposition}
One immediately recognizes that this is an example of the cluster $\mathscr{X}$-mutation formula in the theory of cluster $\mathscr{X}$-varieties, or cluster Poisson varieties, of Fock and Goncharov \cite{FG06, FG06b, FG09}. So the chart of $\mathscr{T}^+(\mathfrak{S})$ formed by the exponentiated shear coordinates for an ideal triangulation~$\Delta$ can be viewed as the cluster chart of the cluster $\mathscr{X}$-variety associated to exchange matrices for ideal triangulations of the surface $\mathfrak{S}$, evaluated at the positive-real semifield $\mathbb{R}_{>0} \subset \mathbb{R}$; see \cite{FG06}.

\subsection{Balanced Chekhov--Fock algebras for triangulations, and naturality}
\label{subsec:balanced_algebras}

A basic question is on the quantization of the Teichm\"uller space $\mathscr{T}(\mathfrak{S})$ (Definition~\ref{def:Teichmuller_space}) of a punctured surface $\mathfrak{S}$, or its enhanced version $\mathscr{T}^+(\mathfrak{S})$ (Definition~\ref{def:enhanced_Teichmuller_space}). In general, by a {\em quantization} of a Poisson manifold $M$ we mean the following package of data:
\begin{enumerate}[label={\rm (Q\arabic*)}]\itemsep=0pt
\item\label{Q1} a Poisson $*$-subalgebra $\mathcal{A}$ of $C^\infty(M;\mathbb{C})$ of classical observables, where the $*$-structure on $C^\infty(M;\mathbb{C})$ is the complex conjugation on the values,

\item a separable complex Hilbert space $\mathscr{H}$ of quantum states,

\item a family of associative $*$-algebras $\mathcal{A}^\hbar$ over~$\mathbb{C}$, parametrized by the real quantum parameter~$\hbar$, such that $\mathcal{A}^0$ is isomorphic to~$\mathcal{A}$,

\item\label{Q4} a family of deformation quantization maps $Q^\hbar \colon \mathcal{A} \to \mathcal{A}^\hbar$ satisfying certain axioms,

\item\label{Q5} a family of representations $\rho^\hbar$ of $\mathcal{A}^\hbar$ on the Hilbert space $\mathscr{H}$.
\end{enumerate}
See, e.g., \cite{KS} for more details on the required conditions. By composing $Q^\hbar$ and $\rho^\hbar$, we get a~quantization map
\[
{\bf Q}^\hbar = \rho^\hbar \circ Q^\hbar \colon \ \mathcal{A} \stackrel{Q^\hbar}{\longrightarrow} \mathcal{A}^\hbar \stackrel{\rho^\hbar}{\longrightarrow} \{\mbox{operators on $\mathscr{H}$}\}
\]
that assigns to each classical observable $f\in \mathcal{A}$ a quantum observable ${\bf Q}^\hbar(f)$, which is an operator on $\mathscr{H}$. The operator version of the most crucial axiom of $Q^\hbar$ states the relationship between the Poisson bracket $\{\cdot,\cdot\}$ on $\mathcal{A}\subset C^\infty(M;\mathbb{C})$ and the operator commutator bracket $[\cdot,\cdot]$:
\[
[{\bf Q}^\hbar(f_1), {\bf Q}^\hbar(f_2)] = {\rm i} 2\pi \hbar \, {\bf Q}^\hbar(\{f_1,f_2\}) + o(\hbar) \quad \mbox{as $\hbar \to 0$}, \qquad \forall f_1,f_2 \in \mathcal{A}.
\]
Moreover, if there is a symmetry group $\mathfrak{G}$ acting on $\mathcal{A}$ as Poisson automorphisms, then one would naturally want to stipulate that the quantization ${\bf Q}^\hbar$ should be equivariant under the action of $\mathfrak{G}$. In our case, we would like to consider the problem of constructing a quantization for $M = \mathscr{T}(\mathfrak{S})$ or $\mathscr{T}^+(\mathfrak{S})$, when the symmetry group $\mathfrak{G}$ is the {\em mapping class group} ${\rm MCG}(\mathfrak{S}) = {\rm Diff}_+(\mathfrak{S}) / {\rm Diff}(\mathfrak{S})_0$, which is the group of orientation-preserving self-diffeomorphisms of $\mathfrak{S}$ considered up to isotopy.

The first results of the program of quantization of Teichm\"uller spaces were obtained in 1990's by Kashaev \cite{Kash98} and independently by Chekhov and Fock \cite{F97, CF}, and the latter approach is based on the shear coordinate systems we investigated in the previous subsection. The basic quantum algebra assigned to each ideal triangulation $\Delta$ is the following.
\begin{Definition}[\cite{BL07,Liu}]
Let $\Delta$ be a labeled ideal triangulation of the triangulable marked surface $\mathfrak{S} = \Sigma \setminus \mathcal{P}$, with the label index set $I$. The {\em Chekhov--Fock algebra} $\mathcal{X}^q_\Delta$ associated to $\Delta$ is the associative $\mathbb{Z}[q^{\pm 1}]$-algebra defined by generators and relations:
\[
\mbox{generators: $X_i^{\pm 1}$, \quad $i \in I$, \qquad relations: $X_i X_j = q^{2\varepsilon_{ij}} X_j X_i$, \quad $\forall i \in I$}.
\]
We equip $\mathcal{X}^q_\Delta$ with a $*$-structure, given by the $*$-map $\mathcal{X}^q_\Delta \to \mathcal{X}^q_\Delta$, $u \mapsto u^*$, defined as the unique ring anti-isomorphism such that
\[
X_i^* = X_i, \quad \forall i \in I, \qquad q^* = q^{-1}.
\]
\end{Definition}
We omit the trivial relations $X_i X_i^{-1} = X_i^{-1} X_i=1$. As the notation suggests, the generator $X_i \in \mathcal{X}^q_\Delta$ is the quantum version of the exponentiated shear coordinate function $X_i$ on $\mathscr{T}^+(\mathfrak{S})$. This algebra, if the $*$-structure is ignored, is an example of the (generalized) quantum torus algebras. In particular, it satisfies the Ore condition \cite{Cohn}, so that the skew field of fractions ${\rm Frac}\bigl(\mathcal{X}^q_\Delta\bigr)$ can be considered; see \cite{BL07, KLS, Liu}.

The following is a key statement on the naturality of the Chekhov--Fock algebras, with respect to change of ideal triangulations.
\begin{Proposition}[\cite{BL07, CF, FG09, Kash98, Liu}]
\label{prop:mu_q_Delta_Delta_prime}
Let $\mathfrak{S}$ be a triangulable marked surface.
For each pair of ideal triangulations $\Delta$ and $\Delta'$ of $\mathfrak{S}$ there exists a quantum coordinate change map
\[
\mu^q_{\Delta,\Delta'} \colon \ {\rm Frac}\bigl(\mathcal{X}^q_{\Delta'}\bigr) \to {\rm Frac}\bigl(\mathcal{X}^q_\Delta\bigr)
\]
which is an $*$-isomorphism of $*$-skew fields, that recovers the classical coordinate change formulas for exponentiated shear coordinates as $q\to 1$, and satisfying the consistency equations:
\[
\mu^q_{\Delta,\Delta'} \circ \mu^q_{\Delta',\Delta''} = \mu^q_{\Delta,\Delta''}
\]
for each triple of ideal triangulations $\Delta$, $\Delta'$ and $\Delta''$ of $\mathfrak{S}$.
\end{Proposition}
One might consider having the system of Chekhov--Fock algebras $\mathcal{X}^q_\Delta$ for ideal triangulations~$\Delta$, together with the isomorphisms \smash{$\mu^q_{\Delta,\Delta'}$}, as having a `quantization' of the enhanced Teichm\"uller space $\mathscr{T}^+(\mathfrak{S})$. However, this is just having a consistent quantum system, but not really a quantization that connects the classical system to the quantum system, in the sense we formulated in \ref{Q1}--\ref{Q5}.

To discuss quantization, we first should choose which classical algebra $\mathcal{A}$ to quantize. One of the most important functions to quantize are the trace-of-monodromy functions along loops. In particular, in \cite{Verlinde}, where the quantization program of Teichm\"uller spaces was proposed for the first time, Verlinde suggested to construct and use a quantized operator for the trace-of-monodromy function as a crucial ingredient for the space of conformal blocks of the Liouville conformal field theory. Also a posteriori, the trace-of-monodromy functions have turned out to play pivotal roles in later developments of quantum Teichm\"uller theory. So, we suggest to set~$\mathcal{A}$ to be the subalgebra generated by the trace-of-monodromy functions $f_{[\gamma]}$ along loops $\gamma$ in $\mathfrak{S}$ that are not null-homotopic.

We would want to quantize these functions $f_{[\gamma]}$ on $\mathscr{T}^+(\mathfrak{S})$ using the framework of Thurston's exponentiated shear coordinates, and their corresponding quantum algebras, the Chekhov--Fock algebras. We first observe the following.
\begin{Proposition}[\cite{F93, F97, FG06, FG07}; see proof of Proposition~\ref{prop:BW_admissible_sum_formula}]
\label{prop:classical_Laurentness}
Let $\mathfrak{S} = \Sigma\setminus \mathcal{P}$ be a triangulable punctured surface. Let $\gamma$ be a loop in $\mathfrak{S}$ that is not null-homotopic. For any labeled ideal triangulation $\Delta$, the function $f_{[\gamma]} \in C^\infty(\mathscr{T}^+(\mathfrak{S}))$ can be written as a Laurent polynomial in the square-roots $\sqrt{X_i}$ of Thurston's exponentiated shear coordinate functions $X_i = X_{i;\Delta}$ $($in~Proposition~{\rm \ref{prop:X_i})}, $i\in I$, with coefficients in $\mathbb{Z}_{\ge 0}$.
\end{Proposition}
An explicit algorithm to compute this square-root Laurent polynomial can be found in the references of this proposition; we do not review it here as it will follow from the quantum counterpart of this algorithm in Section~\ref{subsec:state-sum}, which is mostly based on \cite{BW}.

One thing to notice is that in general, $f_{[\gamma]} \in C^\infty(\mathscr{T}^+(\mathfrak{S}))$ is expressed only by the square-roots of the exponentiated shear coordinates, so to quantize it we need a square-root version of the Chekhov--Fock algebra.
\begin{Definition}
Let $\Delta$ be a labeled ideal triangulation of a triangulable marked surface $\mathfrak{S} = \Sigma\setminus\mathcal{P}$, with $I$ as the labeling index set.

The {\em square-root Chekhov--Fock algebra} $\mathcal{Z}^\omega_\Delta$ is the associative $\mathbb{Z}\bigl[\omega^{\pm 1}\bigr]$-algebra defined by generators and relations:
\[
\mbox{generators: $Z_i^{\pm 1}$, \quad $i \in I$, \qquad relations: $Z_i Z_j = \omega^{2\varepsilon_{ij}} Z_j Z_i$, \quad $\forall i \in I$},
\]
equipped with the $*$-structure given by $Z_i^* = Z_i$, $\forall i \in I$, $\omega^* = \omega^{-1}$, where the Chekhov--Fock algebra $\mathcal{X}^q_\Delta$ is thought of as being embedded into $\mathcal{Z}^\omega_\Delta$ as a subalgebra via the embedding
\[
q \mapsto \omega^4, \qquad
X_i \mapsto Z_i^2, \quad \forall i \in I.
\]
\end{Definition}
To be more precise, one has to work with a certain subalgebra of the square-root version of the Chekhov--Fock algebra, to which the quantum coordinate change isomorphism can be extended algebraically. To describe it more conveniently, and for our later purposes, we recall the following well-known coordinate-free formulation of quantum torus algebras, which we will eventually apply to our $\mathcal{Z}^\omega_\Delta$.

\begin{Definition}
\label{def:QT}
Let $(N,\langle \cdot,\cdot\rangle)$ be a skew symmetric lattice, i.e., a free $\mathbb{Z}$-module $N$ equipped with a skew symmetric $\mathbb{Z}$-bilinear form $\langle \cdot, \cdot \rangle \colon N \times N \to \mathbb{Z}$. We denote $(N,\langle \cdot,\cdot\rangle)$ by $N$ if $\langle \cdot,\cdot\rangle$ is clear. The {\em quantum torus algebra} associated to the skew symmetric lattice $N$, with a formal quantum parameter $t$, is the associative $*$-algebra $\mathbb{T}^t_N$ over $\mathbb{Z}\bigl[t^{\pm 1}\bigr]$, defined by generators and relations as
\[
\mbox{generators}: Z_v, \quad v \in N, \qquad
\mbox{relations}: Z_v Z_w = t^{\langle v,w\rangle} Z_{v+w},\quad \forall v,w\in N,
\]
with the $*$-structure given by $Z_v^*=Z_v$, $\forall v\in N$, $t^* = t^{-1}$.
\end{Definition}
Note in particular that each generator $Z_v$ is invertible, with the inverse being $Z_{-v}$. It is easy to observe that $Z_0 = 1$.

We introduce a logarithmic version of the quantum torus algebra as follows, which can be viewed mostly as a notational gadget; in fact, we believe that this gives a better heuristic understanding of the situation.

\begin{Definition}
\label{def:HA}
Let $N$ be as in Definition~\ref{def:QT}. Let $\mathbb{H}_N$ be the $\mathbb{Z}$-module defined as the quotient of the free $\mathbb{Z}$-module generated by the set of symbols $\{ z_v \mid v\in N\} \sqcup \{{\rm c}\}$ by the relations
\begin{align}
\label{eq:HA_relations}
z_v + z_w = z_{v+w}, \qquad \forall v,w\in N.
\end{align}
Equip $\mathbb{H}_N$ with the unique Lie bracket $[\cdot,\cdot]$ over $\mathbb{Z}$, satisfying
\begin{align}
\label{eq:HA_commutation}
[z_v, z_w] = 2 \langle v,w\rangle \, {\rm c}, \qquad [z_v,{\rm c}]=0, \qquad \forall v,w\in N.
\end{align}
This Lie algebra $\mathbb{H}_N$ is called the {\em Heisenberg algebra} associated to $N$.
\end{Definition}
From equation~\eqref{eq:HA_relations}, we see that each element of $\mathbb{H}_N$ can be uniquely written as $\alpha {\rm c} + z_v$ for some $\alpha \in \mathbb{Z}$ and $v\in N$. It is straightforward to see that $[\alpha {\rm c} + z_v, \beta {\rm c} + z_w ] = 2 \langle v,w\rangle {\rm c}$ yields a well-defined Lie bracket, justifying the above definition. Two elements $z$ and $z'$ of $\mathbb{H}_N$ are said to {\em commute} with each other if $[z,z']=0$.
\begin{Definition}
\label{def:HA_exp}
Let $N$ and $\mathbb{T}^t_N$ be as in Definition~\ref{def:QT}, and $\mathbb{H}_N$ be as in Definition~\ref{def:HA}. Define the {\em exponential} map
\[
\exp \colon \ \mathbb{H}_N \to \mathbb{T}^t_N
\]
as
\[
\exp( \alpha {\rm c} + z_v ) = t^\alpha Z_v, \qquad \forall \alpha \in \mathbb{Z}, \quad \forall v\in N.
\]
For each $z\in \mathbb{H}_N$, we also write
\[
{\rm e}^z := \exp(z).
\]
\end{Definition}
Under the above notation, we see
\[
t^{\pm 1} = {\rm e}^{\pm {\rm c}}, \qquad Z_v = {\rm e}^{z_v}, \qquad \forall v\in N.
\]
What makes this notation particularly useful in practice is the following.
\begin{Proposition}[Baker--Campbell--Hausdorff (BCH) formula]
\label{prop:BCH}
One has
\[
{\rm e}^z {\rm e}^{z'} = {\rm e}^{\frac{1}{2} [z,z']} {\rm e}^{z+z'}, \qquad \forall z,z' \in \mathbb{H}_N.
\]
\end{Proposition}
Notice that $\frac{1}{2} [z,z'] \in \mathbb{Z}{\rm c}$, so that ${\rm e}^{\frac{1}{2} [z,z']}$ is an integer power of $t$. We omit a proof of this proposition as it is a standard exercise; readers can consult, e.g., \cite{Kim_mutation} for a proof.

To obtain a coordinate dependent version of the quantum torus algebra, suppose that one fixes a free $\mathbb{Z}$-module basis $\{v_i \mid i \in I\}$ of the skew symmetric lattice $N$, enumerated by some index set $I$. One can record the skew symmetric form of $N$ by the numbers $n_{ij} = \langle v_i, v_j \rangle$ and then present the quantum torus algebra $\mathbb{T}^t_N$ by
\begin{align*}
\mbox{generators: } Z_{v_i}^{\pm 1}, \quad i\in I, \qquad \mbox{and the relations: } Z_{v_i} Z_{v_j} = t^{2 n_{ij}} Z_{v_j} Z_{v_i}, \quad \forall i,j\in I.
\end{align*}
An element of $\mathbb{T}^t_N$ is a Laurent polynomial in the variables $Z_{v_i}$, $i\in I$, with coefficients in $\mathbb{Z}\bigl[t^{\pm 1}\bigr]$. For any finite sequence of indices $i_1,i_2,\dots,i_n \in I$, define the {\em Weyl-ordered Laurent monomial}~as
\begin{align}
\label{eq:Weyl-ordering}
\bigl[Z_{v_{i_1}} Z_{v_{i_2}} \cdots Z_{v_{i_n}} \bigr]_{\rm Weyl} := t^{-\sum_{j<k} n_{i_j i_k}} Z_{v_{i_1}} Z_{v_{i_2}} \cdots Z_{v_{i_n}}.
\end{align}
It is well known that this is invariant under changing the order of the sequence $i_1,\dots,i_n$; this fact follows from the equality
\[
\bigl[Z_{v_{i_1}} Z_{v_{i_2}} \cdots Z_{v_{i_n}} \bigr]_{\rm Weyl} = Z_{v_{i_1} + v_{i_2} + \cdots + v_{i_n}},
\]
which is a straightforward exercise to check. Thus, for any $v\in N$, the associated element $Z_v \in \mathbb{T}^t_N$ can be called a Weyl-ordered Laurent monomial. By a {\em Laurent monomial} in $\mathbb{T}^t_N$, we mean an element of the form $\pm t^n Z_v$ with $n \in \mathbb{Z}$ and $v\in N$. For a positive Laurent monomial of the form $Z_{v_{i_1}} Z_{v_{i_2}} \cdots Z_{v_{i_n}}$, we say that $Z_{v_{i_1}+\cdots+v_{i_n}}$ is the {\em Weyl-ordering} of it.

We finally come back to our situation of quantum Teichm\"uller theory, and apply the above construction to the following lattice $N_\Delta$.
\begin{Definition}\label{def:N_Delta}
For a labeled ideal triangulation $\Delta$ of a triangulable marked surface $\mathfrak{S}$, the associated skew symmetric lattice $N_\Delta$ is defined as a free $\mathbb{Z}$-module with basis $\{v_i \mid i \in I\}$, with the skew symmetric form $\langle \cdot, \cdot \rangle$ given on the basis by
\begin{align}
\label{eq:skew-form_for_Delta}
\langle v_i, v_j \rangle = \varepsilon_{ij}, \qquad \forall i,j \in I,
\end{align}
where $\varepsilon = (\varepsilon_{ij})_{i,j\in \Delta}$ is as in Definition~\ref{def:varepsilon_Delta}.
\end{Definition}
Then naturally, we can identify $\mathcal{Z}^\omega_\Delta$ as a quantum torus
\[
\mathcal{Z}^\omega_\Delta \cong \mathbb{T}^\omega_{N_\Delta},
\]
with $Z_i$ corresponding to $Z_{v_i}$.

Now we move on to the special subalgebra of $\mathcal{Z}^\omega_\Delta$, on which the square-root version of the quantum coordinate change formula shall be defined.

\begin{Definition}[\cite{BW, Hiatt}]\label{def:balanced_algebra}
Let $\Delta$ be an ideal triangulation of a triangulable marked surface~$\mathfrak{S}$. We say that a Laurent monomial $Z_v$ is {\em balanced} if $v = \sum_{i \in \Delta} a_i v_i$ with $a_i \in \mathbb{Z}$ satisfies the following: for each ideal triangle $T$ of $\Delta$, the sum of $a_i$ for the sides $i$ of $T$, where a self-folded side is counted with multiplicity $2$, is even.

An element of $\mathcal{Z}^\omega_\Delta$ is {\em balanced} if it is a $\mathbb{Z}\bigl[\omega^{\pm 1}\bigr]$-linear combination of balanced Laurent monomials.

The {\em balanced square-root Chekhov--Fock algebra} $\bigl(\mathcal{Z}^\omega_\Delta\bigr)_{\rm bl}$ is the set of all elements of $\mathcal{Z}^\omega_\Delta$ that are balanced.
\end{Definition}
It is easy to check that $\bigl(\mathcal{Z}^\omega_\Delta\bigr)_{\rm bl}$ is a subalgebra.
\begin{Remark}
We suspect that the coefficient ring of the balanced subalgebra \smash{$\bigl(\mathcal{Z}^\omega_\Delta\bigr)_{\rm bl}$} may be taken to be a smaller ring, either $\mathbb{Z}[q^{\pm 1/2}]$ or $\mathbb{Z}[q^{\pm 1}]$, to serve our purposes. But we do not try to justify this.
\end{Remark}
It is known that $\mathcal{Z}^\omega_\Delta$, as well as \smash{$\bigl(\mathcal{Z}^\omega_\Delta\bigr)_{\rm bl}$}, satisfies the Ore condition \cite{Cohn, KLS}, so that the skew field of fractions can be considered. The skew field of fraction \smash{${\rm Frac}\bigl(\bigl(\mathcal{Z}^\omega_\Delta\bigr)_{\rm bl}\bigr)$} is where the quantum coordinate change isomorphism $\mu^q_{\Delta,\Delta'}$ can be extended. Noticing that ${\rm Frac}\bigl(\mathcal{X}^q_\Delta\bigr)$ is naturally a~subalgebra of \smash{${\rm Frac}\bigl(\bigl(\mathcal{Z}^\omega_\Delta\bigr)_{\rm bl}\bigr)$}.
\begin{Proposition}[balanced quantum coordinate change isomorphisms; \cite{BW, Hiatt,KLS,Son}]
\label{prop:Theta}
Let $\mathfrak{S}= \Sigma\setminus\mathcal{P}$ be a triangulated marked surface. For each pair of ideal triangulations $\Delta$ and $\Delta'$ of~$\mathfrak{S}$, there exists a skew field isomorphism
\[
\Theta^\omega_{\Delta,\Delta'} \colon \ {\rm Frac}\bigl(\bigl(\mathcal{Z}^\omega_{\Delta'}\bigr)_{\rm bl}\bigr) \to {\rm Frac}\bigl(\bigl(\mathcal{Z}^\omega_\Delta\bigr)_{\rm bl}\bigr)
\]
that preserves the $*$-structure, extends $\mu^q_{\Delta,\Delta'} \colon {\rm Frac}\bigl(\mathcal{X}^q_{\Delta'}\bigr) \to {\rm Frac}\bigl(\mathcal{X}^q_\Delta\bigr)$ (Proposition~\ref{prop:mu_q_Delta_Delta_prime}) and satisfies the consistency equation
\[
\Theta^\omega_{\Delta,\Delta'} \circ \Theta^\omega_{\Delta',\Delta''} = \Theta^\omega_{\Delta,\Delta''}
\]
for each triple of ideal triangulations $\Delta$, $\Delta'$ and $\Delta''$ of $\mathfrak{S}$.
\end{Proposition}
To describe $\Theta^\omega_{\Delta,\Delta'}$, it suffices to do so for the case when $\Delta' = \mu_k(\Delta)$ or when $\Delta' = P_\sigma(\Delta)$. For the latter case $\Delta' = P_\sigma(\Delta')$, $\Theta^\omega_{\Delta,\Delta'}$ is induced by the isomorphism $\mathcal{Z}^\omega_{\Delta'} \to \mathcal{Z}^\omega_\Delta$ given on the generators by \smash{$Z_{v_{\sigma(i)}'} \mapsto Z_{v_i}$}. For the former case $\Delta' = \mu_k(\Delta)$, $\Theta^\omega_{\Delta,\Delta'} = \Theta^\omega_k$ is given as the composition
\[
\xymatrix{
\Theta^\omega_k = \Theta^{\sharp \omega}_k \circ \Theta'_k \colon \ {\rm Frac}((\mathcal{Z}^\omega_{\Delta'})_{\rm bl}) \ar[r]^-{\Theta'_k} & {\rm Frac}(\bigl(\mathcal{Z}^\omega_\Delta\bigr)_{\rm bl}) \ar[r]^-{\Theta^{\sharp \omega}_k} & {\rm Frac}(\bigl(\mathcal{Z}^\omega_\Delta\bigr)_{\rm bl}),
}
\]
which we describe now. The `monomial transformation part' $\Theta'_k\colon {\rm Frac}((\mathcal{Z}^\omega_{\Delta'})_{\rm bl}) \to {\rm Frac}(\bigl(\mathcal{Z}^\omega_\Delta\bigr)_{\rm bl})$ is given on the generators $Z_{v'}$, $v' \in N_{\Delta'}$, by
\[
\Theta'_k(Z_{v'}) = Z_{C_k(v)},
\]
where $C_k \colon N_{\Delta'} \to N_\Delta$ is a $\mathbb{Z}$-linear map given on the generators $v_i'$, $i\in I$, by
\[
C_k(v_i') =\begin{cases}
-v_k & \mbox{if } i=k, \\
v_i + [\varepsilon_{ik}]_+ v_k & \mbox{if } i\neq k,
\end{cases}
\]
where $\varepsilon_{ik}$ is the value of the exchange matrix for $\Delta$, and $[a]_+ := \max(a,0)$.
The `automorphism part' $\Theta^{\sharp\omega}_k \colon {\rm Frac}\bigl(\bigl(\mathcal{Z}^\omega_\Delta\bigr)_{\rm bl}\bigr) \to {\rm Frac}\bigl(\bigl(\mathcal{Z}^\omega_\Delta\bigr)_{\rm bl}\bigr)$ is given on the generators $Z_v$, $v\in N_\Delta$, the balanced Laurent monomials, by
\[
\Theta^{\sharp\omega}_k(Z_v) = Z_v \, F^q(X_k; \langle v_k,v\rangle),
\]
where $X_k = Z_{v_k}^2$, and $F^q$ is as follows:
\begin{align}
\nonumber
F^q(x;\alpha) = {\textstyle \prod}_{r=1}^{|\alpha|} \bigl( 1 + q^{(2r-1) {\rm sign}(\alpha)} x\bigr)^{{\rm sign}(\alpha)}, \quad \mbox{for } \alpha \in\mathbb{Z} \setminus \{0\}, \qquad F^q(x;0)=1.
\end{align}
For example, we have $F^q(x;1) = 1+qx$, $F^q(x;2) = (1+qx)(1+q^3x)$, $F^q(x;-1) = (1+q^{-1}x)^{-1}$, $F^q(x;-2) = (1+q^{-1} x)^{-1}(1+q^{-3}x)^{-1}$.

In the next subsection, we review how the trace-of-monodromy function $f_{[\gamma]}$ on $\mathscr{T}^+(\mathfrak{S})$ is algebraically quantized by elements of $\bigl(\mathcal{Z}^\omega_\Delta\bigr)_{\rm bl}$, in a manner that is independent of the choice of an ideal triangulation $\Delta$.

\subsection{Quantized trace-of-monodromy via Bonahon--Wong quantum trace maps}
\label{subsec:quantized_trace-of-monodromy}

\begin{figure}
\centering
\scalebox{0.7}{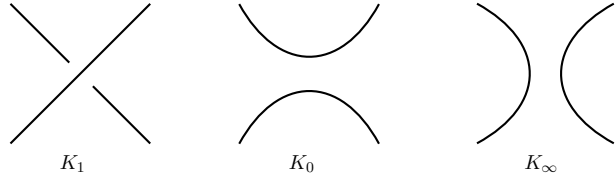}
\vspace{-3mm}
\caption{Kauffman triple of framed links.}
\label{fig:Kauffman_triple}
\end{figure}

One original problem that we discussed is that of quantization of the Teichm\"uller space $\mathscr{T}(\mathfrak{S})$, or the enhanced Teichm\"uller space $\mathscr{T}^+(\mathfrak{S})$. As mentioned, as the algebra $\mathcal{A}$ to be quantized we consider the subalgebra generated by the trace-of-monodromy functions $f_{[\gamma]}$ along loops $\gamma$ in $\mathfrak{S}$. It is well known that the following algebra serves as a quantum algebra $\mathcal{A}^\hbar$ for this $\mathcal{A}$.
\begin{Definition}[skein algebra of a surface; \cite{Pr, Turaev}]
\label{def:skein_algebra}
Let $\mathfrak{S} = \Sigma \setminus \mathcal{P}$ be a marked surface. Denote by
\[
{\bf I} = (-1,1) \subset \mathbb{R}
\]
the open interval, given the usual topology.
\begin{enumerate}[label={\rm (SA\arabic*)}]\itemsep=0pt \setlength{\leftskip}{0.22cm}
\item A {\em framed link} $K$ in $\mathfrak{S} \times {\bf I}$ is a properly embedded closed unoriented 1-submanifold $K$ of $\mathfrak{S}\times {\bf I}$ equipped with a {\em framing}, which is a continuous choice of a (framing) vector in $T_x(\mathfrak{S} \times {\bf I}) \setminus T_x K$ per each $x\in K$. We say that a framing vector is {\em upward vertical} if it is parallel to the ${\bf I}$ factor and points toward $1$ (rather than $-1$). An isotopy of framed links is a homotopy within the class of framed links.

When part of a framed link $K$ is depicted in a picture, we project it to a subset of the plane $\mathbb{R}^2$, where the framing is assumed to be pointing toward the eyes of the reader. We apply the projection only when the projection is at most $2$-to-$1$ in this picture; two points of $K$ projecting to a same point in the plane forms a {\em crossing} in the picture, where we indicate the over- and under-crossing information by broken lines, e.g. as in $K_1$ of Figure~\ref{fig:Kauffman_triple}.

\item An ordered triple of framed links $(K_1,K_0,K_\infty)$ is called the {\em Kauffman triple} if they coincide outside a small ball, in which they look as in Figure~\ref{fig:Kauffman_triple}.

\item The {\em $($Kauffman bracket$)$ skein algebra} $\mathcal{S}^t(\mathfrak{S})$ is a free associative algebra over $\mathbb{Z}\bigl[t^{\pm 1}\bigr]$ defined as follows: as a $\mathbb{Z}\bigl[t^{\pm 1}\bigr]$-module, it is the quotient of the free $\mathbb{Z}\bigl[t^{\pm 1}\bigr]$-module with a~free basis $\{\mbox{isotopy classes $[K]$ of framed links $K$ in $\mathfrak{S} \times {\bf I}$}\}$ by the $\mathbb{Z}\bigl[t^{\pm 1}\bigr]$-submodule generated by the elements
\[
[K_1] - t^2 [K_0] - t^{-2} [K_\infty],
\]
where $(K_1,K_0,K_\infty)$ is a Kauffman triple, and the elements
\[
[K \cup {\rm O}] + \bigl(t^4+t^{-4}\bigr) [K],
\]
where ${\rm O}$ denotes the {\em unknot} in $\mathfrak{S} \times {\bf I}$, i.e., a circle that lives in a constant elevation $\mathfrak{S} \times \{ a\}$ bounding a disc and is equipped with the upward vertical framing, and $K$ is any framed link such that $K \cup {\rm O}$ has no crossings.

Denote the element of $\mathcal{S}^t(\mathfrak{S})$ represented by the isotopy class $[K]$ still as $[K]$.

The multiplication of $\mathcal{S}^t(\mathfrak{S})$ is given by superposition, i.e., when $K_1 \subset \mathfrak{S} \times (-1,0)$ and $K_2 \subset \mathfrak{S} \times (0,1)$, then
\[
[K_1] \cdot [K_2] := [K_1 \cup K_2].
\]
\end{enumerate}
\end{Definition}
As a consequence of the definition, in $\mathcal{S}^t(\mathfrak{S})$ we have the relations $[K_1] = t^2 [K_0] + t^{-2} [K_\infty]$, called the {\em skein relations}, and $[K\cup {\rm O}] = - \bigl(t^4+t^{-4}\bigr) [K]$, called the {\em trivial loop relations} (or the `framing relations' as in \cite{Bullock}). We allow a framed link to be the empty one $\varnothing$, and one may easily notice that $[{\varnothing}]$ represented by the empty skein is the multiplicative identity $1$ of the algebra~$\mathcal{S}^t(\mathfrak{S})$. A traditional definition uses $A$ as a quantum parameter, where we put $A = t^{-2}$ for later notational convenience.

We suggest that the following version of the skein algebra is what would fit Verlinde's suggestion \cite{Verlinde} and the modular functor conjecture \cite{F97, FG09}.
\begin{Definition}[\cite{AS, FKL}]
\label{def:relative_skein_algebra}
Let $\mathfrak{S} = \Sigma \setminus \mathcal{P}$ be a triangulable punctured surface. Let $\vec{l} = (l_p)_{p \in \mathcal{P}} \in (\mathbb{R}_{\ge 0})^\mathcal{P}$.
Define the {\em relative skein algebra} $\mathcal{S}^t(\mathfrak{S})_{\vec{l}}$ as the quotient of $\mathcal{S}^t(\mathfrak{S})$ by the two-sided ideal generated by $[K_{\gamma_p}] + 2\cosh(l_p) [{\varnothing]}$, $p\in \mathcal{P}$, where $\gamma_p$ is a peripheral loop around $p$, and $K_{\gamma_p}$ is a~lift of $\gamma_p$ in $\mathfrak{S} \times \mathbf{I}$ at a~constant elevation, say lying in $\mathfrak{S} \times \{0\}$, with the upward vertical framing.\looseness=-1
\end{Definition}

We might want to formulate the final results of the sequel \cite{sequel} in terms of the relative skein algebra; see Section~\ref{subsec:future_topics}. However, we do not necessarily need to work with it in the present paper. So, we just stick to the non-relative skein algebra $\mathcal{S}^t(\mathfrak{S})$ for our purposes.

In quantum topology, a nice basis of $\mathcal{S}^t(\mathfrak{S})$ has been studied \cite{AK, DThurston}, based on crossingless skeins and their suitable modifications. Using this basis, one could construct a map $\mathcal{S}^1(\mathfrak{S}) \to \mathcal{S}^t(\mathfrak{S})$, say for $t = {\rm e}^{\pi {\rm i} \hbar/4}$, to solve the deformation quantization problem \ref{Q4}. However, in order to solve \ref{Q5}, one should construct representations of the noncommutative skein algebra $\mathcal{S}^t(\mathfrak{S})$, which is quite a nontrivial task to accomplish. Meanwhile, one might instead use Thurston's shear coordinates, together with various versions of the Chekhov--Fock algebras, to solve the quantization problem \ref{Q1}--\ref{Q5}. Now the deformation quantization problem \ref{Q4} is difficult, while the representation problem \ref{Q5} is easier. So one might try to combine the two approaches, one by skein algebras and the other by shear coordinates and Chekhov--Fock algebras. This led to the developments as we describe from now on.

We recall Bonahon and Wong's {\em quantum trace maps}, which relate the skein algebras to the balanced square-root Chekhov--Fock algebras and their fraction skew fields. We only write down a partial version as follows.
\begin{Proposition}[Bonahon--Wong quantum trace; \cite{BW}, see also \cite{Le19}]
\label{prop:BW}
For each triangulable marked surface $\mathfrak{S} = \Sigma\setminus\mathcal{P}$ and its ideal triangulation $\Delta$, there is a $\mathbb{Z}\bigl[\omega^{\pm 1}\bigr]$-algebra homomorphism
\[
{\rm Tr}^\omega_\Delta = {\rm Tr}^\omega_{\Delta;\mathfrak{S}} \colon \ \mathcal{S}^\omega(\mathfrak{S}) \to \bigl(\mathcal{Z}^\omega_\Delta\bigr)_{\rm bl} \subset \mathcal{Z}^\omega_\Delta
\]
called the {\em $($Bonahon--Wong$)$ quantum trace} map $($see Definition~{\rm \ref{def:skein_algebra}} for $\mathcal{S}^\omega(\mathfrak{S})$ and Definition~{\rm \ref{def:balanced_algebra}} for $\bigl(\mathcal{Z}^\omega_\Delta\bigr)_{\rm bl})$, satisfying
\begin{enumerate}[label={\rm (QT\arabic*)}]\itemsep=0pt\setlength{\leftskip}{0.30cm}
\item\label{QT1} $($quantization of trace-of-monodromy with respect to shear coordinates$)$ if $K$ is a connected framed link in $\mathfrak{S} \times {\bf I}$ and if $\gamma$ is a loop in $\mathfrak{S}$ obtained from $K$ by forgetting the framing and applying the projection $\mathfrak{S} \times {\bf I} \to \mathfrak{S}$, then the specialization of the quantum trace ${\rm Tr}^\omega_{\Delta;\mathfrak{S}}([K]) \in \mathcal{Z}^\omega_\Delta$ at $\omega=1$ recovers the trace-of-monodromy function $f_{[\gamma]} \in C^\infty(\mathscr{T}^+(\mathfrak{S}))$ $($as in Proposition~{\rm \ref{prop:classical_Laurentness})} written in terms of the $($square-roots of the$)$ exponentiated shear coordinates for $\Delta$, where each generator $Z_{v_i} \in \mathcal{Z}^\omega_\Delta$, $i\in \Delta$, corresponds to the square-root~$\sqrt{X_i}$ of the exponentiated shear coordinate function $X_i$, as $\omega \to 1$;

\item\label{QT2} $($naturality$)$ for any two ideal triangluations $\Delta$ and $\Delta'$ of $\mathfrak{S}$, the following compatibility equation holds:
\[
{\rm Tr}^\omega_{\Delta;\mathfrak{S}} = \Theta^\omega_{\Delta,\Delta'} \circ {\rm Tr}^\omega_{\Delta';\mathfrak{S}}.
\]
\end{enumerate}
\end{Proposition}

One important immediate consequence is that each image of the quantum trace map is universally balanced (quantum) Laurent. That is, ${\rm Tr}^\omega_{\Delta;\mathfrak{S}}([K])$ lies in
\begin{align*}
\mathbb{BL}^\omega_\Delta := \bigcap_{\Delta'} \Theta^\omega_{\Delta,\Delta'}\bigl( \bigl(\mathcal{Z}^\omega_{\Delta'}\bigr)_{\rm bl} \bigr) \subset \bigl(\mathcal{Z}^\omega_\Delta\bigr)_{\rm bl} \subset {\rm Frac}\bigl( \bigl( \mathcal{Z}^\omega_\Delta\bigr)_{\rm bl} \bigr),
\end{align*}
where the intersection is over all ideal triangulations $\Delta'$ of $\mathfrak{S}$. A characterizing property of the quantum trace maps is the compatibility under cutting and gluing of the marked surfaces along ideal arcs. In order to formulate this property, we need to use the {\em stated} skein algebras $\mathcal{S}^\omega(\mathfrak{S})_{\rm s}$, built not just by framed links in $\mathfrak{S} \times {\bf I}$ but also by framed tangles in $\mathfrak{S} \times {\bf I}$, which may have endpoints at the boundary $\partial \mathfrak{S} \times {\bf I}$, equipped with state values $\in \{+,-\}$ at endpoints. However, in the present paper, we will only deal with framed links, but not framed tangles, so we are content with the above version.

\begin{Remark}
\label{rem:highest_term}
When $\mathfrak{S}$ has empty boundary so that $\mathfrak{S}$ is a punctured surface, from the statement proved in \cite{AK} for the highest term of the image under ${\rm Tr}^\omega_\Delta$ of the elements of a certain natural basis of $\mathcal{S}^\omega(\mathfrak{S})$ \cite{DThurston}, it follows that ${\rm Tr}^\omega_\Delta \colon \mathcal{S}^\omega(\mathfrak{S}) \to \mathbb{BL}^\omega_\Delta$ is injective. This map is not too far from being surjective, and it is not hard to describe the `difference' between $\mathbb{BL}^\omega_\Delta$ and the image of ${\rm Tr}^\omega_\Delta$, although a proof has never appeared in the literature. We do not delve into this matter in the present paper.
\end{Remark}

The following special values of the quantum trace maps are of particular importance.
\begin{Definition}[quantized trace-of-monodromy: algebraic]
\label{def:quantized_trace-of-monodromy}
Let $\mathfrak{S} = \Sigma\setminus\mathcal{P}$ be a triangulable marked surface and $\Delta$ its ideal triangulation. Let $\gamma$ be an unoriented essential simple loop in~$\mathfrak{S}$, and $[\gamma]$ its isotopy class. Let $K_\gamma$ be the framed link in $\mathfrak{S} \times \mathbf{I}$ obtained by lifting~$\gamma$ to~$\mathfrak{S} \times \mathbf{I}$ at a~constant elevation, equipped with the upward vertical framing. The {\em quantized trace-of-monodromy} for $[\gamma]$ with respect to $\Delta$ is the element
\begin{align}
\label{eq:f_omega}
f^\omega_{[\gamma],\Delta} = f^\omega_{[\gamma],\Delta;\mathfrak{S}} := {\rm Tr}^\omega_{\Delta;\mathfrak{S}}([K_\gamma]) \in
\mathbb{BL}^\omega_\Delta \subset \bigl(\mathcal{Z}^\omega_\Delta\bigr)_{\rm bl} \subset \mathcal{Z}^\omega_\Delta.
\end{align}
\end{Definition}
For a fixed $\Delta$, it is easy to observe that \smash{$f^\omega_{[\gamma],\Delta}$} depends on the isotopy class $[\gamma]$, and not on a~particular choice of representative $\gamma$. Also, Proposition~\ref{prop:BW}\,\ref{QT2} guarantees that, for a~fixed~$[\gamma]$, the elements $f^\omega_{[\gamma],\Delta}$ for different ideal triangulations $\Delta$ are compatible with each other in the sense that
\[
f^\omega_{[\gamma],\Delta} = \Theta^\omega_{\Delta,\Delta'}\bigl(f^\omega_{[\gamma],\Delta'}\bigr)
\]
holds for any two ideal triangulations $\Delta$ and $\Delta'$, where $\Theta^\omega_{\Delta,\Delta'}$ is as in Proposition~\ref{prop:Theta}. Therefore, for each $[\gamma]$, the family of elements \smash{$\bigl\{ f^\omega_{[\gamma],\Delta} \bigr\}_{\Delta\colon \text{triangulations}}$}, which is based on the Bonahon--Wong quantum trace maps \cite{BW} and is suggested in \cite{AK}, can be interpreted as providing an answer~to the problem of finding a triangulation-independent (algebraic) quantization of the trace-of-monodromy function $f_{[\gamma]}$ on $\mathscr{T}^+(\mathfrak{S})$ along $[\gamma]$. This problem had sometimes been called a~quantum ordering problem, e.g., as stated in Conjecture~1 of \cite{CF2}.

Before this, only incomplete answers were known. See, e.g., \cite{F97, CF2} and the discussion in the next subsection of the present paper after Proposition~\ref{prop:BW_admissible_sum_formula}. Another attempt was the one by Teschner \cite{T}, who constructed a Hilbert space and a system of operators ${\bf L}_{[\gamma],\sigma}$ on it corresponding to a `marking' $\sigma$, which is a pants decomposition $\mathfrak{S}$ together with some extra data, and to each constituent loop $\gamma$ of $\sigma$, so that ${\bf L}_{[\gamma],\sigma}$ should play the role of the operator version of quantized trace-of-monodromy along $\gamma$. The approach of Teschner provides some nice ideas, which we partly absorb and translate into our language for our purposes, which will constitute a key technical tool of the present paper. However, Teschner's construction do not have a control over the product of such operators for loops that do intersect, and it is not clear which geometric space these operators are quantizing.

Using the Bonahon--Wong quantum trace maps and a basis of $\mathcal{S}^\omega(\mathfrak{S})$ mentioned in Remark~\ref{rem:highest_term}, one may construct a deformation quantization map for $\mathscr{T}^+(\mathfrak{S})$ in the sense of \ref{Q4}, in terms of the shear coordinates and (balanced) Chekhov--Fock algebras; see \cite{AK,CKKO,KS}. Although the main results of the present paper is deeply related to this deformation quantization map, they can be understood without invoking this quantization map.

\subsection{On quantized length operators}
\label{subsec:on_quantized_length_operators}

A missing piece in the quantization problem is \ref{Q5}, and as mentioned, it is much more doable for various Chekhov--Fock algebras, than for skein algebras directly. A general strategy is to construct representations of the Chekhov--Fock algebras and pull them back via Bonahon--Wong quantum trace maps. We will be reviewing and developing this in detail in a sequel to the present paper \cite{sequel}. Here we just give a slight preview, which may serve as a motivation to the present paper. Per each ideal triangulation $\Delta$, one aims to construct a complex Hilbert space $\mathscr{H}_\Delta$ on which the algebra $\mathcal{Z}^\omega_\Delta$ acts as a representation, in a certain sense. Per each element $u \in \mathcal{Z}^\omega_\Delta$, we would like to assign a densely defined operator on $\mathscr{H}_\Delta$. Pullback via the quantum trace map ${\rm Tr}^\omega_\Delta$ yields a representation of $\mathcal{S}^\omega(\mathfrak{S})$ on $\mathscr{H}_\Delta$. The domains of these operators are different from one another, and one can formulate the notion of the {\it Schwartz space} $\mathscr{S}_\Delta$ (see equation~\eqref{eq:Schwartz}, \cite[Definition~5.2]{FG09}), a dense subspace of $\mathscr{H}_\Delta$ playing the role of a common maximal domain, following the ideas of Fock and Goncharov \cite{FG09}. The data so far can be summarized as
\[
\rho^\hbar_\Delta = \rho^\hbar_{\Delta;\mathfrak{S}} \colon \ \mathcal{S}^\omega(\mathfrak{S}) \to \{\mbox{$\mathbb{C}$-linear operators $\mathscr{S}_\Delta \to \mathscr{S}_\Delta$}\},
\]
which is supposed to be an algebra homomorphism; $\omega$ is represented as the number ${\rm e}^{\pi {\rm i} \hbar/4}$, so that $q=\omega^4$ is represented as ${\rm e}^{\pi {\rm i} \hbar}$. For a surface $\mathfrak{S}$, we require the representations $\rho^\hbar_\Delta$ for different $\Delta$ to be compatible with each other, in the following sense. One looks for a unitary intertwining operator ${\bf K}^\hbar_{\Delta,\Delta'} \colon \mathscr{H}_{\Delta'} \to \mathscr{H}_\Delta$ per each pair of ideal triangulations $\Delta$ and $\Delta'$, which preserves the Schwartz subspaces, satisfy the intertwining equations
\[
{\bf K}^\hbar_{\Delta,\Delta'} \circ \rho^\hbar_{\Delta'}(u) = \rho^\hbar_\Delta(u) \circ {\bf K}^\hbar_{\Delta,\Delta'}, \qquad \forall u \in \mathcal{S}^\omega(\mathfrak{S}),
\]
and satisfies the consistency equations \smash{${\bf K}^\hbar_{\Delta,\Delta'} \circ {\bf K}^\hbar_{\Delta',\Delta''} = {\bf K}^\hbar_{\Delta,\Delta''}$} up to multiplicative constants. Hence, dealing with $\rho^\hbar_\Delta$ is equivalent to dealing with $\rho^\hbar_{\Delta'}$ in a suitable sense.

Each element $u \in \mathcal{S}^\omega(\mathfrak{S})$ can now be represented as an operator $\rho^\hbar_\Delta(u)$ on the Schwartz subspace $\mathscr{S}_\Delta$ of the Hilbert space $\mathscr{H}_\Delta$. Of course the most important case is when $u = [K_\gamma]$, the constant elevation lift of an essential simple loop $\gamma$ in $\mathfrak{S}$. This operator $\rho^\hbar_\Delta([K_\gamma])$ is what we should use to decompose the Hilbert space $\mathscr{H}_\Delta$, following the ideas of Verlinde \cite{Verlinde} and Fock and Goncharov \cite{FG09}. However, the issue here is that $\mathscr{S}_\Delta$ is not a maximal domain of $\rho^\hbar_\Delta([K_\gamma])$, hence $\rho^\hbar_\Delta([K_\gamma])$ is not self-adjoint on this domain. So the very first task is to find a canonical self-adjoint extension of $\rho^\hbar_\Delta([K_\gamma])$ to (a maximal domain in) $\mathscr{H}_\Delta$, and then investigate the spectral properties of this (densely defined) self-adjoint operator. The best situation would be that there is a unique self-adjoint extension, i.e., $\rho^\hbar_\Delta([K_\gamma])$ is essentially self-adjoint on $\mathscr{S}_\Delta$.

\begin{Definition}[quantized trace-of-monodromy operator]
If $\rho^\hbar_\Delta([K_\gamma])$ is essentially self\nobreakdash-ad\-joint on the Schwartz space $\mathscr{S}_\Delta$, denote by ${\bf f}^\hbar_{[\gamma],\Delta}$ its unique self-adjoint extension in $\mathscr{H}_\Delta$.
\end{Definition}
We suggest that this operator ${\bf f}^\hbar_{[\gamma],\Delta}$, the operator version of $f^\omega_{[\gamma],\Delta}$, is the operator quantization of the trace-of-monodromy function $f_{[\gamma]}$, which one should use to realize Verlinde's suggestion~\cite{Verlinde} and to answer the modular functor conjecture of Fock and Goncharov \cite{FG09, F97}. For this operator to be a good candidate for the quantized trace-of-monodromy function, there should also exist a corresponding {\em quantized length operator} \smash{${\bf l}^\hbar_{[\gamma],\Delta}$}, which is supposed to be a~self-adjoint opera\-tor~on~$\mathscr{H}_\Delta$ having simple spectrum $[0,\infty)$, satisfying
\begin{align}
\label{eq:l_hbar_operator_definition}
{\bf f}^\hbar_{[\gamma],\Delta} = 2\cosh \bigl({\bf l}^\hbar_{[\gamma],\Delta}/2\bigr),
\end{align}
where the right-hand side is made sense by functional calculus of \smash{${\bf l}^\hbar_{[\gamma],\Delta}$}. This operator \smash{${\bf l}^\hbar_{[\gamma],\Delta}$} can be viewed as a quantization of the length function $l_{[\gamma]}$, as equation~\eqref{eq:l_hbar_operator_definition} is a quantum counterpart of the classical equation~\eqref{eq:trace_of_monodromy_as_length}, by inspection.

The existence and uniqueness of such an operator \smash{${\bf l}^\hbar_{[\gamma],\Delta}$} is guaranteed if \smash{${\bf f}^\hbar_{[\gamma],\Delta}$} has simple spectrum $[2,\infty)$. In this case, the spectral decomposition for \smash{${\bf l}^\hbar_{[\gamma],\Delta}$} gives the direct integral decomposition
\begin{align}
\label{eq:direct_integral_decomposition}
\mathscr{H}_\Delta \to \int^\oplus_{[0,\infty)} (\mathscr{H}_\Delta)_\chi \, {\rm d}\chi,
\end{align}
which appears in Fock and Goncharov's modular functor conjecture \cite{F97, FG09}, which stipulates further equivariance conditions. Importantly, each slice $(\mathscr{H}_\Delta)_\chi$ should be a representation of the skein algebra $\mathcal{S}^\omega(\mathfrak{S}')$ of the surface $\mathfrak{S}'$ obtained by removing $\gamma$ from $\mathfrak{S}$ and shrinking the two new holes to two new punctures, such that for the peripheral loop $\eta$ around each of these two new punctures, the element $[K_\eta]$ should act as the scalar $2\cosh(\chi/2)$. In particular, one sees that the correct algebras that should be used to describe such representations are the relative skein algebras (Definition~\ref{def:relative_skein_algebra}). In order for the above direct integral decomposition to be compatible with the natural map between the skein algebras of $\mathfrak{S}$ and $\mathfrak{S}'$, one should make sure that for each essential simple loop $\xi$ in $\mathfrak{S}$ that is disjoint from $\gamma$, the operators ${\bf f}^\hbar_{[\xi],\Delta}$ and ${\bf f}^\hbar_{[\gamma],\Delta}$ should strongly commute, i.e., their spectral projections should commute. We notice that this strong commutativity has been in the list of sought-for properties, as can be seen in \cite[p.~160, item~6]{CF2} \cite[p.~26, item~3]{F97} \cite[Section~15.1, item~(b)]{T}.

What have been discussed should be the first major steps toward a direct settlement of the modular functor conjecture of Fock and Goncharov \cite{F97, FG09}, going back to Verlinde \cite{Verlinde}. We formulate these steps here.
\begin{Conjecture}[steps toward modular functor conjecture]
\label{conjecture:sequel}
Let $\mathfrak{S} = \Sigma \setminus \mathcal{P}$ be a triangulable marked surface. Let $\Delta$ be an ideal triangulation. Let $\hbar \in \mathbb{R}_{>0}$.
\begin{enumerate}[label={\rm (\arabic*)}]\itemsep=0pt
\item For each essential simple loop $\gamma$ in $\mathfrak{S}$, the operator $\rho^\hbar_\Delta([K_\gamma])$ is essentially self-adjoint on the Schwartz space $\mathscr{S}_\Delta \subset \mathscr{H}_\Delta$.

\item The unique self-adjoint extension \smash{${\bf f}^\hbar_{[\gamma],\Delta}$} of $\rho^\hbar_\Delta([K_\gamma])$ in $\mathscr{H}_\Delta$, called the quantized trace\nobreakdash-of-mo\-no\-dromy operator, has simple spectrum $[2,\infty)$, so that there exists a unique self\nobreakdash-ad\-joint operator ${\bf l}^\hbar_{[\gamma],\Delta}$ having simple spectrum $[0,\infty)$ and satisfying equation~\eqref{eq:l_hbar_operator_definition}, which is called the quantized length operator.

\item For any two disjoint essential simple loops $\gamma$ and $\xi$ in $\mathfrak{S}$, the operators ${\bf f}^\hbar_{[\gamma],\Delta}$ and ${\bf f}^\hbar_{[\xi],\Delta}$ strongly commute.
\end{enumerate}
\end{Conjecture}

The present paper is to establish some algebraic properties of the quantized trace\nobreakdash-of\nobreakdash-mo\-no\-dromy $f^\omega_{[\gamma],\Delta}$, which we will use to tackle the above conjecture in a sequel paper \cite{sequel}; we present some key ideas in Section~\ref{sec:consequences} of the present paper. One general strategy that we take in the present paper is that, for a loop or a pair of loops, we look for a suitable ideal triangulation for which the quantized trace-of-monodromy have nice properties.

\subsection{State-sum formula for the quantized trace-of-monodromy}
\label{subsec:state-sum}

We collect some results on computation of the values of the quantum trace map ${\rm Tr}^\omega_\Delta \colon \mathcal{S}^\omega(\mathfrak{S}) \to \bigl(\mathcal{Z}^\omega_\Delta\bigr)_{\rm bl}$, for any triangulable marked surface $\mathfrak{S}$ and its ideal triangulation $\Delta$. For our purposes, it suffices to deal with the value ${\rm Tr}^\omega_\Delta([K_\gamma])$ when $K_\gamma$ is a constant-elevation lift in $\mathfrak{S} \times {\bf I}$ of an essential simple loop $\gamma$ in $\mathfrak{S}$ with an upward vertical framing, in view of equation~\eqref{eq:f_omega} of Definition~\ref{def:quantized_trace-of-monodromy}. By applying an isotopy if necessary, assume that $\gamma$ meets $\Delta$ tranversally at a~minimum number of points in the following sense.
\begin{Definition}
\label{def:minimal_and_transversal}
Let $\Delta$ be a collection of ideal arcs in a marked surface $\mathfrak{S}$; here in this definition, $\Delta$ is not necessarily an ideal triangulation. Let $\gamma$ be an essential simple loop in $\mathfrak{S}$. Denote by $|\gamma \cap \Delta|$ the cardinality of the intersection $\gamma \cap \Delta$. Let the {\em intersection number} ${\rm int}(\gamma, \Delta)$ between $\gamma$ and $\Delta$ be the minimum of $|\gamma' \cap \Delta|$, where $\gamma'$ runs through all essential simple loops in $\mathfrak{S}$ isotopic to $\gamma$. We say $\gamma$ is in a {\em minimal position} with respect to $\Delta$ if $|\gamma \cap \Delta| = {\rm int}(\gamma,\Delta)$. We say $\gamma$ is in a {\em transversal position} with respect to $\Delta$ if $\gamma$ meets transversally with $\Delta$.

In case $\Delta$ consists of a single element, say $\Delta = \{i\}$, then we write ${\rm int}(\gamma,\{i\}) = {\rm int}(\gamma,i)$.
\end{Definition}
Suppose that $\gamma$ is in a minimal and transversal position with respect to $\Delta$. Each point of intersection of $\gamma$ and $\Delta$ is called a {\em juncture}. The junctures divide $\gamma$ into {\em loop segments}. From the minimality condition, it follows that each loop segment is a simple path in an ideal triangle of $\Delta$ connecting two distinct sides of this triangle (see \cite[Lemma~2.8]{CKKO} for an argument dealing with the case of a self-folded triangle). A {\em juncture-state} is a map
\[
J \colon \ \gamma \cap \Delta \to \{+,-\},
\]
which assigns a sign to each juncture. For a juncture-state $J$, define the element $v_J \in N_\Delta$ as
\begin{align}
\label{eq:v_J}
&v_J := \sum_{i\in \Delta} \biggl( \underset{x\in \gamma \cap i}{\textstyle \sum} J(x) \biggr) \til{v}_i \in N_\Delta, \qquad \mbox{where} \\
\label{eq:til_v_i}
&\til{v}_i := \begin{cases}
v_i & \mbox{if } i \mbox{ is not a self-folded arc}, \\
v_i - v_{\pi(i)} & \mbox{if } i \mbox{ is a self-folded arc, with } \pi \mbox{ as in equation}~\eqref{eq:pi},
\end{cases}
\end{align}
where the signs $+$ and $-$ are viewed as numbers $+1$ and $-1$, respectively. Define the corresponding Weyl-ordered Laurent monomial $Z^J \in \mathcal{Z}^\omega_\Delta$ as
\begin{align}
\label{eq:Z_J}
Z^J := Z_{v_J} \in \mathcal{Z}^\omega_\Delta.
\end{align}
That is, \smash{$Z^J = \bigl[ \prod_{i \in \Delta} \til{Z}_i^{a_i} \bigr]_{\rm Weyl}$} (equation~\eqref{eq:Weyl-ordering}), where for each ideal arc $i\in \Delta$, the power \smash{$a_i=\sum_{x\in \gamma \cap i} J(x)$} is the net sum of juncture-state values of the junctures lying on $i$, and $\til{Z}_i := Z_{\til{v}_i}$ (which equals $Z_i$ if $i$ is not self-folded and equals $Z_i Z_{\pi(i)}$ if $i$ is self-folded).

\begin{figure}[!ht]
\centering\hspace*{-20mm}
\begin{subfigure}[b]{0.32\textwidth}\centering
\scalebox{0.7}{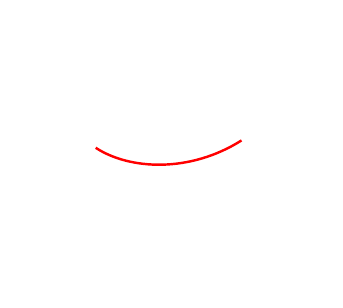}
\caption*{(1) non-self-folded ideal triangle}
\end{subfigure}
\qquad
\begin{subfigure}[b]{0.44\textwidth}\centering
\scalebox{0.7}{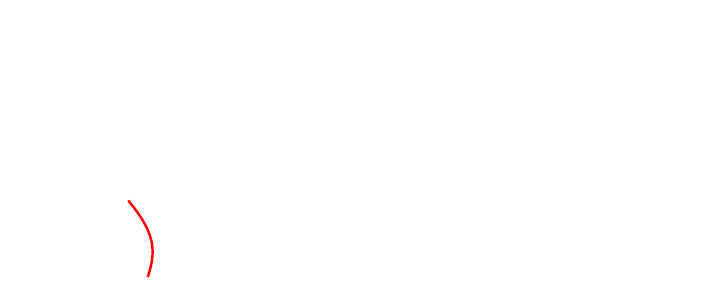}
\caption*{(2) self-folded triangle}
\end{subfigure}
\caption{Non-admissible juncture-states on a loop segment.}
\label{fig:non-admissible_juncture-states}
\end{figure}

A juncture-state $J$ is said to be {\em admissible} if there is no loop segment such that the juncture-state values at its endpoints are as in Figure~\ref{fig:non-admissible_juncture-states}. So, for each loop segment of $\gamma$ in a triangle, there are only three allowed (admissible) possibilities for the juncture-state values at its endpoints. The easiest example of an admissible juncture-state is a constant juncture-state, which either assigns $+$ to all junctures, or assigns $-$ to all junctures. These constant juncture-states are often used, so we give names to them: for each $\epsilon \in \{+,-\}$, denote by $J^\epsilon_{\gamma} = J^\epsilon_{\gamma,\Delta} \colon \gamma \cap \Delta \to \{+,-\}$ the $\epsilon$-valued constant juncture-state:
\begin{align}
\label{eq:constant_juncture-state}
J_\gamma^\epsilon(x) = \epsilon, \qquad \forall x\in \gamma \cap \Delta.
\end{align}

The following is a special case of the {\it state-sum} formula for the quantized trace-of-monodromy written in terms of admissible juncture-states.
\begin{Proposition}[admissible state-sum formula for quantized trace-of-monodromy along a loop, non-redundant case]
\label{prop:BW_admissible_sum_formula}
Let $\Delta$ be an ideal triangulation of a triangulable marked surface $\mathfrak{S} = \Sigma\setminus\mathcal{P}$. Let $\gamma$ be an essential simple loop in $\mathfrak{S}$, in a minimal and transversal position with respect to $\Delta$ $($Definition~{\rm \ref{def:minimal_and_transversal})}. Let $K_\gamma$ be a constant elevation lift of $\gamma$ in $\mathfrak{S} \times {\bf I}$ with upward vertical framing. Suppose that the following {\em non-redundancy} condition holds:
\begin{align}
\label{eq:admissible_non-redundancy_condition}
Z^J \neq Z^{J'}
\end{align}
(i.e., $v_J\neq v_{J'}$) whenever $J$ and $J'$ are distinct admissible juncture-states.
Then the quantized trace-of-monodromy $f^\omega_{[\gamma],\Delta} = {\rm Tr}^\omega_\Delta([K_\gamma]) \in \mathcal{Z}^\omega_\Delta$ $($Definition~{\rm \ref{def:quantized_trace-of-monodromy})} is given by
\begin{align}
\label{eq:admissible_sum_formula}
f^\omega_{[\gamma],\Delta}
= \sum_J Z^J,
\end{align}
where the sum is over all admissible juncture-states $J \colon \gamma\cap \Delta \to \{+,-\}$.
\end{Proposition}

\begin{proof}
For the classical case $\omega=1$, it is known that equation~\eqref{eq:admissible_sum_formula} holds for {\it any} essential simple loop, without the assumption in equation~\eqref{eq:admissible_non-redundancy_condition}; see \cite[Lemma~4, and Section~6 with $\omega=1$]{BW}. We~remark that, in view of \ref{QT1} of Proposition~\ref{prop:BW}, this implies Proposition~\ref{prop:classical_Laurentness}. For the quantum case, in \cite{CKKO} it is proved that, for any essential simple loop $\gamma$ and any ideal triangulation $\Delta$, \smash{$f^\omega_{[\gamma],\Delta}$}~is a~linear combination of Weyl-ordered Laurent monomials $Z_v$, $v\in N_\Delta$, with coefficients in $\mathbb{Z}_{\ge 0}\bigl[\omega^{\pm 1}\bigr]$; here the emphasis is on $\mathbb{Z}_{\ge 0}\bigl[\omega^{\pm 1}\bigr]$, as opposed to $\mathbb{Z}\bigl[\omega^{\pm 1}\bigr]$. Thus,
\[
f^\omega_{[\gamma],\Delta} = \sum_{r=1}^m \omega^{\alpha_r} Z_{v_r}
\]
 for some sequences $v_1,\dots,v_m \in N_\Delta$ and $\alpha_1,\dots,\alpha_m \in \mathbb{Z}$. Note \smash{$\bigl(f^\omega_{[\gamma],\Delta}\bigr)^* = \sum_{r=1}^m \omega^{-\alpha_r} Z_{v_r}$}. Meanwhile, \cite[Lemma~2.30]{AK} says
\[
\bigl(f^\omega_{[\gamma],\Delta}\bigr)^* = f^\omega_{[\gamma],\Delta}.
\]
Now if we assume the condition in equation~\eqref{eq:admissible_non-redundancy_condition}, then $v_1,\dots,v_m$ are mutually distinct, so from \smash{$\bigl(f^\omega_{[\gamma],\Delta}\bigr)^* = f^\omega_{[\gamma],\Delta}$} it follows that $\omega^{\alpha_r} = \omega^{-\alpha_r}$, hence $\alpha_r=0$, $\forall r$, and \smash{$f^\omega_{[\gamma],\Delta} = \sum_{r=1}^m Z_{v_r}$}. Putting $\omega=1$, we obtain $f_{[\gamma]} = \sum_{r=1}^m Z_{v_r}$; since $f_{[\gamma]} = \sum_J Z^J$ is summed over all admissible juncture-states $J$, we see that the set $\{v_J \mid \mbox{admissible juncture-states $J$}\}$ coincides with $\{ v_r \mid r=1,\dots,m\}$. Therefore, we get \smash{$f^\omega_{[\gamma],\Delta} = \sum_J Z^J$} summed over all admissible juncture-states $J$, as desired in equation~\eqref{eq:admissible_sum_formula}.
\end{proof}

The above proposition says that, in case no two Laurent monomial terms of $f_{[\gamma]}$ written with respect to $\Delta$ coincide with each other, then the quantized version $f^\omega_{[\gamma],\Delta}$ is obtained by replacing each term of $f_{[\gamma]}$ by the corresponding Weyl-ordered quantum Laurent monomial, i.e., by the termwise Weyl-ordering. However, the condition equation~\eqref{eq:admissible_non-redundancy_condition} does not hold in general, in which case the construction of \smash{$f^\omega_{[\gamma],\Delta}$} is much more complicated, where in fact the heart of the whole construction of Bonahon and Wong \cite{BW} lies in.

Observe that equation~\eqref{eq:admissible_non-redundancy_condition} is a condition on both $\gamma$ and $\Delta$, not just on $\gamma$. For each essential simple loop $\gamma$, we shall find a suitable $\Delta$ such that equation~\eqref{eq:admissible_non-redundancy_condition} holds, so that we can use Proposition~\ref{prop:BW_admissible_sum_formula} to conveniently compute $f^\omega_{[\gamma],\Delta}$. In fact, among such $\Delta$ we will try to find an even better one so that either the number of terms of $f^\omega_{[\gamma],\Delta}$ is small, or the terms can be organized in a special way that fits our purpose of eventually studying the analytic properties.

We remark that the termwise Weyl-ordered expression in the right-hand side of~\eqref{eq:admissible_sum_formula} was historically the first naive attempt to quantize the trace-of-monodromy function $f_{[\gamma]}$; in particular, the Bonahon--Wong quantum trace was not known then. For example, the construction for a~quantized trace-of-monodromy for $\gamma$ suggested in \cite{F97, CF2} is to first find an ideal triangulation~$\Delta$ such that there exists an ideal arc of $\Delta$ that meets $\gamma$ exactly once, and define \smash{$f^\omega_{[\gamma],\Delta}$} by the term-by-term Weyl-ordered formula in equation~\eqref{eq:admissible_sum_formula}. It turns out that the quantized trace-of-monodromy constructed via the Bonahon--Wong quantum trace indeed satisfies this (see~\cite{AK}), but the problem is that for some loop $\gamma$, there does not exist such an ideal triangulation~$\Delta$, contrary to what was claimed in \cite[below equation~(57)]{F97}. Proposition~\ref{prop:BW_admissible_sum_formula} presents another condition on $\gamma$ and $\Delta$, namely the non-redundancy condition in equation~\eqref{eq:admissible_non-redundancy_condition}, that is a sufficient condition for the termwise Weyl-ordered formula to work. In fact, we shall see that one corollary (Corollary~\ref{cor:AP}) of our arguments and proofs to be established in the following sections implies the following, which might be useful.
\begin{Proposition}[termwise Weyl-ordered quantized trace-of-monodromy; Proposition~\ref{prop:intro_term-by-term}]
\label{prop:term-by-term_Weyl-ordered_f_omega_gamma_Delta}
Let $\mathfrak{S} = \Sigma\setminus\mathcal{P}$ be a triangulable marked surface, and let $\gamma$ be an essential simple loop in $\mathfrak{S}$. Then there exists an ideal triangulation $\Delta$ of $\mathfrak{S}$ such that equation~\eqref{eq:admissible_non-redundancy_condition} holds, and hence equation~\eqref{eq:admissible_sum_formula} holds.
\end{Proposition}

We also remark that there is a more general version of the state-sum formula for the value of the quantum trace ${\rm Tr}^\omega_\Delta([K])$ for any framed link $K$, or even for any framed tangle (which we haven't dealt with in this paper); see \cite{BW}. Although mostly we only need Proposition~\ref{prop:BW_admissible_sum_formula} for the present paper, we find it convenient to state the following weak consequence of the general version of the state-sum formula, which we will be using in later sections.
\begin{Proposition}[support of sum expression of quantized trace-of-monodromy; \cite{BW}]
\label{prop:weak_sum}
Let $\mathfrak{S}$ be a triangulable marked surface, $\Delta$ an ideal triangulation of $\mathfrak{S}$, and $\gamma$ an essential simple loop in~$\mathfrak{S}$ in a transversal position with respect to $\Delta$. Define the $\mathbb{Z}$-submodule $N_\gamma$ of $N$ as the $\mathbb{Z}$-span of $\{\til{v}_i \mid \mbox{$i$ meets $\gamma$}\}$, where $\til{v}_i$ is as in equation~\eqref{eq:til_v_i}. Then
\[
f^\omega_{[\gamma],\Delta} = \sum_{v\in N_\gamma} a_v Z_v
\]
for some $a_v \in \mathbb{Z}\bigl[\omega^{\pm 1}\bigr]$ which are zero for all but finitely many $v\in N_\gamma$.
\end{Proposition}

\section{Simple loops and ideal triangulations}

As mentioned, one strategy to study the algebraic structure and properties of the quantized trace-of-monodromy $f^\omega_{[\gamma],\Delta} \in \mathcal{Z}^\omega_\Delta$ (Definition~\ref{def:quantized_trace-of-monodromy}) along a loop $\gamma$ with respect to an ideal triangulation~$\Delta$ is to find a nice ideal triangulation $\Delta$ per each loop $\gamma$, so that the investigation of~$f^\omega_{[\gamma],\Delta}$ is doable. As a topological step, in the present section we establish certain classifications of single simple loops, and pairs of disjoint simple loops, keeping in mind their relative positions with respect to ideal triangulations.

\subsection{A topological classification of simple loops}
\label{subsec:topological_classification_of_simple_loops}

In order to use throughout the paper, we first clearly set the terminology to refer to some surfaces.
\begin{Definition}
\label{def:surface_terminology}
Let $\Sigma$ be a compact oriented surface with possibly empty boundary. Each boundary component of $\Sigma$ is called a {\em hole}.
\begin{enumerate}[label={\rm (\arabic*)}]\itemsep=0pt
\item For integers $n\ge 1$ and $g\ge 0$, we say $\Sigma$ is an {\em $n$-holed surface} of {\em genus $g$}, if it is obtained from a closed oriented surface of genus $g$ by removing $n$ open discs whose closures are mutually disjoint.

\item If $g=0$, we say $\Sigma$ is an {\em $n$-holed sphere}. If $g=1$, we say $\Sigma$ is an {\em $n$-holed torus}.

\item If $g=0$ and $n=2$, we say $\Sigma$ is a {\em closed annulus}.

\end{enumerate}
\end{Definition}
For small $n$, the number $n$ is written as a word; the adjective $n$-holed is written as one-holed, two-holed, three-holed, etc. Next, we introduce the following, which will become handy.

\begin{Definition}
\label{def:minimal-marked}
Let $(\Sigma,\mathcal{P})$ be a marked surface, and let $\mathfrak{S} = \Sigma\setminus\mathcal{P}$. We say $(\Sigma,\mathcal{P})$ or $\mathfrak{S}$ is {\em minimal-marked} if $\mathcal{P}$ is a finite subset of $\Sigma$ of the smallest cardinality, among all finite sets $\mathcal{P}' \subset \Sigma$ such that $(\Sigma,\mathcal{P}')$ is a marked surface.
\end{Definition}

The notion of minimal-marked surface $(\Sigma,\mathcal{P})$ becomes particularly convenient when $\Sigma$ has a~nonempty boundary, in which case $(\Sigma,\mathcal{P})$ is minimal-marked if $\mathcal{P}$ has exactly one point on each boundary component of $\Sigma$ and no puncture. Hence, if $\Sigma$ is an $n$-holed surface, then a~minimal-marked surface $(\Sigma,\mathcal{P})$ has $n$ marked points in total. For example, a~minimal-marked one-holed surface $(\Sigma,\mathcal{P})$ (or the corresponding $\mathfrak{S} = \Sigma\setminus\mathcal{P}$) has a single marked point, lying in its unique boundary component of $\Sigma$. A minimal-marked annulus $(\Sigma, \mathcal{P})$ has two marked points, one in each of the two boundary components of $\Sigma$.

We also often deal with punctured surfaces $\mathfrak{S} = \Sigma\setminus\mathcal{P}$, i.e., when $\Sigma$ has empty boundary.
\begin{Definition}
Let $n\ge 1$ and $g\ge 0$. Let $(\Sigma,\mathcal{P})$ be a marked surface such that $\Sigma$ is a closed oriented surface of genus $g$, without boundary, and $\mathcal{P}$ consists of $n$ points. We say $\mathfrak{S} = \Sigma\setminus\mathcal{P}$ is an {\em $n$-times-punctured surface of genus $g$}.
\end{Definition}
When $n=1$ or $2$, `$n$-times-punctured' can be replaced by `once-punctured' or `twice\nobreakdash-punc\-tured', respectively. When $g=0$ or $1$, `surface of genus $g$' can be replaced by `sphere' or `torus', respectively. Hence, for example, one can consider once-punctured torus, twice-punctured torus, once-punctured surface of genus $2$, etc.

Below we establish some terminology for subsurfaces.

\begin{Definition}
\label{def:subsurfaces}
Let $(\Sigma, \mathcal{P})$ be a marked surface, and let $\mathfrak{S} = \Sigma\setminus\mathcal{P}$.
\begin{enumerate}[label={\rm (\arabic*)}]\itemsep=0pt
\item A {\em marked subsurface} of $(\Sigma,\mathcal{P})$ is a marked surface $(\Sigma',\mathcal{P}')$, where $\Sigma'$ is an embedded oriented subsurface of $\Sigma$ and $\mathcal{P}' \subset \mathcal{P}$. In this case, we say $\mathfrak{S}' = \Sigma'\setminus \mathcal{P}'$ is a marked subsurface of $\mathfrak{S} = \Sigma \setminus \mathcal{P}$.

\item Suppose $(\Sigma,\mathcal{P})$ is triangulable. A marked subsurface $(\Sigma',\mathcal{P}')$ of $(\Sigma,\mathcal{P})$ is said to be {\em triangle-compatible} if it is triangulable and the closure of its complement in $(\Sigma,\mathcal{P})$ is either empty or also a triangulable marked subsurface of $(\Sigma,\mathcal{P})$.

\end{enumerate}
\end{Definition}

The following notion will become useful.
\begin{Definition}\label{def:cutting}
For a compact oriented surface $\Sigma$ with possibly empty boundary, consider an essential simple loop $\gamma$ in the interior of $\Sigma$. By a surface obtained from $\Sigma$ by {\em cutting along $\gamma$}, we~mean a surface diffeomorphic to $\Sigma_\gamma$ defined as follows. Let $N$ be a tubular neighborhood of~$\gamma$ in~$\Sigma$, so that in particular, $N$ is diffeomorphic to an open annulus. Let $\Sigma_\gamma$ be the closure in~$\Sigma$ of the complement $\Sigma\setminus N$.

Let $\Gamma$ be a finite collection of mutually disjoint essential simple loops in $\Sigma$. Let $\Sigma_\Gamma$ be a~surface obtained from $\Gamma$ by cutting along all loops of $\Gamma$. We say that $\Gamma$ {\em cuts out} a surface $\Sigma'$ if there is a connected component of $\Sigma_\Gamma$ that is diffeomorphic to $\Sigma'$ and whose boundary contains at least one boundary circle coming from a loop of $\Gamma$.
\end{Definition}
Then $\Sigma_\gamma$ has two more boundary components than $\Sigma$, and these new boundary components are said to {\em come from}~$\gamma$.

Meanwhile, when $(\Sigma,\mathcal{P})$ is a marked surface and $i$ is an ideal arc of $(\Sigma,\mathcal{P})$ whose interior lies in the interior of $\Sigma$, {\em cutting $\Sigma$ along $i$} yields a well-defined marked surface $(\Sigma_i,\mathcal{P}_i)$, together with a projection map $\Sigma_i \to \Sigma$, restricting to a projection $\mathcal{P}_i \to \mathcal{P}$. This marked surface $(\Sigma_i,\mathcal{P}_i)$ together with the projection maps is uniquely determined up to isomorphism. This cutting process, which is somewhat different from Definition~\ref{def:cutting} (see, e.g., \cite[Definition and Lemma~3.12]{KLS}), will also be used.

We recall a standard terminology from the literature.
\begin{Definition}
Let $\mathfrak{S}$ be a surface which is not necessarily connected, and let $\gamma$ be a simple loop in $\mathfrak{S}$. We say $\gamma$ is {\em separating} if $\mathfrak{S}_0 \setminus\gamma$ is disconnected, where $\mathfrak{S}_0$ is the connected component of $\mathfrak{S}$ containing $\gamma$. We say $\gamma$ is {\em non-separating} otherwise.
\end{Definition}

We introduce the following terminology, as an analog of a peripheral loop (Definition~\ref{def:essential_and_peripheral}).
\begin{Definition}
\label{def:hole-surrounding}
Let $\mathfrak{S} = \Sigma\setminus\mathcal{P}$ be a marked surface, where $\Sigma$ has a nonempty boundary. An~essential simple loop in $\mathfrak{S}$ or $\Sigma$ is called a {\em hole surrounding loop} if it is a loop that is freely homotopic in $\Sigma$ to a component of $\partial \Sigma$.
\end{Definition}

We now establish a classification of essential simple loops in a surface, suited to our purposes.
\begin{Proposition}[a classification of essential simple loops]
\label{prop:loop_classification}
Let $\gamma$ be an essential simple loop in a triangulable marked surface $\mathfrak{S} = \Sigma\setminus \mathcal{P}$. Assume $\mathfrak{S}$ is connected. Then at least one of the following holds:
\begin{enumerate}[label={\rm (L\arabic*)}]\itemsep=0pt
\item\label{L1} $\gamma$ is a peripheral loop around a puncture.

\item\label{L2} $\gamma$ is a separating loop, and each of the two connected components of $\mathfrak{S} \setminus \gamma$ contains a point of $\mathcal{P}$.

\item\label{L3} $\gamma$ is a non-separating loop and $\mathcal{P}$ has at least two points.

\item\label{L4} $\mathfrak{S}$ is a once-punctured torus and $\gamma$ is a non-separating loop in $\mathfrak{S}$.

\item\label{L5} $\gamma$ is contained in a triangle-compatible marked subsurface $\mathfrak{S}'$ of $\mathfrak{S}$ that is isomorphic to a~minimal-marked one-holed torus, and $\gamma$ is a non-separating loop of $\mathfrak{S}'$.

\item\label{L6} $\gamma$ is a separating loop contained in a triangle-compatible marked subsurface $\mathfrak{S}'$ of $\mathfrak{S}$ that is isomorphic to a minimal-marked one-holed surface of genus $g\ge 1$, and $\gamma$ is a hole surrounding loop of $\mathfrak{S}'$.
\end{enumerate}
\end{Proposition}
This can be proved by simple case-by-case topological arguments, together with the classification of oriented surfaces. An interested reader can find a proof in a previous version of this paper \cite{ver3}.

When investigating the structure of the quantized trace-of-monodromy \smash{$f^\omega_{[\gamma],\Delta}$} for a single essential simple loop $\gamma$, we shall deal with each of the cases \ref{L1}--\ref{L6} in the classification result for loops $\gamma$ obtained in Proposition~\ref{prop:loop_classification}. In fact, the cases \ref{L2}--\ref{L5} can be characterized by the following property, whose proof is straightforward and hence omitted (see \cite{ver3} for a proof).
\begin{Proposition}[non-peripheral loops with intersection number 2 with a triangulation]
\label{prop:2_to_5_have_intersection_2}
If an essential simple loop $\gamma$ in a triangulable marked surface $\mathfrak{S} = \Sigma\setminus\mathcal{P}$ falls into one of the items \ref{L2}--\ref{L5} but not into the item \ref{L1} in Proposition~{\rm \ref{prop:loop_classification}}, then there exists an ideal triangulation~$\Delta$ of~$\mathfrak{S}$ such that ${\rm int}(\gamma,\Delta)=2$.
\end{Proposition}

We shall see in a later section that for a loop $\gamma$ and an ideal triangulation $\Delta$ as appearing in Proposition~\ref{prop:2_to_5_have_intersection_2}, the quantized trace-of-monodromy $f^\omega_{[\gamma],\Delta}\in \mathcal{Z}^\omega_\Delta$ (Definition~\ref{def:quantized_trace-of-monodromy}) takes a~particularly simple form, which is well suited to the study of analytic properties of the corresponding operator, as will be outlined in Section~\ref{subsec:toward_spectral_properties}.

\subsection{On pairs of simple loops}
\label{subsec:pais_of_simple_loops}

For each pair of disjoint essential simple loops $\gamma_1$ and $\gamma_2$ in $\mathfrak{S}$, we would like to show that there exists an ideal triangulation $\Delta$ of $\mathfrak{S}$ such that $f^\omega_{[\gamma_1],\Delta}$ and $f^\omega_{[\gamma_2],\Delta}$ commute in a certain sense that is stronger than a mere algebraic commutativity in the algebra $\mathcal{Z}^\omega_\Delta$.

When one of $\gamma_1$ and $\gamma_2$ is a peripheral loop, the situation will be easy. So let us assume that none of them are peripheral. The next easiest case is the following.
\begin{Definition}
\label{def:triangle-disjoint}
Let $\gamma_1$ and $\gamma_2$ be disjoint essential simple loops in a triangulable marked surface $\mathfrak{S} = \Sigma \setminus \mathcal{P}$. We say that $\gamma_1$ and $\gamma_2$ are {\em triangle-disjoint} with respect to an ideal triangulation~$\Delta$ of $\mathfrak{S}$ if, when they are isotoped into minimal and transversal positions with respect to $\Delta$, the collection of ideal triangles of $\Delta$ meeting $\gamma_1$ and the collection of ideal triangles of $\Delta$ meeting~$\gamma_2$ have no common ideal triangle.
\end{Definition}
In this case, roughly speaking, the variables $Z_i$ involved in $f^\omega_{[\gamma_1],\Delta}$ and those involved in $f^\omega_{[\gamma_2],\Delta}$ are disjoint and commute in a strong sense. We will make this precise in Lemma~\ref{lem:triangle-disjoint_implies_strong_commutativity} algebraically, and in the sequel paper \cite{sequel} analytically.

However, not every pair of disjoint essential simple loops are triangle-disjoint for some ideal triangulation. Here we focus on classifying such pairs, which are not triangle-disjoint with respect to any ideal triangulation. Then in a later section (Section~\ref{subsec:Teschner_recursion}), we will devise a new type of algebraic commutativity to deal with them.

To lessen notational confusion, we name the two loops as $\xi_1$ and $\xi_2$ here.

\begin{Proposition}[a classification of pairs of disjoint loops]
\label{prop:classification_of_pair_of_loops}
Let $\xi_1$ and $\xi_2$ be non-peripheral essential simple loops $($Definition~{\rm \ref{def:essential_and_peripheral})} in a connected triangulable marked surface $\mathfrak{S} = \Sigma\setminus\mathcal{P}$. Suppose that $\xi_1$ and $\xi_2$ are disjoint from and non-homotopic to each other. Then one of the following holds:
\begin{enumerate}[label={\rm (PL\arabic*)}]\itemsep=0pt\setlength{\leftskip}{0.22cm}
\item\label{PL1} $\xi_1$ and $\xi_2$ are triangle-disjoint with respect to some ideal triangulation $\Delta$ of $\mathfrak{S}$ $($Defini\-tion~{\rm \ref{def:triangle-disjoint})}.

\item\label{PL2} There is a triangle-compatible marked subsurface $\mathfrak{S}'$ $($Definition~{\rm \ref{def:subsurfaces})} of $\mathfrak{S}$ isomorphic to a minimal-marked $($Definition~{\rm \ref{def:minimal-marked})} one-holed torus $($Definition~{\rm \ref{def:surface_terminology})} such that $\xi_1$ is a hole surrounding loop $($Definition~{\rm \ref{def:hole-surrounding})} in $\mathfrak{S}'$ and $\xi_2$ is a non-separating loop in $\mathfrak{S}'$.

\item\label{PL3} There is a triangle-compatible marked subsurface $\mathfrak{S}'$ of $\mathfrak{S}$ isomorphic to a minimal-marked one-holed surface of genus $g\ge 2$ such that $\xi_1$ is a hole surrounding loop in $\mathfrak{S}'$ and $\xi_2$ is a~non-separating loop in $\mathfrak{S}'$.

\item\label{PL4} There is a triangle-compatible marked subsurface $\mathfrak{S}' = \Sigma'\setminus \mathcal{P}'$ of $\mathfrak{S}$ isomorphic to a~mini\-mal-marked one-holed surface of genus $g\ge 2$ such that $\xi_1$ and $\xi_2$ are non-separating loops in~$\mathfrak{S}'$, which together with the boundary circle of $\Sigma'$ bounds a three-holed sphere subsurface.

\item\label{PL5} $\mathfrak{S}$ is a once-punctured surface of genus $g\ge 2$ without boundary, where $\xi_1$ and $\xi_2$ are non-separating loops such that $\{\xi_1,\xi_2\}$ cuts out a two-holed genus $0$ surface with one puncture.

\item\label{PL6} There is a triangle-compatible marked subsurface $\mathfrak{S}' = \Sigma'\setminus \mathcal{P}'$ of $\mathfrak{S}$ isomorphic to a~mini\-mal-marked two-holed surface of genus $g\ge 1$ such that $\xi_1$ and $\xi_2$ are hole surrounding loops in~$\mathfrak{S}'$ homotopic to distinct boundary components of $\Sigma'$.

\item\label{PL7} There is a triangle-compatible marked subsurface $\mathfrak{S}' = \Sigma'\setminus \mathcal{P}'$ of $\mathfrak{S}$ isomorphic to a~mini\-mal-marked one-holed surface of genus $g\ge 2$ such that $\xi_1$ is a hole surrounding loop in~$\mathfrak{S}'$, and $\xi_2$ is a separating loop in $\mathfrak{S}'$ that cuts out a one-holed surface of genus~$g'$ with $1\le g'<g$.
\end{enumerate}
\end{Proposition}
Again, one can prove this by a case-by-case argument, together with the classification of oriented surface, so a proof is omitted in the current version, and can be found in a previous version~\cite{ver3}.
In later sections, we~will use the above classification of pairs of disjoint essential simple loops.

\section[Algebraic structure of quantized trace-of-monodromy for a single loop]{Algebraic structure of quantized trace-of-monodromy\\ for a single loop}
\label{sec:algebraic_structure_of_quantized_trace-of-monodromy}

Let $\mathfrak{S}$ be a triangulable marked surface (Definitions~\ref{def:generalized_marked_surface} and~\ref{def:triangulable}) and $\Delta$ its ideal triangulation (Definition~\ref{def:ideal_triangulation}).
Let $\gamma$ be an essential simple loop in $\mathfrak{S}$ (Definition~\ref{def:essential_and_peripheral}). Let
\[
f^\omega_{[\gamma],\Delta} = {\rm Tr}^\omega_{\Delta;\mathfrak{S}}([K_\gamma,{\varnothing}]) \in
\mathbb{BL}^\omega_\Delta \subset \bigl(\mathcal{Z}^\omega_\Delta\bigr)_{\rm bl} \subset \mathcal{Z}^\omega_\Delta
\]
be the quantized trace-of-monodromy for $ [\gamma]$ as defined in Definition~\ref{def:quantized_trace-of-monodromy}, where $K_\gamma$ is a constant elevation lift of $\gamma$ in $\mathfrak{S} \times (-1,1)$ with upward vertical framing, and ${\rm Tr}^\omega_{\Delta;\mathfrak{S}}$ is the quantum trace map in Proposition~\ref{prop:BW}. For the definitions of the codomain algebras $\mathbb{BL}^\omega_\Delta$, $\bigl(\mathcal{Z}^\omega_\Delta\bigr)_{\rm bl}$ and $\mathcal{Z}^\omega_\Delta$, see Sections~\ref{subsec:balanced_algebras}--\ref{subsec:quantized_trace-of-monodromy}. In the present subsection, we investigate the algebraic structure of $f^\omega_{[\gamma],\Delta}$. Using the topological classification of simple loops obtained in Proposition~\ref{prop:loop_classification} of Section~\ref{subsec:topological_classification_of_simple_loops}, we will study $f^\omega_{[\gamma],\Delta}$ for each of the cases classified there.

Before we delve into our investigation of the structure of the quantized trace-of-monodromy $f^\omega_{[\gamma],\Delta} \in \mathcal{Z}^\omega_\Delta$, we first recall from Section~\ref{subsec:balanced_algebras} some basic notations.

An element of $\mathcal{Z}^\omega_\Delta$ is a Laurent polynomial in the variables $Z_i$, $i\in \Delta$, with coefficients in $\mathbb{Z}\bigl[\omega^{\pm 1}\bigr]$. It can also be presented as a $\mathbb{Z}\bigl[\omega^{\pm 1}\bigr]$-linear combination of the set $\{Z_v \mid v\in N_\Delta\}$, where $N_\Delta$ is the skew symmetric lattice associated to $\Delta$, as defined in Definition~\ref{def:N_Delta}. Note that $N_\Delta$ is a free $\mathbb{Z}$-module generated by $v_i$, $i\in \Delta$, with the skew symmetric $\mathbb{Z}$-bilinear form given by $\langle v_i,v_j\rangle = \varepsilon_{ij}$; we have the identification $Z_i = Z_{v_i}$. We also defined in Definition~\ref{def:HA} the Heisenberg algebra $\mathbb{H}_{N_\Delta}$ as the free $\mathbb{Z}$-module generated by $\{z_v \mid v\in N\} \sqcup\{c\}$, modded out by the relations $z_v + z_w = z_{v+w}$. Via the exponential map $\exp \colon \mathbb{H}_{N_\Delta} \to \mathcal{Z}^\omega_\Delta$, $\alpha c + z_v \mapsto {\rm e}^{\alpha c + z_v} := \omega^\alpha Z_v$, defined in Definition~\ref{def:HA_exp}, we have $Z_v = {\rm e}^{z_v}$ and ${\rm e}^c = \omega$.

Here we add one more notation which will be used throughout the remaining part of the paper.
\begin{Definition}[total intersection element for a loop]
\label{def:v_gamma}
Let $\mathfrak{S}$ be a triangulable marked surface and $\Delta$ its ideal triangulation. Let $\gamma$ be an essential simple loop in $\mathfrak{S}$. Define the {\em total intersection} element $v_\gamma = v_{\gamma,\Delta} \in N_\Delta$ for $\gamma$ with respect to $\Delta$ as
\[
v_\gamma = v_{\gamma,\Delta} := \underset{i\in \Delta}{\textstyle \sum} {\rm int}(\gamma,i) \til{v}_i,
\]
where ${\rm int}(\gamma,i) \in \mathbb{Z}_{\ge 0}$ is the intersection number between $\gamma$ and $i$ as defined in Definition~\ref{def:minimal_and_transversal}, and $\til{v}_i$ is as in equation~\eqref{eq:til_v_i}.
\end{Definition}
Note that $v_\gamma = v_{\gamma,\Delta}$ depends only on $\Delta$ and the isotopy class $[\gamma]$ of $\gamma$, so we may write it as $v_{[\gamma]} = v_{[\gamma],\Delta}$. However, for ease of notation, we will mostly stick to $v_\gamma = v_{\gamma,\Delta}$. This element $v_\gamma$ of~$N_\Delta$ yields corresponding elements $z_{v_\gamma} \in \mathbb{H}_{N_\Delta}$ and $Z_{v_\gamma} = {\rm e}^{z_{v_\gamma}} \in \mathcal{Z}^\omega_\Delta$, which we will conveniently use from now on.

\subsection{Peripheral loops}

We now go into the study of the quantized trace-of-monodromy \smash{$f^\omega_{[\gamma],\Delta}$}. We start with a peripheral loop $\gamma$, which appears as the case \ref{L1} in the classification list of Proposition~\ref{prop:loop_classification}.

\begin{Proposition}[quantized trace-of-monodromy for a peripheral loop]
\label{prop:quantized_trace-of-monodromy_peripheral}
Let $\mathfrak{S}$ be a triangulable marked surface. If $\gamma$ is a peripheral loop around a puncture $($Definition~{\rm \ref{def:essential_and_peripheral})}, then for any ideal triangulation $\Delta$ of $(\mathfrak{S})$, one has
\begin{align*}
f^\omega_{[\gamma],\Delta} = Z_{v_\gamma} + Z_{-v_\gamma} = {\rm e}^{z_{v_\gamma}} + {\rm e}^{-z_{v_\gamma}},
\end{align*}
where the total intersection element $v_\gamma = v_{\gamma;\Delta} \in N_\Delta$ is defined in Definition~{\rm \ref{def:v_gamma}}.
\end{Proposition}

\begin{proof} This is well known in the literature; see \cite[Lemma~3.24]{AK}. We still provide a proof here, as a warm-up exercise for the computation of more complicated quantized trace-of-monodromy. Also, our proof here makes use of Proposition~\ref{prop:BW_admissible_sum_formula}, resulting in a much simpler proof than the one in~\cite{AK}.

As shall be done throughout the paper, we find it convenient to give an arbitrary orientation on $\gamma$, so that each loop segment of $\gamma$ becomes oriented. Assume $\gamma$ is in a minimal and transversal position with respect to $\Delta$ (Definition~\ref{def:minimal_and_transversal}), after an isotopy if necessary. An oriented loop segment $k$ is either a {\em left turn} or a {\em right turn}; see Figure~\ref{fig:left_turns} for left turns. There is no U\nobreakdash-turn, due to the minimality condition. Note that a juncture-state $J$ restricts to a state of $k$, i.e., a map $\partial k \to \{+,-\}$ that assigns a sign to the endpoints of $k$. Denote the state of $k$ by $(\varepsilon_{\rm in}, \varepsilon_{\rm ter}) \in \{+,-\}^2$, the pair of state values at the initial and the terminal points of $k$. Then, in view of Figure~\ref{fig:non-admissible_juncture-states} in Section~\ref{subsec:state-sum}, $J$ is admissible if and only if there is no left turn with state $(-,+)$, nor right turn with state $(+,-)$. This observation will be used throughout the paper.

\begin{figure}
\centering\hspace*{-15mm}
\begin{subfigure}[b]{0.33\textwidth}\centering
 \scalebox{0.7}{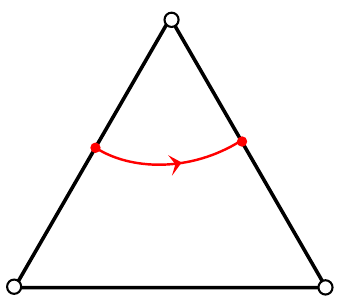}
\caption*{(1) non-self-folded ideal triangle}
\end{subfigure}
\qquad
\begin{subfigure}[b]{0.4\textwidth}\centering
\scalebox{0.7}{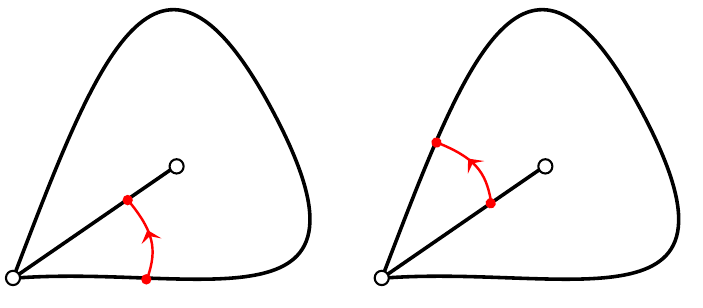}
\caption*{(2) self-folded triangle}
\end{subfigure}
\caption{Oriented loop segments that are left turns.}\label{fig:left_turns}
\end{figure}

Coming back to the case of a peripheral loop, we begin with the following well-known yet straightforward observation on a characterization of peripheral loops.
\begin{Lemma}
\label{lem:peripheral_characterization}
Let $\Delta$ be an ideal triangulation of a triangulable marked surface $\mathfrak{S}$. An oriented essential simple loop $\gamma$ in $\mathfrak{S}$ is a peripheral loop if and only if, when put into a minimal and transversal position with respect to $\Delta$ after an isotopy, it consists solely of left turns or solely of right turns.
\end{Lemma}

Note that the two constant juncture-states $J^+_{\gamma}$ and $J^-_{\gamma}$ in equation~\eqref{eq:constant_juncture-state} are admissible. We claim that these are the only admissible juncture-states. Indeed, if a juncture-state $J$ is not constant, then necessarily there is either a left turn segment with state value $(-,+)$ or a right turn with state value $(+,-)$, so that $J$ is not admissible. Note for each $\epsilon \in \{+,-\}$ that
\[
v_{J_{\gamma,\Delta}^\epsilon} \stackrel{{\rm equation}~\eqref{eq:v_J}}{=} \epsilon \sum_{i\in \Delta} |\gamma \cap i| \til{v}_i
\stackrel{{\rm Definition}~\ref{def:minimal_and_transversal}}{=}
\epsilon \sum_{i\in \Delta} {\rm int}(\gamma,i) \til{v}_i
\stackrel{{\rm Definition}~\ref{def:v_gamma}}{=}
\epsilon v_\gamma,
\]
therefore equation~\eqref{eq:admissible_non-redundancy_condition} is satisfied, and hence the state-sum formula equation~\eqref{eq:admissible_sum_formula} of Proposition~\ref{prop:BW_admissible_sum_formula} yields
\[
f^\omega_{[\gamma],\Delta} = Z^{J^+_{\gamma,\Delta}} + Z^{J^-_{\gamma,\Delta}}
\stackrel{{\rm equation}~\eqref{eq:Z_J}}{=}
Z_{v_{J^+_{\gamma,\Delta}}} + Z_{v_{J^-_{\gamma,\Delta}}}
= Z_{v_\gamma} + Z_{-v_\gamma}.
\]
This completes the proof of Proposition~\ref{prop:quantized_trace-of-monodromy_peripheral}.
\end{proof}

\subsection{Loops with intersection number 2}\label{subsec:loops_with_intersection_number_2}

We move to the next four cases \ref{L2}--\ref{L5} in the classification of Proposition~\ref{prop:loop_classification}; we exclude the case \ref{L1}, i.e., the peripheral loop case. In Proposition~\ref{prop:2_to_5_have_intersection_2}, we characterized these cases by the intersection number with an ideal triangulation; namely, for a loop $\gamma$ falling into these cases, ${\rm int}(\gamma,\Delta)=2$ (Definition~\ref{def:minimal_and_transversal}) holds for some ideal triangulation $\Delta$.

So these cases can be dealt with simultaneously by the following statement.
\begin{Proposition}[quantized trace-of-monodromy for a loop with intersection number 2]
\label{prop:quantized_trace-of-monodromy_for_loop_with_int_2}
Let $\mathfrak{S}$ be a triangulable marked surface, and $\Delta$ its ideal triangulation. If $\gamma$ is a non-peripheral essential simple loop in $\mathfrak{S}$ with ${\rm int}(\gamma,\Delta)=2$ in the sense of Definition~{\rm \ref{def:minimal_and_transversal}}, then
\begin{align*}
f^\omega_{[\gamma],\Delta} = Z_{v_a + v_b} + Z_{\epsilon(v_a - v_b)} + Z_{-v_a-v_b}
 = {\rm e}^{z_{v_a} + z_{v_b}} + {\rm e}^{\epsilon(z_{v_a} - z_{v_b})} + {\rm e}^{-z_{v_a} - z_{v_b}},
\end{align*}
where $a,b \in \Delta$ are the ideal arcs of $\Delta$ that meet $\gamma$ when $\gamma$ is isotoped to a minimal position with respect to $\Delta$, and $\epsilon = {\rm sign}(\varepsilon_{ab}) \in \{+,-\}$; here, $\varepsilon_{ab} \in \{2,-2\}$.
\end{Proposition}

\begin{proof} Suppose $\gamma$ is in a minimal and transversal position with respect to $\Delta$ (Definition~\ref{def:minimal_and_transversal}), and let $a,b\in \Delta$ be the arcs meeting $\gamma$. Let $x_a$, $x_b$ be the junctures of $\gamma$ lying on $a$, $b$, respectively. Give an arbitrary orientation on $\gamma$. Let $k_1$, $k_2$ be the oriented loop segments of $\gamma$, so that $k_1$ goes from $x_a$ to $x_b$, and $k_2$ goes from $x_b$ to $x_a$. One of $k_1$ and~$k_2$ is a left turn and the other is a right turn, because otherwise $\gamma$ is a peripheral loop by Lemma~\ref{lem:peripheral_characterization}. First, assume that $k_1$ is a~left turn and $k_2$ is a right turn. Then it follows that $\varepsilon_{ab} = 2$, thus $\epsilon=+$. Denote a juncture-state $J \colon \gamma \cap \Delta \to \{+,-\}$ of $\gamma$ as a pair of values $(J(x_a),J(x_b))$. In view of the definition of admissibility (in Section~\ref{subsec:state-sum}; see Figure~\ref{fig:non-admissible_juncture-states} for non-admissible configurations), we note that the admissible juncture-states are $(+,+) = J_{\gamma}^+$, $(-,-) = J_{\gamma}^-$ (equation~\eqref{eq:constant_juncture-state}) and $(+,-) =: J_0^+$. On the other hand, if $k_1$ is a right turn and $k_2$ is a left turn, then $\varepsilon_{ab} = -2$, so $\epsilon=-$, and the admissible juncture-states are $J_\gamma^+$, $J_\gamma^-$ and $(-,+) =: J_0^-$. In either case, the elements $v_J$ (equation~\eqref{eq:v_J}) for admissible juncture-states $J$ are
\[
v_{J_\gamma^+} = v_a + v_b, \qquad v_{J_\gamma^-} = -v_a - v_b, \qquad
v_{J_0^\epsilon} = \epsilon(v_a - v_b)
\]
Hence, equation~\eqref{eq:admissible_non-redundancy_condition} is satisfied, so Proposition~\ref{prop:BW_admissible_sum_formula} applies. The formula in equation~\eqref{eq:admissible_sum_formula}, together with equation~\eqref{eq:Z_J}, yields the desired result.
\end{proof}

\subsection[Algebraic strong commutativity, Teschner recursion and Teschner triple of loops]{Algebraic strong commutativity, Teschner recursion\\ and Teschner triple of loops}
\label{subsec:Teschner_recursion}

In order to deal with the remaining case \ref{L6} in the classification of Proposition~\ref{prop:loop_classification}, we develop a new tool.

As a preliminary to this investigation, and also to be used in our later study of strong commutativity of quantized trace-of-monodromy for two loops, we first define some algebraic counterpart of the notion of strong commutativity of self-adjoint operators.

We remark that the representation of the square-root Chekhov--Fock algebra $\mathcal{Z}^\omega_\Delta = \mathbb{T}^\omega_{N_\Delta}$ which we will be using in the sequel paper \cite{sequel} to tackle Conjecture \ref{conjecture:sequel}, as mentioned in Section~\ref{subsec:on_quantized_length_operators}, enjoys a property saying that if $v,w\in N_\Delta$ satisfy $[z_v, z_w] =0$ in $\mathbb{H}_{N_\Delta}$, i.e., if $\langle v,w\rangle =0$, then the self-adjoint operators representing the elements $Z_v, Z_w \in \mathcal{Z}^\omega_\Delta$ {\em strongly commute}, in the sense, e.g., as in \cite[Section~VIII.5]{RS}. See \cite{FG09, Kim_irreducible} for such representations. This inspires one algebraic version of strong commutativity as follows, between an element of the Heisenberg algebra (Definition~\ref{def:HA}) and the corresponding quantum torus algebra (Definition~\ref{def:QT}).

\begin{Definition}[algebraic strong commutativity between Heisenberg algebra and quantum torus algebra]
\label{def:algebraic_SC1}
Let $\mathbb{H}_N$ and $\mathbb{T}^t_N$ be the Heisenberg algebra (Definition~\ref{def:HA}) and the quantum torus algebra (Definition~\ref{def:QT}) associated to a skew symmetric lattice $N$ (Definition~\ref{def:QT}), with the underlying skew symmetric form $\langle \cdot,\cdot \rangle$.

An element $z_v$ of $\mathbb{H}_N$, with $v\in N$, is said to {\em $($algebraically$)$ strongly commute} with an element~$X$ of~$\mathbb{T}^t_N$ if $X = \sum_{u \in N} a_u Z_u$ with $a_u \in \mathbb{Z}\bigl[t^{\pm 1}\bigr]$, where $a_u =0$ for all but finitely many $u\in N$, and $\langle v,u\rangle=0$ holds for all $u\in N$ with $a_u \neq0$.
\end{Definition}

A similar but stronger property of these representations also inspires the following algebraic version of strong commutativity among elements of a quantum torus algebra.
\begin{Definition}[algebraic strong commutativity for quantum torus algebra]
\label{def:algebraic_SC2}
Let $\mathbb{T}^t_N$ be the quantum torus algebra (Definition~\ref{def:QT}) associated to a skew symmetric lattice $N$ (Definition~\ref{def:QT}), with the underlying skew symmetric form $\langle \cdot,\cdot \rangle$.

Denote the real span of $N$ by $N_\mathbb{R} := N \otimes_\mathbb{Z} \mathbb{R}$, equipped with the skew symmetric $\mathbb{R}$-bilinear form $\langle \cdot, \cdot\rangle$ extending the form $\langle \cdot,\cdot\rangle$ of $N$. Let
\begin{align}
\label{eq:N_R_perp}
N_\mathbb{R}^\perp := \{ v\in N_\mathbb{R} \mid \langle v,w\rangle=0, ~ \forall w\in N_\mathbb{R}\}.
\end{align}
For any $\mathbb{Z}$-submodule $M$ of $N$, denote its real span by $M_\mathbb{R} := M \otimes_\mathbb{Z} \mathbb{R} \subset N \otimes_\mathbb{Z} \mathbb{R} = N_\mathbb{R}$.

Elements $X$ and $Y$ of $\mathbb{T}^t_N$ are said to {\em $($algebraically$)$ strongly commute} with each other if there exists nonzero $\mathbb{Z}$-submodules $N_1$ and $N_2$ of $N$ such that
\begin{enumerate}[label={\rm (SC\arabic*)}]\itemsep=0pt\setlength{\leftskip}{0.22cm}
\item\label{SC1} $N_1 \perp N_2$ holds in the sense of
$\langle v,w\rangle=0$, $\forall v\in N_1$, $\forall w\in N_2,$

\item\label{SC2}
$(N_1)_\mathbb{R} \cap (N_2)_\mathbb{R} \subset N_\mathbb{R}^\perp$ holds, and

\item\label{SC3} $X$ and $Y$ can be written as finite sums $X = \sum_{v\in N_1} a_v Z_v$ and $Y = \sum_{v \in N_2} b_v Z_v$ for some $a_v,b_v \in \mathbb{Z}\bigl[t^{\pm 1}\bigr]$ which are zero for all but finitely many $v$.
\end{enumerate}
\end{Definition}
In the sequel paper \cite{sequel}, we will justify that the above two definitions of algebraic strong commutativity indeed imply the strong commutativity of corresponding operators assigned by the representations that we will use.

The following is the motivation of our notion of a triangle-disjoint pair of loops (Definition~\ref{def:triangle-disjoint}).
\begin{Lemma}[triangle-disjoint implies strong commutativity]
\label{lem:triangle-disjoint_implies_strong_commutativity}
Let $\mathfrak{S}$ be a triangulable marked surface and $\Delta$ be its ideal triangulation. Suppose $\gamma_1$ and $\gamma_2$ are triangle-disjoint with respect to~$\Delta$, in the sense of Definition~{\rm \ref{def:triangle-disjoint}}. Then the quantized trace-of-monodromy \smash{$f^\omega_{[\gamma_1],\Delta}$} and \smash{$f^\omega_{[\gamma_2],\Delta}$} $($Definition~{\rm \ref{def:quantized_trace-of-monodromy})}, which are elements in $\mathcal{Z}^\omega_\Delta = \mathbb{T}^\omega_{N_\Delta}$, strongly commute with each other in the sense of Definition~{\rm \ref{def:algebraic_SC2}}.
\end{Lemma}

\begin{proof} The case when $\gamma_1$ or $\gamma_2$ is a peripheral loop is a special case, which is easier; we postpone it to Section~\ref{subsec:easy_part}. Now assume that these loops are not peripheral. Put $\gamma_1$ and $\gamma_2$ in minimal and transversal positions with respect to $\Delta$. For $j=1,2$, let $N_j := N_{\gamma_j}$ be the $\mathbb{Z}$-submodule of~$N_\Delta$ defined in Proposition~\ref{prop:weak_sum}. By hypothesis, we have $(N_1)_\mathbb{R} \cap (N_2)_\mathbb{R}=0$, and it is easy to see that $N_1 \perp N_2$; so we have \ref{SC1}--\ref{SC2}. The only possible subtle situation for proving \ref{SC2} is when a loop meets a self-folded arc, say because the only case when $\varepsilon_{ik} \neq 0$ can happen while~$i$,~$k$ do not appear in a same triangle is when one of $i$ and $k$ is a self-folded arc. Suppose $\gamma_1$ meets a self-folded arc~$i$; since $\gamma_1$ is not peripheral, it must also meet $\pi(i)$ (see equation~\eqref{eq:pi}), thus $\gamma_1$ meets the two triangles $T_1$, $T_2$ having $\pi(i)$ as a side, where $T_1$ is the self-folded one having~$i$. If~$T_2$ is also self-folded, then $\mathfrak{S}$ is a thrice-punctured sphere, being the union of $T_1$ and $T_2$; so there is no $\gamma_2$ that is triangle-disjoint with~$\gamma_1$. So $T_2$ is not self-folded. Note that $\varepsilon_{ik} \neq 0$ can occur if either $k$ is a side of~$T_2$, or $k$ is a self-folded arc of a self-folded triangle $T_3$ such that~$\pi(k)$ is a side of $T_2$. In either case, one can see that~$v_k$ is not involved in any of the generators of~$N_{\gamma_2}$, so that \ref{SC2} holds. One can verify~\ref{SC1} in a similar manner.

By Proposition~\ref{prop:weak_sum}, the sum expression in \smash{$f^\omega_{[\gamma_j],\Delta} = \sum_{v\in N} a_v Z_v$} is supported on $N_j$, i.e., $a_v =0$ if $v \notin N_j$; this gives~\ref{SC3}.
\end{proof}

We are now ready to study the promised new tool, which is one of the most crucial definitions of the present paper.
\begin{Definition}
\label{def:Teschner_triple_level1}
Let $\mathfrak{S}$ be a triangulable marked surface and $\Delta$ be its ideal triangulation. Let $(\gamma,\gamma_1,\gamma_2)$ be a {\em pants triple} in $\mathfrak{S}$, which means an ordered triple of mutually non-intersecting essential non-peripheral simple loops in $\mathfrak{S}$ that cuts out a three-holed sphere in the sense of Definition~\ref{def:cutting}.

We say that a pants triple $(\gamma,\gamma_1,\gamma_2)$ in $\mathfrak{S}$ is a {\em strong Teschner triple} with respect to $\Delta$ if there exist $v_1,v_2 \in N_\Delta$ (Definition~\ref{def:N_Delta}) satisfying the following.
\begin{enumerate}[label={\rm (TR\arabic*)}]\itemsep=0pt\setlength{\leftskip}{0.28cm}
\item\label{TR1} The quantized trace-of-monodromy (Definition~\ref{def:quantized_trace-of-monodromy}) of the three loops satisfy the following equation called the {\em Teschner recursion relation} with respect to $\Delta$:
\begin{align}
f^\omega_{[\gamma],\Delta}
& = Z_{v_\gamma} + Z_{-v_\gamma} + Z_{v_1+v_2} + Z_{v_1} \, f^\omega_{[\gamma_1],\Delta}
+ Z_{v_2} f^\omega_{[\gamma_2],\Delta} \nonumber\\
& = {\rm e}^{z_{v_\gamma}} + {\rm e}^{-z_{v_\gamma}} + {\rm e}^{z_{v_1} + z_{v_2}} + {\rm e}^{z_{v_1}} f^\omega_{[\gamma_1],\Delta} + {\rm e}^{z_{v_2}} f^\omega_{[\gamma_2],\Delta},\label{eq:TR1}
\end{align}
where $v_\gamma = v_{\gamma;\Delta} = \sum_{i\in \Delta} {\rm int}(\gamma,i) \til{v}_i \in N_\Delta$ is the total intersection element for the loop $\gamma$, as defined as in Definition~\ref{def:v_gamma}. Here the term $Z_{v_1+v_2} = {\rm e}^{z_{v_1} + z_{v_2}}$ is called the {\em connecting term}.

\item\label{TR2} $\begin{cases} v_1 - v_2 = v_\gamma \\ \langle v_1, v_2 \rangle = -4 \end{cases}$ or $\begin{cases} v_1 - v_2 = - v_\gamma \\ \langle v_1, v_2 \rangle = 4 \end{cases}$ holds.

\item\label{TR3} $z_{v_i} \in \mathbb{H}_{N_\Delta}$ strongly commutes with $f^\omega_{[\gamma_j],\Delta} \in \mathcal{Z}^\omega_\Delta = \mathbb{T}^\omega_{N_\Delta}$ for all $i,j\in\{1,2\}$, in the sense of Definition~\ref{def:algebraic_SC1}.

\item\label{TR4} The quantized trace-of-monodromy $f^\omega_{[\gamma_1],\Delta}$ and $f^\omega_{[\gamma_2],\Delta}$ (Definition~\ref{def:quantized_trace-of-monodromy}), which are elements of $\mathcal{Z}^\omega_\Delta = \mathbb{T}^\omega_{N_\Delta}$, either strongly commute with each other in the sense of Definition~\ref{def:algebraic_SC2}, or coincide with each other.
\end{enumerate}
\end{Definition}
Note that some of the three curves of a pants triple are allowed to be isotopic to each other.
\begin{Remark}
In principle, one may drop the non-peripheral condition in the definition of a~pants triple, but in this paper we deal only with triples consisting entirely of non-peripheral loops.
\end{Remark}

The above definition of a strong Teschner triple is inspired by Teschner's construction in~\cite{T} of a system of operators ${\bf L}_{[\gamma], \sigma}$ on a Hilbert space $\mathcal{H}(\sigma)$, associated to a `marking' $\sigma$ of a surface, which is a pants decomposition with an extra data, and its constituent loop $\gamma$. We note that Teschner's original recursion relation is not exactly as in equation~\eqref{eq:TR1} of \ref{TR1}; the roles of the elements $z_{v_\gamma}$, $z_{v_1}$ and $z_{v_2}$ are played by some other operators in \cite{T}. Our $z_v$ are linear combinations of quantum version of Thurston's shear coordinates of the Teichm\"uller space associated to an ideal triangulation of the surface. The operators that Teschner uses in place of $z_{v_\gamma}$, $z_{v_1}$ and $z_{v_2}$ are combinatorially built operators. Teschner then uses the resulting operators~${\bf L}_{[\gamma],\sigma}$ to obtain operators~${\bf l}_{[\gamma],\sigma}$ via ${\bf L}_{[\gamma],\sigma} = 2\cosh\bigl({\bf l}_{[\gamma],\sigma}/2\bigr)$, and suggest them as being models of a system of quantum length operators. Roughly speaking, the logic of Teschner goes in the opposite direction to that taken in the present paper.

The sought-for consequences of the notion of a Teschner triple are two-fold. One consequence is that the operator for $\gamma$ would strongly commute with those of $\gamma_1$ and $\gamma_2$, and the other is that the spectral properties of the operator for $\gamma$ could be deduced from those of $\gamma_1$ and $\gamma_2$. We will materialize these two consequences in the sequel \cite{sequel}, with the appropriate modification applied to original arguments of Teschner \cite{T}; see Section~\ref{sec:consequences} for a brief preview. In fact, the latter consequence partially depends on the former. This dependence is translated into an algebraic language as follows, yielding a {\it weak} Teschner triple, which we will need for the whole program, in addition to the strong Teschner triple. The difference between the weak and the strong Teschner triples is at the commutativity requirement \ref{TR4}.
\begin{Definition}
\label{def:Teschner-commute}
Let a pants triple $(\gamma,\gamma_1,\gamma_2)$ in a triangulable marked surface $\mathfrak{S}$ be a strong Teschner triple with respect to an ideal triangulation $\Delta$ of $\mathfrak{S}$. Then, for each $j=1,2$, we say that $f^\omega_{[\gamma],\Delta}$ and $f^\omega_{[\gamma_j],\Delta}$ {\em strongly Teschner-commute} with each other.
\end{Definition}

\begin{Definition}
\label{def:Teschner_triple_level2}
Let $(\gamma,\gamma_1,\gamma_2)$ be a pants triple in a triangulable marked surface $\mathfrak{S}$ (Definition~\ref{def:Teschner_triple_level1}). We say that $(\gamma,\gamma_1,\gamma_2)$ is a {\em weak Teschner triple} with respect to an ideal triangulation~$\Delta$ of~$\mathfrak{S}$ if there exist $v_1, v_2 \in N_\Delta$ (Definition~\ref{def:N_Delta}) satisfying \ref{TR1}, \ref{TR2} and \ref{TR3} of Definition~\ref{def:Teschner_triple_level1}, as well as the following:
\begin{enumerate}[label={\rm (TR\arabic*)}]\itemsep=0pt\setlength{\leftskip}{0.28cm}
\setcounter{enumi}{4}
\item\label{TR5} $f^\omega_{[\gamma_1],\Delta}$ and $f^\omega_{[\gamma_2],\Delta}$ strongly Teschner-commute with each other.
\end{enumerate}
In this case, we say, for each $j=1,2$, that $f^\omega_{[\gamma],\Delta}$ and $f^\omega_{[\gamma_j],\Delta}$ {\em weakly Teschner-commute}.
\end{Definition}
The condition \ref{TR5} is equivalent to:
\begin{enumerate}[label={\rm (TR\arabic*$'$)}]\itemsep=0pt\setlength{\leftskip}{0.38cm}
\setcounter{enumi}{4}
\item\label{TR5_prime} There exists an essential simple loop $\gamma_3$ in $\mathfrak{S}$ such that $(\gamma_1,\gamma_2,\gamma_3)$ is a strong Teschner triple with respect to $\Delta$.
\end{enumerate}

One basic observation is that the role of $\gamma_1$ and that of $\gamma_2$ in the definition of a Teschner triple is symmetric.
\begin{Lemma}
Let $(\gamma,\gamma_1,\gamma_2)$ be a pants triple in a triangulable marked surface $\mathfrak{S}$ $($see Definition~{\rm \ref{def:Teschner_triple_level1})}.
\begin{enumerate}\itemsep=0pt
\item[$(1)$] $(\gamma,\gamma_1,\gamma_2)$ is a strong Teschner triple with respect to $\Delta$ if and only if $(\gamma,\gamma_2,\gamma_1)$ is.

\item[$(2)$] $(\gamma,\gamma_1,\gamma_2)$ is a weak Teschner triple with respect to $\Delta$ if and only if $(\gamma,\gamma_2,\gamma_1)$ is.
\end{enumerate}
\end{Lemma}
A proof is omitted, as it is straightforward.

A few remarks follow.

\begin{Remark}\label{rem:further_weakness}
The notion of a weak Teschner triple was obtained by weakening the commutativity requirement in \ref{TR4} from the algebraic strong commutativity to the strong Teschner-commutativity, as in \ref{TR5}. Inductively, one can define Teschner triples of further levels of weakness, by weakening \ref{TR5} from the strong Teschner-commutativity to the weak Teschner-commutativity, and so on.
\end{Remark}

\begin{Remark}
The focus of Teschner's constructions are on the analytic properties. So, in view of our analytic proofs in \cite{sequel}, we remark that there is no need to use a single ideal triangulation only, when defining the notion of a Teschner triple of loops. For example, \ref{TR4} can be replaced by the condition saying that \smash{$f^\omega_{[\gamma_1],\Delta'}$} and \smash{$f^\omega_{[\gamma_2],\Delta'}$} strongly commute for {\em some} ideal triangulation~$\Delta'$, and \ref{TR5_prime} can be replaced by the condition saying that there exists a loop~$\gamma_3$ such that $(\gamma_1,\gamma_2,\gamma_3)$ is a strong Teschner triple with respect to {\em some} ideal triangulation~$\Delta'$.
\end{Remark}

\begin{Remark}
It is easy to show that if \smash{$f^\omega_{[\gamma_1],\Delta}$} and \smash{$f^\omega_{[\gamma_2],\Delta}$} strongly or weakly Teschner-commute, then they satisfy the following type of commutativity:
\[
f^\omega_{[\gamma_1],\Delta} = \sum_{v \in N_1} a_u Z_v \qquad \text{and} \qquad f^\omega_{[\gamma_2],\Delta} = \sum_{v\in N_2} b_v Z_v
\]
 holds for some $a_v,b_v \in \mathbb{Z}\bigl[t^{\pm 1}\bigr]$ with some subsets $N_1,N_2 \subset N_\Delta$, such that $\langle v,w\rangle =0$ for all ${v\in N_1}$, ${w\in N_2}$. This is weaker than the algebraic strong commutativity defined in Definition~\ref{def:algebraic_SC2}, but will be useful in proving the analytic strong commutativity \cite{sequel}. It easily implies the usual algebraic commutativity
\[
 f^\omega_{[\gamma_1],\Delta}f^\omega_{[\gamma_2],\Delta}=f^\omega_{[\gamma_2],\Delta}f^\omega_{[\gamma_1],\Delta},
\]
 which in fact automatically holds because~$\gamma_1$ and~$\gamma_2$ are disjoint.
\end{Remark}

\begin{Remark}\label{rem:bracket_of_v1_and_v2}
In \ref{TR2}, we specified the value of $\langle v_1,v_2\rangle$. In our analytic arguments in \cite{sequel}, we will only need this value to be any nonzero integer.
\end{Remark}

To remind the reader again, the notion of Teschner triple of loops will be used in upcoming subsections to investigate the property of the quantized trace-of-monodromy for the case \ref{L6}, and also in the next section to study some version of strong commutativity of quantized trace-of-monodromy for disjoint loops.

\subsection{Hole surrounding loop in one-holed torus subsurface}
\label{subsec:loops_in_one-holed_torus}

Here we deal with the item \ref{L6} of Proposition~\ref{prop:loop_classification}, under a further stipulation that the genus~$g$ should be $1$. Namely, we deal with a hole surrounding loop $\gamma$ in a one-holed torus subsurface $\mathfrak{S}' = \Sigma' \setminus \mathcal{P}'$ of $\mathfrak{S} = \Sigma\setminus\mathcal{P}$. To investigate the structure of the quantized trace-of-monodromy along $\gamma$, we make use of the notion of Teschner triple defined in the previous subsection. As~an auxiliary loop, we consider a non-separating loop $\eta$ in the subsurface $\mathfrak{S}'$. We already know the structure of the quantized trace-of-monodromy along $\eta$, because $\eta$ falls into the case \ref{L5} of Proposition~\ref{prop:loop_classification}, so that Proposition~\ref{prop:2_to_5_have_intersection_2} applies, and therefore so does Proposition~\ref{prop:quantized_trace-of-monodromy_for_loop_with_int_2} of Section~\ref{subsec:loops_with_intersection_number_2}. As $\gamma$ is isotopic to the boundary loop of $\Sigma'$, it follows that after an isotopy, we can assume that $\eta$ and $\gamma$ are disjoint. We deal with the quantized trace-of-monodromy of $\gamma$ using that of $\eta$ as follows. For convenience, we formulate the statement in terms only of $\mathfrak{S}' = \Sigma' \setminus \mathcal{P}'$ without referring to $\mathfrak{S} = \Sigma \setminus \mathcal{P}$, where $\mathfrak{S}'$ is a triangle-compatible marked subsurface of a~triangulable marked surface $\mathfrak{S}$ (Definition~\ref{def:subsurfaces}).

\begin{Proposition}[loops in a one-holed torus]
\label{prop:loops_in_a_one-holed_torus}
Let $\mathfrak{S}' = \Sigma'\setminus\mathcal{P}'$ be a triangulable marked surface that is isomorphic to a minimal-marked one-holed torus $($Definitions~{\rm \ref{def:surface_terminology}} and {\rm \ref{def:minimal-marked})}. Let~$\gamma$ be a hole surrounding loop $($Definition~{\rm \ref{def:hole-surrounding})} in $\mathfrak{S}'$ and $\eta$ be a non-separating loop in $\mathfrak{S}'$. Then there exists an ideal triangulation~$\Delta'$ of~$\mathfrak{S}'$ such that
\begin{enumerate}\itemsep=0pt
\item[\rm (1)] ${\rm int}(\eta,\Delta')=2$ $($Definition~{\rm \ref{def:minimal_and_transversal})}, so that Proposition~{\rm \ref{prop:quantized_trace-of-monodromy_for_loop_with_int_2}} applies to $\eta$.
\item[\rm (2)] The triple $(\gamma,\eta,\eta)$ is a strong Teschner triple with respect to~$\Delta'$, in the sense of Definition~{\rm \ref{def:Teschner_triple_level1}}.
\end{enumerate}
\end{Proposition}

\begin{figure}
\centering
\begin{subfigure}[b]{0.6\textwidth}\centering
\scalebox{0.65}{\raisebox{-5mm}{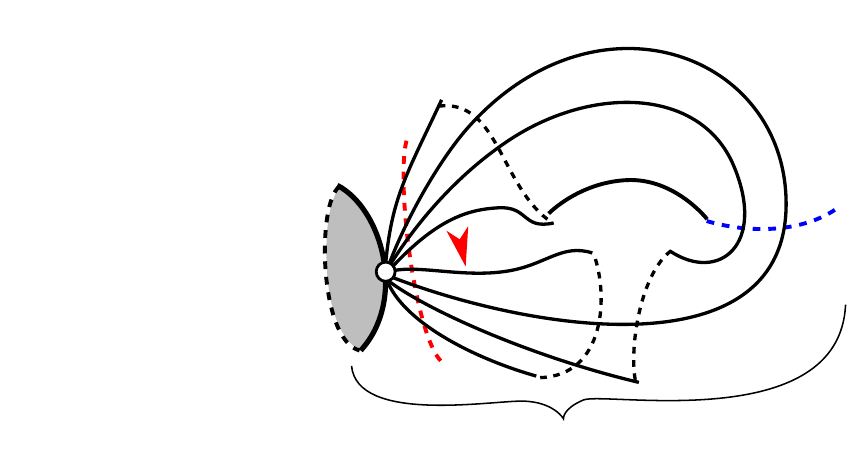}}
\caption*{(1) one-holed torus subsurface $(\Sigma',\mathcal{P}')$ of a surface $(\Sigma,\mathcal{P})$, with hole surrounding loop $\gamma$ and non-separating loop $\eta$ in $(\Sigma',\mathcal{P}')$}
\end{subfigure}
\hspace{5mm}
\begin{subfigure}[b]{0.3\textwidth}\centering
\scalebox{0.65}{\raisebox{0mm}{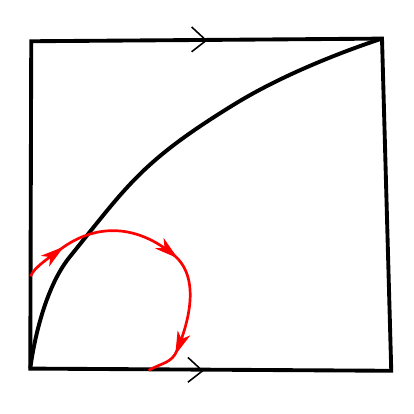}}
\caption*{(2) one-holed torus subsurface $(\Sigma',\mathcal{P}')$ presented as quotient of closed square with hole}
\end{subfigure}
\caption{Loops in one-holed torus subsurface.}\label{fig:loops_in_one-holed_torus_subsurface}
\end{figure}

\begin{proof} We have already touched upon the item (1) when proving Proposition~\ref{prop:2_to_5_have_intersection_2}, to deal with the case \ref{L5} of {Proposition~\ref{prop:loop_classification}}. Namely, by the cutting and gluing idea used in the proof of Proposition~\ref{prop:2_to_5_have_intersection_2}, one can visualize $\eta$ as in Figure~\ref{fig:loops_in_one-holed_torus_subsurface}\,(1); it is easy to see that $\gamma$ can be isotoped to look like in the figure. Also, one can present the minimal-marked one-holed torus as the quotient of a closed square with hole, with each side being identified with its parallel opposite one, as in Figure~\ref{fig:loops_in_one-holed_torus_subsurface}\,(2); in particular, one can consider an ideal triangulation~$\Delta'$ of $\mathfrak{S}'$ as in the figure. We~observe that ${\rm int}(\eta,\Delta')=2$, so that Proposition~\ref{prop:quantized_trace-of-monodromy_for_loop_with_int_2} applies, proving the item (1) of the current proposition. More precisely, since $a$, $b$ are the ideal arcs of $\Delta'$ meeting $\eta$ and $\varepsilon_{ab}=-2$, we have
\begin{align}
\label{eq:oht_eta_f_omega}
f^\omega_{[\eta],\Delta'} = {\rm e}^{ z_{v_a} + z_{v_b} } + {\rm e}^{ - z_{v_a} + z_{v_b}} + {\rm e}^{-z_{v_a} - z_{v_b}}.
\end{align}

Let us now explicitly compute the quantized trace-of-monodromy $f^\omega_{[\gamma],\Delta'}$ for the hole surrounding loop $\gamma$. Give an orientation on $\gamma$ and label its loop segments by $\circled{1}$, \dots, $\circled{8}$, as depicted in~Figure~\ref{fig:loops_in_one-holed_torus_subsurface}. Then $\circled{1}$ is a left turn (indicated as L), and all other segments are right turns (indicated as~R). A juncture-state assigns a sign to each juncture, so there are $2^8$ possible juncture-states in total. Recall that a juncture-state is admissible if there is no occurrence of \smash{$+ \overset{\rm R}{\longrightarrow} -$} or \smash{$- \overset{\rm L}{\longrightarrow} +$}. All admissible juncture-states are as listed in Table \ref{tab:hole-surrounding_loop_in_one-holed_torus}, which one can easily verify; some entries in this table are colored red to facilitate the understanding, and can be ignored if not helpful. In this table, note that the rightmost juncture is the same as the leftmost one, hence must be assigned the same signs. We put parentheses for the entries of this table for the leftmost juncture, as they are redundant; but we still keep them in the table, for convenience.

\begin{table}[htb]
\centering
{\setlength{\tabcolsep}{0,2mm}
\begin{tabular}{ r @{\hspace{1mm} } lllllllll@{\hspace{1mm}} l}
\textcolor{red}{$\gamma$ }: & \textcolor{red}{ \tablearrow{1}{L} } & \textcolor{red}{ \tablearrow{2}{R} } & \textcolor{red}{ \tablearrow{3}{R} } & \textcolor{red}{ \tablearrow{4}{R} } & \textcolor{red}{ \tablearrow{5}{R} } & \textcolor{red}{ \tablearrow{6}{R} } & \textcolor{red}{ \tablearrow{7}{R} } & \textcolor{red}{ \tablearrow{8}{R} } & \hspace{-0,7mm} \textcolor{red}{ $\bullet$ } & \\
ideal arcs : & \hspace*{-2,5mm} ($d$) & $c$ & $a$ & $b$ & $c$ & $d$ & $a$ & $b$ & $d$ &
\\ \hline
$J_0$ : & \hspace*{-2,7mm} ($+$) & \textcolor{red}{$+$} & $+$ & $+$ & $+$ & $+$ & $+$ & $+$ & $+$ & all plus\\ \hline
$J_1$ : & \hspace*{-2,7mm} ($+$) & \textcolor{red}{$-$} \vbar & $+$ & $+$ \vbar & $+$ & $+$ & $+$ & $+$ & $+$ & %
\\
$J_2$ : & \hspace*{-2,7mm} ($+$) & \textcolor{red}{$-$} \vbar & \textcolor{red}{$-$} & $+$ \vbar & $+$ & $+$ & $+$ & $+$ & $+$ & \\
$J_3$ : & \hspace*{-2,7mm} ($+$) & \textcolor{red}{$-$} \vbar & \textcolor{red}{$-$} & \textcolor{red}{$-$} \vbar & $+$ & $+$ & $+$ & $+$ & $+$ & \\ \hline
$J_4$ : & \hspace*{-2,7mm} ($+$) & \textcolor{red}{$-$} & \textcolor{red}{$-$} & \textcolor{red}{$-$} & \textcolor{red}{$-$} & $+$ & $+$ & $+$ & $+$ & connecting term \\ \hline
$J_5$ : & \hspace*{-2,7mm} ($+$) & \textcolor{red}{$-$} & \textcolor{red}{$-$} & \textcolor{red}{$-$} & \textcolor{red}{$-$} & \textcolor{red}{$-$} \vbar & $+$ & $+$ \vbar & $+$ & %
\\
$J_6$ : & \hspace*{-2,7mm} ($+$) & \textcolor{red}{$-$} & \textcolor{red}{$-$} & \textcolor{red}{$-$} & \textcolor{red}{$-$} & \textcolor{red}{$-$} \vbar & \textcolor{red}{$-$} & $+$ \vbar & $+$ & \\
$J_7$ : & \hspace*{-2,7mm} ($+$) & \textcolor{red}{$-$} & \textcolor{red}{$-$} & \textcolor{red}{$-$} & \textcolor{red}{$-$} & \textcolor{red}{$-$} \vbar & \textcolor{red}{$-$} & \textcolor{red}{$-$} \vbar & $+$ & \\ \hline
$J_8$ : & \hspace*{-2,7mm} \textcolor{red}{($-$)} & \red{$-$} & \red{$-$} & \red{$-$} & \red{$-$} & \red{$-$} & \red{$-$} & \red{$-$} & \textcolor{red}{$-$} & all minus \\
\end{tabular}
}
\caption{All admissible juncture-states for hole surrounding loop $\gamma$ in a minimal-marked one-holed torus $(\Sigma',\mathcal{P}')$.}\label{tab:hole-surrounding_loop_in_one-holed_torus}
\end{table}

Let us compute the Weyl-ordered Laurent monomials $Z^{J_r} = Z_{v_{J_r}} = \exp(z_{v_{J_r}})$ for these admissible juncture-states $J_r$, following the definition given in equation~\eqref{eq:v_J}--\eqref{eq:Z_J}. This computation can be done by reading the entries of Table \ref{tab:hole-surrounding_loop_in_one-holed_torus} for each row. For example, the row~$J_4$ reads $-,-,-,-,+,+,+,+$ (without the first column), and since the ideal arcs on which the corresponding junctures lie in are $c$, $a$, $b$, $c$, $d$, $a$, $b$, $d$ as written in the table, we have
\[
v_{J_4} = -v_c - v_a - v_b - v_c + v_d + v_a + v_b + v_d
= -2 v_c + 2v_d.
\]
Note that $\til{v}_i=v_i$ in this case, because there are no self-folded arcs.

Likewise, we compute the others
\begin{align*}
 v_{J_0} & = 2v_a + 2v_b + 2v_c + 2v_d, \ & v_{J_1} & = 2v_a + 2v_b + 2v_d, \ & v_{J_2} & = 2v_b + 2v_d, \\
 v_{J_3} & = 2v_d, \ & v_{J_4} & = -2v_c + 2v_d, \ & v_{J_5} & = -2v_c, \\
 v_{J_6} & = - 2v_a - 2v_c, \ & v_{J_7} & = -2v_a - 2v_b - 2v_c, \ & v_{J_8} & = -2v_a -2 v_b - 2v_c - 2v_d.
\end{align*}

Note from Definition~\ref{def:v_gamma} that $v_\gamma = 2 v_a + 2v_b + 2v_c + 2v_d$, as $\gamma$ meets each of $a$, $b$, $c$, $d$ twice, and does not meet the boundary arc of $\Delta'$. Thus,
\[
v_{J_0} = v_\gamma, \qquad v_{J_8} = - v_\gamma.
\]

We now group the remaining elements $v_{J_r}$ appropriately, as hinted by the horizontal dividing lines in Table~\ref{tab:hole-surrounding_loop_in_one-holed_torus}.
For $J_1$, $J_2$, $J_3$, note that the juncture-state values for the two endpoints of the segment $\circled{3}$ (lying in arcs $a$ and $b$) form the values of all admissible juncture-states of the loop $\eta$ at the endpoints of the (right turn) segment $\circled{$1'$}$ (lying in $a$ and $b$); call these latter juncture-states~$J_1'$, $J_2'$, $J_3'$. In Table \ref{tab:hole-surrounding_loop_in_one-holed_torus}, these values are surrounded by the vertical bars. Taking the signed sum of the elements $v_i$ for the entries in the rows for~$J_1$, $J_2$, $J_3$ outside the vertical bars, one defines the element $v_1 \in N_{\Delta'}$ as
\[
v_1 := -v_c + v_c+v_d+v_a+v_b+v_d= v_a + v_b + 2v_d.
\]
So we have
\[
v_{J_r} = v_1 + v_{J_r'}, \qquad \forall r=1,2,3.
\]
We now show $\langle v_1, v_{J_r'} \rangle=0$ for $r=1,2,3$. Since $v_{J_r'}$ involves only $v_a$ and $v_b$, it suffices to show $\langle v_1, v_a\rangle = \langle v_1 , v_b \rangle=0$. From Figure~\ref{fig:loops_in_one-holed_torus_subsurface}\,(2), one can read
\[
\varepsilon_{ab} = -2, \qquad
\varepsilon_{ac} = 1, \qquad
\varepsilon_{ad}=1, \qquad
\varepsilon_{bc}=-1,\qquad
\varepsilon_{bd}=-1,\qquad
\varepsilon_{cd}=-1,
\]
according to the definition of the exchange matrix $\varepsilon_{ij}$ for $\Delta'$ in Definition~\ref{def:varepsilon_Delta}. Using $v_1 = v_a+v_b+2v_d$ and equation~\eqref{eq:skew-form_for_Delta}, one can now easily verify $\langle v_1, v_a \rangle = \langle v_1, v_b\rangle=0$. Hence, indeed $\langle v_1,v_{J_r'}\rangle=0$ for $r=1,2,3$. Thus, we have
\begin{align*}
\sum_{r=1}^{3} Z^{J_r} & = \sum_{r=1}^{3} {\rm e}^{z_{v_{J_r}}}
= \sum_{r=1}^{3} {\rm e}^{z_{v_1} + z_{v_{J_r'}}}
 \stackrel{{\rm Proposition}~\ref{prop:BCH}}{=}
\sum_{r=1}^{3} {\rm e}^{z_{v_1}} {\rm e}^{z_{v_{J_r'}}}\\
& = {\rm e}^{z_{v_1}} \sum_{r=1}^{3} {\rm e}^{z_{v_{J_r'}}}
 \stackrel{{\rm equation}~\eqref{eq:oht_eta_f_omega}}{=}
{\rm e}^{z_{v_1}} f^\omega_{[\eta],\Delta'}.
\end{align*}
In the last equality, we used $v_{J_1'} = v_a+v_b$, $v_{J_2'} = -v_a+v_b$ and $v_{J_3'} = -v_a - v_b$. We also showed~$z_{v_1}$ strongly commutes with~$f^\omega_{[\eta],\Delta'}$, in the sense of Definition~\ref{def:algebraic_SC1}.

To deal with $J_5$, $J_6$, $J_7$, note that the juncture-state values for the two endpoints of the segment $\circled{7}$ lying at arcs $a$ and $b$, surrounded by the vertical bars in Table \ref{tab:hole-surrounding_loop_in_one-holed_torus}, form the values of the admissible juncture-states $J_1'$, $J_2'$, $J_3'$ of the loop $\eta$ at the endpoints of the segment $\circled{$1'$}$. Taking the signed sum of the elements $v_i$ for the entries in the rows for $J_5$, $J_6$, $J_7$ outside the vertical bars, one defines the element $v_2 \in N_{\Delta'}$ as
\[
v_2 := -v_c - v_a- v_b - v_c - v_d + v_d = - v_a - v_b - 2v_c.
\]
One can immediately check that $v_1 - v_2 = v_\gamma$ holds as stipulated by part of \ref{TR2} of Definition~\ref{def:Teschner_triple_level1},
and that for the `connecting term' $Z^{J_4}$ we have $v_{J_4} = v_1 + v_2$ as required for the connecting term in the Teschner recursion relation (equation~\eqref{eq:TR1}) in Definition~\ref{def:Teschner_triple_level1}\,\ref{TR1}. We~have
\[
v_{J_5} = v_2 + v_{J_1'}, \qquad
v_{J_6} = v_2 + v_{J_2'}, \qquad
v_{J_7} = v_2 + v_{J_3'}.
\]
Let us show $\langle v_2, v_{J_r'} \rangle=0$ for $r=1,2,3$. It suffices to show $\langle v_2,v_a\rangle = \langle v_2,v_b \rangle=0$, which one can easily verify using $v_2 = -v_a-v_b-2v_c$ and equation~\eqref{eq:skew-form_for_Delta}. Thus,
\begin{align*}\begin{split}
\sum_{r=5}^{7} Z^{J_r} & = \sum_{r=5}^{7} {\rm e}^{z_{v_{J_r}}}
= \sum_{r=5}^{7} {\rm e}^{z_{v_2} + z_{v_{J_{r-4}'}}}
 \stackrel{{\rm Proposition}~\ref{prop:BCH}}{=}~
\sum_{r=5}^{7} {\rm e}^{z_{v_2}} {\rm e}^{z_{v_{J_{r-4}'}}} \\
 & = {\rm e}^{z_{v_2}} \sum_{r=1}^{3} {\rm e}^{z_{v_{J_r'}}}
 \stackrel{{\rm equation}~\eqref{eq:oht_eta_f_omega}}{=}
{\rm e}^{z_{v_2}} f^\omega_{[\eta],\Delta'}.
\end{split}
\end{align*}
We also showed $z_{v_2}$ strongly commutes with $f^\omega_{[\eta],\Delta'}$, in the sense of Definition~\ref{def:algebraic_SC1}. Hence, we~have checked the condition \ref{TR3} of Definition~\ref{def:Teschner_triple_level1}.

One now observes
\begin{align*}
& f^\omega_{[\gamma],\Delta'}
\stackrel{\substack{{\rm Proposition}~\ref{prop:BW_admissible_sum_formula} \\ {\rm Table}~\ref{tab:hole-surrounding_loop_in_one-holed_torus}}}{=}
\sum_{r=0}^{8} Z^{J_r} \\
& \hspace{19mm} = Z^{J_0} + Z^{J_8} + Z^{J_4} + \bigl(Z^{J_1}+Z^{J_2}+Z^{J_3}\bigr) + \bigl(Z^{J_5}+Z^{J_6}+Z^{J_7}\bigr) \\
& \hspace{19mm} = {\rm e}^{z_{v_\gamma}} + {\rm e}^{-z_{v_\gamma}} + {\rm e}^{z_{v_1} + z_{v_2}}
+ {\rm e}^{z_{v_1}} f^\omega_{[\eta],\Delta'} + {\rm e}^{z_{v_2}} f^\omega_{[\eta],\Delta'},
\end{align*}
thus the Teschner recursion relation in \ref{TR1} of Definition~\ref{def:Teschner_triple_level1} holds. By inspection, the condition \ref{TR4} also holds, as $f^\omega_{[\eta],\Delta'}$ coincides with itself.

Finally, it remains to check $\langle v_1,v_2\rangle=-4$ for \ref{TR2}. Recall $v_1 = v_a + v_b + 2v_d$ and $v_2 = -v_a-v_b-2v_c$. We computed $\langle v_1, v_a\rangle = 0 = \langle v_1, v_b\rangle$, so $\langle v_1,v_2 \rangle = \langle v_1, -2v_c\rangle = -2( \varepsilon_{ac} + \varepsilon_{bc} + 2\varepsilon_{dc}) = -2(1 - 1 + 2) = -4$, as desired.

This finishes our proof showing that $(\gamma,\eta,\eta)$ is a strong Teschner triple with respect to $\Delta'$, in the sense of Definition~\ref{def:Teschner_triple_level1}.
\end{proof}

In fact, the arguments to prove Table \ref{tab:hole-surrounding_loop_in_one-holed_torus}, specifically the pattern colored red there, applies to any oriented loop whose oriented loop segments with respect to an ideal triangulation consists only of right turns except for one left turn; inspired by Lemma~\ref{lem:peripheral_characterization}, we make the following definition.
\begin{Definition}
\label{def:AP}
Let $\mathfrak{S}$ be a triangulable marked surface and $\Delta$ its ideal triangulation. Let~$\xi$ be an oriented essential simple loop in $\mathfrak{S}$ in a minimal and transversal position with respect to~$\Delta$ (Definition~\ref{def:minimal_and_transversal}).

In case all the oriented loop segments of $\xi$ with respect to $\Delta$ are right turns except for exactly one left turn, or in case all segments are left turns except for one right turn, we say $\xi$ is an {\em almost-peripheral} loop with respect to $\Delta$.
\end{Definition}
Such a loop will appear repeatedly throughout the paper, so we make a lemma about its quantized trace-of-monodromy.
\begin{Lemma}[quantized trace-of-monodromy of an almost-peripheral loop]
\label{lem:AP}
Let $\mathfrak{S}$, $\Delta$ and an almost-peripheral loop $\xi$ be as in Definition~{\rm \ref{def:AP}}. Label the oriented loop segments of $\xi$ formed by the junctures $\xi \cap \Delta$ by {\rm $\circled{1}, \circled{2},\dots,\circledd{$N$}$}, appearing consecutively in this order along the orientation. Suppose that {\rm $\circled{$k$}$} is a left turn for a unique $k \in\{1,\dots,N\}$, and all other oriented loop segments are right turns.

Denote the junctures, i.e., elements of $\xi \cap \Delta$, by $x_0,x_1,\dots,x_N$, where the loop segment {\rm $\circled{$i$}$} goes from $x_{i-1}$ to $x_i$; in particular, $x_0 = x_N$, and $\xi\cap\Delta = \{x_1,x_2,\dots,x_N\}$.

For $r=0,1,2,\dots,N$, define $J_r \colon \xi \cap \Delta \to \{+,-\}$ to be the juncture-state given by
\begin{align*}
J_r(x_{k+j}) = \begin{cases}
- & \mbox{for all $j\in \{0,1,\dots,r-1\}$}, \\
+ & \mbox{for all $j\in \{r,\dots,N-1\}$},
\end{cases}
\end{align*}
where the subscripts of $x_*$ are considered modulo $N$; i.e., $x_{i+N} = x_i$. Particular cases are $J_0$ and $J_N$, which can be understood as being the constant juncture-states $J_0 = J_\xi^+$ and $J_N = J_\xi^-$ by definition $($see equation~\eqref{eq:constant_juncture-state}$)$.

Then, the following hold:
\begin{enumerate}[label={\rm (AP\arabic*)}]\itemsep=0pt\setlength{\leftskip}{0.28cm}
\item\label{AP1} $\{J_0,J_1,\dots,J_N\}$ is the set of all admissible juncture-states.
\item\label{AP2} For $r=1,\dots,N$, let $i_r$ be the ideal arc of $\Delta$ on which the juncture $x_r$ lies in. Then
\begin{align}
\label{eq:AP2}
v_{J_r} = v_{J_{r-1}} - 2 \til{v}_{i_r},
\end{align}
where $v_{J_r}$ and $v_{J_{r-1}}$ are as in equation~\eqref{eq:v_J}, and $\til{v}_i$ is as in equation~\eqref{eq:til_v_i}.

\item\label{AP3} The non-redundancy condition in equation~\eqref{eq:admissible_non-redundancy_condition} is satisfied, hence Proposition~{\rm \ref{prop:BW_admissible_sum_formula}} applies, so that \smash{$f^\omega_{[\xi],\Delta}$} $($Definition~{\rm \ref{def:quantized_trace-of-monodromy})} equals the sum of the Weyl-ordered Laurent monomials~$Z^{J_r}$ $($equation~\eqref{eq:Z_J}$)$:
\begin{align*}
f^\omega_{[\xi],\Delta} = \sum_{r=0}^N Z^{J_r}.
\end{align*}
\end{enumerate}
\end{Lemma}

\begin{proof} Without loss of generality, we assume $k=1$ for convenience. Then one would just apply cyclic shifts of labelings of loop segments and junctures to deduce the lemma for general $k$.

Let $J$ be an admissible juncture-state of $\xi$. We first establish the following claim.

\begin{Claim*} If $J(x_r)=+$ for some $r \in \{1,\dots,N-1\}$, then $J(x_s)=+$ for all $s \in \{r,r+1,\dots,N\}$.
\end{Claim*}

Suppose $J(x_r)=+$ for some $r \in \{1,\dots,N-1\}$. Since $x_r$ is the initial end of the oriented loop segment $\circlee{r+1}$ which is a right turn, for which the juncture-state value $+ \to +$ is admissible but $+ \to -$ is not, it follows that $J(x_{r+1})=+$. By induction, $J(x_s)=+$ for all $s = r,r+1,\dots,N$, proving the claim.

From this claim, it follows that in case there exists $r\in \{1,\dots,N\}$ such that $J(x_r)=+$, the sequence of values $J(x_1), J(x_2),\dots,J(x_N)$ is completely determined by the smallest $r \in \{1,\dots,N\}$ such that $J(x_r)=+$, and the sequence is $-,-,\dots,-,+,+,\dots,+$, with the first $r-1$ values being $-$ and the rest being $+$. This is indeed an admissible juncture-state, which equals $J_{r-1}$ in the statement of Lemma~\ref{lem:AP}.

In case there is no $r\in \{1,\dots,N\}$ such that $J(x_r)=+$, we have $J(x_r)=-$ for all $r=1,\dots,N$, which is the $-$-valued constant juncture-state $J_\xi^-$ (equation~\eqref{eq:constant_juncture-state}). It is admissible, and equals $J_N$ in the statement of Lemma~\ref{lem:AP}. Hence, item \ref{AP1} is established.

For $r=1,\dots,N$, we notice that $J_r(x_j) = J_{r-1}(x_j)$ for all $j\in\{1,\dots,N\}$ except for $j=r$, where $J_r(x_r)=-$ and $J_{r-1}(x_r)=+$. Thus, equation~\eqref{eq:AP2} follows from the definition of~$v_{J_r}$ and~$v_{J_{r-1}}$ in equation~\eqref{eq:v_J}. This settles~\ref{AP2}.

By \ref{AP2}, for any $r,s\in \{0,\dots,N\}$ with $r>s$, we have $v_{J_s} - v_{J_r} = \sum_{k=s+1}^r 2 \til{v}_{i_k}$, hence one deduce that $v_{J_s} - v_{J_r} \neq 0$. Thus, $v_{J_s} \neq v_{J_r}$, and therefore $Z^{J_s} \neq Z^{J_r}$. Therefore, \ref{AP3} indeed holds by Proposition~\ref{prop:BW_admissible_sum_formula}.
\end{proof}

\subsection{Hole surrounding loop in one-holed subsurface of higher genus}
\label{subsec:hole-surrounding_loop_in_one-holed_subsurface}

We now deal with the remaining part of the item \ref{L6} of Proposition~\ref{prop:loop_classification} when the genus $g$ is at least $2$. Namely, we deal the quantized trace-of-monodromy $f^\omega_{[\gamma],\Delta'}$ for a hole surrounding loop~$\gamma$ (Definition~\ref{def:hole-surrounding}) of a minimal-marked (Definition~\ref{def:minimal-marked}) one-holed surface $\mathfrak{S}' = \Sigma'\setminus\mathcal{P}'$ of genus $g\ge 2$ (Definition~\ref{def:surface_terminology}).

\begin{Proposition}[hole surrounding loop in higher genus one-holed subsurface]
\label{prop:hole-surrounding_loop_in_higher_genus_one-holed_subsurface}
Let $\mathfrak{S}'=\Sigma'\setminus\mathcal{P}'$ be a triangulable marked surface that is
isomorphic to a minimal-marked one-holed surface of genus $g\ge 2$, and let $\gamma$ be a hole surrounding loop in $\mathfrak{S}'$. Then there exist non-separating essential simple loops $\eta$ and $\sigma$ in $\mathfrak{S}'$ and an ideal triangulation $\Delta'$ of $\mathfrak{S}'$ such that $(\gamma,\eta,\sigma)$ is a~weak Teschner triple with respect to $\Delta'$, in the sense of Definition~{\rm \ref{def:Teschner_triple_level2}}.
\end{Proposition}

The remainder of the present subsection is wholly devoted to a proof of this proposition. For convenience, the proof is divided into several steps, each forming a paragraph.

\subsubsection{Proof of Proposition~\ref{prop:hole-surrounding_loop_in_higher_genus_one-holed_subsurface}: The topological set-up}

Since $\mathcal{P}'$ consists of just one point, for any chosen ideal triangulation $\Delta'$ of $(\Sigma',\mathcal{P}')$ the two endpoints of each ideal arc of $\Delta'$ coincide with the unique marked point of $(\Sigma',\mathcal{P}')$. So a~hole surrounding loop $\gamma$ in $(\Sigma',\mathcal{P}')$ in a minimal and transversal position with respect to $\Delta'$ (Definition~\ref{def:minimal_and_transversal}) meets each internal ideal arc of $\Delta'$ twice, and does not meet the boundary arc of $\Delta'$. These intersection numbers between $\gamma$ and the ideal arcs of $\Delta'$ uniquely determine the simple loop $\gamma$ in $\Sigma'\setminus\mathcal{P}'$ up to isotopy. Indeed, if we pick two points in each interior arc of $\Delta'$, then in each ideal triangle of $\Delta'$, there is a unique (up to isotopy rel boundary) way of connecting these points by non-intersecting simple paths in this triangle.

\begin{figure}[!ht]
\centering
\scalebox{1.0}{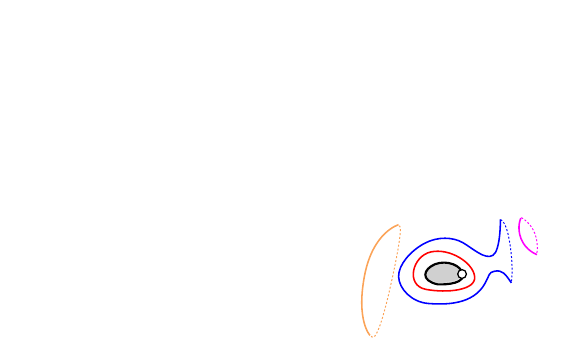}
\vspace{-4mm}
\caption{Minimal-marked one-holed surface $(\Sigma',\mathcal{P}')$ of higher genus $g\ge 2$.}
\label{fig:one-holed_surface_of_higher_genus}
\end{figure}

\begin{figure}[!ht]
\centering
\scalebox{1}{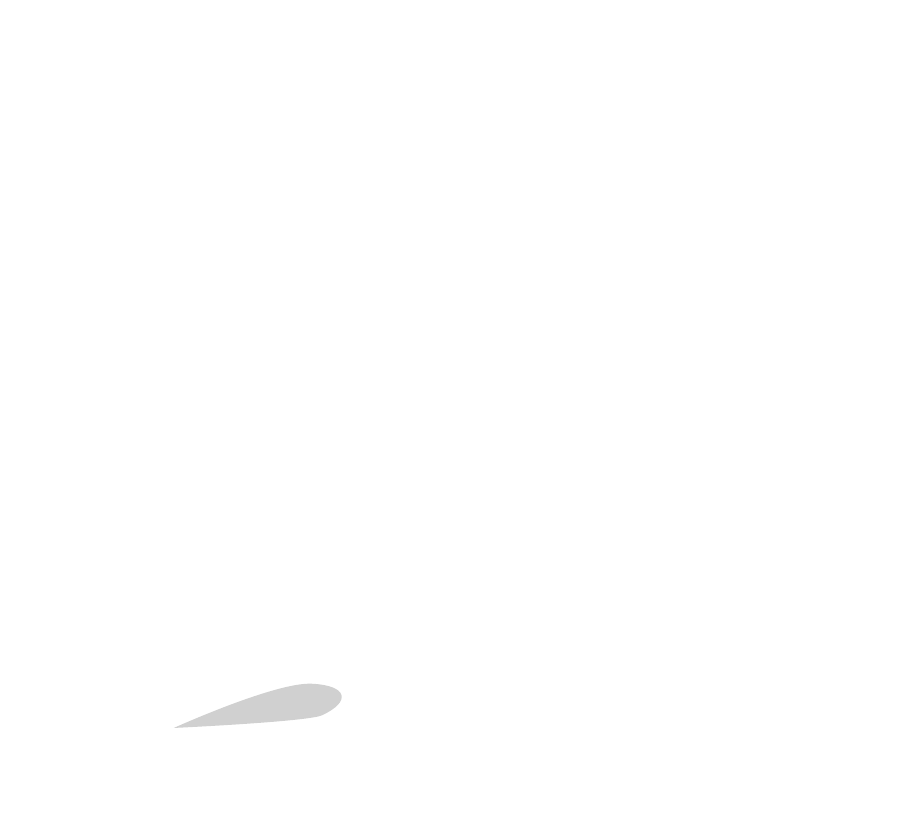}
\vspace{-3mm}
\caption{One-holed surface $(\Sigma',\mathcal{P}')$ of genus $g\ge 2$, as quotient of $4g$-gon, with an ideal triangulation $\Delta'$; this picture shows the case $g=4$.}\label{fig:one-holed_surface_of_higher_genus_flat}\vspace{-2mm}
\end{figure}

A compact surface of genus $g$ can be obtained by attaching $g$ `handles' to the sphere $S^2$; then we remove an open disc to obtain $\Sigma'$ as portrayed in Figure~\ref{fig:one-holed_surface_of_higher_genus}. Recall that one can realize a~compact surface of genus $g$ as quotient of a $4g$-gon, with its edges identified in pairs according to the word of oriented labels of sides being \smash{$a_1b_1a_1^{-1}b_1^{-1} a_2b_2 a_2^{-1} b_2^{-1} \cdots a_gb_ga_g^{-1}b_g^{-1}$}. Here, each subword \smash{$a_i b_i a_i^{-1} b_i^{-1}$} corresponds to the $i$-th handle. Such a $4g$-gon is depicted in Figure~\ref{fig:one-holed_surface_of_higher_genus_flat}.

We choose to work with a specific ideal triangulation $\Delta'$ of $(\Sigma',\mathcal{P}')$ as depicted in Figure~\ref{fig:one-holed_surface_of_higher_genus_flat} for the case $g=4$; the case for other $g\ge 2$ can be deduced from this picture. For convenience, we denote the exchange matrix of $\Delta'$ by $\varepsilon = (\varepsilon_{ij})_{i,j\in \Delta'}$. The advantage of our choice of an ideal triangulation $\Delta'$ is that one can recognize which triangle contributes to which of the $g$ handles relatively easily, as annotated in Figures~\ref{fig:one-holed_surface_of_higher_genus}--\ref{fig:one-holed_surface_of_higher_genus_flat}. We now explain the recipe for the labeling of ideal arcs of $\Delta'$; a reader could skip the following description in words, and just refer to Figure~\ref{fig:one-holed_surface_of_higher_genus_flat}. For convenience, let us call the $4g$ vertices of the $4g$-gon in Figure~\ref{fig:one-holed_surface_of_higher_genus_flat} by $p_1,p_2,\dots,p_{4g}$, arranged counterclockwise in this order, where $p_1$ is the unique vertex lying in the hole in Figure~\ref{fig:one-holed_surface_of_higher_genus_flat}. We remark that these vertices are in fact identical in the actual surface $\Sigma'$ or $\Sigma$, and we named them differently just for this 4$g$-gon picture, for convenience.
\begin{enumerate}[label={\rm (IA\arabic*)}]\itemsep=0pt\setlength{\leftskip}{0.12cm}
\item For each $i = 1,2,\dots,g$ (representing the $i$-th handle), the ideal arc connecting $p_{4(i-1)+1}$ and $p_{4(i-1)+2}$ is denoted by $a_i$, that connecting $p_{4(i-2)+2}$ and $p_{4(i-2)+3}$ is denoted by $b_i$.

\item For each $i=2,\dots,g$, the ideal arc connecting $p_{4(i-1)+1}$ and $p_{4(i-1)+3}$ is denoted by $c_i$, and that connecting $p_{4(i-1)+3}$ and $p_{4i+1}$ by $d_i$, where $p_{4g+1}$ is set to mean $p_1$.

\item For $i=1$, the two ideal arcs connecting $p_1$ and $p_3$ are called $c_1$ and $c_1'$, so that the three arcs $a_1,b_1$ and $c_1$ form an ideal triangle of $\Delta'$.

\item The ideal arc connecting $p_3$ and $p_5$ is called $d_1$.

\item For each $i=2,\dots,g-1$ (for $g \ge 2$), the ideal arc connecting $p_1$ and $p_{4(i-1)+1}$ is called~$f_i$.

\item For each $i=2,\dots,g$, the ideal arc connecting $p_{4(i-1)+1}$ and $p_{4i+1}$ is called $e_i$.
\end{enumerate}
In particular, for example, if $g=2$, then there are only the first and the last handles, and no `middle' handles, so that there is no ideal arc labeled by $f_i$ for some $i$.

In view of the earlier description of a hole surrounding loop $\gamma$ in a minimal and transversal position in terms of the intersection numbers with the ideal arcs of $\Delta'$, we see that our loop $\gamma$ can be isotoped to look as in Figure~\ref{fig:one-holed_surface_of_higher_genus_flat} with respect to $\Delta'$. Now, define two more loops $\eta$ and~$\sigma$ in $(\Sigma',\mathcal{P}')$ as in Figure~\ref{fig:one-holed_surface_of_higher_genus_flat}.
We will show that $\eta$ and~$\sigma$ are non-separating, and that $(\gamma,\eta,\sigma)$ is a weak Teschner triple with respect to $\Delta'$, in the sense of Definition~\ref{def:Teschner_triple_level2}; under the notation there, we could say $\eta = \gamma_1$ and $\sigma = \gamma_2$.

For convenience, give orientations on the loops $\gamma$, $\eta$ and $\sigma$ as in Figure~\ref{fig:one-holed_surface_of_higher_genus_flat}. We label the oriented loop segments of $\gamma$ by $\circled{1},\circled{2},\dots,\circled{$N$}$ arranged in the consecutive order, as in Figure~\ref{fig:one-holed_surface_of_higher_genus_flat}; it is straightforward to verify that $\Delta'$ has $4g-1$ ideal triangles, and that the number of loop segments of $\gamma$ is $N = 12g-4$. The loop $\eta$ consists of $N-6$ oriented loop segments, which we label as $\circled{$7'$},\circled{8},\dots,\circled{$N$}$ in the consecutive order, as in Figure~\ref{fig:one-holed_surface_of_higher_genus_flat}. For $j=8,9,\dots,N$, the loop segment of $\gamma$ labeled by $\circled{$j$}$ is homotopic to that of $\eta$ labeled by the same symbol through a~hotomopy of simple paths in a triangle under the condition that the endpoints of a loop segment should lie in the interior of an ideal arc of $\Delta'$ at all times. Note that the loop $\sigma$ consists of only two oriented loop segments, which we label as $\circled{$1'$},\circled{$2'$}$ as in Figure~\ref{fig:one-holed_surface_of_higher_genus_flat}.

The loop $\sigma$ meets only twice with the triangulation $\Delta'$, one with $a$ and one with $b$; in view of the earlier discussion on the handles and subwords, it is not hard to observe that $\sigma$ is a non-separating loop in the first handle, i.e., as depicted in Figure~\ref{fig:one-holed_surface_of_higher_genus}. One can observe that the loop~$\eta$ which is defined in Figure~\ref{fig:one-holed_surface_of_higher_genus_flat} looks as in Figure~\ref{fig:one-holed_surface_of_higher_genus}. We leave the verification of this fact as an exercise to the reader; see a previous version of the present paper \cite{ver3} for a detailed explanation. Then, it is easy to see from Figure~\ref{fig:one-holed_surface_of_higher_genus} that $\eta$ is non-separating, and that $(\gamma,\eta,\sigma)$ is a pants triple (Definition~\ref{def:Teschner_triple_level1}).

\subsubsection{Proof of Proposition~\ref{prop:hole-surrounding_loop_in_higher_genus_one-holed_subsurface}: the state-sum formulas}

Now we compute the quantized trace-of-monodromy for the loops $\gamma$, $\eta$ and $\sigma$. Let us start with the easy one, namely $\sigma$, which consists only of two segments $\circled{$1'$}$ and $\circled{$2'$}$, meeting two arcs~$a_1$ and~$b_1$ of $\Delta'$, with $\varepsilon_{a_1b_1} = -2$. Hence, Proposition~\ref{prop:quantized_trace-of-monodromy_for_loop_with_int_2} applies, giving us\vspace{-1mm}
\begin{align}
\label{eq:ohs_proof_f_omega_sigma}
f^\omega_{[\sigma],\Delta} = {\rm e}^{z_{v_{a_1}} + z_{v_{b_1}}} + {\rm e}^{-z_{v_{a_1}} + z_{v_{b_1}}} + {\rm e}^{-z_{v_{a_1}}-z_{v_{b_1}}}.\vspace{-1mm}
\end{align}

We move on to the loops $\gamma$ and $\eta$. We write down whether each oriented loop segment is a left turn or a right turn, and indicate by the symbols L and R. It is not hard to see that, among the oriented loop segments $\circled{1},\dots,\circled{$N$}$ of $\gamma$, the only left turn is $\circled{2}$. Hence, $\gamma$ is an almost-peripheral loop (Definition~\ref{def:AP}), and Lemma~\ref{lem:AP} applies. Likewise, among the oriented loop segments $\circled{$7'$}$, $\circled{8}$, \dots, $\circled{$N$}$ of $\eta$, the only left turn is $\circled{$7'$}$, hence $\eta$ is an almost-peripheral loop and Lemma~\ref{lem:AP} applies. This justifies the lists of all admissible juncture-states of $\gamma$ and $\eta$ presented in Table \ref{tab:hole-surrounding_loop_in_one-holed_higher_genus_surface}. In this table, we also record the label of the ideal arc on which each juncture lives~in. The rightmost juncture is the same as the leftmost one, so the information there are put inside parentheses.\looseness=-1

\begin{table}[htb]
\centering
{\setlength{\tabcolsep}{0,2mm}
\begin{tabular}{@{}r @{\hspace{1mm} } llllllllllll@{\hspace{2mm}}lr @{\hspace{4mm}}l@{}}
$\gamma$ : & \tablearrow{1}{R} & \tablearrow{2}{L} & \tablearrow{3}{R} & \tablearrow{4}{R} & \tablearrow{5}{R} & \tablearrow{6}{R} & \tablearrow{7}{R} & \tablearrow{8}{R} & \tablearrow{9}{R} & \tablearrowd{10}{R} & \tablearrowd{11}{R} & $\hspace{0,5mm} \bullet \hspace{-1,3mm} \overset{\cdots}{\longrightarrow} \hspace{-1,3mm}$ & $\cdots$ & \tablearrow{$N$}{R} $\bullet$ & \\
arcs for $\gamma$ : & $f_2$ & $c_1'$ & $c_1$ & $a_1$ & $b_1$ & $c_1$ & $c_1'$ & $d_1$ & $a_1$ & $b_1$ & $d_1$ & $f_2$ & & ($f_2$)\hspace{-2mm} \\
$\eta$ : & \multicolumn{7}{l}{$\hspace{-0,5mm}\xymatrix@C+41,8mm{\bullet \hspace{-0,8mm} \ar[r]^{\overset{\circled{$7'$}}{\rm L}} & } \hspace{-3,0mm}$} & \tablearrow{8}{R} & \tablearrow{9}{R} & \tablearrowd{10}{R} & \tablearrowd{11}{R} & $\hspace{0,5mm} \bullet \hspace{-1,3mm} \overset{\cdots}{\longrightarrow} \hspace{-1,3mm}$ & $\cdots$ & \tablearrow{$N$}{R} $\bullet$ \\
arcs for $\eta$ : & \multicolumn{7}{l}{$f_2$} & $d_1$ & $a_1$ & $b_1$ & $d_1$ & $f_2$ & & ($f_2$) \hspace{-3,7mm} \\ \hline
\multicolumn{10}{l}{all admissible juncture-states for $\gamma$:} \\
$J_1$ : & + & + & \textcolor{red}{+} & + & + & + & + & + & + & + & + & + & $\cdots$ & (+) \hspace{-3,5mm} & all plus \\ \hline
$J_2$ : & + & + & \textcolor{red}{$-$} \vbar & \textcolor{red}{+} & + \vbar & + & + & + & + & + & + & + & $\cdots$ & (+) \hspace{-3,5mm} & \rdelim\}{3}{1pt}[corr. to $\sigma$] \\
$J_3$ : & + & + & \textcolor{red}{$-$} \vbar & \textcolor{red}{$-$} & \textcolor{red}{+} \vbar & + & + & + & + & + & + & + & $\cdots$ & (+) \hspace{-3,5mm} & \\
$J_4$ :& + & + & \textcolor{red}{$-$} \vbar & \textcolor{red}{$-$} & \textcolor{red}{$-$} \vbar & \textcolor{red}{+} & + & + & + & + & + & + & $\cdots$ & (+) \hspace{-3,5mm} & \\ \hline
$J_5$ : & + & + & \textcolor{red}{$-$} & \textcolor{red}{$-$} & \textcolor{red}{$-$} & \textcolor{red}{$-$} & \textcolor{red}{+} & + & + & + & + & + & $\cdots$ & (+) \hspace{-3,5mm} & connecting term \\ \hline
$J_6$ : & + \vbar & + & \textcolor{red}{$-$} & \textcolor{red}{$-$} & \textcolor{red}{$-$} & \textcolor{red}{$-$} & \textcolor{red}{$-$} \vbar & \textcolor{red}{+} & + & + & + & + & $\cdots$ & (+) \hspace{-3,5mm} & \rdelim\}{5}{1pt}[corr. to $\eta$] \\
$J_7$ : & + \vbar & + & \textcolor{red}{$-$} & \textcolor{red}{$-$} & \textcolor{red}{$-$} & \textcolor{red}{$-$} & \textcolor{red}{$-$} \vbar & \textcolor{red}{$-$} & \textcolor{red}{+} & + & + & + & $\cdots$ & (+) \hspace{-3,5mm} & \\
$\vdots$ \hspace{2mm} & & & & & & & & & & & & & $\cdots$& & \\
\hspace{-3mm} $J_{N-1}$ : & \textcolor{red}{+} \vbar & + & \textcolor{red}{$-$} & \textcolor{red}{$-$} & \textcolor{red}{$-$} & \textcolor{red}{$-$} & \textcolor{red}{$-$} \vbar & \textcolor{red}{$-$} & \textcolor{red}{$-$} & \textcolor{red}{$-$} & \textcolor{red}{$-$} & \textcolor{red}{$-$} & $\cdots$ & \textcolor{red}{$-$} \hspace{-1mm} \textcolor{red}{(+)} \hspace{-3,5mm} & \\
$J_N$ : & \textcolor{red}{$-$} \vbar & \textcolor{red}{+} & \textcolor{red}{$-$} & \textcolor{red}{$-$} & \textcolor{red}{$-$} & \textcolor{red}{$-$} & \textcolor{red}{$-$} \vbar & \textcolor{red}{$-$} & \textcolor{red}{$-$} & \textcolor{red}{$-$} & \textcolor{red}{$-$} & \textcolor{red}{$-$} & $\cdots$ & \textcolor{red}{$-$} \hspace{-1mm} \textcolor{red}{($-$)} \hspace{-3,5mm} & \\ \hline
\hspace{-3mm} $J_{N+1}$ : & \textcolor{red}{$-$} & \textcolor{red}{$-$} & \textcolor{red}{$-$} & \textcolor{red}{$-$} & \textcolor{red}{$-$} & \textcolor{red}{$-$} & \textcolor{red}{$-$} & \textcolor{red}{$-$} & \textcolor{red}{$-$} & \textcolor{red}{$-$} & \textcolor{red}{$-$} & \textcolor{red}{$-$} & $\cdots$ & \textcolor{red}{$-$} \hspace{-1mm} \textcolor{red}{($-$)} \hspace{-3,5mm} & all minus \\ \hline
\multicolumn{10}{l}{all admissible juncture-states for $\eta$:} \\
$J_6'$ : & + & \multicolumn{6}{c}{ } & + & + & + & + & + & $\cdots$ & (+) \hspace{-3,5mm} & \\
$J_7'$ : & + & \multicolumn{6}{c}{ } & $-$ & $+$ & + & + & + & $\cdots$ & (+) \hspace{-3,5mm} & \\
$J_8'$ : & + & \multicolumn{6}{c}{ } & $-$ & $-$ & + & + & + & $\cdots$ & (+) \hspace{-3,5mm} & \\
$\vdots$ \hspace{2mm} & & & & & & & & & & & & & $\cdots$& & \\
$J_{N-1}'$ : & + & \multicolumn{6}{c}{ } & $-$ & $-$ & $-$ & $-$ & $-$ & $\cdots$ & $-$ \hspace{-1mm} (+) \hspace{-3,5mm} & \\
$J_{N}'$ : & $-$ & \multicolumn{6}{c}{ } & $-$ & $-$ & $-$ & $-$ & $-$ & $\cdots$ & $-$ \hspace{-1mm} ($-$) \hspace{-3,5mm} & \\ \hline
\end{tabular}}
\caption{All admissible juncture-states for a hole surrounding loop $\gamma$ and a non-separating loop $\eta$ in a minimal-marked one-holed surface $(\Sigma',\mathcal{P}')$ of higher genus $g\ge 2$, with respect to the ideal triangulation $\Delta'$.}\label{tab:hole-surrounding_loop_in_one-holed_higher_genus_surface}
\end{table}

By Lemma~\ref{lem:AP}\,\ref{AP3}, the non-redundancy condition in equation~\eqref{eq:admissible_non-redundancy_condition} is satisfied for both~$\gamma$ and~$\eta$, with respect to~$\Delta'$, so that Proposition~\ref{prop:BW_admissible_sum_formula} applies, yielding\vspace{-1mm}
\begin{align}
\label{eq:ohs_proof_sum_formula_for_two_loops}
f^\omega_{[\gamma],\Delta'} = \sum_{r=1}^{N+1} Z^{J_r}, \qquad
f^\omega_{[\eta],\Delta'} = \sum_{r=6}^{N} Z^{J_r'},\vspace{-1mm}
\end{align}
where $Z^J$ for a juncture-state $J$ is defined as in equation~\eqref{eq:Z_J}.

We now start investigating \smash{$f^\omega_{[\gamma],\Delta'}$} and \smash{$f^\omega_{[\eta],\Delta'}$} in more detail; we shall describe all the juncture-states appearing in the state-sum formula \eqref{eq:ohs_proof_sum_formula_for_two_loops}. We first compute the
the total intersection elements for the relevant loops (Definition~\ref{def:v_gamma}), by counting the intersection numbers of the loops with the ideal arcs of $\Delta'$ through an inspection on Figure~\ref{fig:one-holed_surface_of_higher_genus_flat}:\vspace{-1mm}
\begin{align}
\label{eq:ohs_proof_v_gamma}
&{v_\gamma} = \sum^g_{i=1} 2 {v_{a_i}} + \sum^g_{i=1} 2 {v_{b_i}} + 2{v_{c_1'}} + \sum^g_{i=1} 2 {v_{c_i}} + \sum^g_{i=1} 2 {v_{d_i}} + \sum^g_{i=2} 2 {v_{e_i}} + \sum^{g-1}_{i=2} 2 {v_{f_i}} \\[-1mm]
\nonumber
&{v_\sigma} = {v_{a_1}} + {v_{b_1}}, \\[-1mm]
\label{eq:ohs_proof_v_eta1}
&{v_\eta} = {v_{a_1}} + \sum^g_{i=2} 2 {v_{a_i}} + {v_{b_1}} + \sum^g_{i=1} 2 {v_{b_i}} + \sum^g_{i=2} 2 {v_{c_i}} + \sum^g_{i=1} 2 {v_{d_i}} + \sum^g_{i=2} 2 {v_{e_i}} + \sum^{g-1}_{i=2} 2 {v_{f_i}} \\[-1mm]
\label{eq:ohs_proof_v_eta2}
&\phantom{{v_\eta}}{} = {v_\gamma} - {v_\sigma} - 2 {v_{c_1'}} - 2 {v_{c_1}}.
\end{align}
Keep in mind $\til{v}_i=v_i$ here, for any ideal arc $i \in \Delta'$. Define the elements $v_1, v_2 \in N_{\Delta'}$ as
\begin{align}
\label{eq:ohs_proof_v_1_v_2}
v_1 := - v_{a_1} - v_{b_1} - 2v_{c_1}, \qquad
v_2 := v_1 + v_\gamma,
\end{align}
In particular, $v_2 - v_1 = v_\gamma$, as stipulated by part of \ref{TR2} of Definition~\ref{def:Teschner_triple_level2}. For $r=6,7,\dots,N$, note from Table \ref{tab:hole-surrounding_loop_in_one-holed_higher_genus_surface} that once we delete the entries from the row for $J_r$ for the six junctures for the loop segments $\circled{2},\circled{3},\dots,\circled{6}$ (these entries are bounded by vertical bars in the table), we get the row for $J_r'$. Reading these six entries, we get
\[
v_{J_r} = v_{J_r'} + (v_{c_1'} - v_{c_1} - v_{a_1} - v_{b_1} - v_{c_1} - v_{c_1'})
= v_{J_r'} + v_1, \qquad \forall r = 6,7,\dots,N;
\]
this is how we came up with $v_1$ in equation~\eqref{eq:ohs_proof_v_1_v_2}.

For the remaining admissible juncture-states of $\gamma$, observe from Table \ref{tab:hole-surrounding_loop_in_one-holed_higher_genus_surface} and equation~\eqref{eq:AP2} of Lemma~\ref{lem:AP}\,\ref{AP2}:
\begin{gather*}
v_{J_1} = v_\gamma, \\
v_{J_2} = v_\gamma - 2 v_{c_1}
= (v_{a_1} + v_{b_1}) + (-v_{a_1} - v_{b_1} - 2v_{c_1} + v_\gamma)
= \smash{(\underbrace{ v_{a_1} + v_{b_1} }_{=v_\sigma})} + v_2, \\
v_{J_3} = v_{J_2} - 2v_{a_1} = (-v_{a_1} + v_{b_1}) + v_2, \\
v_{J_4} = v_{J_3} - 2v_{b_1} = ( - v_{a_1} - v_{b_1}) + v_2, \\
v_{J_5} = v_{J_4} - 2v_{c_1} = (- v_{a_1} - v_{b_1} - 2v_{c_1}) + v_2 = v_1 + v_2, \\
v_{J_{N+1}} = - v_{J_1} = - v_\gamma.
\end{gather*}
Notice that $v_{J_2}$, $v_{J_3}$ and $v_{J_4}$ are of the form $v_J + v_2$ for some admissible juncture-state $J$ of $\sigma$, where $J$ can be read from the rows for $J_2$, $J_3$ and $J_4$ in Table \ref{tab:hole-surrounding_loop_in_one-holed_higher_genus_surface}, in between the vertical bars. That is, $v_2$ is obtained from the entries of these rows, outside the vertical bars. This is how we originally found $v_2$, but we expressed it as $v_2 = v_1 + v_\gamma$ as in equation~\eqref{eq:ohs_proof_v_1_v_2} for brevity.

\subsubsection{Proof of Proposition~\ref{prop:hole-surrounding_loop_in_higher_genus_one-holed_subsurface}: The Teschner recursion \ref{TR1}}

Let us first check some of the commutation relations required in \ref{TR3} of Definition~\ref{def:Teschner_triple_level2}, which we also need to obtain the structural result for \smash{$f^\omega_{[\gamma],\Delta'}$}, namely \ref{TR1}. More precisely, we shall first check $[z_{v_1}, z_{J_r'}]=0$ for all $r=6,7,\dots,N$, and $[z_{v_2}, z_{v_{a_1}}+ z_{v_{b_1}}] = [z_{v_2}, -z_{v_{a_1}} + z_{v_{b_1}}] = [z_{v_2}, - z_{v_{a_1}}- z_{v_{b_1}}]=0$. For this it suffices to check $\langle v_1,v_{J_r'}\rangle=0$ for all $r=6,\dots,N$, and $\langle v_2, v_{a_1} \rangle = \langle v_2, v_{b_1}\rangle=0$, in view of equation~\eqref{eq:HA_commutation}. Note from Definition~\ref{def:varepsilon_Delta} and Figure~\ref{fig:one-holed_surface_of_higher_genus_flat}:
\begin{gather}
\varepsilon_{a_1 b_1} = - 2, \qquad
\varepsilon_{a_1 c_1} = 1, \qquad
\varepsilon_{a_1 d_1} = 1, \nonumber\\
\varepsilon_{a_1 h} = 0 \quad \mbox{for all other ideal arcs $h \in \Delta'$}, \label{eq:ohs_proof_some_varepsilon1}\\
\varepsilon_{b_1 a_1} = 2, \qquad
\varepsilon_{b_1 c_1} = -1, \qquad
\varepsilon_{b_1 d_1} = -1, \nonumber\\
\varepsilon_{b_1 h} =0 \quad \mbox{for all other ideal arcs $h\in \Delta'$}.\label{eq:ohs_proof_some_varepsilon2}
\end{gather}
Using equation~\eqref{eq:skew-form_for_Delta}, above values of $\varepsilon_{ij}$ when $i\in \{a_1,b_1\}$, equation~\eqref{eq:ohs_proof_v_gamma} and equation~\eqref{eq:ohs_proof_v_1_v_2}, one can verify that
\begin{align}
\label{eq:ohs_proof_some_commutations1}
0 = \langle v_{a_1}, v_\gamma \rangle
= \langle v_{b_1}, v_\gamma \rangle
= \langle v_{a_1}, v_1\rangle
= \langle v_{b_1}, v_1\rangle
= \langle v_{a_1}, v_2 \rangle
= \langle v_{b_1}, v_2 \rangle.
\end{align}
Hence, $[z_{v_2}, z_{v_{a_1}}+ z_{v_{b_1}}] = [z_{v_2}, -z_{v_{a_1}} + z_{v_{b_1}}] = [z_{v_2}, - z_{v_{a_1}}- z_{v_{b_1}}]=0$, so
\begin{align*}
Z^{J_2} + Z^{J_3} + Z^{J_4}
& = {\rm e}^{(z_{v_{a_1}} + z_{v_{b_1}}) + z_{v_2}}
+ {\rm e}^{(-z_{v_{a_1}} + z_{v_{b_1}}) + z_{v_2}}
+ {\rm e}^{(-z_{v_{a_1}} - z_{v_{b_1}}) + z_{v_2}} \\
& = {\rm e}^{z_{v_2}}( {\rm e}^{z_{v_{a_1}} + z_{v_{b_1}}}
+ {\rm e}^{-z_{v_{a_1}} + z_{v_{b_1}}}
+ {\rm e}^{-z_{v_{a_1}} - z_{v_{b_1}}}) \quad (\because \mbox{Proposition~\ref{prop:BCH}}) \\
& = {\rm e}^{z_{v_2}} f^\omega_{[\sigma],\Delta'} \quad (\because \mbox{equation~\eqref{eq:ohs_proof_f_omega_sigma}}),
\end{align*}
and $z_{v_2}$ strongly commutes with $f^\omega_{[\sigma],\Delta'}$ in the sense of Definition~\ref{def:algebraic_SC1}.

By looking at the ideal arcs of $\Delta'$ intersecting $\eta$ in Figure~\ref{fig:one-holed_surface_of_higher_genus_flat}, one finds that for each $r\in\{6,\dots,N\}$,
\begin{align}
\label{eq:ohs_proof_v_J_r_prime_involves}
\mbox{$v_{J_r'}$ involves only $a_1,\dots,a_g$, $b_1,\dots,b_g$, $c_2,\dots,c_g$, $d_1,\dots,d_g$, $e_2,\dots,e_g$, $f_2,\dots,f_{g-1}$.}\!\!\!
\end{align}
From Definition~\ref{def:varepsilon_Delta} and Figure~\ref{fig:one-holed_surface_of_higher_genus_flat}, note
\[
\varepsilon_{c_1 a_1} = -1, \qquad
\varepsilon_{c_1 b_1} = 1, \qquad
\varepsilon_{c_1 c_1'} = -1, \qquad
\varepsilon_{c_1 e} = 1,
\]
where $e$ is the boundary arc for the hole of $\Sigma'$,
\[
\varepsilon_{c_1 h} =0 \quad \mbox{for all other ideal arcs $h \in \Delta'$}.
\]
Now, by equation~\eqref{eq:skew-form_for_Delta},
\begin{gather*}
 \langle v_{a_1}, v_{b_1}\rangle = -2, \qquad
\langle v_{a_1}, v_{c_1}\rangle = 1, \qquad
\langle v_{a_1}, v_{d_1}\rangle = 1, \\
\langle v_{a_1},v_h\rangle=0 \quad \mbox{for all other ideal arcs $h \in \Delta'$}, \\
\langle v_{b_1}, v_{a_1}\rangle=2,\qquad
\langle v_{b_1}, v_{c_1}\rangle = -1,\qquad
\langle v_{b_1},v_{d_1}\rangle = -1, \\
\langle v_{b_1},v_h\rangle=0 \quad \mbox{for all other ideal arcs $h \in \Delta'$}, \\
\langle v_{c_1}, v_{a_1}\rangle = -1, \qquad
\langle v_{c_1}, v_{b_1}\rangle = 1, \\
\langle v_{c_1}, v_h\rangle=0 \quad \mbox{for all other ideal arcs $h$ involved in $z_{J_r'}$}.
\end{gather*}
Thus, we can verify that
\begin{align}
\label{eq:ohs_proof_some_commutations2}
\langle v_1, v_i \rangle = \langle -v_{a_1}-v_{b_1}-2v_{c_1},v_i \rangle = 0 \quad \mbox{ for all ideal arcs $i$ involved in $z_{J_r'}$.}
\end{align}
Hence, $[z_{v_1},z_{J_r'}]=0$ for all $r=6,\dots,N$, so
\begin{align*}
\hspace{-2mm} \sum_{r=6}^{N} Z^{J_r} & = \sum_{r=6}^{N} \exp(z_{J_r})
= \sum_{r=6}^{N} \exp(z_{J_r'} + z_{v_1})
= \sum_{r=6}^{N} \exp(z_{v_1}) \exp(z_{J_r'}) \\
& = {\rm e}^{z_{v_1}} \sum_{r=6}^{N} Z^{J_r'}
\stackrel{{\rm equation}~\eqref{eq:ohs_proof_sum_formula_for_two_loops}}{=} {\rm e}^{z_{v_1}} f^\omega_{[\eta],\Delta'},
\end{align*}
and $z_{v_1}$ strongly commutes with $f^\omega_{[\eta],\Delta'}$ in the sense of Definition~\ref{def:algebraic_SC1}.

Assembling what we observed so far, we have
\begin{align*}
f^\omega_{[\gamma],\Delta'}
 & \stackrel{{\rm equation}~\eqref{eq:ohs_proof_sum_formula_for_two_loops}}{=} Z_{J_1} + Z_{J_{N+1}} + (Z_{J_2}+Z_{J_3}+Z_{J_4}) + Z_{J_5} + \sum_{r=6}^{N+1} Z^{J_r} \\
 & \hspace{7.9mm} = {\rm e}^{v_\gamma} + {\rm e}^{-v_\gamma} + {\rm e}^{z_{v_2}} f^\omega_{[\sigma],\Delta'} + {\rm e}^{z_{v_1} + z_{v_2}} + {\rm e}^{z_{v_1}} f^\omega_{[\eta],\Delta'},
\end{align*}
so \ref{TR1} of Definition~\ref{def:Teschner_triple_level2}, i.e., the Teschner recursion, holds.

\subsubsection{Proof of Proposition~\ref{prop:hole-surrounding_loop_in_higher_genus_one-holed_subsurface}: Verification of \ref{TR3} and \ref{TR2}}

We still have to check the remaining conditions of Definition~\ref{def:Teschner_triple_level2}. Let us start with \ref{TR3}, which so far we have proved only half of. Namely, we showed that $z_{v_1}$ strongly commutes with $f^\omega_{[\eta],\Delta'}$ and that $z_{v_2}$ strongly commutes with $f^\omega_{[\sigma],\Delta'}$. We still have to show that $z_{v_1}$ strongly commutes with $f^\omega_{[\sigma],\Delta'}$ and that $z_{v_2}$ strongly commutes with $f^\omega_{[\eta],\Delta'}$. For the former, recall from equation~\eqref{eq:ohs_proof_some_commutations1} that we showed $\langle v_{a_1},v_1 \rangle =\langle v_{b_1},v_1\rangle=0$. From equation~\eqref{eq:ohs_proof_f_omega_sigma}, we see that $f^\omega_{[\sigma],\Delta'}$ only involves $z_{v_i}$ with $i \in \{a_1,b_1\}$. Thus, $z_{v_1}$ strongly commutes with $f^\omega_{[\sigma],\Delta'}$. For the latter, note from equation~\eqref{eq:ohs_proof_sum_formula_for_two_loops} that $f^\omega_{[\eta],\Delta'}$ involves only $z_{v_i}$ for the ideal arcs $i$ appearing in equation~\eqref{eq:ohs_proof_v_J_r_prime_involves}, so it suffices to show that $\langle v_i, v_2\rangle=0$ holds for all $i$ in equation~\eqref{eq:ohs_proof_v_J_r_prime_involves}. Recall that we showed in equation~\eqref{eq:ohs_proof_some_commutations2} that $\langle v_i, v_1\rangle=0$ holds for these $i$. Since $v_2 = v_1 +v_\gamma$, it suffices then to show
\begin{align}
\label{eq:ohs_proof_some_commutations3}
\langle v_i, v_\gamma\rangle = 0 \quad \text{for all ideal arcs $i$ appearing in equation~\eqref{eq:ohs_proof_v_J_r_prime_involves}.}
\end{align}

For this we collect the values of $\varepsilon_{ih}$ for $i$ in equation~\eqref{eq:ohs_proof_v_J_r_prime_involves} and $h$ internal ideal arcs of $\Delta'$, by inspection on Figure~\ref{fig:one-holed_surface_of_higher_genus_flat} as follows. For convenience, we set the symbol $f_g$, which was not defined, to mean $e_g$,
\begin{gather*}
\varepsilon_{a_i b_i} =-2, \qquad
\varepsilon_{a_i c_i} = 1, \qquad
\varepsilon_{a_i d_i} = 1, \\
\varepsilon_{a_i h}=0, \quad\mbox{$\forall$ other internal arcs $h\in \Delta'$, $\forall i=1,\dots,g$}, \\
 \varepsilon_{b_i a_i} =2, \qquad
\varepsilon_{b_i c_i} = -1, \qquad
\varepsilon_{b_i d_i}=-1,\\
\varepsilon_{b_i h}=0 \quad\mbox{$\forall$ other internal arcs $h\in \Delta'$, $\forall i=1,\dots,g$}, \\
 \varepsilon_{c_i a_i} = -1, \qquad
\varepsilon_{c_i b_i} = 1, \qquad
\varepsilon_{c_i d_i} = -1, \qquad
\varepsilon_{c_i e_i} = 1, \\
\varepsilon_{c_i h}=0 \quad \mbox{$\forall$ other internal arcs $h\in \Delta'$, $\forall i=2,\dots,g$}, \\
 \varepsilon_{d_1 a_1} = -1,\qquad
\varepsilon_{d_1 b_1}=1,\qquad
\varepsilon_{d_1 c_1'}=1,\qquad
\varepsilon_{d_1 f_2}=-1,\\
\varepsilon_{d_1 h}=0 \quad \mbox{$\forall$ other internal arcs $h\in \Delta'$}, \\
 \varepsilon_{d_i a_i}=-1,\qquad
\varepsilon_{d_i b_i}=1,\qquad
\varepsilon_{d_i c_i}=1,\qquad
\varepsilon_{d_i e_i}=-1,\\
\varepsilon_{d_i h}=0 \quad \mbox{$\forall$ other internal arcs $h\in \Delta'$, $\forall i=2,\dots,g$}, \\
 \varepsilon_{e_i c_i}=-1,\qquad
\varepsilon_{e_i d_i}=1,\qquad
\varepsilon_{e_i f_i}=1,\qquad
\varepsilon_{e_i f_{i+1}}=-1,\\
\varepsilon_{e_i h}=0 \quad \mbox{$\forall$ other internal arcs $h\in \Delta'$,} \
 \mbox{$\forall i=2,\dots,g-1$, if $g\ge 3$}, \\
 \varepsilon_{e_g c_g}=1,\qquad
\varepsilon_{e_g d_g}=1,\qquad
\varepsilon_{e_g e_{g-1}}=1,\qquad
\varepsilon_{e_g f_{g-1}}=-1,\\
\varepsilon_{e_g h}=0\quad \mbox{$\forall$ other internal arcs $h\in \Delta'$, if $g\ge 3$}, \\
 \varepsilon_{f_2 c_1'}=-1,\qquad
\varepsilon_{f_2 d_1}=1,\qquad
\varepsilon_{f_2 c_2}=-1,\qquad
\varepsilon_{f_2 d_2}=1,\\
\varepsilon_{f_2 h}=0 \quad \mbox{$\forall$other internal arcs $h\in \Delta'$, if $g=2$}, \\
 \varepsilon_{f_2 c_1'}=-1,\qquad
\varepsilon_{f_2 d_1}=1,\qquad
\varepsilon_{f_2 e_2}=-1,\qquad
\varepsilon_{f_2 f_3}=1,\\
\varepsilon_{f_2 h}=0 \quad \mbox{$\forall$ other internal arcs $h\in \Delta'$, if $g\ge 3$}, \\
 \varepsilon_{f_i e_{i-1}}=1,\qquad
\varepsilon_{f_i e_i}=-1,\qquad
\varepsilon_{f_i f_{i-1}}=-1,\qquad
\varepsilon_{f_i f_{i+1}}=1,\\
\varepsilon_{f_i h}=0 \quad \mbox{$\forall$ other internal arcs $h\in \Delta$} \
 \mbox{$\forall i=3,\dots,g-1$, if $g\ge 4$}.
\end{gather*}
Note that equation~\eqref{eq:ohs_proof_v_gamma} says $v_\gamma = \sum_h 2 v_h$, where the sum is over all internal ideal arcs $h$ of $\Delta'$. For each ideal arc $i$ in equation~\eqref{eq:ohs_proof_v_J_r_prime_involves}, from the above computation of $\varepsilon_{ih}$ we can easily check that $\sum_h \varepsilon_{ih}=0$, where the sum is over all internal ideal arcs $h$ of $\Delta'$. Thus, by equation~\eqref{eq:skew-form_for_Delta}, we get $\langle v_i, v_\gamma\rangle =0$, as desired in equation~\eqref{eq:ohs_proof_some_commutations3}. This shows that $z_{v_2}$ strongly commutes with~$f^\omega_{[\eta],\Delta}$, finishing the verification of the condition \ref{TR3} of Definition~\ref{def:Teschner_triple_level2}.

For the remaining part of \ref{TR2}, we must verify $\langle v_1, v_2 \rangle=4$. Indeed, note $\langle v_1, v_2 \rangle = \langle v_1, v_1 + v_\gamma \rangle = \langle v_1, v_\gamma \rangle = \langle - v_{a_1} - v_{b_1} - 2v_{c_1}, v_\gamma\rangle = -2 \langle v_{c_1}, v_\gamma\rangle$, using results we established. From equation~\eqref{eq:ohs_proof_v_gamma} and the values of $\varepsilon_{c_1 h}$ which we collected below equation~\eqref{eq:ohs_proof_some_commutations3}, one can verify $-2 \langle v_{c_1}, v_\gamma\rangle = -2 \cdot (-2) = 4$, as desired.

\subsubsection{Proof of Proposition~\ref{prop:hole-surrounding_loop_in_higher_genus_one-holed_subsurface}: Verification of \ref{TR5}}
\label{subsec:hole-surrounding_loop_in_one-holed_subsurface_deeper}

Lastly, for \ref{TR5} or \ref{TR5_prime} of Definition~\ref{def:Teschner_triple_level2}, we consider a loop $\zeta$ as drawn in Figure~\ref{fig:one-holed_surface_of_higher_genus}. A~precise construction of $\zeta$ can be described using the 4$g$-gon picture in Figure~\ref{fig:one-holed_surface_of_higher_genus_flat}; we let $\zeta$ be the oriented loop obtained by concatenating the oriented segments in Figure~\ref{fig:one-holed_surface_of_higher_genus_flat} labeled by $\circledd{$12'$}$, $\circledd{13}$, $\circledd{14}$, \dots, $\circleddd{$N\hspace{-1mm}-\hspace{-1mm}1$}$, in this order; then one can isotope $\zeta$ so that it is disjoint from $\eta$. Note that all these segments except for $\circledd{$12'$}$ are segments of the loop $\eta$. We leave to readers as an exercise to verify that the loop $\zeta$ defined by using the loop segments from Figure~\ref{fig:one-holed_surface_of_higher_genus_flat} indeed looks like Figure~\ref{fig:one-holed_surface_of_higher_genus}, in particular that $(\eta,\zeta,\sigma)$ forms a pants triple. For \ref{TR5_prime}, we shall now prove the following.
\begin{Proposition}
\label{prop:ohs_going_one_level_deeper}
Let $\mathfrak{S}' = \Sigma' \setminus \mathcal{P}'$, $\Delta'$, $\eta$, $\sigma$ and $\zeta$ be as above. Then the triple of loops $(\eta,\zeta,\sigma)$ in $\mathfrak{S}'$ is a strong Teschner triple with respect to $\Delta'$.
\end{Proposition}

\begin{proof} We use the same notations as before. We should study the quantized trace-of-mo\-no\-dromy for the loops $\eta$, $\zeta$ and $\sigma$. We already know from equation~\eqref{eq:ohs_proof_f_omega_sigma} that
 \[f^\omega_{[\sigma],\Delta'} = {\rm e}^{z_{v_{a_1}} + z_{v_{b_1}}} + {\rm e}^{-z_{v_{a_1}} + z_{v_{b_1}}} + {\rm e}^{-z_{v_{a_1}}-z_{v_{b_1}}}.\]
 From equation~\eqref{eq:ohs_proof_sum_formula_for_two_loops}, we know $f^\omega_{[\eta],\Delta'} = \sum_{r=6}^{N} Z^{J_r'}$, where the admissible juncture-states~$J_r'$ for~$\eta$ are as shown in Table~\ref{tab:hole-surrounding_loop_in_one-holed_higher_genus_surface}. We first re-organize the table of all admissible juncture-states for~$\eta$ in Table~\ref{tab:more_loops_in_one-holed_higher_genus_surface}; in this new table, the leftmost juncture coincides with the rightmost juncture, so that the juncture-state values for the former is the same as those for the latter. To indicate this redundancy, we put the values for the former in parentheses.

Since the oriented loop segments $\circledd{$12'$}$, $\circled{13}$, \dots, $\circleddd{$N\hspace{-1mm}-\hspace{-1mm}1$}$ of $\zeta$ are all right turns except for a single left turn $\circledd{$12'$}$, the loop $\zeta$ is almost-peripheral with respect to $\Delta'$ (Definition~\ref{def:AP}), and hence Lemma~\ref{lem:AP} applies. This justifies the list of all admissible juncture-states of $\zeta$ presented in Table \ref{tab:more_loops_in_one-holed_higher_genus_surface} as $J_{11}'', \dots, J_{N-1}''$. Notice that these juncture-states of $\zeta$ are labeled this way because they can be obtained from the juncture-states $J_{11}', \dots, J_{N-1}'$ of $\eta$ by removing the six sign values for the loop segments $\circled{$7'$}$, $\circled{8}, \dots, \circled{11}$ bounded by the vertical bars in the table. Similarly, we see from the table that from the juncture-states $J_7'$, $J_8'$, $J_9'$ of $\eta$, the two sign values for the loop segment $\circled{9}$ comprise the set of all admissible juncture-states for the loop $\sigma$.

\begin{table}[htb]
\centering
{\setlength{\tabcolsep}{0,2mm}
\begin{tabular}{ r @{\hspace{1mm} } lllllllll@{\hspace{2mm}}lr @{\hspace{3mm}}l}
$\eta$ : & \tablearrow{$N$}{R} & \tablearrow{$7'$}{L} & \tablearrow{8}{R} & \tablearrow{9}{R} & \tablearrowd{10}{R} & \tablearrowd{11}{R} & \tablearrowd{12}{R} & \tablearrowd{13}{R} & $\hspace{0,5mm} \bullet \hspace{-1,3mm} \overset{\cdots}{\longrightarrow} \hspace{-1,3mm}$ & $\cdots$ & \tablearrowdd{$N\hspace{-1mm}-\hspace{-1mm}1$}{R} $\bullet$ \\
arcs for $\eta$ : & \hspace{-3mm} ($f_3$) & $f_2$ & $d_1$ & $a_1$ & $b_1$ & $d_1$ & $f_2$ & $e_2$ & $c_2$ & & $*$ \hspace{1,7mm} $f_3$ \\
$\zeta$ : & \multicolumn{7}{l}{$\hspace{-0,5mm}\xymatrix@C+41,8mm{\bullet \hspace{-0,8mm} \ar[r]^{\overset{\circledd{$12'$}}{\rm L}} & } \hspace{-3,0mm}$} & \tablearrowd{13}{R} & $\hspace{0,5mm} \bullet \hspace{-1,3mm} \overset{\cdots}{\longrightarrow} \hspace{-1,3mm}$ & $\cdots$ & \tablearrowdd{$N\hspace{-1mm}-\hspace{-1mm}1$}{R} $\bullet$ \\
arcs for $\zeta$ : & \multicolumn{7}{l}{\hspace{-2mm}($f_3$)} & $e_2$ & $c_2$ & & $*$ \hspace{1,7mm} $f_3$ \\ \hline
\multicolumn{7}{l}{all admissible juncture-states for $\eta$:} \\
$J_6'$ : & \hspace{-2,5mm} (+) & + & \red{+} & + & + & + & + & + & + & $\cdots$ & \hspace{0,5mm} + \hspace{2mm} + \hspace{-2mm} & \mbox{all plus} \\ \hline
$J_7'$ : & \hspace{-2,5mm} (+) & + & \red{$-$} \vbar & \red{$+$} & + \vbar & + & + & + & + & $\cdots$ & + \hspace{2mm} + \hspace{-2mm} & \rdelim\}{3}{1pt}[corr. to $\sigma$] \\
$J_8'$ : & \hspace{-2,5mm} (+) & + & \red{$-$} \vbar & \red{$-$} & \red{+} \vbar & + & + & + & + & $\cdots$ & + \hspace{2,0mm} + \hspace{-2mm} & \\
$J_9'$ : & \hspace{-2,5mm} (+) & + & \red{$-$} \vbar & \red{$-$} & \red{$-$} \vbar & \red{+} & + & + & + & $\cdots$ & + \hspace{2,0mm} + \hspace{-2mm} & \\ \hline
$J_{10}'$ : & \hspace{-2,5mm} (+) & + & \red{$-$} & \red{$-$} & \red{$-$} & \red{$-$} & \red{+} & + & + & $\cdots$ & + \hspace{2,0mm} + \hspace{-2mm} & connecting term \\ \hline
$J_{11}'$ : & \hspace{-2,5mm} (+) \hspace{-2,5mm} \vbar & + & \red{$-$} & \red{$-$} & \red{$-$} & \red{$-$} & \red{$-$} \vbar & \red{$+$} & + & $\cdots$ & + \hspace{2,0mm} + \hspace{-2mm} & \rdelim\}{5}{1pt}[corr. to $\zeta$] \\
$J_{12}'$ : & \hspace{-2,5mm} (+) \hspace{-2,5mm} \vbar & + & \red{$-$} & \red{$-$} & \red{$-$} & \red{$-$} & \red{$-$} \vbar & \red{$-$} & \red{+} & $\cdots$ & + \hspace{2,0mm} + \hspace{-2mm} & \\
$\vdots$ \hspace{2mm} & & & & & & & & & & $\cdots$& & \\
$J_{N-2}'$ : & \hspace{-2,5mm} \red{(+)} \hspace{-2,5mm} \vbar & + & \red{$-$} & \red{$-$} & \red{$-$} & \red{$-$} & \red{$-$} \vbar & \red{$-$} & \red{$-$} & $\cdots$ & \hspace{0,5mm} \red{$-$} \hspace{2,0mm} \red{+} \hspace{-2mm} & \\
$J_{N-1}'$ : & \hspace{-2,5mm} \red{($-$)} \hspace{-2,5mm} \vbar & \red{+} & \red{$-$} & \red{$-$} & \red{$-$} & \red{$-$} & \red{$-$} \vbar & \red{$-$} & \red{$-$} & $\cdots$ & \hspace{0,5mm} \red{$-$} \hspace{2,0mm} \red{$-$} \hspace{-2mm} & \\ \hline
$J_{N}'$ : & \hspace{-2,5mm} \red{($-$)} & \red{$-$} & \red{$-$} & \red{$-$} & \red{$-$} & \red{$-$} & \red{$-$} & \red{$-$} & \red{$-$} & $\cdots$ & \red{$-$} \hspace{2mm} \red{$-$} \hspace{-2mm} & \mbox{all minus} \\ \hline
\multicolumn{7}{l}{all admissible juncture-states for $\zeta$:} \\
$J_{11}''$ : & \hspace{-2,5mm} (+) & & & & & & & + & + & $\cdots$ & + \hspace{2mm} + \hspace{-2mm} & \\
$J_{12}''$ : & \hspace{-2,5mm} (+) & & & & & & & + & $-$ & $\cdots$ & + \hspace{2mm} + \hspace{-2mm} & \\
 & & & & & & & & & & $\cdots$ & & \\
 $J_{N-2}''$ : & \hspace{-2,5mm} (+) & & & & & & & $-$ & $-$ & $\cdots$ & $-$ \hspace{2mm} + \hspace{-2mm} & \\
 $J_{N-1}''$ : & \hspace{-2,5mm} ($-$) & & & & & & & $-$ & $-$ & $\cdots$ & $-$ \hspace{2mm} $-$ \hspace{-2mm} & \\ \hline
\end{tabular}}
\caption{All admissible juncture-states for $\eta$ and $\zeta$ in one-holed surface $(\Sigma',\mathcal{P}')$ of higher genus $g\ge 2$, with respect to the ideal triangulation $\Delta'$.}\label{tab:more_loops_in_one-holed_higher_genus_surface}
\end{table}

By Lemma~\ref{lem:AP}\,\ref{AP3}, the non-redundancy condition in equation~\eqref{eq:admissible_non-redundancy_condition} is satisfied for $\zeta$, so that Proposition~\ref{prop:BW_admissible_sum_formula} applies, yielding the state-sum formula
\begin{align}\label{eq:ohs_proof2_sum_formula_for_zeta}
f^\omega_{[\zeta],\Delta'} = \sum_{r=11}^{N-1} Z^{J_r''}.
\end{align}

For $J_{11}'$, \dots, $J_{N-1}'$ in Table~\ref{tab:more_loops_in_one-holed_higher_genus_surface}, define $v_1' \in N_{\Delta'}$ as the sum of elements~$v_i$ with the ideal arcs $i$ running through the ideal arcs (with multiplicity) on which the junctures of the loop segments $\circled{$7'$}$, $\circled{8}$, \dots, $\circled{11}$ of $\eta$ lie in; that is, read the sign entries in between the vertical bars in Table~\ref{tab:more_loops_in_one-holed_higher_genus_surface}. For $J_7'$, $J_8'$, $J_9'$, define $v_2' \in N_{\Delta'}$ as the sum of elements $v_i$ with $i$ running through the ideal arcs on which the junctures of the loop segments of $\eta$ lie in, {\it except for} the segment $\circled{9}$; that is, read the sign entries outside the vertical bars in Table~\ref{tab:more_loops_in_one-holed_higher_genus_surface}:
\begin{gather*}
v_1' := +v_{f_2} - v_{d_1} - v_{a_1} - v_{b_1} - v_{d_1} - v_{f_2} = - v_{a_1} - v_{b_1} - 2 v_{d_1}, \\
v_2' := v_\eta - v_{a_1} - v_{b_1} - 2 v_{d_1}.
\end{gather*}
In particular, we have $v_1' - v_2' = - v_\eta$, as required by part of \ref{TR2} of Definition~\ref{def:Teschner_triple_level1}. Now, from the table we see
\[
v_{J_r'} = v_{J_r''} + v_1', \qquad \forall r = 11, \dots, N-1.
\]
For the remaining admissible juncture-states of $\eta$, observe the following, keeping in mind the principle in equation~\eqref{eq:AP2} of Lemma~\ref{lem:AP}\,\ref{AP2}:
\begin{align*}
v_{J_6'} & = v_\eta, \\
v_{J_7'} & = v_\eta - 2 v_{d_1} = v_2' + (v_{a_1} + v_{b_1}), \\
v_{J_8'} & = v_{J_7'} - 2 v_{a_1} = v_2' + (-v_{a_1} + v_{b_1}), \\
v_{J_9'} & = v_{J_8'} - 2 v_{b_1} = v_2' + (-v_{a_1} - v_{b_1}), \\
v_{J_{10}'} & = v_{J_9'} - 2 v_{d_1} = v_2' + (-v_{a_1} - v_{b_1} - 2v_{d_1}) = v_1' + v_2', \\
v_{J_N'} & = - v_{J_6'} = - v_\eta.
\end{align*}

For $J_7'$, $J_8'$, $J_9'$, let us show
\begin{align}
\label{eq:ohs_proof2_eqn1}
\mbox{$\langle v_j', \epsilon_1 v_{a_1} + \epsilon_2 v_{b_1} \rangle=0$ for all $j\in\{1,2\}$ and $\epsilon_1,\epsilon_2 \in \{+,-\}$;}
\end{align}
it suffices to show $\langle v_j', v_{a_1}\rangle=\langle v_j',v_{b_1}\rangle=0$. For $j=1$, we can verify this directly, using $v_1' = - v_{a_1} - v_{b_1} - 2v_{d_1}$, equation~\eqref{eq:skew-form_for_Delta} and equation~\eqref{eq:ohs_proof_some_varepsilon1}-- \eqref{eq:ohs_proof_some_varepsilon2}.
For $j=2$, since $v_2' = v_\eta - v_1'$, it suffices to show $\langle v_\eta,v_{a_1}\rangle = \langle v_\eta, v_{b_1}\rangle=0$. Recall from equation~\eqref{eq:ohs_proof_v_eta2} that $v_\eta = {v_\gamma} - {v_\sigma} - 2 {v_{c_1'}} - 2 {v_{c_1}} = v_\gamma - v_{a_1} - v_{b_1} - 2v_{c_1'} - 2v_{c_1}$. Recall that we showed in the previous subsection (equation~\eqref{eq:ohs_proof_some_commutations3}) that $\langle v_\gamma,v_i\rangle=0$ holds for all ideal arcs $i$ involved in $f^\omega_{[\eta],\Delta'}$, i.e., for those $i$ appearing in equation~\eqref{eq:ohs_proof_v_J_r_prime_involves}. This includes $i=a_1$ and $b_1$, so we know $\langle v_\gamma,v_{a_1}\rangle=\langle v_\gamma,v_{b_1}\rangle=0$. So it remains to verify the following
\begin{align*}
\langle v_\eta - v_\gamma,v_{a_1} \rangle & = \langle -v_{a_1} - v_{b_1} -2v_{c_1'} - 2v_{c_1}, v_{a_1}\rangle = - \varepsilon_{a_1 a_1} - \varepsilon_{b_1 a_1} - 2\varepsilon_{c_1' a_1} - 2\varepsilon_{c_1 a_1} =0, \\
\langle v_\eta - v_\gamma,v_{b_1} \rangle & = \langle -v_{a_1} - v_{b_1} -2v_{c_1'} - 2v_{c_1}, v_{b_1}\rangle = - \varepsilon_{a_1 b_1} - \varepsilon_{b_1 b_1} - 2\varepsilon_{c_1' b_1} - 2\varepsilon_{c_1 b_1} = 0,
\end{align*}
where we used the values of $\varepsilon_{i a_1} = - \varepsilon_{a_1 i}$ and $\varepsilon_{i b_1} = - \varepsilon_{b_1 i}$ as computed in the previous subsection. So equation~\eqref{eq:ohs_proof2_eqn1} is proved. In particular,
\begin{align*}
Z^{J_7'} + Z^{J_8'} + Z^{J_9'} & = {\rm e}^{z_{v_2'} + (z_{v_{a_1}} + z_{v_{b_1}})} + {\rm e}^{z_{v_2'} + (- z_{v_{a_1}} + z_{v_{b_1}})} + {\rm e}^{z_{v_2'} + (-z_{v_{a_1}} - z_{v_{b_1}})} \\
& = {\rm e}^{z_{v_2'}}( {\rm e}^{z_{v_{a_1}} + z_{v_{b_1}}} + {\rm e}^{-z_{v_{a_1}} + z_{v_{b_1}}} + {\rm e}^{-z_{v_{a_1}} -z_{v_{b_1}}}) \quad (\because \mbox{Proposition~\ref{prop:BCH}}) \\
& = {\rm e}^{z_{v_2'}} f^\omega_{[\sigma],\Delta'} \quad (\because \mbox{equation~\eqref{eq:ohs_proof_f_omega_sigma}}),
\end{align*}
and $z_{v_2'}$ strongly commutes with $f^\omega_{[\sigma],\Delta'}$ in the sense of Definition~\ref{def:algebraic_SC1}. We also proved that $z_{v_1'}$ strongly commutes with $f^\omega_{[\sigma],\Delta'}$.

For $J_{11}',\dots,J_{N-1}'$, let us show
\begin{align}
\label{eq:ohs_proof2_eqn2}
\langle v_j', v_{J_r''}\rangle=0 \qquad \mbox{for all $j\in \{1,2\}$ and $r=11,12,\dots,N-1$}.
\end{align}
By looking at the ideal arcs of $\Delta'$ intersecting $\zeta$ in Figure~\ref{fig:one-holed_surface_of_higher_genus_flat}, one finds that for each $r\in\{6,\dots,N\}$,
\begin{align}
\label{eq:ohs_proof2_v_J_r_double_prime_involves}
\mbox{$v_{J_r''}$ involves only $a_2,\dots,a_g$, $b_2,\dots,b_g$, $c_2,\dots,c_g$, $d_2,\dots,d_g$, $e_2,\dots,e_g$, $f_3,\dots,f_{g-1}$.}\!\!\!
\end{align}
In particular, the set of ideal arcs $i$ in equation~\eqref{eq:ohs_proof2_v_J_r_double_prime_involves} is a subset of the set of ideal arcs involved in $v_{J_r'}$ for the loop $\eta$ as written in equation~\eqref{eq:ohs_proof_v_J_r_prime_involves}; what are missing in \eqref{eq:ohs_proof2_v_J_r_double_prime_involves} are $a_1$, $b_1$, $d_1$, $f_2$. It then suffices to show
\begin{align*}
\langle v_j', v_i\rangle=0 \qquad \mbox{for all $j\in \{1,2\}$ and $i$ appearing in equation~\eqref{eq:ohs_proof2_v_J_r_double_prime_involves}}.
\end{align*}
From equation~\eqref{eq:ohs_proof_some_commutations2} in the previous subsection, we showed $\langle v_1, v_i \rangle = \langle - v_{a_1} - v_{b_1} - 2v_{c_1}, v_i\rangle =0$ for these $i$. Thus, for $j=1$, in order to show $\langle v_1', v_i \rangle = \langle - v_{a_1} - v_{b_1} - 2v_{d_1}, v_i \rangle =0$ it suffices to show $\langle v_{c_1}, v_i \rangle = \langle v_{d_1}, v_i\rangle$, i.e., $\varepsilon_{c_1 i} = \varepsilon_{d_1 i}$. Recall that we have computed in the previous subsection the values of $\varepsilon_{ih}$, for all $i$ in equation~\eqref{eq:ohs_proof_v_J_r_prime_involves} hence for all $i$ in equation~\eqref{eq:ohs_proof2_v_J_r_double_prime_involves} and for all internal ideal arcs $h$ of $\Delta'$. From these values, we indeed see that $\varepsilon_{i c_1} = \varepsilon_{i d_1}$ holds, as desired. For $j=2$, note $v_2' = v_\eta - v_1'$, so it suffices to show $\langle v_\eta, v_i \rangle=0$ for all $i$ in equation~\eqref{eq:ohs_proof2_v_J_r_double_prime_involves}. Note $v_\eta = v_\gamma - v_{a_1} - v_{b_1} - 2v_{c_1} - 2v_{c_1'} = v_\gamma - v_1 - 2v_{c_1'}$. From equation~\eqref{eq:ohs_proof_some_commutations2} and equation~\eqref{eq:ohs_proof_some_commutations3} of the previous subsection, we showed $\langle v_1, v_i \rangle=0$ and $\langle v_\gamma, v_i \rangle=0$ hold for all $i$ appearing in equation~\eqref{eq:ohs_proof_v_J_r_prime_involves}; hence they hold for all $i$ in equation~\eqref{eq:ohs_proof2_v_J_r_double_prime_involves}. Thus, what remains to show is $\langle v_{c_1'}, v_i \rangle =0$, i.e., $\varepsilon_{c_1' i}=0$ for all $i$ in equation~\eqref{eq:ohs_proof2_v_J_r_double_prime_involves}; this can be verified from the values of $\varepsilon_{ih}$ computed in the previous subsection. So equation~\eqref{eq:ohs_proof2_eqn2} is proved, and therefore
\begin{align*}
\sum_{r=11}^{N-1} Z^{J_r'} = \sum_{r=11}^{N-1} {\rm e}^{z_{v_{J_r'}}}
= \sum_{r=11}^{N-1} {\rm e}^{z_{v_1'} + z_{v_{J_r''}}}
 \stackrel{{\rm Proposition}~\ref{prop:BCH}}{=}
\sum_{r=11}^{N-1} {\rm e}^{z_{v_1'}} {\rm e}^{z_{v_{J_r''}}}
 \stackrel{{\rm equation}~\eqref{eq:ohs_proof2_sum_formula_for_zeta}}{=}
{\rm e}^{z_{v_1'}} f^\omega_{[\zeta],\Delta'}.
\end{align*}
We showed that each of $z_{v_1'}$ and $z_{v_2'}$ strongly commutes with $f^\omega_{[\zeta],\Delta'}$ (Definition~\ref{def:algebraic_SC1}). This finishes the verification of \ref{TR3} of Definition~\ref{def:Teschner_triple_level1}.

Now observe that
\begin{align*}
f^\omega_{[\eta],\Delta'} \stackrel{{\rm equation} \eqref{eq:ohs_proof_sum_formula_for_two_loops}}{=}
\sum_{r=6}^{N} Z^{J_r'}
& = Z^{J_6'} + Z^{J_N'} +Z^{J_{10}'} + \sum_{r=11}^{N-1} Z^{J_r'} + \bigl(Z^{J_7'} + Z^{J_8'} + Z^{J_9'}\bigr)
 \\
& = {\rm e}^{z_{v_\eta}} + {\rm e}^{-z_{v_\eta}} + {\rm e}^{z_{v_1'} + z_{v_2'}} + {\rm e}^{z_{v_1'}} f^\omega_{[\zeta],\Delta'} + {\rm e}^{z_{v_2'}} f^\omega_{[\sigma],\Delta'},
\end{align*}
proving \ref{TR1} of Definition~\ref{def:Teschner_triple_level1}, i.e., the Teschner recursion formula holds.

For \ref{TR4}, observe from Figure~\ref{fig:one-holed_surface_of_higher_genus_flat} that the loop $\sigma$ intersects exactly two ideal triangles of $\Delta'$, and that $\zeta$ does not meet any of these two triangles. Hence, $\sigma$ and $\zeta$ are triangle-disjoint with respect to $\Delta'$ (Definition~\ref{def:triangle-disjoint}), yielding \ref{TR4} in view of Lemma~\ref{lem:triangle-disjoint_implies_strong_commutativity}.

For the remaining part of \ref{TR2}, we need to verify $\langle v_1' , v_2'\rangle=4$. Using the results we collected so far, note $\langle v_1' , v_2'\rangle = \langle v_1', v_\eta + v_1'\rangle = \langle v_1', v_\eta\rangle = \langle v_1', v_\gamma - v_{a_1} - v_{b_1} - 2v_{c_1'} - 2v_{c_1}\rangle = \langle v_1', v_\gamma\rangle + \langle v_1', -2 v_{c_1'} -2 v_{c_1}\rangle$. Note $\langle v_1', v_\gamma \rangle = \langle -v_{a_1}-v_{b_1}-2v_{d_1},v_\gamma \rangle$, which is zero by equation~\eqref{eq:ohs_proof_some_commutations3}. Note $\langle v_1', -2v_{c_1'} - 2v_{c_1}\rangle = \langle v_{a_1}+v_{b_1}+2v_{d_1}, 2v_{c_1'} + 2v_{c_1}\rangle =2\varepsilon_{a_1 c_1'} + 2\varepsilon_{a_1 c_1} + 2\varepsilon_{b_1 c_1'} + 2 \varepsilon_{b_1 c_1} + 4 \varepsilon_{d_1 c_1'} + 4 \varepsilon_{d_1 c_1} =0+2+0+(-2)+4+0=4$, hence $\langle v_1',v_2'\rangle=4$ as desired.

We showed that $(\eta,\zeta,\sigma)$ is a strong Teschner triple with respect to $\Delta'$.
This completes the proof of Proposition~\ref{prop:ohs_going_one_level_deeper}.
\end{proof}

This finally concludes the proof of Proposition~\ref{prop:hole-surrounding_loop_in_higher_genus_one-holed_subsurface}.

\subsection{The first main theorem}

We summarize the results obtained in the present section on the algebraic structure of the quantized trace-of-monodromy, which is our first main theorem.
\begin{Theorem}[algebraic structure of quantized trace-of-monodromy of a single loop; Theorem~\ref{thm:main1}]
\label{thm:AS}
Let $\mathfrak{S} = \Sigma\setminus\mathcal{P}$ be a connected triangulable marked surface. Let $\gamma$ be an essential simple loop in $\mathfrak{S}$, which is classified by \ref{L1}--\ref{L6} in Proposition~{\rm \ref{prop:loop_classification}}.
\begin{enumerate}[label={\rm (AS\arabic*)}]\itemsep=0pt\setlength{\leftskip}{0.22cm}
\item\label{AS1} In the case of \ref{L1}, i.e., when $\gamma$ is a peripheral loop $($Definition~{\rm \ref{def:essential_and_peripheral})} around a puncture, the quantized trace-of-monodromy $f^\omega_{[\gamma],\Delta}$ $($Definition~{\rm \ref{def:quantized_trace-of-monodromy})} for any ideal triangulation $\Delta$ of~$\mathfrak{S}$ is of the form
\[
f^\omega_{[\gamma],\Delta} = Z_{v_\gamma} + Z_{-v_\gamma},
\]
where $v_\gamma \in N_\Delta$ is as defined in Definition~{\rm \ref{def:v_gamma}}.

\item\label{AS2} If $\gamma$ falls into \ref{L2}, \ref{L3}, \ref{L4} or \ref{L5}, but not into \ref{L1}, then there exists an ideal triangulation $\Delta$ of $\mathfrak{S}$ such that ${\rm int}(\gamma,\Delta)=2$ $($Definition~{\rm \ref{def:minimal_and_transversal})}, and if $a,b \in \Delta$ denote the ideal arcs of $\Delta$ that meet $\gamma$ when $\gamma$ is isotoped to a minimal position with respect to~$\Delta$, then $\varepsilon_{ab} \in \{2,-2\}$, $v_\gamma = v_a + v_b$, and
\[
f^\omega_{[\gamma],\Delta} = Z_{v_a+v_b} + Z_{\epsilon(v_a - v_b)} + Z_{-v_a-v_b},
\]
if $\epsilon = {\rm sign}(\varepsilon_{ab}) \in \{+,-\}$.

\item\label{AS3} Suppose that $\gamma$ falls into \ref{L6}, i.e., $\gamma$ is a hole surrounding loop $($Definition~{\rm \ref{def:hole-surrounding})} of a~triangle-compatible marked subsurface $\mathfrak{S}'$ $($Definition~{\rm \ref{def:subsurfaces})} of $\mathfrak{S}$ that is isomorphic to a~minimal-marked one-holed surface of genus $g\ge 1$ $($Definitions~{\rm \ref{def:surface_terminology}} and {\rm \ref{def:minimal-marked})}; then
\begin{enumerate}[label={\rm (\alph*)}]\itemsep=0pt
\item Suppose that $g=1$. If $\eta$ is any non-separating loop in $\mathfrak{S}'$, then there exists an ideal triangulation $\Delta'$ of $\mathfrak{S}'$, which extends to an ideal triangulation~$\Delta$ of $\mathfrak{S}$, such that ${\rm int}(\eta,\Delta)=2$ so that the item \ref{AS2} holds for $f^\omega_{[\eta],\Delta}$ and that $(\gamma,\eta,\eta)$ is a strong Teschner triple with respect to $\Delta$, in the sense of Definition~{\rm \ref{def:Teschner_triple_level1}}.

\item Suppose that $g\ge 2$. Then there exist non-separating essential simple loops~$\eta$ and~$\sigma$ in $\mathfrak{S}'$ and an ideal triangulation $\Delta'$ of $\mathfrak{S}'$, which extends to an ideal triangulation~$\Delta$ of $\mathfrak{S}$, such that $(\gamma,\eta,\sigma)$ is a weak Teschner triple with respect to $\Delta$, in the sense of Definition~{\rm \ref{def:Teschner_triple_level2}}. Further, there exists a separating loop $\zeta$ in $\mathfrak{S}'$ that cuts out $($Definition~{\rm \ref{def:cutting})} a one-holed subsurface of genus $g-1$ without any marked points of $\mathcal{P}$, such that $(\eta,\zeta,\sigma)$ is a strong Teschner triple with respect to the same $\Delta$; see Figure~{\rm \ref{fig:one-holed_surface_of_higher_genus}}. Here $\eta$ and $\sigma$ fall into \ref{L5}; moreover, ${\rm int}(\sigma,\Delta)=2$, so that the item~\ref{AS2} applies to \smash{$f^\omega_{[\sigma],\Delta}$}. The loop $\zeta$ falls into~\ref{L6};
$\zeta$ is a hole surrounding loop of a triangle-compatible marked subsurface $\mathfrak{S}''$ of $\mathfrak{S}'$ isomorphic to a minimal-marked one-holed surface of genus~$g-1$.
\end{enumerate}
\end{enumerate}
\end{Theorem}
\begin{proof} \ref{AS1} This is Proposition~\ref{prop:quantized_trace-of-monodromy_peripheral}.
\ref{AS2} This follows from Propositions~\ref{prop:2_to_5_have_intersection_2} and~\ref{prop:quantized_trace-of-monodromy_for_loop_with_int_2}.
\ref{AS3}\,{\rm (a)} This is Proposition~\ref{prop:loops_in_a_one-holed_torus}.
\ref{AS3}\,{\rm (b)} This follows from Propositions~\ref{prop:hole-surrounding_loop_in_higher_genus_one-holed_subsurface} and~\ref{prop:ohs_going_one_level_deeper}. The fact that ${\rm int}(\sigma,\Delta)=2$ follows from their proofs. The fact that $\eta$ and $\sigma$ fall into \ref{L5} can be easily seen from Figure~\ref{fig:one-holed_surface_of_higher_genus}. For the last assertion about $\zeta$, note that the union of all ideal triangles of $\Delta$ intersecting $\zeta$ forms the sought-for subsurface $\mathfrak{S}''$; see Figure~\ref{fig:one-holed_surface_of_higher_genus_flat}.
\end{proof}

Of course, one can remove the connectedness assumption for the surface $\mathfrak{S}$. We will be using this structural result in the sequel \cite{sequel} to investigate the analytic properties of the operator representing $f^\omega_{[\gamma],\Delta}$ via suitable representations; see Section~\ref{sec:consequences} for a preview.

One useful corollary which follows from the proof arguments in the present section is the following.
\begin{Corollary}\label{cor:AP}
Let $\mathfrak{S}= \Sigma\setminus\mathcal{P}$ be a triangulable marked surface. Let $\gamma$ be an essential non-peripheral simple loop in $\mathfrak{S}$. Then there exists an ideal triangulation $\Delta$ of $\mathfrak{S}$ such that $\gamma$ is almost-peripheral with respect to $\Delta$, in the sense of Definition~{\rm \ref{def:AP}}.
\end{Corollary}
Thus, Lemma~\ref{lem:AP} applies, hence Proposition~\ref{prop:term-by-term_Weyl-ordered_f_omega_gamma_Delta} holds, as promised.

\section[On commutativity of quantized trace-of-monodromy for disjoint loops]{On commutativity of quantized trace-of-monodromy \\ for disjoint loops}
\label{sec:on_commutativity}

\subsection{The second main theorem}

The goal of the present section is to show the following result on a certain kind of algebraic commutativity of quantized trace-of-monodromy for two disjoint essential simple loops, which is our second main theorem.
\begin{Theorem}[algebraic commutativity of quantized trace-of-monodromy of disjoint loops; Theorem~\ref{thm:intro_second_main}]
\label{thm:algebraic_commutativity}
Let $\xi_1$ and $\xi_2$ be disjoint non-homotopic essential simple loops in a triangulable marked surface $\mathfrak{S} = \Sigma\setminus \mathcal{P}$. Consider the corresponding quantized trace-of-monodromy \smash{$f^\omega_{[\xi_i],\Delta} \in \mathcal{Z}^\omega_\Delta$} $($Definition~{\rm \ref{def:quantized_trace-of-monodromy})} with respect to an ideal triangulation $\Delta$ of $\mathfrak{S}$, for $i=1,2$. The following hold:
\begin{enumerate}[label={\rm (AC\arabic*)}]\itemsep=-0.1pt\setlength{\leftskip}{0.25cm}
\item\label{AC1} If one of $\xi_1$ and $\xi_2$ is peripheral $($Definition~{\rm \ref{def:essential_and_peripheral})}, then $f^\omega_{[\xi_1],\Delta}$ and $f^\omega_{[\xi_2],\Delta}$ strongly commute in the sense of Definition~{\rm \ref{def:algebraic_SC2}}, for any ideal triangulation $\Delta$ of $\mathfrak{S}$.

\item\label{AC2} If there exists an ideal triangulation $\Delta$ of $\mathfrak{S}$ such that $\xi_1$ and $\xi_2$ are triangle-disjoint with respect to $\Delta$, in the sense of Definition~{\rm \ref{def:triangle-disjoint}}, then $f^\omega_{[\xi_1],\Delta}$ and $f^\omega_{[\xi_2],\Delta}$ strongly commute in the sense of Definition~{\rm \ref{def:algebraic_SC2}}.

\item\label{AC3} If none of the above holds, then there exists an ideal triangulation $\Delta$ of $\mathfrak{S}$ such that $f^\omega_{[\xi_1],\Delta}$ and $f^\omega_{[\xi_2],\Delta}$ strongly or weakly Teschner-commute in the sense of Definitions~{\rm \ref{def:Teschner-commute}--\ref{def:Teschner_triple_level2}}.
\end{enumerate}
\end{Theorem}

\subsection{Easy part of the proof}\label{subsec:easy_part}

We begin with the item \ref{AC1} of Theorem~\ref{thm:algebraic_commutativity}.

Suppose that $\xi_1$ is a peripheral loop. Let $N_1 = \mathbb{Z} v_{\xi_1}$, where $v_{\xi_1}$ is as defined in Definition~\ref{def:v_gamma}. It is well known in the literature that $\langle v_{\xi_1}, v\rangle=0$ for all $v\in N_\Delta$, so $(N_1)_\mathbb{R} \subset N_\mathbb{R}^\perp$ (equation~\eqref{eq:N_R_perp}); see \cite[Lemma~8]{BL07} or \cite[Lemma~3.9]{AK}. Isotope $\xi_2$ to a minimal and transversal position with respect to $\Delta$, and let $N_2 := N_{\xi_2} \subset N_\Delta$ as defined in Proposition~\ref{prop:weak_sum}. We claim that these $N_1$ and~$N_2$ satisfy the conditions \ref{SC1}--\ref{SC3} of Definition~\ref{def:algebraic_SC2}. First, by \ref{AS1} of Theorem~\ref{thm:AS} we have \smash{$f^\omega_{[\xi_1],\Delta} = Z_{v_{\xi_1}} + Z_{-v_{\xi_1}}$}; here the subscripts $v_{\xi_1}$ and $-v_{\xi_1}$ belong to $N_1$. Meanwhile, by Proposition~\ref{prop:weak_sum}, we have \smash{$f^\omega_{[\xi_2],\Delta} = \sum_{v\in N_2} a_v Z_v$} for some $a_v\in \mathbb{Z}\bigl[\omega^{\pm 1}\bigr]$. Thus, \ref{SC3} is satisfied. We already saw that \ref{SC1} would hold. Note $(N_1)_\mathbb{R} \cap (N_2)_\mathbb{R} \subset (N_1)_\mathbb{R} \subset N_\mathbb{R}^\perp$, hence \ref{SC2} is satisfied. Thus, \smash{$f^\omega_{[\xi_j],\Delta}$}, $j=1,2$, strongly commute with each other in the sense of Definition~\ref{def:algebraic_SC2}. This proves \ref{AC1}, also finishing the proof of Lemma~\ref{lem:triangle-disjoint_implies_strong_commutativity}.

Item \ref{AC2} is proved in Lemma~\ref{lem:triangle-disjoint_implies_strong_commutativity}.

What remains is the item \ref{AC3} of Theorem~\ref{thm:algebraic_commutativity}, for which none of the loops $\xi_1$ and $\xi_2$ is peripheral, and $\xi_1$ and $\xi_2$ are not triangle-disjoint with respect to any ideal triangulation. Such pairs of loops $(\xi_1,\xi_2)$ are classified in Proposition~\ref{prop:classification_of_pair_of_loops} of Section~\ref{subsec:pais_of_simple_loops}, by items \ref{PL2}--\ref{PL7}. Some of them are already dealt with in Section~\ref{sec:algebraic_structure_of_quantized_trace-of-monodromy}, as we observe as follows.

Case {\ref{PL2}}: There is a triangle-compatible minimal-marked one-holed torus subsurface $\mathfrak{S}'$ of $\mathfrak{S}$ such that $\xi_1$ is a hole surrounding loop in $\mathfrak{S}'$ and $\xi_2$ is a non-separating loop in $\mathfrak{S}'$.

This is dealt with by Proposition~\ref{prop:loops_in_a_one-holed_torus} of Section~\ref{subsec:loops_in_one-holed_torus}, where $\xi_1$ and $\xi_2$ are $\gamma$ and $\xi$ of Proposition~\ref{prop:loops_in_a_one-holed_torus}. So we have shown that $(\xi_1,\xi_2,\xi_2)$ is a strong Teschner triple with respect to some ideal triangulation $\Delta$, hence \smash{$f^\omega_{[\xi_1],\Delta}$} and \smash{$f^\omega_{[\xi_2],\Delta}$} strongly Teschner-commute.

Case {\ref{PL3}}: There is a triangle-compatible minimal-marked one-holed marked subsurface $\mathfrak{S}'$ of $\mathfrak{S}$ of genus $g\ge 2$, where $\xi_1$ is a hole surrounding loop in $\mathfrak{S}'$ and $\xi_2$ is a non-separating loop in~$\mathfrak{S}'$.\looseness=-1

This is dealt with by Proposition~\ref{prop:hole-surrounding_loop_in_higher_genus_one-holed_subsurface} of Section~\ref{subsec:hole-surrounding_loop_in_one-holed_subsurface}, where $\xi_1$ and $\xi_2$ are $\gamma$ and $\sigma$ of Proposition~\ref{prop:hole-surrounding_loop_in_higher_genus_one-holed_subsurface}; see Figure~\ref{fig:one-holed_surface_of_higher_genus}. Indeed, one can see that there is a self diffeomorphism of $\mathfrak{S}'$ sending~$\xi_1$ and~$\xi_2$ to~$\gamma$ and~$\sigma$, respectively ($\xi_1$ and $\gamma$ are already equivalent, while cutting along $\xi_2$ and cutting along $\sigma$ yield diffeomorphic surfaces). In Proposition~\ref{prop:hole-surrounding_loop_in_higher_genus_one-holed_subsurface}, we have shown that there exist an essential simple loop $\xi$ in $\mathfrak{S}$ and an ideal triangulation $\Delta$ of $\mathfrak{S}$ such that $(\xi_1,\xi,\xi_2)$ is a weak Teschner triple with respect to $\Delta$. So $f^\omega_{[\xi_1],\Delta}$ and $f^\omega_{[\xi_2],\Delta}$ weakly Teschner-commute. We remark that this is the only case where we have the weak Teschner-commutativity, instead of strong Teschner-commutativity, in Theorem~\ref{thm:algebraic_commutativity}\,\ref{AC3}.

{\samepage Case {\ref{PL4}}: There is a triangle-compatible minimal-marked one-holed marked subsurface $\mathfrak{S}'=\Sigma'\setminus\mathcal{P}'$ of $\mathfrak{S}$ of genus $g\ge 2$, where $\xi_1$ and $\xi_2$ are non-separating loop in $\mathfrak{S}'$, which together with the boundary circle of $\Sigma'$ bounds a three-holed sphere subsurface.

}

This is dealt with by Proposition~\ref{prop:ohs_going_one_level_deeper} of Section~\ref{subsec:hole-surrounding_loop_in_one-holed_subsurface} which is proved in Section~\ref{subsec:hole-surrounding_loop_in_one-holed_subsurface_deeper}, where~$\xi_1$ and~$\xi_2$ are this time $\eta$ and $\sigma$ of Proposition~\ref{prop:hole-surrounding_loop_in_higher_genus_one-holed_subsurface}; see Figure~\ref{fig:one-holed_surface_of_higher_genus}. There we have shown that there exist a loop $\zeta$ in $\mathfrak{S}' \subset \mathfrak{S}$ and an ideal triangulation $\Delta$ of $\mathfrak{S}$ such that $(\eta,\zeta,\sigma) = (\xi_1,\zeta,\xi_2)$ is a~strong Teschner triple with respect to $\Delta$. So \smash{$f^\omega_{[\xi_1],\Delta}$} and \smash{$f^\omega_{[\xi_2],\Delta}$} strongly Teschner-commute.

The remaining cases \ref{PL5}--\ref{PL7} of Proposition~\ref{prop:classification_of_pair_of_loops} will be studied in the remaining parts of the present section. In each case, the loops $\xi_1$ and $\xi_2$ are contained in a triangle-compatible marked subsurface $\mathfrak{S}'$ of $\mathfrak{S}$. So, when studying quantized trace-of-monodromy of $\xi_1$ and $\xi_2$, it suffices to focus on $\mathfrak{S}'$ and forget about $\mathfrak{S}$.

\subsection[Case {\ref{PL5}}: Non-separating loops in once-punctured surface of genus g>=2]{Case {\ref{PL5}}: Non-separating loops in once-punctured surface\\ of genus $\boldsymbol{g\ge 2}$}
\label{subsec:Case5_commutativity}

As in the item \ref{PL5} of Proposition~\ref{prop:classification_of_pair_of_loops}, let $\mathfrak{S}$ be a once-punctured surface of genus $g\ge 2$ without boundary, where $\xi_1$ and $\xi_2$ are non-separating loops such that $\{\xi_1,\xi_2\}$ cuts out a two-holed genus $0$ surface with one puncture. (In this case, $\mathfrak{S}' = \mathfrak{S}$.)

To deal with a once-punctured genus $g\ge 2$ surface, we make use of Figures~\ref{fig:one-holed_surface_of_higher_genus}--\ref{fig:one-holed_surface_of_higher_genus_flat} of Section~\ref{subsec:hole-surrounding_loop_in_one-holed_subsurface}. Namely, we fill in the hole there, i.e., glue a monogon to the hole to remove the hole; then the marked point on the hole now becomes a puncture. Then, from Figure~\ref{fig:one-holed_surface_of_higher_genus}, we see that the loops~$\eta$ and $\sigma$ there (with the hole filled in) satisfy the conditions of our loops $\xi_1$ and $\xi_2$; we also give orientations to these loops as in Figure~\ref{fig:one-holed_surface_of_higher_genus}--\ref{fig:one-holed_surface_of_higher_genus_flat}. We use the ideal triangulation $\Delta$ of $\mathfrak{S}$ obtained from the one presented in Figure~\ref{fig:one-holed_surface_of_higher_genus_flat} by deleting the hole and the ideal arc $c_1'$. We use the same notations for the remaining ideal arcs, and the labels of the loop segments of $\eta$ and $\sigma$. The loops $\eta$ and $\sigma$ in Figure~\ref{fig:one-holed_surface_of_higher_genus_flat} do not meet the ideal triangle formed by $c_1'$, $c_1$ and the hole. One can observe that for each oriented loop segment of $\eta$ and $\sigma$, whether it is a left turn or right turn stays the same, after deleting the hole and $c_1'$. Hence, the computation for the quantized trace-of-monodromy $f^\omega_{[\eta],\Delta}$ and $f^\omega_{[\sigma],\Delta}$ in our current setting is exactly same as that for $f^\omega_{[\eta],\Delta'}$ and $f^\omega_{[\sigma],\Delta'}$ which we performed in Section~\ref{subsec:hole-surrounding_loop_in_one-holed_subsurface}. The loop $\zeta$ can be defined in the same way as in Figure~\ref{fig:one-holed_surface_of_higher_genus}--\ref{fig:one-holed_surface_of_higher_genus_flat} with the hole filled in, and again, the computation of the quantized trace-of-monodromy for $\zeta$ is the same as that done in Section~\ref{subsec:hole-surrounding_loop_in_one-holed_subsurface_deeper}.

The only change is on the values of the exchange matrix $\varepsilon_{ij}$, where at least one of $i$ and $j$ is in $\{c_1,c_1'\}$. We claim that Proposition~\ref{prop:ohs_going_one_level_deeper} still holds, i.e., $(\eta,\zeta,\sigma)$ is a strong Teschner triple with respect to $\Delta$. Indeed, none of \smash{$f^\omega_{[\eta],\Delta'}$, $f^\omega_{[\zeta],\Delta'}$, $f^\omega_{[\sigma],\Delta'}$}, $v_1'$ and $v_2'$ appearing in the proof of this lemma in Section~\ref{subsec:hole-surrounding_loop_in_one-holed_subsurface_deeper} involves $c_1$ or $c_1'$. It may seem like some $c_1$ and $c_1'$ are involved in our proof in Section~\ref{subsec:hole-surrounding_loop_in_one-holed_subsurface_deeper} because we dealt with expressions like $v_\eta = v_\gamma - v_\sigma - 2v_{c_1'} - 2v_{c_1}$ (which is from equation~\eqref{eq:ohs_proof_v_eta2}) there, but that's only because we wanted to shorten our calculation and wanted to use some results obtained in Section~\ref{subsec:hole-surrounding_loop_in_one-holed_subsurface}. In practice, $v_\eta$ does not involve $c_1$ or~$c_1'$, as can be seen from equation~\eqref{eq:ohs_proof_v_eta1}. So the proof in Section~\ref{subsec:hole-surrounding_loop_in_one-holed_subsurface_deeper} can be rewritten in a way that does not involve $c_1$ or~$c_1'$, and by adapting such a proof we deduce that $(\eta,\zeta,\sigma)$ is a strong Teschner triple with respect to $\Delta$ in our current setting. Hence, $f^\omega_{[\xi_1],\Delta} = f^\omega_{[\eta],\Delta}$ and $f^\omega_{[\xi_2],\Delta} = f^\omega_{[\sigma],\Delta}$ strongly Teschner-commute, as desired.

\subsection{Case {\ref{PL6}}: Hole surrounding loops in two-holed subsurface}

As in the item \ref{PL6} of Proposition~\ref{prop:classification_of_pair_of_loops}, suppose that there is a triangle-compatible marked subsurface $\mathfrak{S}'=\Sigma'\setminus\mathcal{P}'$ of $\mathfrak{S}$ isomorphic to a minimal-marked two-holed surface of genus $g\ge 1$, where $\xi_1$ and $\xi_2$ are hole surrounding loops in $\mathfrak{S}'$ for distinct boundary components of $\Sigma'$.

Consider any ideal triangulation $\Delta'$ of $\mathfrak{S}'$. Now glue the two holes together, which yields a~new marked surface $\mathfrak{S}'' = \Sigma'' \setminus \mathcal{P}''$, which has no boundary and one puncture, and is of genus $g+1 \ge 2$. The ideal triangulation $\Delta'$ induces an ideal triangulation $\Delta''$ of $\mathfrak{S}''$. Note that the two boundary arcs $\eta_1$ and $\eta_2$ of $\Delta'$ (corresponding to $\xi_1$ and $\xi_2$, respectively) now get identified with each other and become an internal arc of $\Delta''$, say $\eta''$. Notice conversely that any ideal triangulation of $\mathfrak{S}''$ containing $\eta''$ as one of its constituent ideal arcs can be obtained this way, from an ideal triangulation of $\mathfrak{S}'$.

There is a natural bijection between $\Delta'_0 := \Delta' \setminus\{\eta_1,\eta_2\}$ and $\Delta''_0 := \Delta'' \setminus\{\eta''\}$, and the exchange matrix values $\varepsilon_{ij}'$ for $\Delta_0'$ are the same as those $\varepsilon''_{ij}$ for $\Delta''_0$. Denote by $\xi_1''$ and $\xi_2''$ the loops in $\mathfrak{S}''$ corresponding to $\xi_1$ and $\xi_2$ of $\mathfrak{S}'$. Notice that $\xi_1$ and $\xi_2$ do not intersect $\eta_1$ or $\eta_2$ in~$\mathfrak{S}'$, and hence $\xi_1''$ and $\xi_2''$ do not intersect $\eta''$ in $\mathfrak{S}''$.

Observe that $\xi_1''$ and $\xi_2''$ are non-separating loops in $\mathfrak{S}''$ such that $\{\xi_1'',\xi_2''\}$ cuts out a two-holed genus $0$ surface with one puncture, hence falls into the item \ref{PL5} of Proposition~\ref{prop:classification_of_pair_of_loops}, which was dealt with in Section~\ref{subsec:Case5_commutativity}. In particular, there exists an ideal triangulation $\widetilde{\Delta}''$ of $\mathfrak{S}''$ such that \raisebox{2pt}{\smash{$f^\omega_{[\xi''_1],\widetilde{\Delta}''}$}} and \raisebox{2pt}{\smash{$f^\omega_{[\xi''_2],\widetilde{\Delta}''}$}} strongly Teschner-commute. In our proof in Section~\ref{subsec:Case5_commutativity}, notice that~$\til{\Delta}''$, which was denoted $\Delta'$ there, contains $\eta''$ as a constituent arc, which was $c_1$ there. As none of the~three loops in the (strong Teschner) triple of loops dealt with in Section~\ref{subsec:Case5_commutativity} intersects~$\eta''$, the proof of Section~\ref{subsec:Case5_commutativity} remains valid after we cut~$\mathfrak{S}''$ along the ideal arc $\eta''$; notice that the marked surface obtained by this cutting process is isomorphic to $\mathfrak{S}'$. Through this cutting process, the ideal triangulation $\til{\Delta}''$ of $\mathfrak{S}''$ yields an ideal triangulation $\til{\Delta}'$ of $\mathfrak{S}'$, and the loops~$\xi_1''$ and~$\xi_2''$ in~$\mathfrak{S}''$ yield the loops $\xi_1$ and $\xi_2$ in $\mathfrak{S}'$. Therefore, it follows that there exists a loop $\zeta$ in $\mathfrak{S}'$ such that $(\xi_1,\zeta,\xi_2)$ is a strong Teschner triple with respect to $\til{\Delta}'$. Therefore, \smash{$f^\omega_{[\xi_1],\til{\Delta}'}$} and \smash{$f^\omega_{[\xi_2],\til{\Delta}'}$} strongly Teschner-commute, as desired.

\subsection[Case \ref{PL7}: Hole surrounding loop and separating loop in one-holed subsurface of genus g>= 2]{Case \ref{PL7}: Hole surrounding loop and separating loop\\ in one-holed subsurface of genus $\boldsymbol{g\ge 2}$}

As in the item \ref{PL7} of Proposition~\ref{prop:classification_of_pair_of_loops}, suppose that there is a triangle-compatible marked subsurface $\mathfrak{S}'=\Sigma'\setminus\mathcal{P}'$ of $\mathfrak{S}$ isomorphic to a minimal-marked one-holed surface of genus $g\ge 2$, where $\xi_1$ is a hole surrounding loop in $\mathfrak{S}'$ and $\xi_2$ is a separating loop in $\mathfrak{S'}$ that cuts out a~one-holed surface of genus $g'$ with $1\le g'<g$.

\begin{figure}[t]
\centering
\vspace{-8mm}
\scalebox{1.0}{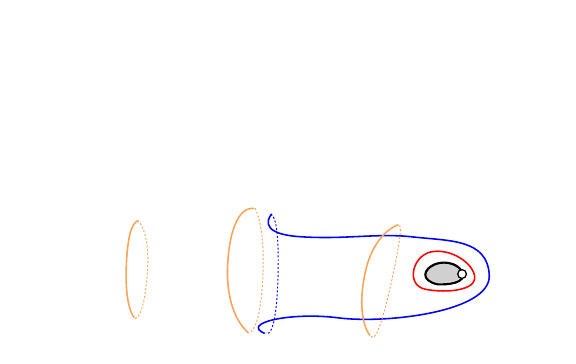}
\vspace{-3mm}

\caption{Separating loops $\zeta_i$ and $\vartheta_i$ ($2\le i \le g$) in a minimal-marked one-holed surface $(\Sigma',\mathcal{P}')$ of higher genus $g\ge 2$.}\label{fig:one-holed_surface_of_higher_genus_zeta_vartheta}
\vspace{-2mm}
\end{figure}

We consider oriented loops $\zeta_2,\zeta_3,\dots,\zeta_g$ in $\mathfrak{S}' = \Sigma'\setminus\mathcal{P}'$ as depicted in Figure~\ref{fig:one-holed_surface_of_higher_genus_zeta_vartheta}, where each $\zeta_i$ is a separating loop in $\mathfrak{S}' = \Sigma' \setminus \mathcal{P}'$ that cuts out a one-holed surface having the $i$-th, $(i+1)$-th, \dots, $g$-th handles, i.e., a one-holed surface of genus $g-i+1$ without a marked point. As $i$ ranges in $2,3,\dots,g$, the genus $g-i+1$ of the one-holed subsurface ranges in $1,2,\dots,g-1$. Note that~$\zeta_2$ coincides with $\zeta$ in Figure~\ref{fig:one-holed_surface_of_higher_genus}. Recall from Figures~\ref{fig:one-holed_surface_of_higher_genus} and~\ref{fig:one-holed_surface_of_higher_genus_flat} an oriented hole surrounding loop $\gamma$, which we also draw in Figure~\ref{fig:one-holed_surface_of_higher_genus_zeta_vartheta}. Thus, $\gamma$ and $\zeta_i$ play the roles of $\xi_1$ and $\xi_2$, respectively, when $g-i+1 = g'$; indeed, one can see that there is a self diffeomorphism of $\mathfrak{S}'$ sending $\xi_1$ and~$\xi_2$ to $\gamma$ and~$\zeta_i$. Hence, it suffices to show that the quantized trace-of-monodromy of $\gamma$ and~$\zeta_i$ strongly or weakly Teschner-commute. We will show that for each $i$, there exists a loop $\vartheta_i$ such that $(\gamma,\zeta_i,\vartheta_i)$ is a strong Teschner triple with respect to some ideal triangulation, yielding the strong Teschner-commutativity.

\begin{figure}[!ht]
\centering
\scalebox{0.85}{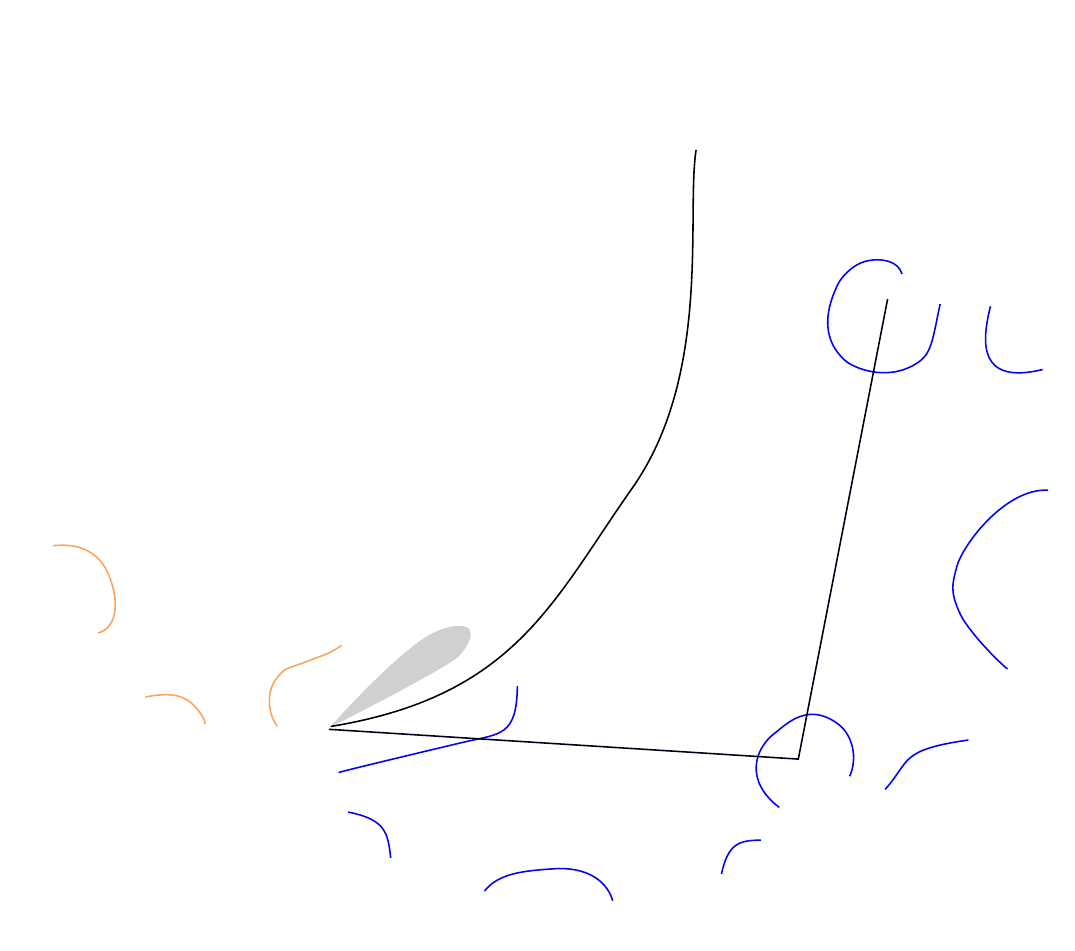}
\vspace{-2mm}
\caption{One-holed surface $(\Sigma',\mathcal{P}')$ of genus $g\ge 2$, as quotient of $4g$-gon, with an ideal triangulation $\Delta'_i$, together with loops $\gamma$, $\zeta_i$, $\vartheta_i$; this picture shows the case $g=5$. Note $k_j = i+9+11(j-2)$ for $j=2,\dots,g$.}\label{fig:one-holed_surface_of_higher_genus_flat5}
\end{figure}

Fix any $i\in \{2,\dots,g\}$. We introduce an oriented loop $\vartheta_i$ in $\mathfrak{S}'$ as depicted in Figure~\ref{fig:one-holed_surface_of_higher_genus_zeta_vartheta}, so that $\{\zeta_i,\vartheta_i,\gamma\}$ cuts out a three-holed sphere. From now on, we will consider only the three loops~$\zeta_i$,~$\vartheta_i$ and~$\gamma$. We introduce an ideal triangulation~$\Delta_i'$ of $\mathfrak{S}'$ as depicted in Figure~\ref{fig:one-holed_surface_of_higher_genus_flat5}, and give a~precise prescription of the construction of $\zeta_i$ and $\vartheta_i$ in terms of loop segments with respect to~$\Delta'_i$. Note that the ideal triangle of $\Delta'_i$ having the hole as one of its sides involves~$f_{i-1}$.
\begin{itemize}\itemsep=0pt
\item Enumerate the oriented loop segments of $\gamma$ by $\circled{1},\circled{2},\dots,\circled{$N$}$ in this order along the orientation of $\gamma$, as depicted in Figure~\ref{fig:one-holed_surface_of_higher_genus_flat5}; note $N = 12g-4$ as mentioned in Section~\ref{subsec:hole-surrounding_loop_in_one-holed_subsurface}.
\end{itemize}
The first segment $\circled{1}$ is a left turn from arc $f_{i-1}$ to arc $f_{i-1}'$.
\begin{itemize}\itemsep=0pt
\item For each $j=2,3,\dots,g$, define $k_j := i+9+11(j-2)$.
\end{itemize}
Then as can be seen in Figure~\ref{fig:one-holed_surface_of_higher_genus_flat5}, for $j=2,\dots,i-1$, the segment $\circledd{$k_j$}$ is the segment that meets~$e_j$ for the first time. For $j=i,i+1,\dots,g$, $\circleddd{$k_j\hspace{-0,7mm}+\hspace{-0,7mm}1$}$ is the segment meeting $e_j$ for the first time. This number $k_j$ depends not only on $j$ but also on $i$, but we are not incorporating $i$ into its notation by writing it for example as $k_{i,j}$, for economy of notations.
\begin{itemize}\itemsep=0pt
\item Consider a new oriented segment in a triangle formed by $f_i$, $f_{i-1}$ and $e_i$, starting from the initial juncture of $\circled{$N$}$ ending at the terminal juncture of $\circleddd{$k_i\hspace{-0,7mm}+\hspace{-0,7mm}1$}$; see Figure~\ref{fig:one-holed_surface_of_higher_genus_flat5}. Label this oriented segment by $\circleddd{$k_i\hspace{-0,7mm}+\hspace{-0,7mm}1'$}$.
\end{itemize}
In particular, the segment $\circleddd{$k_i\hspace{-0,7mm}+\hspace{-0,7mm}1'$}$ goes from the arc $f_i$ to the arc $e_i$.
\begin{itemize}\itemsep=0pt
\item Define an oriented simple loop $\zeta_i$ as the concatenation of the oriented segments $\circleddd{$k_i\hspace{-0,7mm}+1\hspace{-0,5mm}'$}$, $\circleddd{$k_i\hspace{-0,7mm}+\hspace{-0,7mm}2$}$, $\circleddd{$k_i\hspace{-0,7mm}+\hspace{-0,7mm}3$}$, \dots, $\circleddd{$N\hspace{-1mm}-\hspace{-1mm}2$}$, $\circleddd{$N\hspace{-1mm}-\hspace{-1mm}1$}$.
\end{itemize}
All these segments except for the first one are (isotopic to) loop segments of $\gamma$ of the same labels. One can see that $\zeta_i$ is a separating loop cutting out a one-holed subsurface `consisting of' the $i$-th, $(i+1)$-th, \dots, and $g$-th handles, hence exactly the loop $\zeta_i$ which we drew in Figure~\ref{fig:one-holed_surface_of_higher_genus_zeta_vartheta}.
\begin{itemize}\itemsep=0pt
\item Consider a new oriented segment in a triangle formed by $f'_{i-1}$, $f_{i-2}$ and $e_{i-1}$, starting from the initial juncture of $\circleddd{$k_{i-1}\hspace{-0,7mm}+\hspace{-0,7mm}10$}$ and ending at the terminal juncture of $\circled{2}$. Label this oriented segment as $\circledd{$2'$}$.
\end{itemize}
The segment $\circledd{$2'$}$ goes from the arc $e_{i-1}$ to the arc $f_{i-2}$.
\begin{itemize}\itemsep=0pt
\item Define an oriented simple loop $\vartheta_i$ as the concatenation of the oriented segments $\circledd{$2'$}$, $\circled{3}$, \dots, $\circleddd{$k_{i-1}\hspace{-0,7mm}+\hspace{-0,7mm}8$}$, $\circleddd{$k_{i-1}\hspace{-0,7mm}+\hspace{-0,7mm}9$}$.
\end{itemize}
All these segments except for the first one are (isotopic to) loop segments of $\gamma$ of the same labels. Then one can see that $\vartheta_i$ is a separating loop disjoint from $\zeta_i$ and cuts out a one-holed subsurface `consisting of' the $1$st, $2$nd, \dots, $(i-1)$-th handles, hence coincides with the loop $\vartheta_i$ which we drew in Figure~\ref{fig:one-holed_surface_of_higher_genus_zeta_vartheta}.

We now show the following statement which is being sought for in the present subsection.
\begin{Proposition}
\label{prop:zeta_i_and_vartheta_i}
Let $\mathfrak{S}' = \Sigma'\setminus\mathcal{P}'$, $\Delta'_i$, $\zeta_i$, $\vartheta_i$ and $\gamma$ as above, with $i \in \{2,\dots,g\}$. Then the triple of loops $(\gamma,\zeta_i,\vartheta_i)$ is a strong Teschner triple with respect to $\Delta'_i$.
\end{Proposition}

The arguments are similar to those in Sections~\ref{subsec:hole-surrounding_loop_in_one-holed_subsurface} and~\ref{subsec:hole-surrounding_loop_in_one-holed_subsurface_deeper}.

\begin{proof}[Proof of Proposition~\ref{prop:zeta_i_and_vartheta_i}] Among the oriented loop segments $\circled{1},\circled{2},\dots,\circled{$N$}$ of $\gamma$, the only left turn is $\circled{1}$ and all others are right turns. Among the segments $\circleddd{$k_i\hspace{-0,7mm}+\hspace{-0,7mm}1'$}$, $\circleddd{$k_i\hspace{-0,7mm}+\hspace{-0,7mm}2$}$, $\circleddd{$k_i\hspace{-0,7mm}+\hspace{-0,7mm}3$}$, \dots, $\circleddd{$N\hspace{-1mm}-\hspace{-1mm}2$}$, $\circleddd{$N\hspace{-1mm}-\hspace{-1mm}1$}$ of $\zeta_i$, the only left turn is $\circleddd{$k_i\hspace{-0,7mm}+\hspace{-0,7mm}1'$}$. Among the segments $\circledd{$2'$}$, $\circled{3}$, \dots, $\circleddd{$k_{i-1}\hspace{-0,7mm}+\hspace{-0,7mm}8$}$, $\circleddd{$k_{i-1}\hspace{-0,7mm}+\hspace{-0,7mm}9$}$ of $\vartheta_i$, the only left turn is $\circledd{$2'$}$. So these three loops are almost-peripheral and therefore Lemma~\ref{lem:AP} applies. Denote the junctures of $\gamma$ by $x_1, \dots,x_N$, so that the oriented loop segment $\circled{$j$}$, for $j=1,\dots,N$, starts from $x_{j-1}$ and ends at $x_j$; we set $x_0=x_N$. For each $r=0,\dots,N$, denote by $J_r : \gamma \cap \Delta_i' = \{x_1,x_2,\dots,x_N\} \to \{+,-\}$ the juncture-state of $\gamma$ given by
\begin{align*}
J_r(x_j) = \begin{cases}
- & \mbox{if $1\le j \le r$}, \\
+ & \mbox{if $r+1 \le j \le N$}.
\end{cases}
\end{align*}
Then by Lemma~\ref{lem:AP}\,\ref{AP1}, $\{J_0,\dots,J_N\}$ is the set of all admissible juncture-states of $\gamma$. Similarly, for $r=k_i,k_i+1,\dots,N-1$ denote by $J_r' \colon \zeta_i \cap \Delta_i' = \{x_{k_i+1},x_{k_i+2},\dots,x_{N-1}\} \to \{+,-\}$ the juncture-state of $\zeta_i$ given by
\begin{align*}
J'_r(x_j) = \begin{cases}
- & \mbox{if $k_i +1 \le j \le r$}, \\
+ & \mbox{if $r+1 \le j \le N-1$}.
\end{cases}
\end{align*}
Then $\{J'_{k_i}, J'_{k_i+1},\dots,J'_{N-1}\}$ is the set of all admissible juncture-states of $\zeta_i$. For $r=1$, $2$, $\dots$, $k_{i-1}+9$, denote by $J''_r \colon \vartheta_i \cap \Delta'_i = \{x_2,x_3,\dots,x_{k_{i-1}+9}\} \to \{+,-\}$ the juncture-state of $\vartheta_i$ given~by
\begin{align*}
J''_r(x_j) = \begin{cases}
- & \mbox{if $2 \le j \le r$}, \\
+ & \mbox{if $r+1 \le j \le k_{i-1}+9$}.
\end{cases}
\end{align*}
Then $\{J''_1,J''_2,\dots,J''_{k_{i-1}+9}\}$ is the set of all admissible juncture-states of $\vartheta_i$. According to Lem\-ma~\ref{lem:AP}\,\ref{AP3}, we have the state-sum formulae
\begin{align*}
f^\omega_{[\gamma],\Delta'_i} = \sum_{r=0}^{N} Z^{J_r}, \qquad
f^\omega_{[\zeta_i],\Delta'_i} = \sum_{r=k_i}^{N-1} Z^{J'_r}, \qquad
f^\omega_{[\vartheta_i],\Delta'_i} = \sum_{r=1}^{k_{i-1}+9} Z^{J''_r}.
\end{align*}
For each $j=1,2,\dots,N$, denote by $h_j \in \Delta'_i$ the ideal arc on which the juncture $x_j$ lies in. Define $v_1,v_2 \in N_{\Delta'_i}$ by
\begin{align*}
v_2 :={}& \sum_{j \notin \{2,\dots,k_{i-1}+9\}} J_r(x_j) v_{h_j} \qquad \mbox{for any $r\in\{1,2,\dots,k_{i-1}+9\}$} \\
={}& -v_{h_1} + \sum_{j=k_{i-1}+10}^{N} v_{h_j} \\
v_1 :={}& \sum_{j \notin \{k_i+1,k_i+2,\dots,N-1\}} J_r(x_j) v_{h_j} \qquad \mbox{for any $r\in \{k_i,k_i+1,\dots,N-1\}$} \\
={} & \biggl( - \sum_{j=1}^{k_i} v_{h_j} \biggr) + v_{h_N}.
\end{align*}
Then we get
\[
v_2 - v_1 = -v_{h_1} + v_{h_{k_{i-1}+10}} + \sum_{j=k_i}^{N} v_{h_j} + v_{h_{k_i}} + \sum_{j=1}^{k_{i-1}+10} v_{h_j} -v_{h_N}
= \sum_{j=1}^{N} v_{h_j} = v_\gamma,
\]
where we used $k_{i-1}+11= k_i$, $h_1 = f'_{i-1}=h_{k_{i-1}+10}$ and $h_N = f_{i-1} = h_{k_i}$. So part of \ref{TR2} of Definition~\ref{def:Teschner_triple_level1} is satisfied. Note
\begin{align*}
&v_{J_0} = \sum_{j=1}^{N} v_{h_j} = v_j, \qquad
 v_{J_N} = -\sum_{j=1}^{N} v_{h_j} = -v_j, \\
&v_{J_r} = v_2 + v_{J_r''}, \qquad \forall r=1,2,\dots,k_{i-1}+9, \\
&v_{J_r} = v_1 + v_{J_r'}, \qquad \forall r = k_i,k_i+1,\dots,N-1,
\end{align*}
and
\begin{align*}
v_{J_{k_{i-1}+10}} & = v_{J_{k_i-1}} = v_2 + v_{J_{k_{i-1}+9}''} - 2 v_{h_{k_{i-1}+10}} \\
& = v_2 + \sum_{j=2}^{k_{i-1}+9} (-v_{h_j}) - 2 v_{h_{k_{i-1}+10}} \\
& = v_2 + \biggl( - \sum_{j=1}^{k_i} v_{h_j} + v_{h_N} \biggr) + \underbrace{ v_{h_1} + v_{h_{k_i}} - v_{h_N} - v_{h_{k_{i-1}+10}} }_{=0} = v_2 + v_1.
\end{align*}
We claim
\begin{align}
\label{eq:Case7proof_eq1}
\langle v_2, v_{h_j} \rangle=0, \qquad \forall j=2,3,\dots,k_{i-1}+9.
\end{align}
To prove this, note $v_2 = -v_{h_1} + \sum_{j=k_{i-1}+10}^N v_{h_j} = -v_{h_1} + v_{h_{k_{i-1}+10}} + \sum_{j=k_i}^N v_{h_j} = \sum_{j=k_i}^N v_{h_j}$. It~is easy to see by inspection on Figure~\ref{fig:one-holed_surface_of_higher_genus_flat5} that
\begin{align*}
\varepsilon_{h_{j'} h_j} = 0, \qquad \forall j' \in \{k_i,k_i+1,\dots,N\}, \quad \forall j \in \{2,3,\dots,k_{i-1}+9\}
\end{align*}
holds, for $h_{j'}$ and $h_j$ do not appear in a same triangle. Hence, equation~\eqref{eq:Case7proof_eq1} holds. We can show
\begin{align}
\label{eq:Case7proof_eq2}
\langle v_1, v_{h_j} \rangle=0, \qquad \forall j=k_i+1,k_i+2,\dots,N-1,
\end{align}
in a similar manner, using $v_1 = - \sum_{j=1}^{k_i} v_{h_j} + v_{h_N} = - \sum_{j=1}^{k_{i-1}+10} v_{h_j} - v_{h_{k_i}} + v_{h_N} = - \sum_{j=1}^{k_{i-1}+10} v_{h_j}$. Meanwhile, we claim
\begin{align}
\label{eq:Case7proof_eq3}
\langle v_\gamma, v_h \rangle=0 \qquad \mbox{for all internal arcs $h$ of $\Delta'_i$ not equal to $f_{i-1}$ or $f_{i-1}'$}.
\end{align}
First, note $v_\gamma = 2 \sum_h v_h$ where the sum is over all internal arcs $h$ of $\Delta'_i$. For each internal arc $h$ of $\Delta'_i$ not equal to $f_{i-1}$ or $f'_{i-1}$, let $h_1',h_2',h_3',h_4'$ be the sides of the ideal quadrilateral of $\Delta'_i$ containing $h$ as the diagonal. Then $\langle v_\gamma, v_h\rangle = 2\bigl(\varepsilon_{h_1' h} +\varepsilon_{h_2' h} +\varepsilon_{h_3' h} +\varepsilon_{h_4' h} \bigr)$, which one can easily check to be zero by inspection on Figure~\ref{fig:one-holed_surface_of_higher_genus_flat5}. Since $v_2 - v_1 = v_\gamma$, by combining equation~\eqref{eq:Case7proof_eq1} and equation~\eqref{eq:Case7proof_eq3}, we get
\begin{align}
\nonumber
\langle v_1, v_{h_j} \rangle=0, \qquad \forall j=2,3,\dots,k_{i-1}+9,
\end{align}
and by combining equation~\eqref{eq:Case7proof_eq2} and equation~\eqref{eq:Case7proof_eq3} we get
\begin{align}
\nonumber
\langle v_2, v_{h_j} \rangle=0, \qquad \forall j=k_i+1,k_i+2,\dots,N-1.
\end{align}
Thus,
\begin{align*}
&\langle v_1, v_{J_r''}\rangle = \langle v_2, v_{J_r''}\rangle = 0 \qquad \forall r = 1,2,\dots,k_{i-1}+9, \\
& \langle v_1, v_{J_r'}\rangle = \langle v_2, v_{J_r'}\rangle = 0 \qquad \forall r = k_i,k_i+1,\dots,N-1.
\end{align*}
Thus, for $r=1,2,\dots,k_{i-1}+9$ we get $Z^{J_r} = \exp \bigl(z_{v_1} + z_{v_{J_r''}}\bigr) = \exp(z_{v_1}) \exp\bigl(z_{v_{J_r''}}\bigr)$, and for $r=k_i,k_i+1,\dots,N-1$, we get $Z^{J_r} = \exp\bigl(z_{v_2} + z_{v_{J_r'}}\bigr) + \exp(z_{v_2}) \exp\bigl(z_{v_{J_r'}}\bigr)$. Thus,
\begin{align*}
f^\omega_{[\gamma],\Delta'_i} & = Z^{J_0} + Z^{J_N} + Z^{J_{k_{i-1}+10}} + \sum_{r=k_i}^{N-1} Z^{J_r} + \sum_{r=1}^{k_{i-1}+9} Z^{J_r} \\
& = {\rm e}^{v_\gamma} + {\rm e}^{-v_\gamma} + {\rm e}^{v_1 + v_2} + {\rm e}^{v_1} f^\omega_{[\zeta_i],\Delta'_i} + {\rm e}^{v_2} f^\omega_{[\vartheta_i],\Delta'_i},
\end{align*}
establishing \ref{TR1} of Definition~\ref{def:Teschner_triple_level1}. We showed that each of $z_{v_1}$ and $z_{v_2}$ strongly commutes with each of \smash{$f^\omega_{[\zeta_i],\Delta'_i}$} and \smash{$f^\omega_{[\vartheta_i],\Delta'_i}$}, establishing~\ref{TR3}. To show~\ref{TR4}, note that the set of all ideal arcs meeting $\zeta_i$ is $\{h_{j'} \mid j' \in \{k_i+1,k_i+2,\dots,N-1\}\}$, and the set of ideal arcs meeting $\vartheta_i$ is $\{h_j \mid j\in \{2,3,\dots,k_{i-1}+9\}\}$. We see from Figure~\ref{fig:one-holed_surface_of_higher_genus_flat5} that the set of ideal triangles having one of $h_{j'}$ as a side is disjoint from that the set of ideal triangles having one of $h_j$ as a side. Thus, $\zeta_i$ and $\vartheta_i$ are triangle-disjoint with respect to $\Delta'_i$, which establishes \ref{TR4} in view of Lemma~\ref{lem:triangle-disjoint_implies_strong_commutativity}.

For the remaining part of \ref{TR2}, we should verify $\langle v_1,v_2 \rangle =4$. Note $\langle v_1, v_2 \rangle = \langle v_1, v_\gamma + v_1\rangle = \langle v_1, v_\gamma \rangle$. Observe $v_1 = v_{h_N} - \sum_{j=1}^{k_i} v_{h_j} = v_{h_N} - v_{h_1} - \dots - v_{h_{k_i-1}} - v_{h_{k_i}} =
f_{i-1} - f'_{i-1} - \dots - f'_{i-1} - f_{i-1} = - 2f'_{i-1} - \cdots$, where $\cdots$ do not involve $f_{i-1}$ nor $f'_{i-1}$. In view of equation~\eqref{eq:Case7proof_eq3}, we then have $\langle v_1, v_\gamma \rangle = \langle -2f'_{i-1}, v_\gamma\rangle$. Note $v_\gamma = 2\sum_h v_h$, with the sum being over all internal ideal arcs of $\Delta_i'$. Note from Figure~\ref{fig:one-holed_surface_of_higher_genus_flat5} that $\varepsilon_{f'_{i-1} f_{i-1}}=-1$, $\varepsilon_{f'_{i-1} f_{i-2}} = -1$, $\varepsilon_{f'_{i-1} e_{i-1}}=1$, and $\varepsilon_{f'_{i-1} h}=0$ for all other internal arcs $h$ of $\Delta'_i$. Hence, it follows that $\langle -2f'_{i-1}, v_\gamma\rangle
= -2 \cdot 2 \bigl(\varepsilon_{f'_{i-1} f_{i-1}} + \varepsilon_{f'_{i-1} f_{i-2}} + \varepsilon_{f'_{i-1} e_{i-1}}\bigr) = -4 (-1 -1 +1) = 4$, as desired.

Therefore, $(\gamma,\zeta_i,\vartheta_i)$ is a strong Teschner triple with respect to $\Delta'_i$, hence Proposition~\ref{prop:zeta_i_and_vartheta_i} follows.
\end{proof}

This finishes our proof of Theorem~\ref{thm:algebraic_commutativity}, the second main theorem of the present paper.

\section{Consequences and future}
\label{sec:consequences}

\subsection{Toward spectral properties of length operators}
\label{subsec:toward_spectral_properties}

For any given essential simple loop $\gamma$ in a triangulable marked surface $\mathfrak{S}$, the first main theorem of the present paper, Theorem~\ref{thm:AS}, studies the algebraic structure of the quantized trace-of-monodromy $f^\omega_{[\gamma],\Delta} \in \mathcal{Z}^\omega_\Delta$ (Definition~\ref{def:quantized_trace-of-monodromy}), with respect to a suitably found ideal triangulation~$\Delta$ of~$\mathfrak{S}$. As mentioned in Section~\ref{subsec:on_quantized_length_operators}, the original motivation is to study the analytic properties of the densely defined operator corresponding to \smash{$f^\omega_{[\gamma],\Delta}$} via a representation of $\mathcal{Z}^\omega_\Delta$ on a Hilbert space~$\mathscr{H}_\Delta$. This will be the topic of the sequel~\cite{sequel}, and here we present some main ideas, which support the usefulness of the results of the present paper. We build on the rough description given in Section~\ref{subsec:on_quantized_length_operators}.

For the kind of representations we will study \cite{FG09, Kim_irreducible}, we have $\mathscr{H}_\Delta = L^2(\mathbb{R}^n, {\rm d}t_1 {\rm d}t_2\cdots {\rm d}t_n)$, and each Weyl-ordered Laurent monomial $Z_v \in \mathcal{Z}^\omega_\Delta$, is represented by $\exp({\bf z}_v)$, where ${\bf z}_v$ is a self-adjoint operator given as an $\mathbb{R}$-linear combination of the position operators ${\bf q}_j =t_j$, the momentum operators \smash{${\bf p}_j = {\rm i} \frac{\partial}{\partial t_j}$}, and the identity operator ${\rm Id}$; we call such an operator ${\bf z}_v$ a~{\em standard self-adjoint operator} in \cite{sequel}. We first represent elements of $\mathcal{Z}^\omega_\Delta$ on a nice dense subspace of $\mathscr{H}_\Delta$, say the `Hermite' subspace
\begin{align}
\label{eq:our_Hermite_subspace}
\mathscr{D}_\Delta = {\rm span}_\mathbb{C}\bigl\{ {\rm e}^{- \sum_{j=1}^n a_j t_j^2 + b_j t_j} P(t_1,\dots,t_n) \mid \mbox{$P$ a polynomial, $a_j>0$, $b_j \in \mathbb{C}$}\bigr\},
\end{align}
which is a variant of Fock and Goncharov's subspace $W_{\bf i}$ \cite{FG09}. For each $u\in \mathcal{Z}^\omega_\Delta$, denote by $\rho^\hbar_{\rm D}(u) = \rho^\hbar_{\Delta;{\rm D}}(u)$ the operator defined on $\mathscr{D}_\Delta$ representing $u \in \mathcal{Z}^\omega_\Delta$; in particular, \[\rho^\hbar_{\rm D}(Z_i) \rho^\hbar_{\rm D}(Z_j) = {\rm e}^{\pi {\rm i} \hbar \varepsilon_{ij}/2} \rho^\hbar_{\rm D}(Z_j) \rho^\hbar_{\rm D}(Z_i)\]
 for all $i,j \in \Delta$. The actual common maximal domain of operators for the skein algebra $\mathcal{S}^\omega(\mathfrak{S})$, called the {\em Schwartz space} $\mathscr{S}_\Delta$, is defined as the intersection of domains of the adjoints, following Fock and Goncharov \cite{FG09},
\begin{align}
\label{eq:Schwartz}
\mathscr{S}_\Delta = \bigcap_{u \in \mathcal{S}^\omega(\mathfrak{S})} {\rm Dom}\bigl( \bigl(\rho^\hbar_{\rm D}( {\rm Tr}^\omega_\Delta(u))\bigr)^* \bigr),
\end{align}
on which $\mathcal{S}^\omega(\mathfrak{S})$ acts via a representation $\rho^\hbar_\Delta$, defined as
\[
\rho^\hbar_\Delta(u) := \bigl(\rho^\hbar_{\rm D}(({\rm Tr}^\omega_\Delta(u))^*)\bigr)^* \restriction \mathscr{S}_\Delta,
\]
where in the right-hand side, the inner $*$ is the $*$-map of the algebra $\mathcal{Z}^\omega_\Delta$, while the outer $*$ is the operator adjoint.

Now, one goal is to show for each essential simple loop $\gamma$ in $\mathfrak{S}$ that $\rho^\hbar_\Delta([K_\gamma])$ is essentially self-adjoint on $\mathscr{S}_\Delta$ (where $K_\gamma$ is a constant elevation lift of $\gamma$), and that the unique self-adjoint extension \smash{${\bf f}^\hbar_{[\gamma],\Delta}$} has simple spectrum $[2,\infty)$. One approach is to prove such a statement for the~restriction $\rho^\hbar_\Delta([K_\gamma])\restriction \mathscr{D}_\Delta = \rho^\hbar_{\rm D}\bigl(f^\omega_{[\gamma],\Delta}\bigr)$ to the nice subspace $\mathscr{D}_\Delta$.

Indeed, for the case \ref{AS1} of Theorem~\ref{thm:AS}, i.e., when $\gamma$ is a peripheral loop, then
\[
\rho^\hbar_{\rm D}\bigl(f^\omega_{[\gamma],\Delta}\bigr) = \bigl( {\rm e}^{{\bf z}_{v_\gamma}} + {\rm e}^{-{\bf z}_{v_\gamma}} \bigr) \restriction \mathscr{D}_\Delta.
\]
It is well known that one can find a (Segal--)Shale--Weil type unitary operator $U$ on $\mathscr{H}_\Delta$ (this operator is reviewed in \cite{Kim_irreducible}), which generalizes Fourier transforms, that preserves $\mathscr{D}_\Delta$ and $U {\bf z}_{v_\gamma} U^{-1} = {\bf q}_1 = t_1$. So the problem boils down to the operator $\bigl({\rm e}^{{\bf q}_1} + {\rm e}^{-{\bf q}_1}\bigr) \restriction \mathscr{D}_\Delta$. We expect that it should be well-known and straightforward to prove directly that this is essentially self-adjoint, whose self-adjoint extension equals ${\rm e}^{{\bf q}_1} + {\rm e}^{-{\bf q}_1}$ which is the self-adjoint operator constructed by functional calculus of ${\bf q}_1$. In particular, it has simple spectrum $[2,\infty)$.

More interesting is the case \ref{AS2}, when
\[
\rho^\hbar_{\rm D}\bigl(f^\omega_{[\gamma],\Delta}\bigr) = \bigl({\rm e}^{{\bf z}_1} + {\rm e}^{-{\bf z}_1} + {\rm e}^{{\bf z}_2} \bigr)\restriction \mathscr{D}_\Delta,
\]
where ${\bf z}_1$ and ${\bf z}_2$ are standard self-adjoint operators such that $[{\bf z}_1, {\bf z}_2] = 8 \pi {\rm i} \hbar$. One can find a~Shale--Weil unitary operator $U$ on $\mathscr{H}_\Delta$ preserving $\mathscr{D}_\Delta$ such that $U {\bf z}_1 U^{-1} = 8 \pi \hbar {\bf p}_1 = 8 \pi \hbar {\rm i} \frac{\partial}{\partial t_1}$ and $U {\bf z}_2 U^{-1} = {\bf q}_1 = t_1$. So the problem boils down to the operator $\bigl({\rm e}^{8\pi \hbar {\bf p}_1} + {\rm e}^{-8\pi \hbar {\bf p}_1} + {\rm e}^{{\bf q}_1}\bigr)\restriction \mathscr{D}_\Delta$, which sends $\varphi(t_1,\dots,t_n) \in \mathscr{D}_\Delta$ to $\varphi(t_1 + 8\pi \hbar {\rm i}, t_2,\dots) + \varphi(t_1 - 8\pi \hbar {\rm i}, \dots) + {\rm e}^{t_1} \varphi(t_1,\dots)$. It is proved by Takhtajan and Faddeev \cite{FT} that this operator is essentially self-adjoint and its unique self-adjoint extension has simple spectrum $[2,\infty)$, exactly as we want. They actually proved it for the space smaller than $\mathscr{D}_\Delta$, which implies what we need.

The most intriguing and difficult is the remaining \ref{AS3}. Suppose that $(\gamma,\gamma_1,\gamma_2)$ is a~strong or weak Teschner triple, with respect to an ideal triangulation $\Delta$. Then the operator on the Schwartz space $\mathscr{S}_\Delta$ corresponding to the loop $\gamma$, or to the algebraic quantized trace-of-monodromy would be
\begin{align*}
\rho^\hbar_\Delta([K_\gamma]) = \bigl( {\rm e}^{{\bf z}_{v_\gamma}} + {\rm e}^{-{\bf z}_{v_\gamma}} + {\rm e}^{{\bf z}_1 + {\bf z}_2} + {\rm e}^{{\bf z}_1} {\bf f}^\hbar_{[\gamma_1],\Delta} + {\rm e}^{{\bf z}_2} {\bf f}^\hbar_{[\gamma_2],\Delta}\bigr) \restriction \mathscr{S}_\Delta,
\end{align*}
where ${\bf z}_1$ and ${\bf z}_2$ are standard self-adjoint operators such that
\[
[{\bf z}_1,{\bf z}_2]=8\pi {\rm i} \hbar, \qquad \text{with ${\bf z}_{v_\gamma} = {\bf z}_2 - {\bf z}_1$};
\]
 similarly for the operator \smash{$\rho^\hbar_{\rm D}\bigl(f^\omega_{[\gamma],\Delta}\bigr)$} on the Hermite subspace $\mathscr{D}_\Delta$. Here, where we wrote \smash{${\bf f}^\hbar_{[\gamma_j],\Delta}$}, we~should a priori write \smash{$\rho^\hbar_\Delta([K_{\gamma_j}])$} or \smash{$\rho^\hbar_{\rm D}(f^\omega_{[\gamma_j],\Delta})$} instead. In order to prove properties of \smash{$\rho^\hbar_{\rm D}\bigl(f^\omega_{[\gamma],\Delta}\bigr)$}, we~use~an induction, in the sense that we assume that \smash{$\rho^\hbar_\Delta([K_{\gamma_j}])$} or \smash{$\rho_{\rm D}(f^\omega_{[\gamma_j],\Delta})$} satisfies the sought-for property, i.e., has a unique self-adjoint extension \smash{${\bf f}^\hbar_{[\gamma_j],\Delta}$} such that it has simple spectrum $[2,\infty)$, hence
\[
{\bf f}^\hbar_{[\gamma_j],\Delta} = 2\cosh\bigl({\bf l}^\hbar_{[\gamma_j],\Delta}/2\bigr)
\]
for a unique self-adjoint operator \smash{${\bf l}^\hbar_{[\gamma_j],\Delta}$} having simple spectrum $[0,\infty)$, which we call the~length ope\-rator. From~\ref{TR3}, we would show that ${\bf z}_j$ strongly commutes with~\smash{${\bf f}^\hbar_{[\gamma_l],\Delta}$} and hence with \smash{${\bf l}^\hbar_{[\gamma_l],\Delta}$}, and from \ref{TR4} or \ref{TR5} we would show that \smash{${\bf f}^\hbar_{[\gamma_j],\Delta}$}, $j=1,2$, strongly commute with each other; hence so do \smash{${\bf l}^\hbar_{[\gamma_j],\Delta}$}, $j=1,2$.

Now, inspired by Teschner's suggestion \cite[equation~(15.11)]{T}, we consider the following unitary operator
\begin{align*}
{\bf U} :=
{\rm e}^{ \frac{1}{16\pi {\rm i} \hbar} ({\bf z}_1 + {\bf z}_2) {\bf l}_2}
\Phi^{4\hbar} \bigl(-{\bf z}_2 + \tfrac{1}{2} {\bf l}_1\bigr)
\Phi^{4\hbar}\bigl(-{\bf z}_2 - \tfrac{1}{2} {\bf l}_1\bigr)
\bigl(\Phi^{4\hbar}\bigl(- {\bf z}_1 - \tfrac{1}{2} {\bf l}_2\bigr)\bigr)^{-1},
\end{align*}
where each factor in the right-hand side is defined via functional calculus, where $\Phi^\hbar(z)$ stands for Faddeev's (non-compact) {\it quantum dilogarithm} function \cite{F95, FK94}, defined for a parameter $\hbar \in \mathbb{R}_{>0}$ and complex variable $z$ living in the strip $|{\operatorname{Im} z}| < \pi(1+\hbar)$ as the formula
\[
\Phi^\hbar(z) := \exp\biggl( - \frac{1}{4} \int_\Omega \frac{{\rm e}^{-{\rm i} pz}}{\sinh(\pi p) \sinh(\pi \hbar p)} \frac{{\rm d}p}{p} \biggr),
\]
going back at least to Barnes \cite{B01}, where $\Omega$ is the contour along the real line that avoids the origin along a small half-circle above the origin. This function analytically continues to a meromorphic function, satisfies $|\Phi^\hbar(x)|=1$ for all $x\in \mathbb{R}$, satisfies the difference equation $\Phi^\hbar(x+2\pi {\rm i} \hbar) = \bigl(1+{\rm e}^{\pi {\rm i} \hbar} {\rm e}^z\bigr) \Phi^\hbar(z)$, and was crucially used to construct unitary intertwiners for representations of quantum Teichm\"uller space \cite{CF, Kash98} and quantum cluster variety \cite{FG09}.

The difference equation gives an operator equation like
\begin{align*}
\Phi^{4\hbar}({\bf A}) {\rm e}^{{\bf B}} \bigl(\Phi^{4 \hbar}({\bf A})\bigr)^{-1} = {\rm e}^{{\bf B}} + {\rm e}^{{\bf A} + {\bf B}}
\end{align*}
whenever ${\bf A}$ and ${\bf B}$ are self-adjoint operators on a separable complex Hilbert space satisfying the `Weyl-relations version' (\cite[Definition~14.2]{Hall}, \cite[Definition~3.5]{Kim_phase}) of $[{\bf A}, {\bf B}] = 8 \pi {\rm i} \hbar$. The strategy is to show using equations like above to show
\[
{\bf U} \, \rho^\hbar_\Delta([K_\gamma]) {\bf U}^{-1} = \bigl( {\rm e}^{ {\bf z}_2 } + {\rm e}^{ - {\bf z}_2 } + {\rm e}^{{\bf z}_1 + {\bf z}_2}\bigr) \restriction \mathfrak{D}
\]
for a suitable domain $\mathfrak{D}$. If one shows that $\mathfrak{D}$ contains a nice subspace like in equation~\eqref{eq:our_Hermite_subspace}, then Takhtajan and Faddeev's result \cite{FT} should apply and we would obtain that the above operator on $\mathfrak{D}$ is essentially self-adjoint and has a unique self-adjoint extension with simple spectrum $[2,\infty)$. Formal computations work easily, and what matters is a careful investigation of the domains. As an example of the formal computations, we exhibit a result of conjugation of the last factor of ${\bf U}$: first, replace ${\bf f}^\hbar_{[\gamma_j],\Delta}$ by
\[
{\rm e}^{ \frac{1}{2} {\bf l}^\hbar_{[\gamma_j],\Delta}} + {\rm e}^{ - \frac{1}{2} {\bf l}^\hbar_{[\gamma_j],\Delta}},
\]
 so that
\[
\rho^\hbar_\Delta([K_\gamma])
= ( {\rm e}^{ {\bf z}_1 - {\bf z}_2 } + {\rm e}^{{\bf z}_2 - {\bf z}_1} + {\rm e}^{{\bf z}_1 + {\bf z}_2} + {\rm e}^{{\bf z}_1 + \frac{1}{2} {\bf l}_1} + {\rm e}^{{\bf z}_1 - \frac{1}{2} {\bf l}_1} + {\rm e}^{{\bf z}_2 + \frac{1}{2} {\bf l}_2} + {\rm e}^{{\bf z}_2 - \frac{1}{2} {\bf l}_2} )\restriction \mathscr{S}_\Delta,
\]
and after conjugation by $(\Phi^{4\hbar}(- {\bf z}_1 - \frac{1}{2} {\bf l}_2))^{-1}$ it becomes
\begin{align*}
& \textstyle \bigl(\Phi^{4\hbar}\bigl(- {\bf z}_1 - \tfrac{1}{2} {\bf l}_2\bigr)\bigr)^{-1} \rho^\hbar_\Delta([K_\gamma]) \bigl(\Phi^{4\hbar}\bigl(- {\bf z}_1 - \tfrac{1}{2} {\bf l}_2\bigr)\bigr) \\
& \qquad{}= \bigl({\rm e}^{{\bf z}_2 + \frac{1}{2} {\bf l}_2} + {\rm e}^{{\bf z}_1+{\bf z}_2} + \bigl({\rm e}^{{\bf z}_1-{\bf z}_2} + {\rm e}^{-{\bf z}_2 - \frac{1}{2} {\bf l}_2}\bigr)
+ {\rm e}^{{\bf z}_1 + \frac{1}{2} {\bf l}_1} + {\rm e}^{{\bf z}_1 - \frac{1}{2} {\bf l}_1}\bigr)\restriction{\mbox{(some domain)}}.
\end{align*}
This is because conjugation by \smash{$\bigl(\Phi^{4\hbar}\bigl(- {\bf z}_1 - \frac{1}{2} {\bf l}_2\bigr)\bigr)^{-1}$} sends ${\rm e}^{{\bf z}_2 + \frac{1}{2} {\bf l}_2} + {\rm e}^{{\bf z}_2 - {\bf z}_1}$ to ${\rm e}^{{\bf z}_2 + \frac{1}{2} {\bf l}_2}$, ${\rm e}^{{\bf z}_1 + {\bf z}_2} + {\rm e}^{{\bf z}_2 - \frac{1}{2} {\bf l}_2}$ to ${\rm e}^{{\bf z}_1+{\bf z}_2}$, ${\rm e}^{{\bf z}_1-{\bf z}_2}$ to ${\rm e}^{{\bf z}_1 - {\bf z}_2} + {\rm e}^{-{\bf z}_2 - \frac{1}{2} {\bf l}_2}$, and fixes the remaining two terms. So, modulo the domain issues, which actually might be quite subtle as investigated by Ruijsenaars in \cite{Ruijsenaars} (which was pointed out to the author by Ivan Ip), we obtained the desired analytic properties of $\rho^\hbar_\Delta([K_\gamma])$, yielding the self-adjoint operator ${\bf f}^\hbar_{[\gamma],\Delta}$ with simple spectrum $[2,\infty)$.

As such domain issues are topics of the sequel \cite{sequel}, for the moment let us assume that those issues are settled. We can note that the above conjugator ${\bf U}$ strongly commutes with ${\bf f}^\hbar_{[\gamma_j],\Delta}$, $j=1,2$. Hence, the situation after conjugation, modulo issues of domains, is
\[
{\bf U} \, {\bf f}^\hbar_{[\gamma],\Delta} \, {\bf U}^{-1}
= {\rm e}^{{\bf z}_2} + {\rm e}^{-{\bf z}_2} + {\rm e}^{{\bf z}_1+{\bf z}_2}, \qquad
{\bf U} \, {\bf f}^\hbar_{[\gamma_j],\Delta} \, {\bf U}^{-1}
= {\bf f}^\hbar_{[\gamma_j],\Delta} , \quad j=1,2.
\]
In particular, we see that \smash{${\bf U} \, {\bf f}^\hbar_{[\gamma],\Delta} \, {\bf U}^{-1}$} strongly commutes with \smash{${\bf U} \, {\bf f}^\hbar_{[\gamma_j],\Delta} \, {\bf U}^{-1}$} for $j=1,2$. Therefore, it follows that \smash{${\bf f}^\hbar_{[\gamma],\Delta}$} strongly commutes with \smash{${\bf f}^\hbar_{[\gamma_j],\Delta}$} for $j=1,2$. This says that the algebraic Teschner-commutativity of the present paper would imply the strong commutativity of corresponding operators, which would be one of the crucial parts of the proof in the sequel \cite{sequel} for item (3) of Conjecture \ref{conjecture:sequel}. Meanwhile, we would also need to establish that the algebraic strong commutativity implies strong commutativity of the corresponding operators.

All these will be carefully dealt with in \cite{sequel}, in order to settle Conjecture \ref{conjecture:sequel}.

We note that, the main difficulty of settling Conjecture \ref{conjecture:sequel} in \cite{sequel} shall lie in dealing with a~simple loop for the case of \ref{AS3} of Theorem~\ref{thm:AS}, and with a pair of simple loops for the case of \ref{AC3} of Theorem~\ref{thm:algebraic_commutativity}. These `difficult' cases involve the Teschner recursion relations. Using our classification results, i.e., Proposition~\ref{prop:loop_classification} for a simple loop and Proposition~\ref{prop:classification_of_pair_of_loops} for a pair of simple loops, we see that for each punctured surface of genus zero, those difficult cases never arise. This means that, the results of the current paper and the arguments we just gave in the present subsection essentially settle Conjecture \ref{conjecture:sequel} for genus zero surfaces.

\subsection{Future topics}\label{subsec:future_topics} We address some possible future directions.

{\bf Toward the modular functor conjecture.}
This is the very motivation of the present work and the sequel~\cite{sequel}. As suggested by Fock and Goncharov in \cite{FG09} and outlined in Section~\ref{subsec:on_quantized_length_operators}, the core statement to aim is to establish the decomposition
\begin{align}
\label{eq:direct_integral_decomposition2}
\mathscr{H}_\Delta \stackrel{\sim}{\to} \int_{[0,\infty)}^\oplus (\mathscr{H}_\Delta)_\chi \, {\rm d}\chi
\end{align}
as in equation~\eqref{eq:direct_integral_decomposition}, with suitable equivariance. As suggested already by Verlinde \cite{Verlinde}, this decomposition should be obtained as the spectral decomposition for the length operator \smash{${\bf l}^\hbar_{[\gamma],\Delta}$} of any chosen essential simple loop $\gamma$ in the triangulable punctured surface $\mathfrak{S} = \Sigma\setminus \mathcal{P}$. But this is only one necessary step. We then have to interpret each slice $(\mathscr{H}_\Delta)_\chi$ as the Hilbert space of representation of the skein algebra of the surface $\mathfrak{S}'$, obtained from $\mathfrak{S}$ by removing $\gamma$ and shrinking the two new holes to two new punctures. From an ideal triangulation $\Delta$ of $\mathfrak{S}$, one should associate an ideal triangulation~$\Delta'$ of the cut surface $\mathfrak{S}'$; one could use the suggestion by Roger~\cite{Roger}. For this to work, we should really be working with representations so that the peripheral loops are represented by scalars. That is, let us choose $\vec{l} = (l_p)_{p \in \mathcal{P}} \in (\mathbb{R}_{\ge 0})^\mathcal{P}$, i.e., a~nonzero real number per puncture of~$\mathfrak{S}$. We should consider a suitable representation of~$\mathcal{S}^\omega(\mathfrak{S})$ on a Hilbert space $\mathscr{H}_{\Delta;\vec{l}}$, such that for each peripheral loop $\gamma_p$ around a puncture $p\in \mathcal{P}$, the operator corresponding to $[K_{\gamma_p}]$ should be the scalar $2\cosh(l_p/2)$; we shall classify and study such representations in \cite{sequel}, using the results of \cite{Kim_irreducible}. Another way of saying this is to work with a representation of the relative skein algebra $\mathcal{S}^\omega(\mathfrak{S})_{\vec{l}}$ on a Hilbert space $\mathscr{H}_{\Delta;\vec{l}}$, which can be thought of as a formulation of quantization of the relative Teichm\"uller space $\mathscr{T}(\mathfrak{S})_{\vec{l}}$. For the cut surface $\mathfrak{S}'$, the set of punctures $\mathcal{P}'$ can be identified as $\mathcal{P} \sqcup \{p_1,p_2\}$, so that an element of~$(\mathbb{R}_{\ge 0})^{\mathcal{P}'}$ can be written as $\vec{l} \sqcup (\chi_1,\chi_2)$ for $\vec{l} \in (\mathbb{R}_{\ge 0})^\mathcal{P}$ and $\chi_1,\chi_2 \in \mathbb{R}_{\ge 0}$. What we need is to interpret the decomposition in equation~\eqref{eq:direct_integral_decomposition2} as the following:
\begin{align}
\label{eq:direct_integral_decomposition3}
\mathscr{H}_{\Delta;\vec{l}} \stackrel{\sim}{\to} \int_{ \mathbb{R}_{\ge 0} }^\oplus \mathscr{H}_{\Delta'; \vec{l}\sqcup(\chi,\chi)} \, {\rm d}\chi,
\end{align}
as suggested in \cite{FG09} (see equation~\eqref{eq:intro_decomposition}). We have a natural map $\mathcal{S}^\omega(\mathfrak{S}') \to \mathcal{S}^\omega(\mathfrak{S})$, so the above decomposition should be equivariant under the action of $\mathcal{S}^\omega(\mathfrak{S}')$. One crucial necessary tool for establishing this equivariance is the strong commutativity as in the item (3) of Conjecture~\ref{conjecture:sequel}. More subtle is the mapping class group equivariance. The representation of $\mathcal{S}^\omega(\mathfrak{S})$ on $\mathscr{H}_\Delta$ should have a suitable equivariance under the action of the mapping class group ${\rm MCG}(\mathfrak{S}) = {\rm Diff}_+(\mathfrak{S})/{\rm Diff}(S)_0$. This indeed is one of the major features of known representations of quantum Teichm\"uller space \cite{F97,CF,Kash98} and quantum cluster variety \cite{FG09}; in particular, ${\rm MCG}(\mathfrak{S})$ acts as unitary maps on $\mathscr{H}_\Delta$ via a projective unitary representation
\[
\rho^\hbar_{\Delta;{\rm MCG}} \colon \ {\rm MCG}(\mathfrak{S}) \to {\rm U}(\mathscr{H}_\Delta) = \mbox{(unitary maps on $\mathscr{H}_\Delta$)},
\]
which satisfies $\rho^\hbar_{\Delta;{\rm MCG}}(g_1g_2) = \rho^\hbar_{\Delta;{\rm MCG}}(g_1)\rho^\hbar_{\Delta;{\rm MCG}}(g_2)$ up to constants. We need such an equivariance for the representation of the relative skein algebra \smash{$\mathcal{S}^\omega(\mathfrak{S})_{\vec{l}}$} on \smash{$\mathscr{H}_{\Delta;\vec{l}}$}, under a subgroup of ${\rm MCG}(\mathfrak{S})$ preserving the data \smash{$\vec{l}$}; for this we will need results from \cite{Kim_irreducible}, as we shall investigate in \cite{sequel}. To establish an equivariance of the decomposition in equation~\eqref{eq:direct_integral_decomposition3} under the action of mapping class groups, one would need to use a natural map from ${\rm MCG}(\mathfrak{S}')$ to ${\rm MCG}(\mathfrak{S})$. A~subtle point is that the map ${\rm MCG}(\mathfrak{S}') \to {\rm MCG}(\mathfrak{S})$ is not uniquely determined, and there is an ambiguity by the Dehn twist $D_{[\gamma]} \in {\rm MCG}(\mathfrak{S})$ along $\gamma$. It seems that one crucial step of resolving this issue is to verify the following conjecture that relates the operator representing the Dehn twist and the operator representing the length function.
\begin{Conjecture}
\label{conjecture:Dehn}
Let $\mathfrak{S}$ be a triangulable marked surface and $\Delta$ its ideal triangulation. Let $\gamma$ be an essential non-peripheral simple loop in $\mathfrak{S}$. In the above setting, the following equality holds up to a constant:
\[
\rho^\hbar_{\Delta;{\rm MCG}} (D_{[\gamma]}) = \exp\biggl( {\rm i} \frac{ ({\bf l}^\hbar_{[\gamma],\Delta})^2 }{8\pi \hbar} \biggr).
\]
\end{Conjecture}
A version of this is conjecture is written in \cite[Proposition~12]{T}, using the operator ${\bf l}_{[\gamma],\sigma}$ defined via the equation $\mathbf{L}_{[\gamma],\sigma} = 2 \cosh(\mathbf{l}_{[\gamma],\sigma}/2)$, where ${\bf L}_{[\gamma],\sigma}$ are the operators Teschner constructs in \cite{T} using a suitably modified operator version of \ref{TR1}, per each `marking' $\sigma$, which is a~pants decomposition of the surface with some extra data, and a consituent loop~$\gamma$ of~$\sigma$. However, as mentioned in Section~\ref{subsec:on_quantized_length_operators}, it is not clear what is being quantized by the operators~${\bf L}_{[\gamma],\sigma}$ and~${\bf l}_{[\gamma],\sigma}$ constructed in \cite{T}, whereas it is clear what our operators ${\bf f}^\hbar_{[\gamma],\Delta}$ and ${\bf l}^\hbar_{[\gamma],\Delta}$ are quantizing.
Teschner also constructs a mapping class group action on the Hilbert space, by which he formulates the statement \cite[Proposition~12]{T}, a counterpart of the above Conjecture~\ref{conjecture:Dehn}.
Teschner refers this statement to Kashaev~\cite{Kash00}, who studied a quantization of certain space related to the Teichm\"uller space. However, Kashaev showed this only for a special case of~$\mathfrak{S}$,~$\Delta$ and~$\gamma$, for which he could define the length operator. Indeed, only by using the present paper and the sequel~\cite{sequel} one can finally make sense of the length operator ${\bf l}^\hbar_{[\gamma],\Delta}$ for general~$\mathfrak{S}$,~$\Delta$ and~$\gamma$. We believe that this conjecture is interesting enough on its own, independent of its connection to the modular functor conjecture.

{\bf On Teschner recursion.}
In this paper, we investigated only some of the Teschner triples $(\gamma,\gamma_1,\gamma_2)$, in order to investigate properties of the quantized trace-of-monodromy along some special loops $\gamma$, and some special pairs of loops $(\gamma,\gamma_1)$. One might further want to delve into the study of Teschner triples, e.g., to classify them.
\begin{Conjecture}
Let $\mathfrak{S}$ be a triangulable marked surface. For {\em each} pants triple $(\gamma,\gamma_1,\gamma_2)$ in $\mathfrak{S}$ $($Definition~{\rm \ref{def:Teschner_triple_level1})}, there exists an ideal triangulation $\Delta$ such that $(\gamma,\gamma_1,\gamma_2)$ is a~Teschner triple of some `weakness level' $($see Remark~{\rm \ref{rem:further_weakness})} with respect to~$\Delta$.
\end{Conjecture}
All pants triples were classified in a previous version of the present paper \cite{ver3}. For some of the cases in the classification list, we already proved the conjecture, and it remains to deal with the remaining cases. Along the way, perhaps it will be useful and meaningful if one can characterize the ideal triangulations with respect to which a given pants triple is a Teschner triple.

Another interesting question is on the meaning of the Teschner recursion formula~\ref{TR1} at the classical limit $\omega \to 1$. This gives an equality
\[
f_{[\gamma]} = {\rm e}^{z_{v_\gamma}} + {\rm e}^{-z_{v_\gamma}} + {\rm e}^{z_{v_1} + z_{v_2}} + {\rm e}^{z_{v_1}} f_{[\gamma_1]} + {\rm e}^{z_{v_2}} f_{[\gamma_2]}
\]
of functions on the enhanced Teichm\"uller space $\mathscr{T}^+(\mathfrak{S})$, where each $z_{v_i}$ is now the function $x_i/2$, where~$x_i$ is the unexponentiated shear coordinate function in equation~\eqref{eq:log_shear}. Note that $(\gamma,\gamma_1,\gamma_2)$ is a pants triple, so that the corresponding length functions $l_{[\gamma]}$, $l_{[\gamma_1]}$, $l_{[\gamma_2]}$ are independent, as they are part of the Fenchel--Nielsen coordinate system associated to any pants decomposition of $\mathfrak{S}$ containing that special pair of pants bounded by $\gamma$, $\gamma_1$, $\gamma_2$. So the functions $f_{[\gamma]}$, $f_{[\gamma_1]}$, $f_{[\gamma_2]}$ are also independent. Yet, the above gives an equation involving these three functions, with the help of some shear coordinate functions. It will be nice if one can come up with a geometric interpretation of this classical equation, which might also help approaching the above mentioned problem of characterization of ideal triangulations with respect to which the triple $(\gamma,\gamma_1,\gamma_2)$ is a Teschner triple.

\subsection*{Acknowledgements}

This research was supported by Basic Science Research Program through the National Research Foundation of Korea (NRF) funded by the Ministry of Education (Grant Number 2017R1D1A1\-B03030230). This work was supported by the National Research Foundation of Korea (NRF) grant funded by the Korea government (MSIT) (No. 2020R1C1C1A01011151). H.K. has been supported by KIAS Individual Grant (MG047203, MG047204) at Korea Institute for Advanced Study. H.K. thanks Dylan Allegretti and Ivan Ip for helpful discussions, Igor Frenkel and Alexander Goncharov for encouragements, and anonymous referees for their effort put into reviewing this paper and improving it.


\pdfbookmark[1]{References}{ref}
\LastPageEnding

\end{document}

%% file: spiraling.pdf_tex
\begingroup%
  \makeatletter%
  \providecommand\color[2][]{%
    \errmessage{(Inkscape) Color is used for the text in Inkscape, but the package 'color.sty' is not loaded}%
    \renewcommand\color[2][]{}%
  }%
  \providecommand\transparent[1]{%
    \errmessage{(Inkscape) Transparency is used (non-zero) for the text in Inkscape, but the package 'transparent.sty' is not loaded}%
    \renewcommand\transparent[1]{}%
  }%
  \providecommand\rotatebox[2]{#2}%
  \newcommand*\fsize{\dimexpr\f@size pt\relax}%
  \newcommand*\lineheight[1]{\fontsize{\fsize}{#1\fsize}\selectfont}%
  \ifx\svgwidth\undefined%
    \setlength{\unitlength}{181.41732283bp}%
    \ifx\svgscale\undefined%
      \relax%
    \else%
      \setlength{\unitlength}{\unitlength * \real{\svgscale}}%
    \fi%
  \else%
    \setlength{\unitlength}{\svgwidth}%
  \fi%
  \global\let\svgwidth\undefined%
  \global\let\svgscale\undefined%
  \makeatother%
  \begin{picture}(1,0.984375)%
    \lineheight{1}%
    \setlength\tabcolsep{0pt}%
    \put(0,0){\includegraphics[width=\unitlength,page=1]{spiraling.pdf}}%
    \put(0.13918605,0.15044904){\color[rgb]{0,0,0}\makebox(0,0)[lt]{\lineheight{1.25}\smash{\begin{tabular}[t]{l}surface (interior)\end{tabular}}}}%
    \put(0.24762911,0.873657){\color[rgb]{0,0,0}\makebox(0,0)[lt]{\lineheight{1.25}\smash{\begin{tabular}[t]{l}funnel end for $p$\end{tabular}}}}%
    \put(0.73994493,0.70907032){\color[rgb]{0,0,0}\makebox(0,0)[lt]{\lineheight{1.25}\smash{\begin{tabular}[t]{l}\red{geodesic loop}\end{tabular}}}}%
    \put(0.72452635,0.54023674){\color[rgb]{0,0,0}\makebox(0,0)[lt]{\lineheight{1.25}\smash{\begin{tabular}[t]{l}\blue{geodesic arc $i$}\end{tabular}}}}%
  \end{picture}%
\endgroup%

%% file: spiraling2.pdf_tex
\begingroup%
  \makeatletter%
  \providecommand\color[2][]{%
    \errmessage{(Inkscape) Color is used for the text in Inkscape, but the package 'color.sty' is not loaded}%
    \renewcommand\color[2][]{}%
  }%
  \providecommand\transparent[1]{%
    \errmessage{(Inkscape) Transparency is used (non-zero) for the text in Inkscape, but the package 'transparent.sty' is not loaded}%
    \renewcommand\transparent[1]{}%
  }%
  \providecommand\rotatebox[2]{#2}%
  \newcommand*\fsize{\dimexpr\f@size pt\relax}%
  \newcommand*\lineheight[1]{\fontsize{\fsize}{#1\fsize}\selectfont}%
  \ifx\svgwidth\undefined%
    \setlength{\unitlength}{181.41732283bp}%
    \ifx\svgscale\undefined%
      \relax%
    \else%
      \setlength{\unitlength}{\unitlength * \real{\svgscale}}%
    \fi%
  \else%
    \setlength{\unitlength}{\svgwidth}%
  \fi%
  \global\let\svgwidth\undefined%
  \global\let\svgscale\undefined%
  \makeatother%
  \begin{picture}(1,0.984375)%
    \lineheight{1}%
    \setlength\tabcolsep{0pt}%
    \put(0,0){\includegraphics[width=\unitlength,page=1]{spiraling2.pdf}}%
    \put(0.30790278,0.12394339){\color[rgb]{0,0,0}\makebox(0,0)[lt]{\lineheight{1.25}\smash{\begin{tabular}[t]{l}surface (interior)\end{tabular}}}}%
    \put(0.24513833,0.84594074){\color[rgb]{0,0,0}\makebox(0,0)[lt]{\lineheight{1.25}\smash{\begin{tabular}[t]{l}funnel end for $p$\end{tabular}}}}%
    \put(0.73994493,0.70907032){\color[rgb]{0,0,0}\makebox(0,0)[lt]{\lineheight{1.25}\smash{\begin{tabular}[t]{l}\red{geodesic loop}\end{tabular}}}}%
    \put(0.72164804,0.46118677){\color[rgb]{0,0,0}\makebox(0,0)[lt]{\lineheight{1.25}\smash{\begin{tabular}[t]{l}\blue{geodesic arc $i$}\end{tabular}}}}%
  \end{picture}%
\endgroup%

%% file: ideal_quadrilateral_in_H.pdf_tex
\begingroup%
  \makeatletter%
  \providecommand\color[2][]{%
    \errmessage{(Inkscape) Color is used for the text in Inkscape, but the package 'color.sty' is not loaded}%
    \renewcommand\color[2][]{}%
  }%
  \providecommand\transparent[1]{%
    \errmessage{(Inkscape) Transparency is used (non-zero) for the text in Inkscape, but the package 'transparent.sty' is not loaded}%
    \renewcommand\transparent[1]{}%
  }%
  \providecommand\rotatebox[2]{#2}%
  \newcommand*\fsize{\dimexpr\f@size pt\relax}%
  \newcommand*\lineheight[1]{\fontsize{\fsize}{#1\fsize}\selectfont}%
  \ifx\svgwidth\undefined%
    \setlength{\unitlength}{167.24409449bp}%
    \ifx\svgscale\undefined%
      \relax%
    \else%
      \setlength{\unitlength}{\unitlength * \real{\svgscale}}%
    \fi%
  \else%
    \setlength{\unitlength}{\svgwidth}%
  \fi%
  \global\let\svgwidth\undefined%
  \global\let\svgscale\undefined%
  \makeatother%
  \begin{picture}(1,0.54237288)%
    \lineheight{1}%
    \setlength\tabcolsep{0pt}%
    \put(0,0){\includegraphics[width=\unitlength,page=1]{ideal_quadrilateral_in_H.pdf}}%
    \put(0.46160816,0.36244975){\color[rgb]{0,0,0}\makebox(0,0)[lt]{\lineheight{1.25}\smash{\begin{tabular}[t]{l}$\til{i}$\end{tabular}}}}%
    \put(0.83193155,0.39528856){\color[rgb]{0,0,0}\makebox(0,0)[lt]{\lineheight{1.25}\smash{\begin{tabular}[t]{l}$\mathbb{H}^2$\end{tabular}}}}%
    \put(0.03549457,0.02753324){\color[rgb]{0,0,0}\makebox(0,0)[lt]{\lineheight{1.25}\smash{\begin{tabular}[t]{l}$p_1$\end{tabular}}}}%
    \put(0.27765566,0.02753324){\color[rgb]{0,0,0}\makebox(0,0)[lt]{\lineheight{1.25}\smash{\begin{tabular}[t]{l}$p_2$\end{tabular}}}}%
    \put(0.65435001,0.02753324){\color[rgb]{0,0,0}\makebox(0,0)[lt]{\lineheight{1.25}\smash{\begin{tabular}[t]{l}$p_3$\end{tabular}}}}%
    \put(0.85166605,0.02753324){\color[rgb]{0,0,0}\makebox(0,0)[lt]{\lineheight{1.25}\smash{\begin{tabular}[t]{l}$p_4$\end{tabular}}}}%
  \end{picture}%
\endgroup%

%% file: kauffmantriple.pdf_tex
\begingroup%
  \makeatletter%
  \providecommand\color[2][]{%
    \errmessage{(Inkscape) Color is used for the text in Inkscape, but the package 'color.sty' is not loaded}%
    \renewcommand\color[2][]{}%
  }%
  \providecommand\transparent[1]{%
    \errmessage{(Inkscape) Transparency is used (non-zero) for the text in Inkscape, but the package 'transparent.sty' is not loaded}%
    \renewcommand\transparent[1]{}%
  }%
  \providecommand\rotatebox[2]{#2}%
  \ifx\svgwidth\undefined%
    \setlength{\unitlength}{368.50393701bp}%
    \ifx\svgscale\undefined%
      \relax%
    \else%
      \setlength{\unitlength}{\unitlength * \real{\svgscale}}%
    \fi%
  \else%
    \setlength{\unitlength}{\svgwidth}%
  \fi%
  \global\let\svgwidth\undefined%
  \global\let\svgscale\undefined%
  \makeatother%
  \begin{picture}(1,0.30769231)%
    \put(0,0){\includegraphics[width=\unitlength,page=1]{kauffmantriple.pdf}}%
    \put(0.12369109,0.02563069){\color[rgb]{0,0,0}\makebox(0,0)[lb]{\smash{$K_1$}}}%
    \put(0.45855945,0.02563069){\color[rgb]{0,0,0}\makebox(0,0)[lb]{\smash{$K_0$}}}%
    \put(0.80251782,0.02563069){\color[rgb]{0,0,0}\makebox(0,0)[lb]{\smash{$K_\infty$}}}%
  \end{picture}%
\endgroup%

%% file: JS_nonadmissible.pdf_tex
\begingroup%
  \makeatletter%
  \providecommand\color[2][]{%
    \errmessage{(Inkscape) Color is used for the text in Inkscape, but the package 'color.sty' is not loaded}%
    \renewcommand\color[2][]{}%
  }%
  \providecommand\transparent[1]{%
    \errmessage{(Inkscape) Transparency is used (non-zero) for the text in Inkscape, but the package 'transparent.sty' is not loaded}%
    \renewcommand\transparent[1]{}%
  }%
  \providecommand\rotatebox[2]{#2}%
  \newcommand*\fsize{\dimexpr\f@size pt\relax}%
  \newcommand*\lineheight[1]{\fontsize{\fsize}{#1\fsize}\selectfont}%
  \ifx\svgwidth\undefined%
    \setlength{\unitlength}{164.40944882bp}%
    \ifx\svgscale\undefined%
      \relax%
    \else%
      \setlength{\unitlength}{\unitlength * \real{\svgscale}}%
    \fi%
  \else%
    \setlength{\unitlength}{\svgwidth}%
  \fi%
  \global\let\svgwidth\undefined%
  \global\let\svgscale\undefined%
  \makeatother%
  \begin{picture}(1,0.89655172)%
    \lineheight{1}%
    \setlength\tabcolsep{0pt}%
    \put(0,0){\includegraphics[width=\unitlength,page=1]{JS_nonadmissible.pdf}}%
    \put(0.33963609,0.3634241){\color[rgb]{0,0,0}\makebox(0,0)[lt]{\lineheight{1.25}\smash{\begin{tabular}[t]{l}\red{loop segment}\end{tabular}}}}%
    \put(0.21084376,0.49002856){\color[rgb]{0,0,0}\makebox(0,0)[lt]{\lineheight{1.25}\smash{\begin{tabular}[t]{l}$-$\end{tabular}}}}%
    \put(0.7345471,0.50162485){\color[rgb]{0,0,0}\makebox(0,0)[lt]{\lineheight{1.25}\smash{\begin{tabular}[t]{l}$+$\end{tabular}}}}%
    \put(0,0){\includegraphics[width=\unitlength,page=2]{JS_nonadmissible.pdf}}%
  \end{picture}%
\endgroup%

%% file: JS_nonadmissible2.pdf_tex
\begingroup%
  \makeatletter%
  \providecommand\color[2][]{%
    \errmessage{(Inkscape) Color is used for the text in Inkscape, but the package 'color.sty' is not loaded}%
    \renewcommand\color[2][]{}%
  }%
  \providecommand\transparent[1]{%
    \errmessage{(Inkscape) Transparency is used (non-zero) for the text in Inkscape, but the package 'transparent.sty' is not loaded}%
    \renewcommand\transparent[1]{}%
  }%
  \providecommand\rotatebox[2]{#2}%
  \newcommand*\fsize{\dimexpr\f@size pt\relax}%
  \newcommand*\lineheight[1]{\fontsize{\fsize}{#1\fsize}\selectfont}%
  \ifx\svgwidth\undefined%
    \setlength{\unitlength}{337.32283465bp}%
    \ifx\svgscale\undefined%
      \relax%
    \else%
      \setlength{\unitlength}{\unitlength * \real{\svgscale}}%
    \fi%
  \else%
    \setlength{\unitlength}{\svgwidth}%
  \fi%
  \global\let\svgwidth\undefined%
  \global\let\svgscale\undefined%
  \makeatother%
  \begin{picture}(1,0.43697479)%
    \lineheight{1}%
    \setlength\tabcolsep{0pt}%
    \put(0,0){\includegraphics[width=\unitlength,page=1]{JS_nonadmissible2.pdf}}%
    \put(0.22545343,0.10319886){\color[rgb]{0,0,0}\makebox(0,0)[lt]{\lineheight{1.25}\smash{\begin{tabular}[t]{l}\red{loop segment}\end{tabular}}}}%
    \put(0.19163277,0.00793754){\color[rgb]{0,0,0}\makebox(0,0)[lt]{\lineheight{1.25}\smash{\begin{tabular}[t]{l}$-$\end{tabular}}}}%
    \put(0.16064254,0.16967608){\color[rgb]{0,0,0}\makebox(0,0)[lt]{\lineheight{1.25}\smash{\begin{tabular}[t]{l}$+$\end{tabular}}}}%
    \put(0,0){\includegraphics[width=\unitlength,page=2]{JS_nonadmissible2.pdf}}%
    \put(0.65430683,0.22693338){\color[rgb]{0,0,0}\makebox(0,0)[lt]{\lineheight{1.25}\smash{\begin{tabular}[t]{l}\red{loop segment}\end{tabular}}}}%
    \put(0.6962875,0.12052479){\color[rgb]{0,0,0}\makebox(0,0)[lt]{\lineheight{1.25}\smash{\begin{tabular}[t]{l}$-$\end{tabular}}}}%
    \put(0.58392223,0.24101852){\color[rgb]{0,0,0}\makebox(0,0)[lt]{\lineheight{1.25}\smash{\begin{tabular}[t]{l}$+$\end{tabular}}}}%
    \put(0,0){\includegraphics[width=\unitlength,page=3]{JS_nonadmissible2.pdf}}%
  \end{picture}%
\endgroup%

%% file: leftturn_NSF.pdf_tex
\begingroup%
  \makeatletter%
  \providecommand\color[2][]{%
    \errmessage{(Inkscape) Color is used for the text in Inkscape, but the package 'color.sty' is not loaded}%
    \renewcommand\color[2][]{}%
  }%
  \providecommand\transparent[1]{%
    \errmessage{(Inkscape) Transparency is used (non-zero) for the text in Inkscape, but the package 'transparent.sty' is not loaded}%
    \renewcommand\transparent[1]{}%
  }%
  \providecommand\rotatebox[2]{#2}%
  \newcommand*\fsize{\dimexpr\f@size pt\relax}%
  \newcommand*\lineheight[1]{\fontsize{\fsize}{#1\fsize}\selectfont}%
  \ifx\svgwidth\undefined%
    \setlength{\unitlength}{164.40944882bp}%
    \ifx\svgscale\undefined%
      \relax%
    \else%
      \setlength{\unitlength}{\unitlength * \real{\svgscale}}%
    \fi%
  \else%
    \setlength{\unitlength}{\svgwidth}%
  \fi%
  \global\let\svgwidth\undefined%
  \global\let\svgscale\undefined%
  \makeatother%
  \begin{picture}(1,0.89655172)%
    \lineheight{1}%
    \setlength\tabcolsep{0pt}%
    \put(0,0){\includegraphics[width=\unitlength,page=1]{leftturn_NSF.pdf}}%
  \end{picture}%
\endgroup%

%% file: leftturn_SF.pdf_tex
\begingroup%
  \makeatletter%
  \providecommand\color[2][]{%
    \errmessage{(Inkscape) Color is used for the text in Inkscape, but the package 'color.sty' is not loaded}%
    \renewcommand\color[2][]{}%
  }%
  \providecommand\transparent[1]{%
    \errmessage{(Inkscape) Transparency is used (non-zero) for the text in Inkscape, but the package 'transparent.sty' is not loaded}%
    \renewcommand\transparent[1]{}%
  }%
  \providecommand\rotatebox[2]{#2}%
  \newcommand*\fsize{\dimexpr\f@size pt\relax}%
  \newcommand*\lineheight[1]{\fontsize{\fsize}{#1\fsize}\selectfont}%
  \ifx\svgwidth\undefined%
    \setlength{\unitlength}{337.32283465bp}%
    \ifx\svgscale\undefined%
      \relax%
    \else%
      \setlength{\unitlength}{\unitlength * \real{\svgscale}}%
    \fi%
  \else%
    \setlength{\unitlength}{\svgwidth}%
  \fi%
  \global\let\svgwidth\undefined%
  \global\let\svgscale\undefined%
  \makeatother%
  \begin{picture}(1,0.43697479)%
    \lineheight{1}%
    \setlength\tabcolsep{0pt}%
    \put(0,0){\includegraphics[width=\unitlength,page=1]{leftturn_SF.pdf}}%
  \end{picture}%
\endgroup%

%% file: one-holed_torus2.pdf_tex
\begingroup%
  \makeatletter%
  \providecommand\color[2][]{%
    \errmessage{(Inkscape) Color is used for the text in Inkscape, but the package 'color.sty' is not loaded}%
    \renewcommand\color[2][]{}%
  }%
  \providecommand\transparent[1]{%
    \errmessage{(Inkscape) Transparency is used (non-zero) for the text in Inkscape, but the package 'transparent.sty' is not loaded}%
    \renewcommand\transparent[1]{}%
  }%
  \providecommand\rotatebox[2]{#2}%
  \newcommand*\fsize{\dimexpr\f@size pt\relax}%
  \newcommand*\lineheight[1]{\fontsize{\fsize}{#1\fsize}\selectfont}%
  \ifx\svgwidth\undefined%
    \setlength{\unitlength}{412.48199223bp}%
    \ifx\svgscale\undefined%
      \relax%
    \else%
      \setlength{\unitlength}{\unitlength * \real{\svgscale}}%
    \fi%
  \else%
    \setlength{\unitlength}{\svgwidth}%
  \fi%
  \global\let\svgwidth\undefined%
  \global\let\svgscale\undefined%
  \makeatother%
  \begin{picture}(1,0.53242852)%
    \lineheight{1}%
    \setlength\tabcolsep{0pt}%
    \put(0,0){\includegraphics[width=\unitlength,page=1]{one-holed_torus2.pdf}}%
    \put(0.54814215,0.26052161){\color[rgb]{0,0,0}\makebox(0,0)[lt]{\lineheight{1.25}\smash{\begin{tabular}[t]{l}$\textcolor{red}{\gamma}$\end{tabular}}}}%
    \put(0.42640197,0.0199637){\color[rgb]{0,0,0}\makebox(0,0)[lt]{\lineheight{1.25}\smash{\begin{tabular}[t]{l}$\mbox{one-holed-torus subsurface $(\Sigma',\mathcal{P}')$}$\end{tabular}}}}%
    \put(0,0){\includegraphics[width=\unitlength,page=2]{one-holed_torus2.pdf}}%
    \put(0.9308869,0.32568466){\color[rgb]{0,0,0}\makebox(0,0)[lt]{\lineheight{1.25}\smash{\begin{tabular}[t]{l}$\textcolor{blue}{\eta}$\end{tabular}}}}%
    \put(0,0){\includegraphics[width=\unitlength,page=3]{one-holed_torus2.pdf}}%
    \put(0.12335584,0.2291707){\color[rgb]{0,0,0}\makebox(0,0)[lt]{\lineheight{1.25}\smash{\begin{tabular}[t]{l}$\mbox{rest of $(\Sigma,\mathcal{P})$}$\end{tabular}}}}%
    \put(0,0){\includegraphics[width=\unitlength,page=4]{one-holed_torus2.pdf}}%
    \put(0.55163156,0.22881847){\color[rgb]{0,0,0}\makebox(0,0)[lt]{\lineheight{1.25}\smash{\begin{tabular}[t]{l}$\textcolor{red}{\circled{1}}$\end{tabular}}}}%
    \put(0.55171393,0.19233448){\color[rgb]{0,0,0}\makebox(0,0)[lt]{\lineheight{1.25}\smash{\begin{tabular}[t]{l}$\textcolor{red}{\circled{2}}$\end{tabular}}}}%
    \put(0.52863067,0.2988661){\color[rgb]{0,0,0}\makebox(0,0)[lt]{\lineheight{1.25}\smash{\begin{tabular}[t]{l}$\textcolor{red}{\circled{8}}$\end{tabular}}}}%
    \put(0.55065955,0.15782636){\color[rgb]{0,0,0}\makebox(0,0)[lt]{\lineheight{1.25}\smash{\begin{tabular}[t]{l}$\textcolor{red}{\circled{3}}$\end{tabular}}}}%
    \put(0.54253466,0.13026633){\color[rgb]{0,0,0}\makebox(0,0)[lt]{\lineheight{1.25}\smash{\begin{tabular}[t]{l}$\textcolor{red}{\circled{4}}$\end{tabular}}}}%
    \put(0.4743614,0.12783549){\color[rgb]{0,0,0}\makebox(0,0)[lt]{\lineheight{1.25}\smash{\begin{tabular}[t]{l}$\textcolor{red}{\circled{5}}$\end{tabular}}}}%
    \put(0.49828375,0.36267602){\color[rgb]{0,0,0}\makebox(0,0)[lt]{\lineheight{1.25}\smash{\begin{tabular}[t]{l}$\textcolor{red}{\circled{6}}$\end{tabular}}}}%
    \put(0.51357421,0.3295087){\color[rgb]{0,0,0}\makebox(0,0)[lt]{\lineheight{1.25}\smash{\begin{tabular}[t]{l}$\textcolor{red}{\circled{7}}$\end{tabular}}}}%
    \put(0.86834813,0.32949153){\color[rgb]{0,0,0}\makebox(0,0)[lt]{\lineheight{1.25}\smash{\begin{tabular}[t]{l}$\textcolor{blue}{\circled{$1'$}}$\end{tabular}}}}%
    \put(0.92149063,0.28645743){\color[rgb]{0,0,0}\makebox(0,0)[lt]{\lineheight{1.25}\smash{\begin{tabular}[t]{l}$\textcolor{blue}{\circled{$2'$}}$\end{tabular}}}}%
    \put(0.42112951,0.20515788){\color[rgb]{0,0,0}\makebox(0,0)[lt]{\lineheight{1.25}\smash{\begin{tabular}[t]{l}$p$\end{tabular}}}}%
    \put(0.57877792,0.2981391){\color[rgb]{0,0,0}\makebox(0,0)[lt]{\lineheight{1.25}\smash{\begin{tabular}[t]{l}$d$\end{tabular}}}}%
    \put(0.6073507,0.23300806){\color[rgb]{0,0,0}\makebox(0,0)[lt]{\lineheight{1.25}\smash{\begin{tabular}[t]{l}$c$\end{tabular}}}}%
    \put(0.64580253,0.16936888){\color[rgb]{0,0,0}\makebox(0,0)[lt]{\lineheight{1.25}\smash{\begin{tabular}[t]{l}$a$\end{tabular}}}}%
    \put(0.7147946,0.38730047){\color[rgb]{0,0,0}\makebox(0,0)[lt]{\lineheight{1.25}\smash{\begin{tabular}[t]{l}$b$\end{tabular}}}}%
    \put(0,0){\includegraphics[width=\unitlength,page=5]{one-holed_torus2.pdf}}%
  \end{picture}%
\endgroup%

%% file: one-holed_torus_flat2.pdf_tex
\begingroup%
  \makeatletter%
  \providecommand\color[2][]{%
    \errmessage{(Inkscape) Color is used for the text in Inkscape, but the package 'color.sty' is not loaded}%
    \renewcommand\color[2][]{}%
  }%
  \providecommand\transparent[1]{%
    \errmessage{(Inkscape) Transparency is used (non-zero) for the text in Inkscape, but the package 'transparent.sty' is not loaded}%
    \renewcommand\transparent[1]{}%
  }%
  \providecommand\rotatebox[2]{#2}%
  \newcommand*\fsize{\dimexpr\f@size pt\relax}%
  \newcommand*\lineheight[1]{\fontsize{\fsize}{#1\fsize}\selectfont}%
  \ifx\svgwidth\undefined%
    \setlength{\unitlength}{198.42519685bp}%
    \ifx\svgscale\undefined%
      \relax%
    \else%
      \setlength{\unitlength}{\unitlength * \real{\svgscale}}%
    \fi%
  \else%
    \setlength{\unitlength}{\svgwidth}%
  \fi%
  \global\let\svgwidth\undefined%
  \global\let\svgscale\undefined%
  \makeatother%
  \begin{picture}(1,1)%
    \lineheight{1}%
    \setlength\tabcolsep{0pt}%
    \put(0,0){\includegraphics[width=\unitlength,page=1]{one-holed_torus_flat2.pdf}}%
    \put(0.41558612,0.05404439){\color[rgb]{0,0,0}\makebox(0,0)[lt]{\lineheight{1.25}\smash{\begin{tabular}[t]{l}$a$\end{tabular}}}}%
    \put(0.41637133,0.92193097){\color[rgb]{0,0,0}\makebox(0,0)[lt]{\lineheight{1.25}\smash{\begin{tabular}[t]{l}$a$\end{tabular}}}}%
    \put(0.90884689,0.9374533){\color[rgb]{0,0,0}\makebox(0,0)[lt]{\lineheight{1.25}\smash{\begin{tabular}[t]{l}$p$\end{tabular}}}}%
    \put(0.92345238,0.04859359){\color[rgb]{0,0,0}\makebox(0,0)[lt]{\lineheight{1.25}\smash{\begin{tabular}[t]{l}$p$\end{tabular}}}}%
    \put(0.03898899,0.94061632){\color[rgb]{0,0,0}\makebox(0,0)[lt]{\lineheight{1.25}\smash{\begin{tabular}[t]{l}$p$\end{tabular}}}}%
    \put(0.04001508,0.04067181){\color[rgb]{0,0,0}\makebox(0,0)[lt]{\lineheight{1.25}\smash{\begin{tabular}[t]{l}$p$\end{tabular}}}}%
    \put(0.70104921,0.46547732){\color[rgb]{0,0,0}\makebox(0,0)[lt]{\lineheight{1.25}\smash{\begin{tabular}[t]{l}$c$\end{tabular}}}}%
    \put(0.2978426,0.45654363){\color[rgb]{0,0,0}\makebox(0,0)[lt]{\lineheight{1.25}\smash{\begin{tabular}[t]{l}$\textcolor{red}{\gamma}$\end{tabular}}}}%
    \put(0,0){\includegraphics[width=\unitlength,page=2]{one-holed_torus_flat2.pdf}}%
    \put(0.23309451,0.23231517){\color[rgb]{0,0,0}\rotatebox{31.291276}{\makebox(0,0)[lt]{\lineheight{1.25}\smash{\begin{tabular}[t]{l}$\mbox{hole}$\end{tabular}}}}}%
    \put(0,0){\includegraphics[width=\unitlength,page=3]{one-holed_torus_flat2.pdf}}%
    \put(0.41364684,0.41494959){\color[rgb]{0,0,0}\makebox(0,0)[lt]{\lineheight{1.25}\smash{\begin{tabular}[t]{l}$\textcolor{red}{\circled{1}}$\end{tabular}}}}%
    \put(0.45984402,0.1569966){\color[rgb]{0,0,0}\makebox(0,0)[lt]{\lineheight{1.25}\smash{\begin{tabular}[t]{l}$\textcolor{red}{\circled{2}}$\end{tabular}}}}%
    \put(0.10840703,0.42235704){\color[rgb]{0,0,0}\makebox(0,0)[lt]{\lineheight{1.25}\smash{\begin{tabular}[t]{l}$\textcolor{red}{\circled{8}}$\end{tabular}}}}%
    \put(0.35849493,0.72971697){\color[rgb]{0,0,0}\makebox(0,0)[lt]{\lineheight{1.25}\smash{\begin{tabular}[t]{l}$\textcolor{blue}{\circled{$2'$}}$\end{tabular}}}}%
    \put(0.8033078,0.19899606){\color[rgb]{0,0,0}\makebox(0,0)[lt]{\lineheight{1.25}\smash{\begin{tabular}[t]{l}$\textcolor{red}{\circled{7}}$\end{tabular}}}}%
    \put(0.57421961,0.83471575){\color[rgb]{0,0,0}\makebox(0,0)[lt]{\lineheight{1.25}\smash{\begin{tabular}[t]{l}$\textcolor{red}{\circled{6}}$\end{tabular}}}}%
    \put(0.7173998,0.74498953){\color[rgb]{0,0,0}\makebox(0,0)[lt]{\lineheight{1.25}\smash{\begin{tabular}[t]{l}$\textcolor{red}{\circled{5}}$\end{tabular}}}}%
    \put(0.82430753,0.61708197){\color[rgb]{0,0,0}\makebox(0,0)[lt]{\lineheight{1.25}\smash{\begin{tabular}[t]{l}$\textcolor{red}{\circled{4}}$\end{tabular}}}}%
    \put(0.20958764,0.79080715){\color[rgb]{0,0,0}\makebox(0,0)[lt]{\lineheight{1.25}\smash{\begin{tabular}[t]{l}$\textcolor{red}{\circled{3}}$\end{tabular}}}}%
    \put(0.71204455,0.35172149){\color[rgb]{0,0,0}\makebox(0,0)[lt]{\lineheight{1.25}\smash{\begin{tabular}[t]{l}$\textcolor{blue}{\circled{$1'$}}$\end{tabular}}}}%
    \put(0.22083097,0.60776546){\color[rgb]{0,0,0}\makebox(0,0)[lt]{\lineheight{1.25}\smash{\begin{tabular}[t]{l}$\textcolor{blue}{\eta}$\end{tabular}}}}%
    \put(0,0){\includegraphics[width=\unitlength,page=4]{one-holed_torus_flat2.pdf}}%
    \put(0.00973629,0.56520931){\color[rgb]{0,0,0}\makebox(0,0)[lt]{\lineheight{1.25}\smash{\begin{tabular}[t]{l}$b$\end{tabular}}}}%
    \put(0.95307819,0.5537549){\color[rgb]{0,0,0}\makebox(0,0)[lt]{\lineheight{1.25}\smash{\begin{tabular}[t]{l}$b$\end{tabular}}}}%
    \put(0.45099334,0.63011762){\color[rgb]{0,0,0}\makebox(0,0)[lt]{\lineheight{1.25}\smash{\begin{tabular}[t]{l}$d$\end{tabular}}}}%
    \put(0,0){\includegraphics[width=\unitlength,page=5]{one-holed_torus_flat2.pdf}}%
  \end{picture}%
\endgroup%

%% file: one-holed_surface2mod2.pdf_tex
\begingroup%
  \makeatletter%
  \providecommand\color[2][]{%
    \errmessage{(Inkscape) Color is used for the text in Inkscape, but the package 'color.sty' is not loaded}%
    \renewcommand\color[2][]{}%
  }%
  \providecommand\transparent[1]{%
    \errmessage{(Inkscape) Transparency is used (non-zero) for the text in Inkscape, but the package 'transparent.sty' is not loaded}%
    \renewcommand\transparent[1]{}%
  }%
  \providecommand\rotatebox[2]{#2}%
  \newcommand*\fsize{\dimexpr\f@size pt\relax}%
  \newcommand*\lineheight[1]{\fontsize{\fsize}{#1\fsize}\selectfont}%
  \ifx\svgwidth\undefined%
    \setlength{\unitlength}{276.4190001bp}%
    \ifx\svgscale\undefined%
      \relax%
    \else%
      \setlength{\unitlength}{\unitlength * \real{\svgscale}}%
    \fi%
  \else%
    \setlength{\unitlength}{\svgwidth}%
  \fi%
  \global\let\svgwidth\undefined%
  \global\let\svgscale\undefined%
  \makeatother%
  \begin{picture}(1,0.63043005)%
    \lineheight{1}%
    \setlength\tabcolsep{0pt}%
    \put(0,0){\includegraphics[width=\unitlength,page=1]{one-holed_surface2mod2.pdf}}%
    \put(0.81455175,0.1686768){\color[rgb]{0,0,0}\makebox(0,0)[lt]{\lineheight{1.25}\smash{\begin{tabular}[t]{l}$\textcolor{red}{\gamma}$\end{tabular}}}}%
    \put(0.78597276,0.23071397){\color[rgb]{0,0,0}\makebox(0,0)[lt]{\lineheight{1.25}\smash{\begin{tabular}[t]{l}$\textcolor{blue}{\eta}$\end{tabular}}}}%
    \put(0.6199298,0.20463169){\color[rgb]{0,0,0}\makebox(0,0)[lt]{\lineheight{1.25}\smash{\begin{tabular}[t]{l}$\textcolor{yellowi}{\zeta}$\end{tabular}}}}%
    \put(0.9107719,0.21799637){\color[rgb]{0,0,0}\makebox(0,0)[lt]{\lineheight{1.25}\smash{\begin{tabular}[t]{l}$\textcolor{realpurple}{\sigma}$\end{tabular}}}}%
    \put(0,0){\includegraphics[width=\unitlength,page=2]{one-holed_surface2mod2.pdf}}%
    \put(0.41668107,0.60174683){\color[rgb]{0,0,0}\makebox(0,0)[lt]{\lineheight{1.25}\smash{\begin{tabular}[t]{l}$\mbox{$g-2$ handles}$\end{tabular}}}}%
    \put(0.80455137,0.37430046){\color[rgb]{0,0,0}\makebox(0,0)[lt]{\lineheight{1.25}\smash{\begin{tabular}[t]{l}$\mbox{first handle}$\end{tabular}}}}%
    \put(0.72409113,0.4931594){\color[rgb]{0,0,0}\makebox(0,0)[lt]{\lineheight{1.25}\smash{\begin{tabular}[t]{l}$\mbox{second handle}$\end{tabular}}}}%
    \put(0.01305109,0.42081373){\color[rgb]{0,0,0}\makebox(0,0)[lt]{\lineheight{1.25}\smash{\begin{tabular}[t]{l}$\mbox{last handle}$\end{tabular}}}}%
    \put(0.01063679,0.46954754){\color[rgb]{0,0,0}\makebox(0,0)[lt]{\lineheight{1.25}\smash{\begin{tabular}[t]{l}$\mbox{($g$-th)}$\end{tabular}}}}%
    \put(0,0){\includegraphics[width=\unitlength,page=3]{one-holed_surface2mod2.pdf}}%
    \put(0.47003904,0.4067789){\color[rgb]{0,0,0}\makebox(0,0)[lt]{\lineheight{1.25}\smash{\begin{tabular}[t]{l}...\end{tabular}}}}%
    \put(0.46866861,0.2432972){\color[rgb]{0,0,0}\makebox(0,0)[lt]{\lineheight{1.25}\smash{\begin{tabular}[t]{l}...\end{tabular}}}}%
    \put(0,0){\includegraphics[width=\unitlength,page=4]{one-holed_surface2mod2.pdf}}%
    \put(0.4337902,0.11625112){\color[rgb]{0,0,0}\makebox(0,0)[lt]{\lineheight{1.25}\smash{\begin{tabular}[t]{l}$\mbox{sphere-body}$\end{tabular}}}}%
    \put(0,0){\includegraphics[width=\unitlength,page=5]{one-holed_surface2mod2.pdf}}%
  \end{picture}%
\endgroup%

%% file: one-holed_surface_flat2mod2.pdf_tex
\begingroup%
  \makeatletter%
  \providecommand\color[2][]{%
    \errmessage{(Inkscape) Color is used for the text in Inkscape, but the package 'color.sty' is not loaded}%
    \renewcommand\color[2][]{}%
  }%
  \providecommand\transparent[1]{%
    \errmessage{(Inkscape) Transparency is used (non-zero) for the text in Inkscape, but the package 'transparent.sty' is not loaded}%
    \renewcommand\transparent[1]{}%
  }%
  \providecommand\rotatebox[2]{#2}%
  \newcommand*\fsize{\dimexpr\f@size pt\relax}%
  \newcommand*\lineheight[1]{\fontsize{\fsize}{#1\fsize}\selectfont}%
  \ifx\svgwidth\undefined%
    \setlength{\unitlength}{430.86614173bp}%
    \ifx\svgscale\undefined%
      \relax%
    \else%
      \setlength{\unitlength}{\unitlength * \real{\svgscale}}%
    \fi%
  \else%
    \setlength{\unitlength}{\svgwidth}%
  \fi%
  \global\let\svgwidth\undefined%
  \global\let\svgscale\undefined%
  \makeatother%
  \begin{picture}(1,0.92105263)%
    \lineheight{1}%
    \setlength\tabcolsep{0pt}%
    \put(0,0){\includegraphics[width=\unitlength,page=1]{one-holed_surface_flat2mod2.pdf}}%
    \put(0.28815974,0.12464834){\color[rgb]{0,0,0}\rotatebox{8.9006151}{\makebox(0,0)[lt]{\lineheight{1.25}\smash{\begin{tabular}[t]{l}$\mbox{hole}$\end{tabular}}}}}%
    \put(0,0){\includegraphics[width=\unitlength,page=2]{one-holed_surface_flat2mod2.pdf}}%
    \put(0.34174624,0.19936808){\color[rgb]{0,0,0}\makebox(0,0)[lt]{\lineheight{1.25}\smash{\begin{tabular}[t]{l}$\textcolor{red}{\circled{1}}$\end{tabular}}}}%
    \put(0.4030431,0.12612369){\color[rgb]{0,0,0}\makebox(0,0)[lt]{\lineheight{1.25}\smash{\begin{tabular}[t]{l}$\textcolor{red}{\circled{2}}$\end{tabular}}}}%
    \put(0.38542057,0.08589202){\color[rgb]{0,0,0}\makebox(0,0)[lt]{\lineheight{1.25}\smash{\begin{tabular}[t]{l}$\textcolor{red}{\circled{3}}$\end{tabular}}}}%
    \put(0.93490664,0.45755803){\color[rgb]{0,0,0}\makebox(0,0)[lt]{\lineheight{1.25}\smash{\begin{tabular}[t]{l}$\textcolor{red}{\circled{4}}$\end{tabular}}}}%
    \put(0.75781977,0.14525828){\color[rgb]{0,0,0}\makebox(0,0)[lt]{\lineheight{1.25}\smash{\begin{tabular}[t]{l}$\textcolor{red}{\circled{5}}$\end{tabular}}}}%
    \put(0.76030642,0.18318025){\color[rgb]{0,0,0}\makebox(0,0)[lt]{\lineheight{1.25}\smash{\begin{tabular}[t]{l}$\textcolor{red}{\circled{6}}$\end{tabular}}}}%
    \put(0.80133676,0.23415725){\color[rgb]{0,0,0}\makebox(0,0)[lt]{\lineheight{1.25}\smash{\begin{tabular}[t]{l}$\textcolor{red}{\circled{7}}$\end{tabular}}}}%
    \put(0.88277566,0.2577808){\color[rgb]{0,0,0}\makebox(0,0)[lt]{\lineheight{1.25}\smash{\begin{tabular}[t]{l}$\textcolor{red}{\circled{8}}$\end{tabular}}}}%
    \put(0.56622072,0.04653693){\color[rgb]{0,0,0}\makebox(0,0)[lt]{\lineheight{1.25}\smash{\begin{tabular}[t]{l}$\textcolor{red}{\circled{9}}$\end{tabular}}}}%
    \put(0.82620357,0.54312794){\color[rgb]{0,0,0}\makebox(0,0)[lt]{\lineheight{1.25}\smash{\begin{tabular}[t]{l}$\textcolor{red}{\circledd{10}}$\end{tabular}}}}%
    \put(0.79424578,0.55120968){\color[rgb]{0,0,0}\makebox(0,0)[lt]{\lineheight{1.25}\smash{\begin{tabular}[t]{l}$\textcolor{red}{\circledd{11}}$\end{tabular}}}}%
    \put(0.77994736,0.58477992){\color[rgb]{0,0,0}\makebox(0,0)[lt]{\lineheight{1.25}\smash{\begin{tabular}[t]{l}$\textcolor{red}{\circledd{12}}$\end{tabular}}}}%
    \put(0.27452856,0.20556047){\color[rgb]{0,0,0}\makebox(0,0)[lt]{\lineheight{1.25}\smash{\begin{tabular}[t]{l}$\textcolor{red}{\circled{$N$}}$\end{tabular}}}}%
    \put(0.2122547,0.19754064){\color[rgb]{0,0,0}\makebox(0,0)[lt]{\lineheight{1.25}\smash{\begin{tabular}[t]{l}$\textcolor{red}{\circleddd{$N\hspace{-1mm}-\hspace{-1mm}1$}}$\end{tabular}}}}%
    \put(0.65027061,0.33976587){\color[rgb]{0,0,0}\makebox(0,0)[lt]{\lineheight{1.25}\smash{\begin{tabular}[t]{l}$\textcolor{blue}{\circled{$7'$}}$\end{tabular}}}}%
    \put(0.87308533,0.32554301){\color[rgb]{0,0,0}\makebox(0,0)[lt]{\lineheight{1.25}\smash{\begin{tabular}[t]{l}$\textcolor{blue}{\circled{8}}$\end{tabular}}}}%
    \put(0.6065327,0.06828311){\color[rgb]{0,0,0}\makebox(0,0)[lt]{\lineheight{1.25}\smash{\begin{tabular}[t]{l}$\textcolor{blue}{\circled{9}}$\end{tabular}}}}%
    \put(0.8538454,0.50788416){\color[rgb]{0,0,0}\makebox(0,0)[lt]{\lineheight{1.25}\smash{\begin{tabular}[t]{l}$\textcolor{blue}{\circledd{10}}$\end{tabular}}}}%
    \put(0.77303492,0.50304869){\color[rgb]{0,0,0}\makebox(0,0)[lt]{\lineheight{1.25}\smash{\begin{tabular}[t]{l}$\textcolor{blue}{\circledd{11}}$\end{tabular}}}}%
    \put(0.71846755,0.56480114){\color[rgb]{0,0,0}\makebox(0,0)[lt]{\lineheight{1.25}\smash{\begin{tabular}[t]{l}$\textcolor{blue}{\circledd{12}}$\end{tabular}}}}%
    \put(0.76906108,0.65265989){\color[rgb]{0,0,0}\makebox(0,0)[lt]{\lineheight{1.25}\smash{\begin{tabular}[t]{l}$\textcolor{blue}{\circledd{13}}$\end{tabular}}}}%
    \put(0.54768949,0.45748766){\color[rgb]{0,0,0}\makebox(0,0)[lt]{\lineheight{1.25}\smash{\begin{tabular}[t]{l}$\textcolor{yellowi}{\circledd{$12'$}}$\end{tabular}}}}%
    \put(0.33120867,0.26756617){\color[rgb]{0,0,0}\makebox(0,0)[lt]{\lineheight{1.25}\smash{\begin{tabular}[t]{l}$\textcolor{blue}{\circled{$N$}}$\end{tabular}}}}%
    \put(0.23917626,0.25867304){\color[rgb]{0,0,0}\makebox(0,0)[lt]{\lineheight{1.25}\smash{\begin{tabular}[t]{l}$\textcolor{blue}{\circleddd{$N\hspace{-1mm}-\hspace{-1mm}1$}}$\end{tabular}}}}%
    \put(0.16217772,0.24440761){\color[rgb]{0,0,0}\makebox(0,0)[lt]{\lineheight{1.25}\smash{\begin{tabular}[t]{l}$\textcolor{blue}{\circleddd{$N\hspace{-1mm}-\hspace{-1mm}2$}}$\end{tabular}}}}%
    \put(0,0){\includegraphics[width=\unitlength,page=3]{one-holed_surface_flat2mod2.pdf}}%
    \put(0.3877552,0.17916471){\color[rgb]{0,0,0}\makebox(0,0)[lt]{\lineheight{1.25}\smash{\begin{tabular}[t]{l}$\textcolor{red}{\gamma}$\end{tabular}}}}%
    \put(0.6355768,0.11624657){\color[rgb]{0,0,0}\makebox(0,0)[lt]{\lineheight{1.25}\smash{\begin{tabular}[t]{l}$\textcolor{realpurple}{\circled{$1'$}}$\end{tabular}}}}%
    \put(0.88258319,0.3861989){\color[rgb]{0,0,0}\makebox(0,0)[lt]{\lineheight{1.25}\smash{\begin{tabular}[t]{l}$\textcolor{realpurple}{\circled{$2'$}}$\end{tabular}}}}%
    \put(0.58620993,0.32110598){\color[rgb]{0,0,0}\makebox(0,0)[lt]{\lineheight{1.25}\smash{\begin{tabular}[t]{l}$\textcolor{blue}{\eta}$\end{tabular}}}}%
    \put(0.62881655,0.2501933){\color[rgb]{0,0,0}\makebox(0,0)[lt]{\lineheight{1.25}\smash{\begin{tabular}[t]{l}{\small $\textcolor{blue}{(\eta \, \circled{$7'$})}$}\end{tabular}}}}%
    \put(0.63238356,0.53716748){\color[rgb]{0,0,0}\makebox(0,0)[lt]{\lineheight{1.25}\smash{\begin{tabular}[t]{l}$\textcolor{yellowi}{\zeta}$\end{tabular}}}}%
    \put(0.87522316,0.4123273){\color[rgb]{0,0,0}\makebox(0,0)[lt]{\lineheight{1.25}\smash{\begin{tabular}[t]{l}$\textcolor{realpurple}{\sigma}$\end{tabular}}}}%
    \put(0,0){\includegraphics[width=\unitlength,page=4]{one-holed_surface_flat2mod2.pdf}}%
    \put(0.38243463,0.05293277){\color[rgb]{0,0,0}\makebox(0,0)[lt]{\lineheight{1.25}\smash{\begin{tabular}[t]{l}$a_1$\end{tabular}}}}%
    \put(0,0){\includegraphics[width=\unitlength,page=5]{one-holed_surface_flat2mod2.pdf}}%
    \put(0.73734321,0.08983635){\color[rgb]{0,0,0}\makebox(0,0)[lt]{\lineheight{1.25}\smash{\begin{tabular}[t]{l}$b_1$\end{tabular}}}}%
    \put(0,0){\includegraphics[width=\unitlength,page=6]{one-holed_surface_flat2mod2.pdf}}%
    \put(0.93870829,0.32300757){\color[rgb]{0,0,0}\makebox(0,0)[lt]{\lineheight{1.25}\smash{\begin{tabular}[t]{l}$a_1$\end{tabular}}}}%
    \put(0,0){\includegraphics[width=\unitlength,page=7]{one-holed_surface_flat2mod2.pdf}}%
    \put(0.89768189,0.53687843){\color[rgb]{0,0,0}\makebox(0,0)[lt]{\lineheight{1.25}\smash{\begin{tabular}[t]{l}$b_1$\end{tabular}}}}%
    \put(0,0){\includegraphics[width=\unitlength,page=8]{one-holed_surface_flat2mod2.pdf}}%
    \put(0.85067184,0.62527363){\color[rgb]{0,0,0}\makebox(0,0)[lt]{\lineheight{1.25}\smash{\begin{tabular}[t]{l}$a_2$\end{tabular}}}}%
    \put(0,0){\includegraphics[width=\unitlength,page=9]{one-holed_surface_flat2mod2.pdf}}%
    \put(0.84757424,0.72771722){\color[rgb]{0,0,0}\makebox(0,0)[lt]{\lineheight{1.25}\smash{\begin{tabular}[t]{l}$b_2$\end{tabular}}}}%
    \put(0,0){\includegraphics[width=\unitlength,page=10]{one-holed_surface_flat2mod2.pdf}}%
    \put(0.64952493,0.80799357){\color[rgb]{0,0,0}\makebox(0,0)[lt]{\lineheight{1.25}\smash{\begin{tabular}[t]{l}$b_2$\end{tabular}}}}%
    \put(0,0){\includegraphics[width=\unitlength,page=11]{one-holed_surface_flat2mod2.pdf}}%
    \put(0.75502607,0.80843317){\color[rgb]{0,0,0}\makebox(0,0)[lt]{\lineheight{1.25}\smash{\begin{tabular}[t]{l}$a_2$\end{tabular}}}}%
    \put(0,0){\includegraphics[width=\unitlength,page=12]{one-holed_surface_flat2mod2.pdf}}%
    \put(0.59197421,0.8053914){\color[rgb]{0,0,0}\makebox(0,0)[lt]{\lineheight{1.25}\smash{\begin{tabular}[t]{l}$a_{g-1}$\end{tabular}}}}%
    \put(0,0){\includegraphics[width=\unitlength,page=13]{one-holed_surface_flat2mod2.pdf}}%
    \put(0.53130871,0.85702135){\color[rgb]{0,0,0}\makebox(0,0)[lt]{\lineheight{1.25}\smash{\begin{tabular}[t]{l}$b_{g-1}$\end{tabular}}}}%
    \put(0,0){\includegraphics[width=\unitlength,page=14]{one-holed_surface_flat2mod2.pdf}}%
    \put(0.36395975,0.84248347){\color[rgb]{0,0,0}\makebox(0,0)[lt]{\lineheight{1.25}\smash{\begin{tabular}[t]{l}$a_{g-1}$\end{tabular}}}}%
    \put(0,0){\includegraphics[width=\unitlength,page=15]{one-holed_surface_flat2mod2.pdf}}%
    \put(0.25253061,0.7783856){\color[rgb]{0,0,0}\makebox(0,0)[lt]{\lineheight{1.25}\smash{\begin{tabular}[t]{l}$b_{g-1}$\end{tabular}}}}%
    \put(0,0){\includegraphics[width=\unitlength,page=16]{one-holed_surface_flat2mod2.pdf}}%
    \put(0.18928178,0.70515202){\color[rgb]{0,0,0}\makebox(0,0)[lt]{\lineheight{1.25}\smash{\begin{tabular}[t]{l}$a_g$\end{tabular}}}}%
    \put(0,0){\includegraphics[width=\unitlength,page=17]{one-holed_surface_flat2mod2.pdf}}%
    \put(0.05509832,0.60281259){\color[rgb]{0,0,0}\makebox(0,0)[lt]{\lineheight{1.25}\smash{\begin{tabular}[t]{l}$b_g$\end{tabular}}}}%
    \put(0,0){\includegraphics[width=\unitlength,page=18]{one-holed_surface_flat2mod2.pdf}}%
    \put(0.01172617,0.36893933){\color[rgb]{0,0,0}\makebox(0,0)[lt]{\lineheight{1.25}\smash{\begin{tabular}[t]{l}$a_g$\end{tabular}}}}%
    \put(0,0){\includegraphics[width=\unitlength,page=19]{one-holed_surface_flat2mod2.pdf}}%
    \put(0.09800926,0.15800117){\color[rgb]{0,0,0}\makebox(0,0)[lt]{\lineheight{1.25}\smash{\begin{tabular}[t]{l}$b_g$\end{tabular}}}}%
    \put(0.18063884,0.08907896){\color[rgb]{0,0,0}\makebox(0,0)[lt]{\lineheight{1.25}\smash{\begin{tabular}[t]{l}{\small $p_1$}\end{tabular}}}}%
    \put(0.5720129,0.02097217){\color[rgb]{0,0,0}\makebox(0,0)[lt]{\lineheight{1.25}\smash{\begin{tabular}[t]{l}{\small $p_2$}\end{tabular}}}}%
    \put(0.88461883,0.16470045){\color[rgb]{0,0,0}\makebox(0,0)[lt]{\lineheight{1.25}\smash{\begin{tabular}[t]{l}{\small $p_3$}\end{tabular}}}}%
    \put(0.97366979,0.48245853){\color[rgb]{0,0,0}\makebox(0,0)[lt]{\lineheight{1.25}\smash{\begin{tabular}[t]{l}{\small $p_4$}\end{tabular}}}}%
    \put(0.86877124,0.68947065){\color[rgb]{0,0,0}\makebox(0,0)[lt]{\lineheight{1.25}\smash{\begin{tabular}[t]{l}{\small $p_6$}\end{tabular}}}}%
    \put(0.81042478,0.78201893){\color[rgb]{0,0,0}\makebox(0,0)[lt]{\lineheight{1.25}\smash{\begin{tabular}[t]{l}{\small $p_7$}\end{tabular}}}}%
    \put(0.69706342,0.83548734){\color[rgb]{0,0,0}\makebox(0,0)[lt]{\lineheight{1.25}\smash{\begin{tabular}[t]{l}{\small $p_8$}\end{tabular}}}}%
    \put(0.42132323,0.87013378){\color[rgb]{0,0,0}\makebox(0,0)[lt]{\lineheight{1.25}\smash{\begin{tabular}[t]{l}{\small $p_{4(g-2)+3}$}\end{tabular}}}}%
    \put(0.08640805,0.65421785){\color[rgb]{0,0,0}\rotatebox{42.260456}{\makebox(0,0)[lt]{\lineheight{1.25}\smash{\begin{tabular}[t]{l}{\small $p_{4(g-1)+2}$}\end{tabular}}}}}%
    \put(0.82375219,0.58965326){\color[rgb]{0,0,0}\makebox(0,0)[lt]{\lineheight{1.25}\smash{\begin{tabular}[t]{l}{\small $p_5$}\end{tabular}}}}%
    \put(0,0){\includegraphics[width=\unitlength,page=20]{one-holed_surface_flat2mod2.pdf}}%
    \put(0.35314214,0.90428353){\color[rgb]{0,0,0}\makebox(0,0)[lt]{\lineheight{1.25}\smash{\begin{tabular}[t]{l}$\mbox{$(g-1)$-th handle}$\end{tabular}}}}%
    \put(0,0){\includegraphics[width=\unitlength,page=21]{one-holed_surface_flat2mod2.pdf}}%
    \put(0.79979015,0.85646492){\color[rgb]{0,0,0}\rotatebox{-40.165541}{\makebox(0,0)[lt]{\lineheight{1.25}\smash{\begin{tabular}[t]{l}$\mbox{$2$nd handle}$\end{tabular}}}}}%
    \put(0.64592455,0.88778026){\color[rgb]{0,0,0}\rotatebox{-10.889808}{\makebox(0,0)[lt]{\lineheight{1.25}\smash{\begin{tabular}[t]{l}$\cdots$\end{tabular}}}}}%
    \put(0,0){\includegraphics[width=\unitlength,page=22]{one-holed_surface_flat2mod2.pdf}}%
    \put(0.90392699,0.07309379){\color[rgb]{0,0,0}\rotatebox{41.456037}{\makebox(0,0)[lt]{\lineheight{1.25}\smash{\begin{tabular}[t]{l}$\mbox{first handle}$\end{tabular}}}}}%
    \put(0,0){\includegraphics[width=\unitlength,page=23]{one-holed_surface_flat2mod2.pdf}}%
    \put(0.03636024,0.69367093){\color[rgb]{0,0,0}\rotatebox{25.322429}{\makebox(0,0)[lt]{\lineheight{1.25}\smash{\begin{tabular}[t]{l}$\mbox{last handle}$\end{tabular}}}}}%
    \put(0.41969569,0.53599774){\color[rgb]{0,0,0}\makebox(0,0)[lt]{\lineheight{1.25}\smash{\begin{tabular}[t]{l}$\mbox{sphere-body}$\end{tabular}}}}%
    \put(0,0){\includegraphics[width=\unitlength,page=24]{one-holed_surface_flat2mod2.pdf}}%
    \put(0.78198956,0.67899849){\color[rgb]{0,0,0}\makebox(0,0)[lt]{\lineheight{1.25}\smash{\begin{tabular}[t]{l}$c_2$\end{tabular}}}}%
    \put(0.6949559,0.74307999){\color[rgb]{0,0,0}\makebox(0,0)[lt]{\lineheight{1.25}\smash{\begin{tabular}[t]{l}$d_2$\end{tabular}}}}%
    \put(0.71311487,0.6702513){\color[rgb]{0,0,0}\makebox(0,0)[lt]{\lineheight{1.25}\smash{\begin{tabular}[t]{l}$e_2$\end{tabular}}}}%
    \put(0.51561772,0.80380845){\color[rgb]{0,0,0}\rotatebox{-36.532084}{\makebox(0,0)[lt]{\lineheight{1.25}\smash{\begin{tabular}[t]{l}$c_{g-1}$\end{tabular}}}}}%
    \put(0.44832924,0.74566976){\color[rgb]{0,0,0}\makebox(0,0)[lt]{\lineheight{1.25}\smash{\begin{tabular}[t]{l}$e_{g-1}$\end{tabular}}}}%
    \put(0.37228975,0.76256871){\color[rgb]{0,0,0}\makebox(0,0)[lt]{\lineheight{1.25}\smash{\begin{tabular}[t]{l}$d_{g-1}$\end{tabular}}}}%
    \put(0.15064221,0.60087801){\color[rgb]{0,0,0}\makebox(0,0)[lt]{\lineheight{1.25}\smash{\begin{tabular}[t]{l}$c_g$\end{tabular}}}}%
    \put(0.12124741,0.32882646){\color[rgb]{0,0,0}\makebox(0,0)[lt]{\lineheight{1.25}\smash{\begin{tabular}[t]{l}$d_g$\end{tabular}}}}%
    \put(0.20842113,0.45099846){\color[rgb]{0,0,0}\makebox(0,0)[lt]{\lineheight{1.25}\smash{\begin{tabular}[t]{l}$e_g$\end{tabular}}}}%
    \put(0.55278492,0.14519376){\color[rgb]{0,0,0}\makebox(0,0)[lt]{\lineheight{1.25}\smash{\begin{tabular}[t]{l}$c_1$\end{tabular}}}}%
    \put(0.8250012,0.3844601){\color[rgb]{0,0,0}\makebox(0,0)[lt]{\lineheight{1.25}\smash{\begin{tabular}[t]{l}$d_1$\end{tabular}}}}%
    \put(0.55296257,0.23440092){\color[rgb]{0,0,0}\makebox(0,0)[lt]{\lineheight{1.25}\smash{\begin{tabular}[t]{l}$c_1'$\end{tabular}}}}%
    \put(0.58549745,0.3997316){\color[rgb]{0,0,0}\makebox(0,0)[lt]{\lineheight{1.25}\smash{\begin{tabular}[t]{l}$f_2$\end{tabular}}}}%
    \put(0.36108524,0.45035016){\color[rgb]{0,0,0}\makebox(0,0)[lt]{\lineheight{1.25}\smash{\begin{tabular}[t]{l}$f_{g-1}$\end{tabular}}}}%
    \put(0.42057222,0.48898135){\color[rgb]{0,0,0}\rotatebox{-1.0250199}{\makebox(0,0)[lt]{\lineheight{1.25}\smash{\begin{tabular}[t]{l}$\cdots$\end{tabular}}}}}%
  \end{picture}%
\endgroup%

%% file: one-holed_surface3.pdf_tex
\begingroup%
  \makeatletter%
  \providecommand\color[2][]{%
    \errmessage{(Inkscape) Color is used for the text in Inkscape, but the package 'color.sty' is not loaded}%
    \renewcommand\color[2][]{}%
  }%
  \providecommand\transparent[1]{%
    \errmessage{(Inkscape) Transparency is used (non-zero) for the text in Inkscape, but the package 'transparent.sty' is not loaded}%
    \renewcommand\transparent[1]{}%
  }%
  \providecommand\rotatebox[2]{#2}%
  \newcommand*\fsize{\dimexpr\f@size pt\relax}%
  \newcommand*\lineheight[1]{\fontsize{\fsize}{#1\fsize}\selectfont}%
  \ifx\svgwidth\undefined%
    \setlength{\unitlength}{276.4190001bp}%
    \ifx\svgscale\undefined%
      \relax%
    \else%
      \setlength{\unitlength}{\unitlength * \real{\svgscale}}%
    \fi%
  \else%
    \setlength{\unitlength}{\svgwidth}%
  \fi%
  \global\let\svgwidth\undefined%
  \global\let\svgscale\undefined%
  \makeatother%
  \begin{picture}(1,0.63043005)%
    \lineheight{1}%
    \setlength\tabcolsep{0pt}%
    \put(0,0){\includegraphics[width=\unitlength,page=1]{one-holed_surface3.pdf}}%
    \put(0.81455175,0.1686768){\color[rgb]{0,0,0}\makebox(0,0)[lt]{\lineheight{1.25}\smash{\begin{tabular}[t]{l}$\textcolor{red}{\gamma}$\end{tabular}}}}%
    \put(0.50387229,0.1880013){\color[rgb]{0,0,0}\makebox(0,0)[lt]{\lineheight{1.25}\smash{\begin{tabular}[t]{l}$\textcolor{blue}{\vartheta_i}$\end{tabular}}}}%
    \put(0.57301717,0.16633715){\color[rgb]{0,0,0}\makebox(0,0)[lt]{\lineheight{1.25}\smash{\begin{tabular}[t]{l}$\textcolor{yellowi}{\zeta=\zeta_2}$\end{tabular}}}}%
    \put(0.40221933,0.19146727){\color[rgb]{0,0,0}\makebox(0,0)[lt]{\lineheight{1.25}\smash{\begin{tabular}[t]{l}$\textcolor{yellowi}{\zeta_{i}}$\end{tabular}}}}%
    \put(0,0){\includegraphics[width=\unitlength,page=2]{one-holed_surface3.pdf}}%
    \put(0.18487666,0.22154507){\color[rgb]{0,0,0}\makebox(0,0)[lt]{\lineheight{1.25}\smash{\begin{tabular}[t]{l}$\textcolor{yellowi}{\zeta_{g}}$\end{tabular}}}}%
    \put(0,0){\includegraphics[width=\unitlength,page=3]{one-holed_surface3.pdf}}%
    \put(0.80455137,0.37430046){\color[rgb]{0,0,0}\makebox(0,0)[lt]{\lineheight{1.25}\smash{\begin{tabular}[t]{l}$\mbox{first handle}$\end{tabular}}}}%
    \put(0.72409113,0.4931594){\color[rgb]{0,0,0}\makebox(0,0)[lt]{\lineheight{1.25}\smash{\begin{tabular}[t]{l}$\mbox{second handle}$\end{tabular}}}}%
    \put(0.01305109,0.42081373){\color[rgb]{0,0,0}\makebox(0,0)[lt]{\lineheight{1.25}\smash{\begin{tabular}[t]{l}$\mbox{last handle}$\end{tabular}}}}%
    \put(0.29677896,0.52129511){\color[rgb]{0,0,0}\makebox(0,0)[lt]{\lineheight{1.25}\smash{\begin{tabular}[t]{l}$\mbox{$i$-th handle}$\end{tabular}}}}%
    \put(0.01063679,0.46954754){\color[rgb]{0,0,0}\makebox(0,0)[lt]{\lineheight{1.25}\smash{\begin{tabular}[t]{l}$\mbox{($g$-th)}$\end{tabular}}}}%
    \put(0,0){\includegraphics[width=\unitlength,page=4]{one-holed_surface3.pdf}}%
    \put(0.50375725,0.3955395){\color[rgb]{0,0,0}\makebox(0,0)[lt]{\lineheight{1.25}\smash{\begin{tabular}[t]{l}...\end{tabular}}}}%
    \put(0.23368975,0.33188096){\color[rgb]{0,0,0}\rotatebox{18.255822}{\makebox(0,0)[lt]{\lineheight{1.25}\smash{\begin{tabular}[t]{l}...\end{tabular}}}}}%
    \put(0.28424817,0.2488646){\color[rgb]{0,0,0}\rotatebox{9.3581074}{\makebox(0,0)[lt]{\lineheight{1.25}\smash{\begin{tabular}[t]{l}...\end{tabular}}}}}%
    \put(0.53433184,0.26622212){\color[rgb]{0,0,0}\rotatebox{2.2018796}{\makebox(0,0)[lt]{\lineheight{1.25}\smash{\begin{tabular}[t]{l}...\end{tabular}}}}}%
    \put(0,0){\includegraphics[width=\unitlength,page=5]{one-holed_surface3.pdf}}%
  \end{picture}%
\endgroup%

%% file: one-holed_surface_flat5more.pdf_tex
\begingroup%
  \makeatletter%
  \providecommand\color[2][]{%
    \errmessage{(Inkscape) Color is used for the text in Inkscape, but the package 'color.sty' is not loaded}%
    \renewcommand\color[2][]{}%
  }%
  \providecommand\transparent[1]{%
    \errmessage{(Inkscape) Transparency is used (non-zero) for the text in Inkscape, but the package 'transparent.sty' is not loaded}%
    \renewcommand\transparent[1]{}%
  }%
  \providecommand\rotatebox[2]{#2}%
  \newcommand*\fsize{\dimexpr\f@size pt\relax}%
  \newcommand*\lineheight[1]{\fontsize{\fsize}{#1\fsize}\selectfont}%
  \ifx\svgwidth\undefined%
    \setlength{\unitlength}{515.90551181bp}%
    \ifx\svgscale\undefined%
      \relax%
    \else%
      \setlength{\unitlength}{\unitlength * \real{\svgscale}}%
    \fi%
  \else%
    \setlength{\unitlength}{\svgwidth}%
  \fi%
  \global\let\svgwidth\undefined%
  \global\let\svgscale\undefined%
  \makeatother%
  \begin{picture}(1,0.88461538)%
    \lineheight{1}%
    \setlength\tabcolsep{0pt}%
    \put(0,0){\includegraphics[width=\unitlength,page=1]{one-holed_surface_flat5more.pdf}}%
    \put(0.37358691,0.25092622){\color[rgb]{0,0,0}\rotatebox{33.019654}{\makebox(0,0)[lt]{\lineheight{1.25}\smash{\begin{tabular}[t]{l}$\mbox{hole}$\end{tabular}}}}}%
    \put(0,0){\includegraphics[width=\unitlength,page=2]{one-holed_surface_flat5more.pdf}}%
    \put(0.44544777,0.29467263){\color[rgb]{0,0,0}\makebox(0,0)[lt]{\lineheight{1.25}\smash{\begin{tabular}[t]{l}$\textcolor{red}{\circled{1}}$\end{tabular}}}}%
    \put(0.4570048,0.25532424){\color[rgb]{0,0,0}\makebox(0,0)[lt]{\lineheight{1.25}\smash{\begin{tabular}[t]{l}$\textcolor{red}{\circled{2}}$\end{tabular}}}}%
    \put(0.67896035,0.52538128){\color[rgb]{0,0,0}\makebox(0,0)[lt]{\lineheight{1.25}\smash{\begin{tabular}[t]{l}$\textcolor{blue}{\circledd{$2'$}}$\end{tabular}}}}%
    \put(0.47719679,0.20502992){\color[rgb]{0,0,0}\makebox(0,0)[lt]{\lineheight{1.25}\smash{\begin{tabular}[t]{l}$\textcolor{red}{\circledd{$i\hspace{-1mm}-\hspace{-1mm}1$}}$\end{tabular}}}}%
    \put(0.4565754,0.21681178){\color[rgb]{0,0,0}\makebox(0,0)[lt]{\lineheight{1.25}\smash{\begin{tabular}[t]{l}$\textcolor{red}{\circled{3}}$\end{tabular}}}}%
    \put(0.40024023,0.1707617){\color[rgb]{0,0,0}\makebox(0,0)[lt]{\lineheight{1.25}\smash{\begin{tabular}[t]{l}$\textcolor{red}{\circledd{$i$}}$\end{tabular}}}}%
    \put(0.32636323,0.15912673){\color[rgb]{0,0,0}\makebox(0,0)[lt]{\lineheight{1.25}\smash{\begin{tabular}[t]{l}$\textcolor{red}{\circleddd{$i\hspace{-0,7mm}+\hspace{-0,7mm}1$}}$\end{tabular}}}}%
    \put(0.65932727,0.10021663){\color[rgb]{0,0,0}\makebox(0,0)[lt]{\lineheight{1.25}\smash{\begin{tabular}[t]{l}$\textcolor{red}{\circleddd{$i\hspace{-0,7mm}+\hspace{-0,7mm}2$}}$\end{tabular}}}}%
    \put(0.43570294,0.07194992){\color[rgb]{0,0,0}\makebox(0,0)[lt]{\lineheight{1.25}\smash{\begin{tabular}[t]{l}$\textcolor{red}{\circleddd{$i\hspace{-0,7mm}+\hspace{-0,7mm}3$}}$\end{tabular}}}}%
    \put(0.50861285,0.08758487){\color[rgb]{0,0,0}\makebox(0,0)[lt]{\lineheight{1.25}\smash{\begin{tabular}[t]{l}$\textcolor{red}{\circleddd{$i\hspace{-0,7mm}+\hspace{-0,7mm}4$}}$\end{tabular}}}}%
    \put(0.57406506,0.05756736){\color[rgb]{0,0,0}\makebox(0,0)[lt]{\lineheight{1.25}\smash{\begin{tabular}[t]{l}$\textcolor{red}{\circleddd{$i\hspace{-0,7mm}+\hspace{-0,7mm}5$}}$\end{tabular}}}}%
    \put(0.35851124,0.118645){\color[rgb]{0,0,0}\makebox(0,0)[lt]{\lineheight{1.25}\smash{\begin{tabular}[t]{l}$\textcolor{red}{\circleddd{$i\hspace{-0,7mm}+\hspace{-0,7mm}6$}}$\end{tabular}}}}%
    \put(0.69842481,0.13093028){\color[rgb]{0,0,0}\makebox(0,0)[lt]{\lineheight{1.25}\smash{\begin{tabular}[t]{l}$\textcolor{red}{\circleddd{$i\hspace{-0,7mm}+\hspace{-0,7mm}7$}}$\end{tabular}}}}%
    \put(0.68169123,0.16335169){\color[rgb]{0,0,0}\makebox(0,0)[lt]{\lineheight{1.25}\smash{\begin{tabular}[t]{l}$\textcolor{red}{\circleddd{$i\hspace{-0,7mm}+\hspace{-0,7mm}8$}}$\end{tabular}}}}%
    \put(0.71806019,0.18960485){\color[rgb]{0,0,0}\makebox(0,0)[lt]{\lineheight{1.25}\smash{\begin{tabular}[t]{l}$\textcolor{red}{\circleddd{$i\hspace{-0,7mm}+\hspace{-0,7mm}9$}}$\end{tabular}}}}%
    \put(0.75095434,0.20437441){\color[rgb]{0,0,0}\makebox(0,0)[lt]{\lineheight{1.25}\smash{\begin{tabular}[t]{l}$\textcolor{red}{\circleddd{$i\hspace{-0,7mm}+\hspace{-0,7mm}10$}}$\end{tabular}}}}%
    \put(0.70160624,0.21295667){\color[rgb]{0,0,0}\makebox(0,0)[lt]{\lineheight{1.25}\smash{\begin{tabular}[t]{l}$\textcolor{red}{\circleddd{$k_{i-2}$}}$\end{tabular}}}}%
    \put(0.76577901,0.51475996){\color[rgb]{0,0,0}\makebox(0,0)[lt]{\lineheight{1.25}\smash{\begin{tabular}[t]{l}$\textcolor{red}{\circleddd{$k_{i-2}\hspace{-0,7mm}+\hspace{-0,7mm}10$}}$\end{tabular}}}}%
    \put(0.76355911,0.22883777){\color[rgb]{0,0,0}\makebox(0,0)[lt]{\lineheight{1.25}\smash{\begin{tabular}[t]{l}$\textcolor{red}{\circleddd{$k_{i-2}\hspace{-0,7mm}+\hspace{-0,7mm}1$}}$\end{tabular}}}}%
    \put(0.79062686,0.55421869){\color[rgb]{0,0,0}\makebox(0,0)[lt]{\lineheight{1.25}\smash{\begin{tabular}[t]{l}$\textcolor{red}{\circleddd{$i\hspace{-0,7mm}+\hspace{-0,7mm}19$}}$\end{tabular}}}}%
    \put(0.74389572,0.5947815){\color[rgb]{0,0,0}\makebox(0,0)[lt]{\lineheight{1.25}\smash{\begin{tabular}[t]{l}$\textcolor{red}{\circleddd{$k_{i-1}$}}$\end{tabular}}}}%
    \put(0.62572439,0.70515212){\color[rgb]{0,0,0}\makebox(0,0)[lt]{\lineheight{1.25}\smash{\begin{tabular}[t]{l}$\textcolor{red}{\circleddd{$k_{i}$}}$\end{tabular}}}}%
    \put(0.59916432,0.71680053){\color[rgb]{0,0,0}\makebox(0,0)[lt]{\lineheight{1.25}\smash{\begin{tabular}[t]{l}$\textcolor{red}{\circleddd{$k_{i}\hspace{-0,7mm}+\hspace{-0,7mm}1$}}$\end{tabular}}}}%
    \put(0.39408107,0.67761684){\color[rgb]{0,0,0}\makebox(0,0)[lt]{\lineheight{1.25}\smash{\begin{tabular}[t]{l}$\textcolor{red}{\circleddd{$k_{i+1}$}}$\end{tabular}}}}%
    \put(0.31753189,0.66916354){\color[rgb]{0,0,0}\makebox(0,0)[lt]{\lineheight{1.25}\smash{\begin{tabular}[t]{l}$\textcolor{red}{\circleddd{$k_{g-1}\hspace{-0,7mm}+\hspace{-0,7mm}1$}}$\end{tabular}}}}%
    \put(0.11384817,0.65540154){\color[rgb]{0,0,0}\makebox(0,0)[lt]{\lineheight{1.25}\smash{\begin{tabular}[t]{l}$\textcolor{red}{\circleddd{$k_{g}$}}$\end{tabular}}}}%
    \put(0.59352723,0.75340978){\color[rgb]{0,0,0}\makebox(0,0)[lt]{\lineheight{1.25}\smash{\begin{tabular}[t]{l}$\textcolor{red}{\circleddd{$k_{i}+2$}}$\end{tabular}}}}%
    \put(0.42977272,0.77482962){\color[rgb]{0,0,0}\makebox(0,0)[lt]{\lineheight{1.25}\smash{\begin{tabular}[t]{l}$\textcolor{red}{\circleddd{$k_{i}\hspace{-0,7mm}+\hspace{-0,7mm}4$}}$\end{tabular}}}}%
    \put(0.41004839,0.7222314){\color[rgb]{0,0,0}\makebox(0,0)[lt]{\lineheight{1.25}\smash{\begin{tabular}[t]{l}$\textcolor{red}{\circleddd{$k_{i}\hspace{-0,7mm}+\hspace{-0,7mm}10$}}$\end{tabular}}}}%
    \put(0.77314123,0.6520958){\color[rgb]{0,0,0}\makebox(0,0)[lt]{\lineheight{1.25}\smash{\begin{tabular}[t]{l}$\textcolor{red}{\circleddd{$k_{i-1}\hspace{-0,7mm}+\hspace{-0,7mm}1$}}$\end{tabular}}}}%
    \put(0.76008564,0.71974258){\color[rgb]{0,0,0}\makebox(0,0)[lt]{\lineheight{1.25}\smash{\begin{tabular}[t]{l}$\textcolor{red}{\circleddd{$k_{i-1}+5$}}$\end{tabular}}}}%
    \put(0.81251176,0.65053621){\color[rgb]{0,0,0}\makebox(0,0)[lt]{\lineheight{1.25}\smash{\begin{tabular}[t]{l}$\textcolor{red}{\circleddd{$k_{i-1}\hspace{-0,7mm}+\hspace{-0,7mm}2$}}$\end{tabular}}}}%
    \put(0.32252087,0.26994708){\color[rgb]{0,0,0}\makebox(0,0)[lt]{\lineheight{1.25}\smash{\begin{tabular}[t]{l}$\textcolor{red}{\circledd{$N$}}$\end{tabular}}}}%
    \put(0.28649596,0.25622758){\color[rgb]{0,0,0}\makebox(0,0)[lt]{\lineheight{1.25}\smash{\begin{tabular}[t]{l}$\textcolor{red}{\circleddd{$N\hspace{-1mm}-\hspace{-1mm}1$}}$\end{tabular}}}}%
    \put(0.5979345,0.39389439){\color[rgb]{0,0,0}\makebox(0,0)[lt]{\lineheight{1.25}\smash{\begin{tabular}[t]{l}$\textcolor{blue}{\vartheta_i}$\end{tabular}}}}%
    \put(0.64343403,0.62588768){\color[rgb]{0,0,0}\makebox(0,0)[lt]{\lineheight{1.25}\smash{\begin{tabular}[t]{l}{\tiny $k_j=i\hspace{-0,7mm}+\hspace{-0,7mm}9\hspace{-0,7mm}+\hspace{-0,7mm}11(j-2)$}\end{tabular}}}}%
    \put(0.42326432,0.53870905){\color[rgb]{0,0,0}\makebox(0,0)[lt]{\lineheight{1.25}\smash{\begin{tabular}[t]{l}$\textcolor{yellowi}{\circleddd{$k_i+1'$}}$\end{tabular}}}}%
    \put(0.40263944,0.32707861){\color[rgb]{0,0,0}\makebox(0,0)[lt]{\lineheight{1.25}\smash{\begin{tabular}[t]{l}$\textcolor{red}{\gamma}$\end{tabular}}}}%
    \put(0.43468194,0.5993289){\color[rgb]{0,0,0}\makebox(0,0)[lt]{\lineheight{1.25}\smash{\begin{tabular}[t]{l}$\textcolor{yellowi}{\zeta_{i}}$\end{tabular}}}}%
    \put(0,0){\includegraphics[width=\unitlength,page=3]{one-holed_surface_flat5more.pdf}}%
    \put(0.68507331,0.73108898){\color[rgb]{0,0,0}\makebox(0,0)[lt]{\lineheight{1.25}\smash{\begin{tabular}[t]{l}$\textcolor{red}{\circleddd{$k_{i-1}\hspace{-0,7mm}+\hspace{-0,7mm}9$}}$\end{tabular}}}}%
    \put(0.64768129,0.69875005){\color[rgb]{0,0,0}\makebox(0,0)[lt]{\lineheight{1.25}\smash{\begin{tabular}[t]{l}$\textcolor{red}{\circleddd{$k_{i-1}\hspace{-0,7mm}+\hspace{-0,7mm}10$}}$\end{tabular}}}}%
    \put(0.71911254,0.76355226){\color[rgb]{0,0,0}\makebox(0,0)[lt]{\lineheight{1.25}\smash{\begin{tabular}[t]{l}$\textcolor{red}{\circleddd{$k_{i-1}\hspace{-0,7mm}+\hspace{-0,7mm}3$}}$\end{tabular}}}}%
    \put(0,0){\includegraphics[width=\unitlength,page=4]{one-holed_surface_flat5more.pdf}}%
    \put(0.27543233,0.29338541){\color[rgb]{0,0,0}\makebox(0,0)[lt]{\lineheight{1.25}\smash{\begin{tabular}[t]{l}$\textcolor{red}{\circleddd{$k_g\hspace{-0,7mm}+\hspace{-0,7mm}10$}}$\end{tabular}}}}%
    \put(0,0){\includegraphics[width=\unitlength,page=5]{one-holed_surface_flat5more.pdf}}%
    \put(0.29463708,0.14830004){\color[rgb]{0,0,0}\makebox(0,0)[lt]{\lineheight{1.25}\smash{\begin{tabular}[t]{l}$a_1$\end{tabular}}}}%
    \put(0,0){\includegraphics[width=\unitlength,page=6]{one-holed_surface_flat5more.pdf}}%
    \put(0.39606726,0.05490747){\color[rgb]{0,0,0}\makebox(0,0)[lt]{\lineheight{1.25}\smash{\begin{tabular}[t]{l}$b_1$\end{tabular}}}}%
    \put(0,0){\includegraphics[width=\unitlength,page=7]{one-holed_surface_flat5more.pdf}}%
    \put(0.61091431,0.03885796){\color[rgb]{0,0,0}\makebox(0,0)[lt]{\lineheight{1.25}\smash{\begin{tabular}[t]{l}$a_1$\end{tabular}}}}%
    \put(0.79796142,0.13868549){\color[rgb]{0,0,0}\makebox(0,0)[lt]{\lineheight{1.25}\smash{\begin{tabular}[t]{l}$a_2$\end{tabular}}}}%
    \put(0.97561693,0.49270626){\color[rgb]{0,0,0}\makebox(0,0)[lt]{\lineheight{1.25}\smash{\begin{tabular}[t]{l}$a_2$\end{tabular}}}}%
    \put(0.91912666,0.21971721){\color[rgb]{0,0,0}\makebox(0,0)[lt]{\lineheight{1.25}\smash{\begin{tabular}[t]{l}$b_2$\end{tabular}}}}%
    \put(0.79849363,0.12400022){\color[rgb]{0,0,0}\makebox(0,0)[lt]{\lineheight{1.25}\smash{\begin{tabular}[t]{l}($a_{i-2}$)\end{tabular}}}}%
    \put(0.91625555,0.20019311){\color[rgb]{0,0,0}\makebox(0,0)[lt]{\lineheight{1.25}\smash{\begin{tabular}[t]{l}($b_{i-2}$)\end{tabular}}}}%
    \put(0,0){\includegraphics[width=\unitlength,page=8]{one-holed_surface_flat5more.pdf}}%
    \put(0.72796477,0.11434079){\color[rgb]{0,0,0}\makebox(0,0)[lt]{\lineheight{1.25}\smash{\begin{tabular}[t]{l}$b_1$\end{tabular}}}}%
    \put(0,0){\includegraphics[width=\unitlength,page=9]{one-holed_surface_flat5more.pdf}}%
    \put(0.85582671,0.63759115){\color[rgb]{0,0,0}\makebox(0,0)[lt]{\lineheight{1.25}\smash{\begin{tabular}[t]{l}$a_{i-1}$\end{tabular}}}}%
    \put(0,0){\includegraphics[width=\unitlength,page=10]{one-holed_surface_flat5more.pdf}}%
    \put(0.85323971,0.72314844){\color[rgb]{0,0,0}\makebox(0,0)[lt]{\lineheight{1.25}\smash{\begin{tabular}[t]{l}$b_{i-1}$\end{tabular}}}}%
    \put(0,0){\includegraphics[width=\unitlength,page=11]{one-holed_surface_flat5more.pdf}}%
    \put(0.68036541,0.77485666){\color[rgb]{0,0,0}\rotatebox{31.178138}{\makebox(0,0)[lt]{\lineheight{1.25}\smash{\begin{tabular}[t]{l}$b_{i-1}$\end{tabular}}}}}%
    \put(0,0){\includegraphics[width=\unitlength,page=12]{one-holed_surface_flat5more.pdf}}%
    \put(0.76460672,0.79357506){\color[rgb]{0,0,0}\rotatebox{-26.868564}{\makebox(0,0)[lt]{\lineheight{1.25}\smash{\begin{tabular}[t]{l}$a_{i-1}$\end{tabular}}}}}%
    \put(0,0){\includegraphics[width=\unitlength,page=13]{one-holed_surface_flat5more.pdf}}%
    \put(0.63977153,0.78801918){\color[rgb]{0,0,0}\makebox(0,0)[lt]{\lineheight{1.25}\smash{\begin{tabular}[t]{l}$a_{i}$\end{tabular}}}}%
    \put(0,0){\includegraphics[width=\unitlength,page=14]{one-holed_surface_flat5more.pdf}}%
    \put(0.58910587,0.8311387){\color[rgb]{0,0,0}\makebox(0,0)[lt]{\lineheight{1.25}\smash{\begin{tabular}[t]{l}$b_{i}$\end{tabular}}}}%
    \put(0.26321937,0.80183395){\color[rgb]{0,0,0}\makebox(0,0)[lt]{\lineheight{1.25}\smash{\begin{tabular}[t]{l}$b_{g-1}$\end{tabular}}}}%
    \put(0,0){\includegraphics[width=\unitlength,page=15]{one-holed_surface_flat5more.pdf}}%
    \put(0.44934187,0.81899718){\color[rgb]{0,0,0}\makebox(0,0)[lt]{\lineheight{1.25}\smash{\begin{tabular}[t]{l}$a_{i}$\end{tabular}}}}%
    \put(0.14633732,0.79683839){\color[rgb]{0,0,0}\makebox(0,0)[lt]{\lineheight{1.25}\smash{\begin{tabular}[t]{l}$a_{g-1}$\end{tabular}}}}%
    \put(0.33457709,0.78689606){\color[rgb]{0,0,0}\rotatebox{-59.153621}{\makebox(0,0)[lt]{\lineheight{1.25}\smash{\begin{tabular}[t]{l}$a_{g-1}$\end{tabular}}}}}%
    \put(0,0){\includegraphics[width=\unitlength,page=16]{one-holed_surface_flat5more.pdf}}%
    \put(0.38117268,0.75849757){\color[rgb]{0,0,0}\rotatebox{57.93463}{\makebox(0,0)[lt]{\lineheight{1.25}\smash{\begin{tabular}[t]{l}$b_i$\end{tabular}}}}}%
    \put(0.06003478,0.75005727){\color[rgb]{0,0,0}\makebox(0,0)[lt]{\lineheight{1.25}\smash{\begin{tabular}[t]{l}$b_{g-1}$\end{tabular}}}}%
    \put(0,0){\includegraphics[width=\unitlength,page=17]{one-holed_surface_flat5more.pdf}}%
    \put(0.02730494,0.60641026){\color[rgb]{0,0,0}\makebox(0,0)[lt]{\lineheight{1.25}\smash{\begin{tabular}[t]{l}$a_g$\end{tabular}}}}%
    \put(0,0){\includegraphics[width=\unitlength,page=18]{one-holed_surface_flat5more.pdf}}%
    \put(0.00964149,0.39956613){\color[rgb]{0,0,0}\makebox(0,0)[lt]{\lineheight{1.25}\smash{\begin{tabular}[t]{l}$b_g$\end{tabular}}}}%
    \put(0,0){\includegraphics[width=\unitlength,page=19]{one-holed_surface_flat5more.pdf}}%
    \put(0.08669997,0.25158052){\color[rgb]{0,0,0}\makebox(0,0)[lt]{\lineheight{1.25}\smash{\begin{tabular}[t]{l}$a_g$\end{tabular}}}}%
    \put(0,0){\includegraphics[width=\unitlength,page=20]{one-holed_surface_flat5more.pdf}}%
    \put(0.21147848,0.19188221){\color[rgb]{0,0,0}\makebox(0,0)[lt]{\lineheight{1.25}\smash{\begin{tabular}[t]{l}$b_g$\end{tabular}}}}%
    \put(0.2962387,0.18978019){\color[rgb]{0,0,0}\makebox(0,0)[lt]{\lineheight{1.25}\smash{\begin{tabular}[t]{l}{\small $p_1$}\end{tabular}}}}%
    \put(0,0){\includegraphics[width=\unitlength,page=21]{one-holed_surface_flat5more.pdf}}%
    \put(0.44030738,0.87061042){\color[rgb]{0,0,0}\makebox(0,0)[lt]{\lineheight{1.25}\smash{\begin{tabular}[t]{l}$\mbox{$i$-th handle}$\end{tabular}}}}%
    \put(0.12413759,0.85246812){\color[rgb]{0,0,0}\makebox(0,0)[lt]{\lineheight{1.25}\smash{\begin{tabular}[t]{l}$\mbox{$(g-1)$-th handle}$\end{tabular}}}}%
    \put(0,0){\includegraphics[width=\unitlength,page=22]{one-holed_surface_flat5more.pdf}}%
    \put(0.81333202,0.83067407){\color[rgb]{0,0,0}\rotatebox{-40.165541}{\makebox(0,0)[lt]{\lineheight{1.25}\smash{\begin{tabular}[t]{l}$\mbox{$(i-1)$-th handle}$\end{tabular}}}}}%
    \put(0.94744117,0.69819473){\color[rgb]{0,0,0}\rotatebox{-53.828438}{\makebox(0,0)[lt]{\lineheight{1.25}\smash{\begin{tabular}[t]{l}$\cdots$\end{tabular}}}}}%
    \put(0.32803901,0.8579103){\color[rgb]{0,0,0}\rotatebox{2.5916418}{\makebox(0,0)[lt]{\lineheight{1.25}\smash{\begin{tabular}[t]{l}$\cdots$\end{tabular}}}}}%
    \put(0.48147921,0.24937698){\color[rgb]{0,0,0}\rotatebox{-82.528728}{\makebox(0,0)[lt]{\lineheight{1.25}\smash{\begin{tabular}[t]{l}\red{$\cdots$}\end{tabular}}}}}%
    \put(0,0){\includegraphics[width=\unitlength,page=23]{one-holed_surface_flat5more.pdf}}%
    \put(0.16803457,0.03118133){\color[rgb]{0,0,0}\rotatebox{0.62622378}{\makebox(0,0)[lt]{\lineheight{1.25}\smash{\begin{tabular}[t]{l}$\mbox{first handle}$\end{tabular}}}}}%
    \put(0.87859934,0.06745959){\color[rgb]{0,0,0}\rotatebox{0.62622378}{\makebox(0,0)[lt]{\lineheight{1.25}\smash{\begin{tabular}[t]{l}$\mbox{2nd handle}$\end{tabular}}}}}%
    \put(0.87834931,0.04181833){\color[rgb]{0,0,0}\rotatebox{0.62622378}{\makebox(0,0)[lt]{\lineheight{1.25}\smash{\begin{tabular}[t]{l}$\mbox{(($i-2$)-th)}$\end{tabular}}}}}%
    \put(0,0){\includegraphics[width=\unitlength,page=24]{one-holed_surface_flat5more.pdf}}%
    \put(0.02701912,0.12786106){\color[rgb]{0,0,0}\rotatebox{-1.0688603}{\makebox(0,0)[lt]{\lineheight{1.25}\smash{\begin{tabular}[t]{l}$\mbox{last handle}$\end{tabular}}}}}%
    \put(0.02784924,0.14958033){\color[rgb]{0,0,0}\rotatebox{-1.0688603}{\makebox(0,0)[lt]{\lineheight{1.25}\smash{\begin{tabular}[t]{l}$\mbox{($g$-th)}$\end{tabular}}}}}%
    \put(0.49589059,0.56303106){\color[rgb]{0,0,0}\makebox(0,0)[lt]{\lineheight{1.25}\smash{\begin{tabular}[t]{l}$\mbox{sphere-body}$\end{tabular}}}}%
    \put(0,0){\includegraphics[width=\unitlength,page=25]{one-holed_surface_flat5more.pdf}}%
    \put(0.78806483,0.68687273){\color[rgb]{0,0,0}\makebox(0,0)[lt]{\lineheight{1.25}\smash{\begin{tabular}[t]{l}$c_{i-1}$\end{tabular}}}}%
    \put(0.72577823,0.73597889){\color[rgb]{0,0,0}\makebox(0,0)[lt]{\lineheight{1.25}\smash{\begin{tabular}[t]{l}$d_{i-1}$\end{tabular}}}}%
    \put(0.74094396,0.67515492){\color[rgb]{0,0,0}\makebox(0,0)[lt]{\lineheight{1.25}\smash{\begin{tabular}[t]{l}$e_{i-1}$\end{tabular}}}}%
    \put(0.57600131,0.78669715){\color[rgb]{0,0,0}\rotatebox{-36.532084}{\makebox(0,0)[lt]{\lineheight{1.25}\smash{\begin{tabular}[t]{l}$c_{i}$\end{tabular}}}}}%
    \put(0.2622204,0.75423046){\color[rgb]{0,0,0}\rotatebox{-29.474133}{\makebox(0,0)[lt]{\lineheight{1.25}\smash{\begin{tabular}[t]{l}$c_{g-1}$\end{tabular}}}}}%
    \put(0.5198043,0.73814176){\color[rgb]{0,0,0}\makebox(0,0)[lt]{\lineheight{1.25}\smash{\begin{tabular}[t]{l}$e_{i}$\end{tabular}}}}%
    \put(0.2090619,0.69831985){\color[rgb]{0,0,0}\makebox(0,0)[lt]{\lineheight{1.25}\smash{\begin{tabular}[t]{l}$e_{g-1}$\end{tabular}}}}%
    \put(0.45629879,0.75225518){\color[rgb]{0,0,0}\makebox(0,0)[lt]{\lineheight{1.25}\smash{\begin{tabular}[t]{l}$d_{i}$\end{tabular}}}}%
    \put(0.15431048,0.7269036){\color[rgb]{0,0,0}\makebox(0,0)[lt]{\lineheight{1.25}\smash{\begin{tabular}[t]{l}$d_{g-1}$\end{tabular}}}}%
    \put(0.0820804,0.46783236){\color[rgb]{0,0,0}\makebox(0,0)[lt]{\lineheight{1.25}\smash{\begin{tabular}[t]{l}$c_g$\end{tabular}}}}%
    \put(0.16100431,0.29163954){\color[rgb]{0,0,0}\makebox(0,0)[lt]{\lineheight{1.25}\smash{\begin{tabular}[t]{l}$d_g$\end{tabular}}}}%
    \put(0.18408898,0.47270657){\color[rgb]{0,0,0}\makebox(0,0)[lt]{\lineheight{1.25}\smash{\begin{tabular}[t]{l}$e_g$\end{tabular}}}}%
    \put(0.41507627,0.11920739){\color[rgb]{0,0,0}\makebox(0,0)[lt]{\lineheight{1.25}\smash{\begin{tabular}[t]{l}$c_1$\end{tabular}}}}%
    \put(0.60327575,0.10749552){\color[rgb]{0,0,0}\makebox(0,0)[lt]{\lineheight{1.25}\smash{\begin{tabular}[t]{l}$d_1$\end{tabular}}}}%
    \put(0.71334543,0.4297382){\color[rgb]{0,0,0}\rotatebox{52.066616}{\makebox(0,0)[lt]{\lineheight{1.25}\smash{\begin{tabular}[t]{l}$f_3$\end{tabular}}}}}%
    \put(0.68784768,0.44023304){\color[rgb]{0,0,0}\rotatebox{51.742243}{\makebox(0,0)[lt]{\lineheight{1.25}\smash{\begin{tabular}[t]{l}$f_{i-2}$\end{tabular}}}}}%
    \put(0.44694118,0.49150121){\color[rgb]{0,0,0}\makebox(0,0)[lt]{\lineheight{1.25}\smash{\begin{tabular}[t]{l}$f_{i-1}$\end{tabular}}}}%
    \put(0.56551919,0.43327135){\color[rgb]{0,0,0}\makebox(0,0)[lt]{\lineheight{1.25}\smash{\begin{tabular}[t]{l}$f_{i-1}'$\end{tabular}}}}%
    \put(0.53541874,0.17933429){\color[rgb]{0,0,0}\makebox(0,0)[lt]{\lineheight{1.25}\smash{\begin{tabular}[t]{l}$f_2$\end{tabular}}}}%
    \put(0.84717979,0.26172812){\color[rgb]{0,0,0}\makebox(0,0)[lt]{\lineheight{1.25}\smash{\begin{tabular}[t]{l}$c_2$\end{tabular}}}}%
    \put(0.88444259,0.46966398){\color[rgb]{0,0,0}\makebox(0,0)[lt]{\lineheight{1.25}\smash{\begin{tabular}[t]{l}$d_2$\end{tabular}}}}%
    \put(0.78829528,0.38915785){\color[rgb]{0,0,0}\makebox(0,0)[lt]{\lineheight{1.25}\smash{\begin{tabular}[t]{l}$e_2$\end{tabular}}}}%
    \put(0.35575912,0.58231126){\color[rgb]{0,0,0}\makebox(0,0)[lt]{\lineheight{1.25}\smash{\begin{tabular}[t]{l}$f_{i}$\end{tabular}}}}%
    \put(0.29329682,0.46860366){\color[rgb]{0,0,0}\makebox(0,0)[lt]{\lineheight{1.25}\smash{\begin{tabular}[t]{l}$f_{g-1}$\end{tabular}}}}%
    \put(0.63542399,0.3695805){\color[rgb]{0,0,0}\rotatebox{-42.548696}{\makebox(0,0)[lt]{\lineheight{1.25}\smash{\begin{tabular}[t]{l}$\cdots$\end{tabular}}}}}%
    \put(0.33615037,0.52133){\color[rgb]{0,0,0}\rotatebox{-1.0250199}{\makebox(0,0)[lt]{\lineheight{1.25}\smash{\begin{tabular}[t]{l}$\cdots$\end{tabular}}}}}%
    \put(0,0){\includegraphics[width=\unitlength,page=26]{one-holed_surface_flat5more.pdf}}%
  \end{picture}%
\endgroup%